\documentclass[reqno]{amsart}
\usepackage{amsmath,amssymb,stmaryrd}
\usepackage{amsrefs}
\usepackage{enumerate}
\usepackage{amsfonts}
\usepackage[colorlinks=true, pdfborder={0 0 0}]{hyperref}
\usepackage{cleveref}
\usepackage{upref}

\newtheorem{theorem}{Theorem}
\newtheorem{defi}{Definition}
\newtheorem{proposition}{Proposition}
\newtheorem{coro}{Corollary}
\newtheorem{lemma}{Lemma}
\newtheorem{step}{Step}

\newtheorem{thdefi}{Theorem-Definition}

\def \rr {\mathbb{R}}
\def \nn {\mathbb{N}}
\def \rn {\mathbb{R}^n}
\def \rnm {\mathbb{R}_{-}^n}
\def \rnmp {\rnm \setminus\{0\}}
\def \rnmpbar  {\overline{\rnm} \setminus\{0\}}

\def \R {\mathbb{R}}

\def \eps {\epsilon}
\def \crit {2^\star(s)}
\def \ue {u_\epsilon}
\def \ve {v_\epsilon}
\def \re {r_\eps}
\def \he {h_\eps}
\def \pe {p_\eps}
\def \tge {\tilde{g}_\eps}
\def \ye {y_\eps}
\def \ze {z_\eps}

\def \elle {\ell_\eps}
\def \tv {\tilde{v}}
\def \tu {\tilde{u}}
\def \tw {\tilde{w}}
\def \tve {\tilde{v}_\eps}
\def \tye {\tilde{y}_\eps}
\def \tue {\tilde{u}_\eps}
\def \twe {\tilde{w}_{\eps}}
\def \tuei {\tilde{u}_{i,\eps}}
\def \tueun {\tilde{u}_{1,\eps}}
\def \tui {\tilde{u}_i}
\def \bue {\overline{u}_\eps}
\def \tze {\tilde{z}_\eps}
\def \mei {\mu_{i,\eps}}
\def \kei {k_{i,\eps}}
\def \meun {\mu_{1,\eps}}
\def \keun {k_{1,\eps}}
\def \xeun {x_{1,\eps}}
\def \meN {\mu_{N,\eps}}
\def \keN {k_{N,\eps}}
\def \Ono {\Omega}
\def \bdry {\partial \Omega}
\def \bono {\partial \Omega \setminus \{ 0 \} }
\def \T {\mathcal{T}}
\def \huno {H_{1,0}^2(\Omega)}
\def \ds {\displaystyle}
\def \beqn {\begin{eqnarray}}
\def \eeqn {\end{eqnarray}}
\def \bequa {\begin{equation}}
\def \eequa {\end{equation}}
\def \bp {\beta_+(\gamma)}
\def \bm {\beta_-(\gamma)}
\def \eucl {\hbox{Eucl}}
\def \Omegabar {\overline{\Omega}}
\def \I {\mathcal{I}}
\def \rnp {\rn\setminus\{0\}}
\def \crits {2^\star(s)}
\def \txi {\tilde{x}_i}
\def \tyi {\tilde{y}_i}

\addtocontents{toc}{\protect\setcounter{tocdepth}{1}}

\title[Multiplicity and Pohozev stability for Hardy-Schr\"odinger equations]{Multiplicity and stability of the Pohozaev obstruction for Hardy-Schr\"odinger equations with boundary singularity}
\date{March 11th, 2020}
 
\author{Nassif Ghoussoub}
\address{Nassif Ghoussoub: Department of Mathematics, University of British Columbia, Vancouver, V6T 1Z2 Canada}
\email{nassif@math.ubc.ca}
\author{Saikat Mazumdar}
\address{Saikat Mazumdar:
Department of Mathematics,
Indian Institute of Technology Bombay,
Mumbai 400076, India.}
\email{saikat@math.iitb.ac.in, saikat.mazumdar@iitb.ac.in}
\author{Fr\'ed\'eric Robert}
\address{Fr\'ed\'eric Robert, UniversitŽ de Lorraine, CNRS, IECL, F-54000 Nancy, France}
\email{frederic.robert@univ-lorraine.fr} 

\thanks{This work was initiated when the second-named author held a postdoctoral position at the University of British Columbia under the supervision of the first-named author, that was partially supported by the Natural Sciences and Engineering Research Council of Canada (NSERC)}

\begin{document}

\begin{abstract}
Let $\Omega$ be a smooth bounded domain in $\rn$ ($n\geq 3$) such that $0\in\partial \Omega$. We consider issues of non-existence, existence, and multiplicity of variational solutions in $\huno$ for the borderline Dirichlet problem,
$$\left\{ \begin{array}{llll}
-\Delta u-\gamma \frac{u}{|x|^2}- h(x) u
&=& \frac{|u|^{\crits-2}u}{|x|^s}  \ \ &\text{in } \Omega,\\
\hfill u&=&0 &\text{on }\bono, 
\end{array} \right.\eqno{(E)}$$
where  $0<s<2$, ${\crits}:=\frac{2(n-s)}{n-2}$, $\gamma\in\rr$ and  
$h\in C^0(\overline{\Omega})$.   
We use sharp blow-up analysis on  --possibly high energy-- solutions of corresponding subcritical problems to establish, for example, that if $\gamma<\frac{n^2}{4}-1$ and the principal curvatures of $\partial\Omega$ at $0$ are  non-positive but not all of them vanishing, then Equation (E) has an infinite number of high energy (possibly sign-changing) solutions in $\huno$. This complements results of the first and third authors, who showed in \cite{gr4} that if $\gamma\leq \frac{n^2}{4}-\frac{1}{4}$ and the mean curvature of $\partial\Omega$ at $0$ is negative, then $(E)$ has a positive least energy solution. 

On the other hand, the sharp blow-up analysis also allows us to show that if the mean curvature at $0$ is nonzero and the mass, when defined, is also nonzero, then there is a surprising stability of regimes where there are no variational positive solutions under $C^1$-perturbations of the potential $h$. In particular, and in sharp contrast with the non-singular case (i.e., when $\gamma=s=0$), we prove non-existence of such solutions  for $(E)$ in any dimension, whenever $\Omega$ is star-shaped and $h$ is close to $0$, which include situations not covered by the classical Pohozaev obstruction. 

\end{abstract}

\maketitle
\tableofcontents

\section{\, Introduction}
This manuscript is the continuation of a long-time project initiated by the first and the third author in \cite{gr1} about nonlinear critical equations involving the Hardy potential when the singularity is located on the boundary of the domain under study.
Let $\Omega$ be such a smooth bounded domain in $\rn$, $n\geq 3$, with $0\in\partial\Omega$. We fix $s\in (0,2)$ and define the critical Sobolev exponent ${\crits}:=\frac{2(n-s)}{n-2}$. For $\gamma\in\rr$ and $h_0\in C^1(\overline{\Omega})$, we consider in the sequel issues of non-existence, existence, and multiplicity of variational solutions in $\huno$ for the borderline Dirichlet problem,
 \begin{eqnarray} \label{one}
\left\{ \begin{array}{llll}
-\Delta u-\gamma \frac{u}{|x|^2}-h_0(x) u
&=& \frac{|u|^{\crits-2}u}{|x|^s}  \ \ &\text{in } \Omega,\\
\hfill u&=&0 &\text{on }\bono. 
\end{array} \right.
\end{eqnarray}
By solutions, we mean here functions $u\in\huno$, i.e., the completion of $C^{\infty}_c(\Omega)$ for the $L_2$-norm of the gradient $\Vert\nabla u\Vert_2$. This problem has by now a long history starting with the fact that when $\gamma=s=0$ and $h_0$ is a constant, it is the counterpart of the Yamabe problem \cites{aubin, LeeParker, schoen1} in Euclidian space, as initiated by Brezis-Nirenberg \cite{bn}, with important contributions in the critical dimension $n=3$, by Druet \cite{d2}, and for multiplicity results for $n\geq 7$, by Devillanova-Solimini \cite{ds}, among many others. 

The case dealing with least energy solutions for $s>0$ but $\gamma =0$, when the singularity $0$ is on the boundary of the domain was initiated by Ghoussoub-Kang \cite{gk} and developed by Ghoussoub-Robert \cite{gr1}. The case involving the Hardy potential, i.e., when $\gamma >0$,  was introduced by Lin-Wadade \cite{LW3} with a follow-up contribution by Ghoussoub-Robert \cite{gr4}. This paper addresses remaining issues about the multiplicity of solutions, but also about obstructions to the existence of solutions and their stability under small perturbations. 

 The existence of solutions is related to the coercivity of the operator$-\Delta -\frac{\gamma}{|x|^2} -h_0(x)$. It is clear that the operator $-\Delta - \frac{\gamma}{|x|^2}$ is coercive on $\huno$ whenever $\gamma <\gamma_H(\Omega)$, where $\gamma_H(\Omega)$ is the Hardy constant associated to the domain $\Omega$, that is  
\begin{align}\label{Hardy inequality}
\gamma_{H}(\Omega):=  \inf_{u \in\huno \setminus \{ 0\}}\frac{ \int _{\Omega} |\nabla u|^{2} ~ dx}{ \int_{\Omega} \frac{u^{2}}{|x|^{2}}~ dx},
\end{align}
which  has been extensively studied (see for example \cite{GM.book} and \cite{gr4}). We recall that if $0\in \Omega$, then 
 \begin{equation}
 \gamma_H (\Omega)=\gamma_{H}(\rn)=\frac{(n-2)^2}{4}.
 \end{equation}
 When $0\in \partial \Omega$, the situation is extremely different. For non-smooth domains modeled on cones, we refer to Egnell \cite{eg}, and the more recent works of Cheikh-Ali \cites{HCA,HCA2}. If $\Omega$ is smooth, then, around $0$, the domain is modeled on the half-space $\rnm:=\{x\in \rn; \, x_1<0\}$. We then get that (see \cite{gr4})
\begin{equation}
\frac{(n-2)^2}{4} < \gamma_{H}(\Omega)\leq\gamma_H(\rnm)= \frac{n^2}{4}.
\end{equation}

Note that when $h_{0}\equiv 0$, (\ref{one}) is the Euler-Lagrange equation for the following Hardy-Sobolev variational problem: For $\gamma <\gamma_H(\Omega)$ and $0\leq s\leq 2$, there exists $\mu_{\gamma, s}(\Omega)>0$ such that
\begin{equation} \label{general}
\mu_{\gamma, s}(\Omega)=\inf\left\{\frac{\int_{\Omega} |\nabla u|^2\, dx-\gamma \int_{\Omega}\frac{u^2}{|x|^2}\, dx}{\left(\int_{\Omega}\frac{|u|^{\crit}}{|x|^s}\, dx\right)^{\frac{2}{\crit}}};\, u \in\huno \setminus \{ 0\} \right\}.
\end{equation}
Note that when $s=2$ and $\gamma=0$, this is the Hardy inequality mentioned above, while if $s=0$ and $\gamma=0$, it is the Sobolev inequality.  If $\Omega =\R^n$, $s\in [0, 2]$ and $\gamma \in (-\infty , \frac{(n-2)^2}{4})$, (\ref{general})  contains -- after a suitable change of variables -- 
the Caffarelli-Kohn-Nirenberg inequalities \cite{ckn}. The latter state that there is a constant $C:=C(a,b,n)>0$ such that,  
\begin{equation} \label{CKN}
\left(\int_{\rn}|x|^{-bq}|u|^q \right)^{\frac{2}{q}}\leq C\int_{\rn}|x|^{-2a}|\nabla u|^2 dx \hbox{ for all }u\in C^\infty_c(\rn),
\end{equation}
where
\begin{equation}\label{cond1}
-\infty<a<\frac{n-2}{2}, \ \ 0 \leq b-a\leq 1, \ \ {\rm and}\ \ q=\frac{2n}{n-2+2(b-a)}.
\end{equation}
The first difficulty in these problems is due to the fact that  $\crits$ is critical from the viewpoint of the Sobolev embeddings, in such a way that if $\Omega$ is bounded, then $\huno$ is embedded in the weighted space $L^p(\Omega, |x|^{-s})$ for $1\leq p\leq\crits $, and the embedding is compact if and only if  $p<\crits$.
This lack of compactness defeats the classical minimization strategy to get extremals for (\ref{general}). In fact, when $s=0$ and $\gamma=0$, this is the setting of the critical case in the classical Sobolev inequalities, which started this whole line of inquiry, due to its connection with the Yamabe problem on compact Riemannian manifolds \cite{aubin}, \cite{schoen1}, \cite{LeeParker}.  Another complicating feature of the problem is that the  term $\frac{u}{|x|^2}$ is as critical as $\frac{u^{2^*(s)-1}}{|x|^s}$ in the sense that they have the same homogeneity as the Laplacian. These difficulties are summarized by the invariance of the problem under conformal transformation. Indeed, for a function $u:\Omega\to\rr$ and $r>0$, let
\begin{equation}\label{conf:trans}
u_r: x\mapsto r^{\frac{n-2}{2}}u(r \cdot x)
\end{equation}
and note that Equation \eqref{one} is then "essentially" invariant under the transformation $u\mapsto u_r$ in the sense that
\begin{eqnarray}\label{eq:trans} 
\left\{ \begin{array}{llll}
-\Delta u_r-\gamma \frac{u_r}{|x|^2}-r^2h_{0}(rx) u_r
&=& \frac{|u_r|^{\crits-2}u_r}{|x|^s}  \ \ &\text{in } r^{-1}\Omega,\\
\hfill u_r&=&0 &\text{on }r^{-1}\bono. 
\end{array} \right.
\end{eqnarray}
This   "invariance" is behind the lack of compactness in the embeddings associated to the variational formulation of  \eqref{one}, which prohibits the use of general abstract topological or variational methods. However, as one notices, the invariance is not complete, since the potential $h$ has  changed, and the domain itself was transformed. As we shall see, both the geometry of the domain and -to a lesser extent- the potential $h$ break the invariance enough that one will be able to recover compactness for  \eqref{one}.

Another important aspect of this problem is the singularity at $0$ and its location within the domain since the Hardy potential does not belong to the Kato class. Elliptic problems with singular potential arise in quantum mechanics, astrophysics, as well as in Riemannian geometry, in particular  in the study of the scalar curvature problem on the standard sphere. Indeed, if the latter is equipped with its standard metric whose scalar curvature is singular at the north and south poles, then by considering its stereographic projection of $\rn$, the problem of finding a conformal metric with prescribed scalar curvature $K(x)$ leads to finding solutions of the form $-\Delta u-\gamma \frac{u}{|x|^2}=K(x)u^{2^\star(0)-1}$ on $\rn$. The latter is a simplified version of the nonlinear Wheeler-DeWitt equation, which appears in quantum cosmology (see \cites{BB,BE,LZ,SmetsTAMS} and the references cited therein).

This paper deals specifically with the case where $0$ belongs to the boundary of a smooth domain $\Omega$. We shall see that the boundary at $0$ plays an important role, and our starting point is the existence Theorem \ref{gro} below for least energy solutions. It was first established by Ghoussoub-Robert \cite{gr1} when $\gamma=0$, then by Lin-Wadade \cite{LW3} when $0<\gamma <\frac{(n-2)^2}{4}$ under the assumption that the mean curvature at $0$ is negative. The result was extended to the range $\gamma\leq \frac{n^2-1}{4}$ in \cite{gr4}, but more importantly, it was shown there that in the remaining range $(\frac{n^2-1}{4}, \frac{n^2}{4})$, the curvature condition does not suffice anymore and a more global condition is needed: the boundary mass $m_{\gamma,h}(\Omega)$ of a domain associated to $\gamma$ and $h$, that we now recall. 

\subsection{The models and the definition of the mass} Letting formally $r\to 0$ in \eqref{eq:trans}, we get that $u$ should behave like solutions to 
\begin{eqnarray}\label{eq:U:sec2} 
\left\{ \begin{array}{llll}
-\Delta U-\gamma \frac{U}{|x|^2}
&=& \frac{|U|^{\crits-2}U}{|x|^s}  \ \ &\text{ in }\rnm,\\
\hfill U&=&0 &\text{on }\partial\rnm. 
\end{array} \right.
\end{eqnarray}
To the best of our knowledge, no explicit positive solution of \eqref{eq:U:sec2} is known. This was the reason why a specific blowup analysis was carried out in \cite{gr1}, which relied on the symmetry properties and a precise description of the asymptotic behavior of such solutions --also established in \cite{gr1}. On the other hand, the asymptotic behavior of such nonlinear problems is governed by the solutions to the linear problem 
\begin{eqnarray}\label{eq:U:lin} 
\left\{ \begin{array}{llll}
-\Delta U-\gamma \frac{U}{|x|^2}
&=& 0  \ \ &\text{ in }\rnm,\\
\hfill U&=&0 &\text{on }\partial\rnm. 
\end{array} \right.
\end{eqnarray}
One can then easily see that a function of the form $u(x)=x_1|x|^{-\beta}$ is a solution to (\ref{eq:U:lin}) if and only if $\beta\in \{\beta_-(\gamma),\beta_+(\gamma)\}$, where 
\begin{align}\label{def:beta}
 \beta_{\pm}(\gamma):=\frac{n}{2}\pm\sqrt{\frac{n^2}{4}-\gamma} \quad \hbox{ for }\gamma<\frac{n^2}{4}. 
 \end{align}

\begin{thdefi}[\cite{gr4}]\label{thdefi:mass} Let $\Omega$ be a smooth bounded domain of $\mathbb{R}^{n}$ $(n \geq 3)$  such that $0 \in \partial \Omega$.  Suppose $\gamma<\frac{n^2}{4}$ and let $h\in C^{1}(\overline\Omega)$ be such that the operator  $-\Delta-\gamma|x|^{-2}-h$ is coercive. Assuming that
$$\gamma>\frac{n^2-1}{4},$$
then there exists $\mathcal{H}\in C^2(\overline{\Omega}\setminus\{0\})$ such that
$$\qquad\qquad \left\{\begin{array}{ll}
-\Delta \mathcal{H}-\frac{\gamma}{|x|^2}\mathcal{H}+h(x) \mathcal{H}=0 &\hbox{ in } \Omega\\
\hfill \mathcal{H}>0&\hbox{ in } \Omega\\
\hfill \mathcal{H}=0&\hbox{ on }\partial\Omega \setminus\{0\}.
\end{array}\right.$$
Then, there exist constants $c_1,c_2\in\rr$ with $c_1>0$ such that
$$\mathcal{H}(x)=c_1\frac{d(x,\bdry)}{|x|^{\beta_+(\gamma)}}+c_2\frac{d(x,\bdry)}{|x|^{\beta_-(\gamma)}}+o\left(\frac{d(x,\bdry)}{|x|^{\beta_-(\gamma)}}\right)$$
as $x\to 0$.  In the spirit of Schoen-Yau \cite{SY}, we define the boundary mass as
$$m_{\gamma,h}(\Omega):=\frac{c_2}{c_1},$$
which is independent of the choice of $\mathcal{H}$.
\end{thdefi}
The problem of existence of least energy solutions can now be summarized in the following theorem, whose proof can also be deduced from the refined blow-up techniques developed in this paper.

\begin{theorem}[G.-R.\cite{gr1}, Lin-Wadade \cite{LW3}, G.-R. \cite{gr4}] \label{gro} Let $\Omega$ be a smooth bounded domain in $\rn$ $(n\geq 3)$ such that the singularity $0$ belongs to the boundary $\partial \Omega$. Suppose that $0< s <2$  
and fix $h_0\in C^1(\overline{\Omega})$ such that $-\Delta-\gamma|x|^{-2}-h_0$ is coercive. Assume one of the following two conditions:
\begin{itemize}
\item $\gamma\leq \frac{n^2-1}{4}$ and the mean curvature of $\partial\Omega$ at $0$ is negative.
\item $\frac{n^2-1}{4}<\gamma<\frac{n^2}{4}$ and the boundary mass $m_{\gamma,h_{0}}(\Omega)$ is positive.
\end{itemize}
Then, there is a positive solution to \eqref{one} that is a minimizer for the associated variational problem,
\begin{equation} 
\inf\left\{\frac{\int_{\Omega} |\nabla u|^2\, dx-\gamma \int_{\Omega}\frac{u^2}{|x|^2}\, dx -  \int_{\Omega} h_0(x) u^2\, dx}{\left(\int_{\Omega}\frac{|u|^{\crit}}{|x|^s}\, dx\right)^{\frac{2}{\crit}}};\, u \in\huno \setminus \{ 0\} \right\}.
\end{equation}
\end{theorem}
Our focus in this project, is to investigate the extent to which the above local curvature condition at $0$ and the global (mass) condition on the domain are necessary for the existence of positive solutions. Most importantly, we give results pertaining to the persistence of the lack of positive solutions for \eqref{one} under $C^1$-perturbations of the potential $h$. We will also show that, under suitable curvature conditions, this equation has an infinite number of non-necessarily positive solutions. 

Both existence and non-existence results will follow from a sharp blow-up analysis of solutions to perturbations of Equation \eqref{one}. More precisely, we consider
\begin{equation}\label{subexpos}
\hbox{$\pe\in [0,\crits-2)$ such that
$\lim_{\eps\to 0}\pe=0$, }
\end{equation}
and a family $(\he)_{\eps>0}\in  C^1(\overline{\Omega})$ such that 
\begin{align}\label{hyp:he}
\lim_{\eps\to
0}\he=h_0 \hbox{ in }C^1(\overline{\Omega})\hbox{ and }-\Delta - \frac{\gamma}{|x|^2}-h_{0}\hbox{ is coercive in }\Omega.
\end{align}
We then perform a blow-up analysis, as $\eps$ go to zero, on a  sequence  of functions $(\ue)_{\eps>0}$ in $\huno$ such that $\ue$ is a solution to the Dirichlet boundary value problems:
$$
\left\{ \begin{array}{llll}
-\Delta \ue-\gamma \frac{\ue}{|x|^2}-\he \ue&=& \frac{|\ue|^{\crits-2-\pe}\ue}{|x|^s}  \ \ &\text{in } \Omega ,\\
\hfill \ue&=&0 &\text{on }\partial \Omega.
\end{array}\right.\eqno{(E_\eps)}
$$
The novelty and delicacy of our analysis stem from the fact that the sequence $(\ue)_{\eps>0}$  might blow up along excited states, as opposed to a unique ground state in \cite{gr1}. Moreover, the sequence $(\ue)_{\eps>0}$ is not assumed to have a fixed sign.

\subsection{Non-existence: stability of the Pohozaev obstruction.} Starting with issues of non-existence of solutions, we shall prove the following surprising stability of regimes where variational positive solutions do not exist. 

\begin{theorem}\label{thm:non:ter} Let $\Omega$ be a smooth bounded domain in $\rn$ $(n\geq 3)$ such that the singularity $0$ belongs to the boundary $\partial \Omega$. Assume that $0< s <2$ and $\gamma<n^2/4$. Fix $h_0\in C^1(\overline{\Omega})$ such that $-\Delta-\gamma|x|^{-2}-h_0$ is coercive, and assume that one of the following conditions hold:
\begin{itemize}
\item  $\gamma\leq \frac{n^2-1}{4}$ and the mean curvature at $0$ is non-zero; 
\item $\gamma> \frac{n^2-1}{4}$ and the boundary mass $m_{\gamma,h_0}(\Omega)$ is non-zero. 
\end{itemize}
If there is no positive variational solution to \eqref{one} with $h=h_0$, then for all $\Lambda>0$, there exists $\epsilon:=\epsilon(\Lambda,h_0)>0$ such that for any $h\in C^1(\overline{\Omega})$ with
$\Vert h-h_0\Vert_{C^1(\Omega)}<\epsilon,$
there is no positive solution to \eqref{one} such that $\Vert \nabla u\Vert_2\leq\Lambda$. 
\end{theorem}
The above result is surprising for the following reason: Assuming $\Omega$ is starshaped with respect to $0$, then the classical Pohozaev obstruction (see Section \ref{sec:app:poho}) yields that \eqref{one} has no positive variational solution whenever
\begin{equation}\label{ineq:h:poho}
\hbox{ $h_0(x)+\frac{1}{2}( \nabla h_0(x), x) \leq 0$ for all $x\in \Omega$.}
 \end{equation}
 We then get the following corollaries.

 \begin{coro}\label{thm:non} Let $\Omega$ be a smooth bounded domain in $\rn$ $(n\geq 3)$ such that $0 \in \partial \Omega$. Assume $\Omega$ is starshaped with respect to $0$, $0< s <2$ and $\gamma<\gamma_H(\Omega)$.   If $\gamma\leq\frac{n^2-1}{4}$, we shall also assume that the mean curvature at $0$ is non-vanishing. If $h_0$ is a potential satisfying \eqref{ineq:h:poho}, then for all $\Lambda>0$, there exists $\epsilon(\Lambda,h_0)>0$ such that for all $h\in C^1(\overline{\Omega})$ satisfying
$\Vert h-h_0\Vert_{C^1(\Omega)}<\epsilon(\Lambda,h_0),$
there is no positive solution to \eqref{one} such that $\Vert \nabla u\Vert_2\leq\Lambda$. 
\end{coro}
\begin{coro}\label{thm:non:bis} Let $\Omega$ be a smooth bounded domain in $\rn$ $(n\geq 3)$, such that $0\in\partial \Omega$. We fix $0< s <2$ and $\gamma<\gamma_H(\Omega)$, the Hardy constant defined in \eqref{Hardy inequality}. Assume that 
$$\Omega\hbox{ is starshaped with respect to }0.$$
When $\gamma\leq\frac{n^2-1}{4}$, we assume that the mean curvature at $0$ is positive. Then for all $\Lambda>0$, there exists $\epsilon(\Lambda)>0$ such that for all $\lambda\in [0,\epsilon(\Lambda))$, there is no positive solution to 
\begin{equation}\label{one:bis}\left\{ \begin{array}{cl}
-\Delta u-\gamma \frac{u}{|x|^2}-\lambda u
= \frac{u^{\crits-1}}{|x|^s}  &\text{ in } \Omega,\\
u>0 &\hbox{ in }\Omega\\
u=0 &\hbox{ on }\bono
\end{array}\right.\end{equation}
with $\Vert \nabla u\Vert_2\leq\Lambda$. 
\end{coro}
It is worth comparing these results to what happens in the nonsingular case. Indeed, in contrast to the singular case, a celebrated result of Brezis-Nirenberg \cite{bn} shows that, for $\gamma=s=0$, a variational solution to \eqref{one:bis} always exists whenever $n\geq 4$ and $0<\lambda <\lambda_1(\Omega)$, with the geometry of the domain playing no role whatsoever. On the other hand, Druet-Laurain \cite{dl} showed that the geometry plays a role in dimension $n=3$, still for $\gamma=s=0$, by proving that when $\Omega$ is star-shaped, then there is no solution to \eqref{one:bis} for all small values of $\lambda>0$ (with no apriori bound on $\Vert \nabla u\Vert_2$). Another point of view is that for $n=3$, the nonexistence of solutions persists under small perturbations, but it does not for $n\geq 4$: the Pohozaev obstruction is stable only for $n=3$ in the nonsingular case.

\smallskip\noindent This is in stark contrast with the situation here, i.e. when $0\in\partial\Omega$ and $s>0$. In this case, for both the existence and non-existence results, the geometry plays a role in all dimensions: it is either the local geometry at $0$ (i.e., depending on whether the mean curvature at $0$ is 
vanishing or not) in high dimensions, or the global geometry of the domain (i.e., depending on whether the mass is positive or the domain is star-shaped) in low dimensions. Corollaries \ref{thm:non} and \ref{thm:non:bis} show that the Pohozaev obstruction is stable in all dimensions in the singular case.

\smallskip\noindent Let us discuss some extensions related to this absence or not of low/large dimension phenomenon.

\begin{itemize}
\item Our stability result still holds under an additional smooth perturbation of the domain $\Omega$, as was done by Druet-Hebey-Laurain \cite{dhl} when $n=3$, $\gamma=s=0$.
\item In the forthcoming paper \cite{GMR2}, we tackle the case of the interior singularity $0\in \Omega$, where the results are much more in the spirit of Brezis-Nirenberg and Druet-Laurain concerning the dichotomy between low and high dimensions. 
\item On of the main features of the stability result of Druet-Laurain \cite{dl} is the absence of any apriori control on $\Vert \nabla u\Vert_2$. In the interor case $0\in\Omega$, we expect to get rid also of the apriori bound in the singular case $s>0$. In the boundary case $0\in\partial\Omega$, bypassing the apriori bound by $\Lambda$ is more delicate and will require extra care. These issues are projects in progress.
\end{itemize}

\medskip\noindent The proof of Theorem \ref{thm:non:ter} (and Corollaries \ref{thm:non} and \ref{thm:non:bis}) relies on the blow-up analysis. Namely, arguing by contradiction, we assume the existence of solutions $(\ue)_\eps$ to \eqref{one:bis} with $p_\eps\equiv 0$ and  $(h_\eps)_\eps\to h_0$ in $C^1$ with a control on the Dirichlet energy. Due to the "invariance" under the conformal transformation \eqref{conf:trans}, the $\ue$'s might concentrate on some peaks at $0$. The formation of these peaks is described via blow-up analysis in Proposition \ref{prop:fund:est}. Then Proposition \ref{prop:rate:sc:2} applies which yields vanishing of the mean curvature or the mass, depending on the dimension, contradicting the hypothesis of Theorem \ref{thm:non:ter}. Concerning Corollaries \ref{thm:non} and \ref{thm:non:bis}, the hypothesis imply that the mass is negative when defined.

\subsection{Multiplicity of sign-changing solutions} As to the question of multiplicity, we shall prove the following result, which uses that in the subcritical case, i.e., when $p_\epsilon>0$, there is an infinite number of higher energy solutions for such $\epsilon$. Again, the core of the proof is a sharp blow-up analysis of such solutions as  $p_\epsilon\to 0$.

\begin{theorem}[{\textit{\rm The general case}}]\label{th:cpct:sc}
Let $\Omega$ be a smooth bounded domain in $\rn$, $n\geq 3$, such that $0\in  \bdry$ and  assume that  $0< s < 2$. Let  $h_{0} \in C^1(\overline{\Omega})$ and  $(\he)_{\eps>0}\in C^1(\overline{\Omega})$ be such that  \eqref{hyp:he} holds, and  
let $(\pe)_{\eps>0}$  be such that \eqref{subexpos} holds. 
Consider a sequence of functions $(\ue)_{\eps>0}$ that  is uniformly bounded in $\huno$ such that  for each $\eps >0$, $u_\epsilon$ satisfies 
Equation  $(E_\eps)$. Then,  \begin{enumerate}
\item If $\gamma<\frac{n^2}{4}-1$ and the principal curvatures of $\partial\Omega$ at $0$ are  non-positive but  not all of them vanish, then the sequence $(\ue)_{\eps>0}$ is pre-compact in $\huno$.
\item  In particular, Equation \eqref{one} has an infinite number of (possibly sign-changing) solutions in $\huno$. 
\end{enumerate}
\end{theorem}

The above result was established by Ghoussoub-Robert \cite{gr2} in the case when $\gamma=0$. The main challenge here is to prove the compactness of the subcritical solutions at high energy levels, as the  nonlinearities approach the critical exponent. The multiplicity result then follows from standard min-max methods. The proof 
relies heavily on pointwise blow-up analysis techniques in the spirit of  Druet-Hebey-Robert \cite{dhr} and Druet \cite{druet.jdg}, though our situation adds considerable difficulties to carrying out the program. 

\subsection{Compactness Theorems and blow-up analysis}
As mentioned above, the central tool is an analysis of the formation of peaks on families $(\ue)_\eps$ of solutions to equations like \eqref{one} when blow-up occurs. This long analysis yields Propositions \ref{prop:rate:sc} and \ref{prop:rate:sc:2} that describe the blow-up rate. When blowup does not occur, there is compactness. The following theorems are immediate consequences of these propositions. 

We note that the restrictions on both $\gamma$ and on the curvature at $0$ are more stringent than for the existence of a ground state solution in Theorem \ref{gro}. The stronger assumptions turned out to be due to the potentially sign-changing approximate solutions  -actually solutions of subcritical problems- and not because they are not necessarily minimizing. Indeed, the following theorem  does not assume any smallness of the energy bound as long as the approximate solutions are positive. It therefore yields another proof for Theorem \ref{gro}, which does not rely on the existence of minimizing sequence below the energy level of a single bubble.

\begin{theorem}[{\textit{\rm The non-changing sign case}}]\label{th:cpct:sc:3} Assume in addition  to the hypothesis of Theorem \ref{th:cpct:sc},  
   that the solutions $(\ue)_{\eps>0}$ satisfy for all $\eps>0$,
 \begin{equation}
 \hbox{$\ue >0$ \quad on $\Omega$.}
 \end{equation}
Then, the sequence $(\ue)_{\eps>0}$ is pre-compact in $\huno$, provided one of the following conditions is satisfied:
\begin{itemize}
\item $\gamma\leq \frac{n^2-1}{4}$ and the mean curvature of $\partial\Omega$ at $0$ is negative.
\item $\frac{n^2-1}{4}<\gamma<\frac{n^2}{4}$ and the boundary mass $m_{\gamma,h_{0}}(\Omega)$ is positive.
\end{itemize}
\end{theorem}
Our method also shows that if the --possibly sign-changing-- sequence is weakly null, then the compactness result in Theorem \ref{th:cpct:sc} will still hold for $\gamma$ up to $\frac{n^2}{4}-\frac{1}{4}$:

\begin{theorem}[{\textit{\rm The case of a weak null limit}}]\label{th:cpct:sc:2} Assume in addition  to the hypothesis of Theorem \ref{th:cpct:sc},  
that the solutions $(\ue)_{\eps>0}$ satisfy,
\begin{equation}\lim_{\eps\to 0}\Vert \ue\Vert_2=0.
\end{equation}
If $\gamma<\frac{n^2-1}{4}$ and the principal curvatures of $\partial\Omega$ at $0$ are  non-positive but  not all of them vanishing, then the sequence  $(\ue)_{\eps>0}$ converges strongly to $0$ in $\huno$. 
\end{theorem}

\subsection{Structure of the manuscript} This paper is organized as follows. Section \ref{sec:setup} consists in  preliminary material in order to introduce the sequence of functions that will be thoroughly analyzed in Sections \ref{sec:blowuplemma:1} to \ref{pf blow-up rates} in the case where they "blow-up". Section \ref{sec:proof:th}  contains the proof of the multiplicity result and Section \ref{sec:nonex} will have the applications to non-existence regimes and their stability under perturbations. We then have five relevant appendices. The first (Appendix A, Section \ref{sec:app:poho}) introduces the Pohozaev identity in our setting. The second (Appendix B, Section \ref{sec:app:lemma}) contains a technical lemma on the continuity of the first eigenvalue $\lambda_1(\Delta +V)$ with respect to variations of the potential $V$. Appendix C (Section \ref{sec:app:regul}) recalls regularity results established in \cite{gr4} about the regularity and behavior at $0$ of solutions of equations involving the Hardy-Schr\"odinger operator on bounded domains having $0$ on their boundary. In  Appendix D (Section \ref{sec:app:c}), we construct the Green functions associated to the operators $-\Delta -\frac{\gamma}{|x|^2} -h$ on such domains, and exhibit some of their properties needed throughout the memoir. The last Appendix E (Section \ref{sec:G:rnm}) does the same but for the Hardy-Schr\"odinger operator $-\Delta -\frac{\gamma}{|x|^2}$ on $\rnm$.

\section{\, Setting up  the blow-up }\label{sec:setup}
Throughout this paper, $\Omega$ will always  be a smooth bounded domain of $\rn$, $n\geq 3$, such that $0\in \bdry$.  We will always assume that  $\gamma <\frac{n^2}{4}$ and  $s\in (0,2)$. We set  $\crits:=\frac{2(n-s)}{n-2}$.  When $\gamma < \gamma_{H}(\Omega)$, then the 
following Hardy-Sobolev inequality holds on $\Omega$:  there exists $C>0$ such that
 \begin{equation}\label{HS-ineq}
C\left(\int_\Omega\frac{|u|^{\crits}}{|x|^s}\,dx\right)^{{2}/{\crits}}\leq \int_\Omega |\nabla u|^2\,dx-\gamma \int_\Omega\frac{u^2}{|x|^2}\,dx\hbox{ for all }u\in \huno.
\end{equation}
For each $\eps>0$, we consider $\pe\in [0,\crits-2)$ such that
\begin{align}\label{lim:pe}
\lim_{\eps\to 0}\pe=0.
\end{align}
Let $h_{0}\in C^1(\overline{\Omega})$ and consider a family $(\he)_{\eps>0}\in  C^1(\overline{\Omega})$ such that (\ref{hyp:he})  holds.  Consider a  sequence  of functions $(\ue)_{\eps>0}$ in $\huno$   such that for all $\eps >0$ the function $\ue$ is a solution to the Dirichlet boundary value problem:
$$
\left\{ \begin{array}{llll}
-\Delta \ue-\gamma \frac{\ue}{|x|^2}-\he \ue&=& \frac{|\ue|^{\crits-2-\pe}\ue}{|x|^s}  \ \ &\text{in } \huno ,\\
\hfill \ue&=&0 &\text{on }\partial \Omega.
\end{array}\right.\eqno{(E_\eps)}
$$
By the regularity result Theorem \ref{th:hopf} in Appendix B, we have $\ue \in C^{2}(\overline{\Omega}\setminus\{0\})$ and there exists $K_\eps\in\rr$ such that $ \lim_{x\to 0}~\frac{ |x|^{\bm} \ue (x)}{ d(x, \bdry) } =K_{\epsilon}$. In  addition, we assume that the sequence $(\ue)_{\eps>0}$ is bounded in $\huno$ and we let  $\Lambda>0$ be such that
\begin{align}\label{bnd:ue}
\int \limits_{\Omega}  \frac{|\ue|^{\crits-\pe}}{|x|^{s}} dx \leq \Lambda  \qquad \hbox{ for all } \eps>0.
\end{align}
It then follows from the weak compactness of the unit ball of $\huno$ that there exists $u_0\in\huno$ such that as $\eps\to 0$
\bequa\label{weak:lim:ue}
\ue\rightharpoonup u_0 \qquad \hbox{  weakly in } \huno. 
\eequa
Note that $u_0$ is a solution to the Dirichlet boundary value problem:
$$
\left\{ \begin{array}{llll}
-\Delta u-\gamma \frac{u}{|x|^2}-h_{0} u&=& \frac{|u|^{\crits-2-\pe}u}{|x|^s}  \ \ & \hbox{ in } \Ono ,\\
\hfill u&=&0 & \hbox{ on } \bono.
\end{array}\right.
$$
From the regularity Theorem \ref{th:hopf} we have  $u_{0} \in C^{2}(\overline{\Omega}\setminus\{0\})$ and $ \ds \lim_{x\to 0}\frac{|x|^{\bm} u_{0} (x)}{d\left( x, \bdry\right)}=K_{0} \in \R$.   It then follows that  
$\displaystyle \sup \limits_{\Omega} \frac{|x|^{\bm}u_{0}(x)}{d(x,\bdry)} $ and hence $\Vert |x|^{\bm-1} u_{0} (x) \Vert_{L^{\infty}(\Omega)}$ is finite.
\medskip

\noindent
We fix $\tau\in\rr$ such that
\begin{align}\label{tau}
\bm-1<\tau<\frac{n-2}{2}.
\end{align}
The following proposition  shows  that the sequence $(\ue)_{\eps}$ is  pre-compact in $\huno$  if $\left(|x|^{\tau}\ue \right)_{\eps>0}$ is uniformly bounded in $L^\infty(\Omega)$. 

\begin{proposition}\label{prop:bounded}
Let $\Omega$ be a smooth bounded domain of $\rn$, $n\geq 3$,  such that $0\in \bdry$ and  assume that  $0< s < 2$,   $\gamma<\frac{n^2}{4}$.  We let $(\ue)$, $(\he)$ and $(\pe)$ be such that $(E_\eps)$, \eqref{hyp:he} and \eqref{lim:pe} holds. Suppose  that there exists $C>0$ such that $|x|^{\tau}|\ue(x)|\leq C$ for all $x\in\Omega$ and  for all $\eps >0$. Then up to a subsequence, $\lim \limits_{\eps\to 0}\ue=u_0$ in $\huno$, where $u_0$ is as  in \eqref{weak:lim:ue}. 
\end{proposition}
\noindent{\it Proof of Proposition \ref{prop:bounded}:} The sequence $(\ue)$ is clearly uniformly bounded in $L^{\infty}(\Omega')$ for any $\Omega' \subset \subset \Omegabar\setminus\{0\}$. Then by standard elliptic estimates and from  \eqref{weak:lim:ue} it follows that $\ue \to u_{0}$ in $C^2_{loc}(\Omegabar\setminus\{0\})$. Now since  $|x|^{\tau}|\ue(x)|\leq C$ for all $x\in\Omega$ and for all $\eps >0$, and  since $\tau<\frac{n-2}{2}$, we have 
\begin{align}\label{lim:int:0}
\lim \limits_{\delta\to 0}\lim \limits_{\eps \to 0} \int \limits_{\Omega \cap B_{\delta}(0)}  \frac{|\ue|^{\crits-\pe}}{|x|^{s}} dx=0 ~\hbox{ and }  ~\lim \limits_{\delta\to 0}\lim \limits_{\eps \to 0}  \int \limits_{\Omega \cap  B_{\delta}(0)}  \frac{|\ue|^{2}}{|x|^{2}} dx=0.
\end{align}
Therefore    
\begin{align*}
\lim \limits_{\eps \to 0} \int \limits_{\Omega}  \frac{|\ue|^{\crits-\pe}}{|x|^{s}} dx =  \int \limits_{\Omega }  \frac{|u_{0}|^{\crits}}{|x|^{s}} dx ~\hbox{ and }  ~
\lim \limits_{\eps \to 0} \int \limits_{\Omega}  \frac{|\ue|^{2}}{|x|^{2}} dx =  \int \limits_{\Omega }  \frac{|u_{0}|^{2}}{|x|^{2}} dx.
\end{align*}
From $(E_{\eps})$  and \eqref{weak:lim:ue} we  then obtain 
\begin{align*}
\lim \limits_{\eps \to 0} \int \limits_{\Omega} \left(  |\nabla \ue|^2-\gamma  \frac{\ue^2}{|x|^2}- \he \ue^2 \right)dx &= \int \limits_{\Omega} \left(  |\nabla u_{0}|^2-\gamma  \frac{u_{0}^2}{|x|^2}- h_{0} u_{0}^2 \right)dx  \\
\hbox{  so then } \lim \limits_{\eps \to 0} \int \limits_{\Omega}   |\nabla \ue|^2 &= \lim \limits_{\eps \to 0} \int \limits_{\Omega}   |\nabla u_{0}|^2.
\end{align*}
And hence  $\lim \limits_{\eps\to 0}\ue=u_0$ in    $\huno$. This proves Proposition \ref{prop:bounded}. \qed

\medskip\noindent From now on, we assume that

\bequa\label{hyp:blowup}
\lim_{\eps\to 0}\Vert |x|^{\tau }\ue\Vert_{L^\infty(\Omega)}=+\infty.
\eequa
We shall say that blow-up occurs whenever  (\ref{hyp:blowup}) holds.

\section{\, Scaling Lemmas}\label{sec:blowuplemma:1}

\noindent
In this section we state and prove two scaling lemmas which we shall use many times in our analysis. We start by describing a  parametrization around  a point of the boundary $\partial \Omega$. Let $ p  \in \partial \Omega $. Then there exists  $U$,$V$ open in ${\R^{n}}$, there exists $I \subset \R$ an open interval, there exists $U' \subset \R^{n-1}$ an open subset, and there exist a smooth diffeomorphism $\mathcal{T}: U \longrightarrow  V$ and $\T_0\in C^\infty(U')$,  such that upto a rotation of coordinates if necessary 
\begin{align}\label{def:T:bdry}
\left\{\begin{array}{ll}
\bullet & 0 \in U=I\times U' \hbox{ and }p \in V.\\
\bullet & \mathcal{T}(0)=p.\\
\bullet &\mathcal{T} \left( U \cap \{x_{1} <0 \} \right)= V \cap   \Omega ~ \hbox{ and } \mathcal{T} \left( U \cap \{x_{1} =0 \} \right)= V \cap \partial  \Omega. \\
\bullet & D_{0} \mathcal{T} = \mathbb{I}_{\R^{n}}. \hbox{ Here $D_{x} \mathcal{T} $ denotes the differential of $ \mathcal{T} $ at the point $x$}\\
& \hbox{ and $ \mathbb{I}_{\R^{n}}$ is the identity map on $\R^{n}$.}\\
\bullet &{\mathcal{T}}_{*}(0)~ (e_{1})= \nu_{p}~\hbox{ where }\nu_{p}\hbox{ denotes the outer unit normal vector to }\\
& \partial \Omega\hbox{ at the  point }p.\\ 
\bullet & \{ { \mathcal{T}}_{*} (0)(e_{2}) , \cdots, {\mathcal{T}}_{*} (0)(e_{n}) \} \hbox{ forms an orthonormal basis of }\\
&T_{p} \partial \Omega .\\
\bullet & \T(x_1,y)=p+(x_1+\T_0(y),y)\hbox{ for all }(x_1,y)\in I\times U'=U \\
\bullet & \T_0(0)=0 \hbox{ and }\nabla\T_0(0)=0.
\end{array}\right.
\end{align}
This boundary  parametrization will be throughout  useful during our analysis. An important remark is that 
\begin{align}\label{rem:T:bdry}
\left( \T(x_{1},y), \bdry \right)=(1+o(1))|x_{1}| \qquad  \hbox{ for all } (x_1,y)\in I\times U'=U ~ \hbox{ close to } 0. 
 \end{align} 
\begin{lemma}{\label{scaling lemma 1}}
Let $\Omega$ be a smooth bounded domain of $\rn$, $n\geq 3$,  such that $0\in  \bdry $ and  assume that  $0< s < 2$,   $\gamma<\frac{n^2}{4}$. Let $(\ue)$, $(\he)$ and $(\pe)$ be such that $(E_\eps)$, \eqref{hyp:he},  \eqref{lim:pe} and \eqref{bnd:ue} holds. 
Let $(y_\eps)_\eps \in \Omega $ and  let
\begin{align*}
\nu_\eps^{-\frac{n-2}{2}}:= |u_{\eps}(y_{\eps})|, \quad \elle:=\nu_\eps^{1-\frac{\pe}{\crits-2}} ~ \hbox{ and } \quad \kappa_{\eps}:= \left| y_{\eps} \right|^{s/2} \elle^{\frac{2-s}{2}} \qquad \hbox{ for } \eps >0
\end{align*}
Suppose $\lim \limits_{\eps \to 0} y_{\eps}=0$ and $\lim \limits_{\eps \to 0} \nu_{\eps}=0$. Assume that for any $R > 0$ there exists  $C(R)>0$ such that  for all  $\eps >0$ 
\begin{align}\label{scale lem:hyp 1 on u}
|u_{\eps}(x)| \leq& C(R) ~\frac{|y_{\eps}|^{\tau }}{| x|^{\tau }} | u_{\eps}(y_{\eps})| \qquad \hbox{ for all } x\in B_{ R \kappa_{\eps}}( y_{\eps})  \cap \Omega.
\end{align}
Then
\begin{align*}
|\ye|=O(\elle) \qquad \hbox{ as } \eps \to 0.  
\end{align*}
\end{lemma}

\noindent{\it Proof of Lemma \ref{scaling lemma 1}:} We proceed by contradiction and assume that
\begin{align}
\label{scaling lem 1: contradiction}
\lim_{\eps\to 0}\frac{|\ye|}{\elle}=+\infty.
\end{align}
Then it follows from the definition of $\kappa_{\eps}$ that
\begin{align}
\label{ppty:ke}
\lim \limits_{\eps\to 0}\kappa_{\eps}=0,\; \lim \limits_{\eps\to 0}\frac{\kappa_{\eps}}{\elle}=+\infty \hbox{ and }\lim \limits_{\eps\to 0}\frac{\kappa_{\eps}}{|\ye|}=0.
\end{align}
\medskip

\noindent{\it Case 1:} We assume that there exists $\rho>0$ such for all $\eps>0$
that
$$\frac{d(\ye,\partial\Omega)}{\kappa_{\eps}}\geq 3\rho.$$

\noindent
We define for all $\eps>0$
$$\ve(x):=\nu_\epsilon^{\frac{n-2}{2}} \ue(\ye+\kappa_{\eps} x) \qquad \hbox{ for } x \in B_{2\rho}(0)$$
Note that this is well defined for $\eps >0$ small enough. It follows from \eqref{scale lem:hyp 1 on u} that there exists 
$C(\rho)>0$ such that all $\eps>0$
\begin{align}\label{bnd:lem:1:bis}
|\ve(x)|\leq C(\rho) \frac{1}{\left| \frac{\ye}{|y_{\epsilon}|}+\frac{\kappa_{\eps}}{|y_{\epsilon}|} x \right|^{\tau }} \qquad  \forall x \in B_{2\rho}(0)
\end{align}
using \eqref{ppty:ke} we  then get as $\eps \to 0$
\begin{align*}
|\ve(x)|\leq C(\rho)\left( 1+o(1)\right)   \qquad  \forall x \in B_{2\rho}(0).
\end{align*}
From  equation $(E_\eps)$ we obtain that  $v_{\eps}$ satisfies 
$$-\Delta\ve-\frac{\kappa_{\eps}^2}{|y_{\eps}|^{2}} \frac{\gamma}{\left|\frac{\ye}{|\ye|}+\frac{\kappa_{\eps}}{|\ye|}x\right|^{2}}~\ve-\kappa_{\eps}^2~\he(\ye+\kappa_{\eps} x)~\ve=\frac{|\ve|^{\crits-2-\pe}\ve}{\left|\frac{\ye}{|\ye|}+\frac{\kappa_{\eps}}{|\ye|}x\right|^s}$$
weakly in $B_{2\rho}(0)$ for all $\eps>0$.  With the help of \eqref{ppty:ke} and  standard elliptic theory it then  follows that there exists $v\in C^1(B_{2\rho}(0))$ such that
$$\lim \limits_{\eps\to 0} \ve = v \qquad \hbox{ in } C^1(B_{\rho}(0)).$$
In particular,
\begin{align}
\label{lim:v:case1:nonzero}
|v(0)|=\lim \limits_{\eps\to 0} |\ve(0)|=1
\end{align}
and therefore $v\not\equiv 0$. 
\medskip

\noindent
On the other hand, a change of variables and the definition of $\kappa_{\eps}$ yields
\begin{align*}
\int \limits_{ B_{\rho \kappa_{\eps}}(\ye)}\frac{|\ue|^{\crits-\pe}}{|x|^s}~dx &= \frac{|\ue(\ye)|^{\crits-\pe}\kappa_{\eps}^n}{|\ye|^s}\int \limits_{B_{\rho}(0)}\frac{|\ve|^{\crits-\pe}}{\left|\frac{\ye}{|\ye|}+\frac{\kappa_{\eps}}{|\ye|} x\right|^s}~ dx\\
&= \elle^{-\left(1+\frac{2(2-s)}{\crits-2-\pe} \right)}\left(\frac{|\ye|}{\elle}\right)^{s\left(\frac{n-2}{2}\right)}\int \limits_{B_{\rho}(0)}\frac{|\ve|^{\crits-\pe}}{\left|\frac{\ye}{|\ye|}+\frac{\kappa_{\eps}}{|\ye|}x\right|^s}~dx \\
&\geq  \left(\frac{|\ye|}{\elle}\right)^{s \left( \frac{n-2}{2} \right)}\int \limits_{B_{\rho}(0)}\frac{|\ve|^{\crits-\pe}}{\left|\frac{\ye}{|\ye|}+\frac{\kappa_{\eps}}{|\ye|}x\right|^s}~dx.
\end{align*}
Using the equation $(E_\eps)$, \eqref{bnd:ue}, \eqref{scaling lem 1: contradiction}, \eqref{ppty:ke} and  passing to the limit $\eps\to 0$ we get that
$$\int_{B_{\rho}(0)}|v|^{\crits}\, dx=0$$
and so then  $v\equiv 0$ in $B_\rho(0)$, a contradiction with \eqref{lim:v:case1:nonzero}. Thus \eqref{scaling lem 1: contradiction} cannot hold  in that case.  
\medskip

\noindent{\it Case 2:} We assume that, up to a subsequence,
\bequa\label{lim:d:be:0}
\lim_{\eps\to 0}\frac{d(\ye,\partial\Omega)}{\kappa_{\eps}}=0.
\eequa
Note that $\lim \limits_{\eps \to 0} y_{\eps}= 0 $. Consider the boundary map  $\T:U\to V$ as in \eqref{def:T:bdry}, where $U,V$ are both  open neighbourhoods of $0$. We let $\tue=\ue\circ\T$, which is defined in $U\cap \rnm$. For any $i,j=1,...,n$, we let $g_{ij}=(\partial_i\T,\partial_j\T)$, where $(\cdot,\cdot)$ denotes the Euclidean scalar product on $\rn$, and we consider $g$ as a metric on $\rn$. We let $\Delta_g=div_g(\nabla)$ the Laplace-Beltrami
operator with respect to the metric $g$. As easily checked, using $(E_{\eps})$ we get  that  for all $\eps >0$
$$-\Delta_g\tue-\frac{\gamma }{|\T(x)|^{2}} \tue-\he\circ\T(x)\cdot \tue=\frac{|\tue|^{\crit-2-\pe}\tue}{|\T(x)|^s}$$weakly in $U\cap \rnm$. We let $\ze\in\partial\Omega$ be such that
\begin{align}
\label{def:ze}
|\ze-\ye|=d(\ye,\partial\Omega).
\end{align}
We let $\tye,\tze\in U$ such that
\begin{align}
\label{def:tye:tze}
\T(\tye)=\ye\hbox{ and }\T(\tze)=\ze.
\end{align}
It follows from the properties \eqref{def:T:bdry} of the boundary map $\T$ that
\begin{align}
\label{ppty:tye:tze}
\lim_{\eps\to 0}\tye=\lim_{\eps\to 0}\tze=0,\; (\tye)_1<0\hbox{ and
}(\tze)_1=0
\end{align}
\medskip

\noindent
We rescale and define for all $\eps >0$
$$\tve(x):=\nu_{\eps}^{\frac{n-2}{2}}\tue(\tze+\kappa_{\eps} x) \qquad \hbox{ for } x\in \frac{U-\tze}{\kappa_{\eps}}\cap \rnm. $$ With \ref{ppty:tye:tze}), we get that $\tve$ is defined on $B_R(0)\cap\{x_1<0\}$ for all $R>0$, for $\eps$ is small enough. Then for all $\eps>0$ the functions $\tve$ satisfies the equation:
$$-\Delta_{\tge}\tve- \frac{\kappa_{\eps}^2}{|y_{\eps}|^{2}} \frac{\gamma}{\left|\frac{ \T (\tze+\kappa_{\eps} x)}{|\ye|}\right|^2}-\kappa_{\eps}^2\he\circ\T(\tze+\kappa_{\eps} x)\tve=\frac{|\tve|^{\crit-2-\pe}\tve}{\left|\frac{ \T (\tze+\kappa_{\eps} x)}{|\ye|}\right|^s}$$
weakly in $B_R(0)\cap\{x_1<0\}$. In this expression, $\tge=g(\tze+\kappa_{\eps} x)$ and $\Delta_{\tge}$ is the Laplace-Beltrami operator with respect to the metric $\tge$. With \eqref{lim:d:be:0}, \eqref{def:ze} and \eqref{def:tye:tze}, we get  for all $\eps>0$
$$\T(\tze+\kappa_{\eps} x)=\ye+O_R(1)\kappa_{\eps} \qquad \hbox{ for all } x\in B_{R}(0)\cap\{x_1\leq 0\}$$
where,  there exists $C_R>0$ such that $|O_R(1)|\leq C_R$ for all $x\in B_{R}(0)\cap\{x_1\leq0\}$. With (\ref{ppty:ke}), we then get that
$$\lim_{\eps\to 0}\frac{|\T(\tze+\kappa_{\eps} x)|}{|\ye|}=1 \qquad \hbox{ in } C^0(B_{R}(0)\cap\{x_1\leq 0\}).$$
It follows from \eqref{scale lem:hyp 1 on u} that there exists  $C'(R)>0$ such that all $\eps>0$
\begin{align}\label{bnd:lem:1:ter}
|\tve(x)|\leq C(R) \frac{1}{{\left|\frac{ \T (\tze+\kappa_{\eps} x)}{|\ye|}\right|^{\tau}}} \qquad  \forall x \in B_{R}(0)\cap\{x_1\leq0 \}.
\end{align}
Using \eqref{ppty:ke}  and the propoerties of the boundary map $\T$ we  then get as $\eps \to 0$
$$|\tve(x)|\leq C(R) \left( 1+o(1)\right) \qquad  \forall x \in B_{R}(0)\cap\{x_1\leq0 \}.$$
With the help of \eqref{ppty:ke} and  standard elliptic theory it then  follows that  there exists $\tv\in C^1(B_R(0)\cap\{x_1\leq 0\})$ such that
$$\lim_{\eps\to 0}\tve=\tv \qquad \hbox{ in } C^{0}(B_{R/2}(0)\cap \{x_1\leq 0\}).$$
Since $\tve$ vanishes on $B_R(0)\cap\{x_1=0\}$ and (\ref{bnd:lem:1:ter}) holds, it follows that 
\begin{align}\label{eq:ter:vanish}
\tv\equiv 0\hbox{ on }B_{R/2}(0)\cap \{x_1=0\}.
\end{align}
Moreover, from \eqref{lim:d:be:0}, \eqref{def:ze} and \eqref{def:tye:tze} we have that 
$$ \left| \tve\left(\frac{\tye-\tze}{\kappa_{\eps}}\right) \right|=1\hbox{ and }\lim_{\eps\to 0}\frac{\tye-\tze}{\kappa_{\eps}}=0.$$
In particular, $\tv(0)=1$,  contradiction with \eqref{eq:ter:vanish}. Thus \eqref{scaling lem 1: contradiction} cannot hold  in {\it Case 2} also.
\medskip

\noindent In both cases, we have contradicted \eqref{scaling lem 1: contradiction} . This proves that $\ye=O(\elle)$ when $\eps\to 0$,
which proves the Lemma. \qed

\begin{lemma}{\label{scaling lemma 2}}
Let $\Omega$ be a smooth bounded domain of $\rn$, $n\geq 3$,  such that $0\in  \bdry $ and  assume that  $0< s < 2$,   $\gamma<\frac{n^2}{4}$. Let $(\ue)$, $(\he)$ and $(\pe)$ such that $(E_\eps)$, \eqref{hyp:he},  \eqref{lim:pe} and \eqref{bnd:ue} holds. Let $(y_\eps)_\eps \in \Omega $ and  let
\begin{align*}
\nu_\eps^{-\frac{n-2}{2}}:= |u_{\eps}(y_{\eps})| \quad  \hbox{ and }  \quad \elle:=\nu_\eps^{1-\frac{\pe}{\crits-2}} \qquad \hbox{ for } \eps >0
\end{align*}
Suppose $\nu_{\epsilon} \to 0$ and $|\ye|=O(\elle) $ as  $\eps \to 0$. 
\medskip

\noindent
Since $0\in\partial\Omega$, we let $\T :U\to V$ as in \eqref{def:T:bdry} with $y_0=0$, where $U,V$ are open neighborhoods of $0$.
For   $\eps >0$ we rescale and define 
\begin{align*}
\twe(x):= \nu_{\eps}^{\frac{n-2}{2}} u_{\eps} \circ \T( \elle x) \qquad \text{ for } x \in \elle^{-1} U \cap \rnmpbar.
\end{align*}
Assume that for any $R >\delta> 0$ there exists  $C(R,\delta)>0$ such that  for all  $\eps>0$ 
\begin{align}\label{scale lem:hyp 2 on u}
|\twe(x)| \leq& C(R, \delta)  \qquad \text{ for all } x\in B_{R }(0) \setminus \overline{ B_{\delta}(0)} \cap \rnm.
\end{align}
\medskip

\noindent
Then there exists $\tw \in   H_{1,0}^2(\rnm) \cap C^1(\rnmpbar)$ such that  
\begin{align*}
 \twe \rightharpoonup  & ~\tw \qquad \hbox{ weakly in } H_{1,0}^2(\rnm) \quad \text{ as } \eps \rightarrow 0 \notag\\
 \twe \rightarrow  & ~\tw  \qquad \hbox{in } C^{1}_{loc}(\rnmpbar)  \qquad \text{ as } \eps \rightarrow 0 
\end{align*}
And  $\tw$  satisfies weakly the equation 
$$-\Delta \tw-\frac{\gamma}{|x|^2}\tw= \frac{|\tw|^{\crits-2} \tw}{|x|^s}\hbox{ in }\rnm.$$

\noindent
Moreover if $\tw \not\equiv 0$, then 
$$\int \limits_{\rnm} \frac{|\tw|^{\crits}}{|x|^{s}} \geq \mu_{\gamma, s}(\rnm)^{\frac{\crits}{\crits-2}} $$ 
and  there exists $t \in (0,1]$ such that $\lim \limits_{\eps\to0} \nu_{\eps}^{\pe}=t$,
where $\mu_{\gamma,s}(\rnm)$ is as in \eqref{general}.
\end{lemma}

\noindent{\it Proof of Lemma \ref{scaling lemma 2}:}  
The proof proceeds in four steps.\\

 \noindent {\bf Step \ref{scaling lemma 2}.1:}
Let $\eta \in C^{\infty}_{c}(\R^{n})$. One has  that  $\eta \twe \in H_{0,1}^2(\rnm)$ for $\eps >0 $ sufficiently small. We claim that there exists $\tw_{\eta} \in H_{1,0}^2(\rnm)$ such that upto a subsequence 
\begin{eqnarray*}
 \left \{ \begin{array} {lc}
          \eta \twe \rightharpoonup  \tw_{\eta} \qquad  \quad \text{ weakly in }  H_{1,0}^2(\rnm) ~  \text{ as } \eps \to 0, \\
          \eta \twe \rightarrow  \tw_{\eta}(x) \qquad  a.e ~  ~ \text{ in } \rnm ~  \text{ as } \eps \to 0.
            \end{array} \right. 
\end{eqnarray*} 
\medskip

\noindent
We prove the claim. Let $x \in \R_{-}^{n}$, then  
\begin{align*}
\nabla \left( \eta \twe \right)(x)=   \twe(x) \nabla \eta(x)  +  \nu_{\eps}^{\frac{n-2}{2}} \elle~ \eta(x) D_{(\elle x)} \T \left[\nabla \ue  \left( \T (\elle x) \right) \right]
\end{align*}
In this expression, $D_x\T$ is the differential of the function $\T$ at $x$.

\medskip\noindent Now for any $\theta >0$, there exists  $C({\theta}) >0$ such that for any $a,b >0$ 
\begin{align*}
(a+b)^{2} \leq C({\theta}) a^{2} + (1+ \theta) b^{2}
\end{align*}
With this inequality we then obtain 
\begin{eqnarray*}
\int \limits_{\R_{-}^{n}} \left|  \nabla \left( \eta  \twe \right)\right|^{2}~ dx& \leq& C(\theta)  \int \limits_{\R_{-}^{n}} | \nabla \eta |^{2} \twe^{2} ~ dx \\
&&+ (1 + \theta) \nu_{\eps}^{\frac{n-2}{2}} \elle \int \limits_{\R_{-}^{n}} \eta^{2} \left|  D_{(\elle x)} 
\T \left[\nabla \ue \left( \T (\elle x) \right) \right] \right|^{2} ~ dx
\end{eqnarray*}
Since $D_{0} \T = \mathbb{I}_{\R^{n}} $ we have as $\eps \to0$ 
\begin{align*}
\int \limits_{\R_{-}^{n}} \left|  \nabla \left( \eta  \twe \right)\right|^{2}~ dx & \leq C(\theta)  \int \limits_{\R_{-}^{n}} | \nabla \eta |^{2} \twe^{2} ~ dx  \\
& + (1 + \theta) \left( 1 + O(\elle )\right) \nu_{\eps}^{\frac{n-2}{2}} \elle \int \limits_{\R_{-}^{n}} \eta^{2} \left| \nabla \ue\left( \T (\elle x) \right)  \right|^{2} (1+ o(1)) ~dx 
\end{align*}
With H\"{o}lder inequality and a change of variables  this becomes 
\begin{align}{\label{b'dd on the integral 1}}
\int \limits_{\R_{-}^{n}} \left|  \nabla \left( \eta  \twe \right)\right|^{2}~ dx  &\leq~ C(\theta) \left\| \nabla \eta \right\|^{2}_{L^{n}}  \left( \frac{\nu_{\eps}}{\elle} \right)^{n-2}   \left( ~\int \limits_{\Omega} |\ue|^{\crit} ~ dx \right)^{\frac{n-2}{n}}  \notag \\ & + (1+ \theta)  \ \left( \frac{\nu_{\eps}}{\elle} \right)^{n-2}   \int \limits_{\Omega}  \left| \nabla u_{\epsilon} \right|^{2} ~ dx
\end{align}
Since $\left\| u_{\epsilon} \right\|_{\huno}= O(1)$,  so for $\eps >0$ small enough
\begin{align*} 
\left\| \eta \twe \right\|_{H_{1,0}^2(\rnm)} \leq C_{\eta}
\end{align*}
Where $C_{\eta}$ is a constant depending on the function $\eta$. The claim then follows from the reflexivity of $H_{1,0}^2(\rnm)$.
\medskip 

\noindent {\bf Step \ref{scaling lemma 2}.2:}
Let  ${\eta}_{1} \in C_{c}^{\infty}(\R^{n})$, $0 \leq  {\eta}_{1} \leq 1$ be a smooth cut-off function, such that
 \begin{eqnarray}
{\eta} _{1} = \left \{ \begin{array} {lc}
                 1 \quad \text{for } \ \   x \in B_1(0)\\
                 0 \quad \text{for } \ \   x \in {\R}^n \backslash B_2(0 )
          \end{array} \right.
\end{eqnarray}
For any $R>0$ we let $\eta_{R} = \eta_{1}(x/R)$. Then with a diagonal argument we can assume that upto a subsequence for any $R >0$ there exists  $\tw_{R} \in H_{1,0}^2(\rnm)$ such  that 
\begin{eqnarray*}
 \left \{ \begin{array} {lc}
          \eta_{R} \twe \rightharpoonup  \tw_{R} \qquad  \qquad \text{ weakly in }   H_{1,0}^2(\rnm) ~  \hbox{ as } \eps \to 0 \\
          \eta_{R} \twe(x) \rightarrow  \tw_{R}(x) \qquad a.e ~ x~\text{ in } \rnm ~  \text{ as } \eps \to 0 
            \end{array} \right. 
\end{eqnarray*}
Since $\displaystyle{\left\| \nabla  \eta_{R} \right\|^{2}_{n} =\left\| \nabla  \eta_{1} \right\|^{2}_{n} }$ for all $R >0$,   letting $\eps \to 0$ in  $(\ref{b'dd on the integral 1})$ we obtain   that 
\begin{align*}
\int_{\R^{n}_{-}} \left|\nabla w_{R} \right|^{2} dx \leq C \qquad \text{for all } R>0
\end{align*}
where $C$ is a constant independent of $R$. So there exists  $ \tw \in  H_{1,0}^2(\rnm)  $ such  that 
\begin{eqnarray*}
 \left \{ \begin{array} {lc}
           \tw_{R} \rightharpoonup  \tw \qquad \qquad \text{ weakly in }  D^{1,2} (\R^{n}) ~  \text{ as } R \to +\infty \\
           \tw_{R}(x) \rightarrow  \tw(x)\qquad  a.e ~ x~\text{ in } \rnm ~  \text{ as } R \to +\infty 
            \end{array} \right. 
\end{eqnarray*}
\medskip 

\noindent {\bf Step \ref{scaling lemma 2}.3:}
We claim that $ \tw \in  C^{1}( \rnmpbar )$ and  it satisfies weakly the equation 
$$
\left\{ \begin{array}{llll}
 -\Delta \tw-\frac{\gamma}{|x|^2} \tw&=& \frac{|\tw|^{\crits-2} \tw}{|x|^s}  \ \ & \hbox{ in } \R^{n}_{-}\\
\hfill \tw&=&0 & \hbox{ on } \partial \rnm \setminus \{ 0\}.
\end{array}\right.
$$
\noindent
We prove the claim.  For any $i,j=1,...,n$, we let $(\tge)_{ij}=(\partial_i\T(\elle x),\partial_j\T(\elle x))$,
where $(\cdot,\cdot)$ denotes the Euclidean scalar product on $\rn$. We
consider $\tge$ as a metric on $\rn$. We let $\Delta_g=div_g(\nabla)$ the Laplace-Beltrami
operator with respect to the metric $g$.  From   $(E_{\eps})$ it follows that   for any $\eps >0$ and $R>0$, $ \eta_{R} \twe $ satisfies  weakly the equation 
\begin{align}{\label{blowup eqn 1}}
-\Delta_{\tge} \left(\eta_{R} \twe \right)- \frac{\gamma}{\left| \frac{\T (\elle x)}{\elle}\right|^{2}}\eta_{R} \twe  -\elle^2~\he \circ \T(\elle x) \eta_{R} \twe =\frac{|\left(\eta_{R} \twe \right)|^{\crits-2-\pe} \left(\eta_{R} \twe \right)}{\left| \frac{\T (\elle x)}{\elle}\right|^{s}}.
\end{align}
and note that $ \eta_{R} \twe \equiv 0$ on $ B_{R}(0) \setminus \{ 0\} \cap \partial \rnm$.
 From  \eqref{def:T:bdry}, \eqref{scale lem:hyp 2 on u} and  using the  standard elliptic estimates  it follows that $\tw_{R} \in C^{1} \left(  B_{R}(0) \setminus \{ 0\} \cap \overline{\rnm} \right) $ and  that  up to a subsequence 
\begin{align*}
\lim \limits_{\epsilon \to 0} \eta_{R} \twe = \tw_{R}  \qquad \hbox{ in } C_{loc}^{1} \left(  B_{R/2}(0) \setminus \{ 0\} \cap \overline{\rnm} \right).
\end{align*}
Letting $\eps \to 0$  in eqn $\eqref{blowup eqn 1}$ gives that $ w_{R} $  satisfies  weakly the equation 
$$-\Delta \tw_{R} - \frac{\gamma}{|x|^{2}} \tw_{R} =\frac{|\tw_{R}|^{\crits-2-\pe} \tw_{R} }{|x|^s}.$$
Again  we have that $|\tw_{R}(x)| \leq C(R, \delta)$ for all $x \in \overline{B_{R/2}(0)} \setminus  \overline{B_{2 \delta}(0)}$ and  then  again from standard elliptic estimates it follows that   $ \tw \in C^{1}( \rnmpbar)$
and  $\lim \limits_{R\to + \infty} \tilde{w}_{R}= \tilde{w}$ in $C^{1}_{loc}(\rnmpbar )$, up to a subsequence.   Letting $R \to + \infty$  we obtain that $ \tw $  satisfies  weakly the equation 
$$
\left\{ \begin{array}{llll}
 -\Delta \tw-\frac{\gamma}{|x|^2} \tw&=& \frac{|\tw|^{\crits-2} \tw}{|x|^s}  \ \ & \hbox{ in } \R^{n}_{-}\\
\hfill \tw&=&0 & \hbox{ on } \partial \rnm \setminus \{ 0\}.
\end{array}\right.
 $$
This proves our claim.
\medskip

 \noindent {\bf Step \ref{scaling lemma 2}.4:}
Coming back to  equation $\eqref{b'dd on the integral 1}$ we have  for $R>0$
\begin{align}{\label{b'dd on the integral 1b}}
\int \limits_{\R_{-}^{n}} \left|  \nabla \left( \eta_{R}  \twe \right)\right|^{2}~ dx  &\leq~ C(\theta) \left( ~\int \limits_{  \{ x \in \rnm: R <|x|<2R  \}} (\eta_{2R}  \twe)^{2^{*}} ~ dx \right)^{\frac{n-2}{n}} \notag \\ & + (1+ \theta)  \ \left( \frac{\nu_{\eps}}{\elle} \right)^{n-2}   \int \limits_{\Omega}  \left| \nabla u_{\epsilon} \right|^{2} ~ dx.
\end{align}
Since  the sequence $(\ue)_{\eps}$ is bounded in $\huno$,  letting $\eps \to 0$    and then  $R \to + \infty$  we  obtain  for some constant $C$
\begin{align*}
\int \limits_{\rnm} \left|  \nabla  \tw \right|^{2}~ dx  \leq C\left(  \lim \limits_{\epsilon \to 0} \left( \frac{\nu_{\eps}}{\elle} \right)   \right)^{n-2}  .
\end{align*}

\noindent
Now if  $w \not\equiv 0$ weakly satisfies  the equation
$$
\left\{ \begin{array}{llll}
 -\Delta \tw-\frac{\gamma}{|x|^2} \tw&=& \frac{|\tw|^{\crits-2} \tw}{|x|^s}  \ \ & \hbox{ in } \R^{n}_{-}\\
\hfill \tw&=&0 & \hbox{ on } \partial \rnmp .
\end{array}\right.
$$
using the definition of  $\mu_{\gamma, s}(\rnm)$ it then  follows that 
$$\int \limits_{\rnm} \frac{|w|^{\crits}}{|x|^{s}} \geq \mu_{\gamma, s}(\rnm)^{\frac{\crits}{\crits-2}}. $$
Hence $\displaystyle \lim \limits_{\epsilon \to 0} \left( \frac{\nu_{\eps}}{\elle} \right)  >0$ which implies that 
\begin{align*}
t:= \lim\limits_{\epsilon \to 0} \nu_{\eps} ^{\pe}  >0.
\end{align*}
Since $\lim \limits_{\eps \to 0} \nu_{\eps}=0 $, therefore  we have that  $ 0 <t \leq 1.$ This completes the lemma. \qed

\section{\, Construction and exhaustion of the blow-up scales}\label{sec:exh}

\noindent
In this section we prove the following proposition in the spirit of Druet-Hebey-Robert \cite{dhr}:

\begin{proposition}\label{prop:exhaust} 
 Let $\Omega$ be a smooth bounded domain of $\rn$, $n\geq 3$,  such that $0\in  \bdry$ and  assume that  $0< s < 2$,   $\gamma<\frac{n^2}{4}$.  Let $(\ue)$, $(\he)$ and $(\pe)$ be such that $(E_\eps)$, \eqref{hyp:he},  \eqref{lim:pe} and \eqref{bnd:ue} holds. Assume that blow-up occurs, that is
$$\lim_{\eps\to 0}\Vert |x|^{\tau} \ue\Vert_{L^\infty(\Omega)}=+\infty ~\hbox{ where }  ~ \bm-1<\tau<\frac{n-2}{2} .$$ 

\noindent
Then there exists $N\in\nn^\star$ families of scales $(\mei)_{\eps>0}$ such that we have: \par

\begin{enumerate}
\item[{\bf(A1)}]
$\lim \limits_{\eps\to 0}\ue=u_0$ in $C^2_{loc}(\overline{\Omega}\setminus\{0\})$ where $u_0$ is as in\eqref{weak:lim:ue}. \\
\item[{\bf(A2)}]
$0< \meun< ...<\meN$, for all $\eps>0$. \\
\item[{\bf(A3)}]
$\lim\limits_{\eps\to 0}\meN=0\hbox{ and } \lim\limits_{\eps\to 0}\frac{\mu_{i+1,\eps}}{\mei}=+\infty\hbox{ for all } 1 \leq  i \leq N-1.$ \\
\item[{\bf(A4)}]
For any $1 \leq  i \leq N$ and for   $\eps >0$ we  rescale and define 
$$ \tuei(x):= \mei^{\frac{n-2}{2}}\ue( \T (\kei x)) \qquad \hbox{ for } x \in k_{i,\eps}^{-1} U \cap \rnmpbar,$$
where $\kei=\mei^{1-\frac{\pe}{\crit-2}}$. 
Then there exists $\tui \in H_{1,0}^2(\rnm) \cap C^1(\rnmpbar)$, $\tui \not\equiv 0$  such that $\tui $ weakly solves the equation 
$$\left\{ \begin{array}{llll}
-\Delta \tui-\frac{\gamma}{|x|^2}\tui &=& \frac{|\tui|^{\crits-2} \tui}{|x|^s}  &\hbox{ in } \rnm \\
\hfill \tui&=&0 & \hbox{ on } \partial \rnmp .
\end{array}\right.$$ 
and
\begin{align*}
\tuei \longrightarrow &~ \tui \qquad \hbox{in } C^{1}_{loc}(\rnmpbar)  \qquad \hbox{ as } \eps \to0, \notag\\
\tuei \rightharpoonup  &~ \tui \qquad \hbox{ weakly in } H_{1,0}^2(\rnm)  \quad \hbox{ as } \eps \to0.
\end{align*}
\item[{\bf(A5)}]
There exists
$C>0$ such that
$$|x|^{\frac{n-2}{2}}|\ue(x)|^{1-\frac{\pe}{\crits-2}}\leq C\quad \hbox{for all $\eps>0$ and all $x\in\Omega$.}$$

\item[{\bf(A6)}] $\lim \limits_{R\to +\infty}\lim \limits_{\eps\to 0}~\sup \limits_{ \Omega \setminus  B_{R k_{N,\eps}}(0)   } |x|^{\frac{n-2}{2}}|\ue(x)-u_0(x)|^{1-\frac{\pe}{\crits-2}}=0.$ \\
\item[{\bf(A7)}]
$\lim \limits_{\delta \to 0}  \lim \limits_{\eps\to 0}~\sup \limits_{  B_{\delta k_{1,\eps}}(0) \cap \Omega } |x|^{\frac{n-2}{2}}\left|\ue(x)-\mu_{1,\eps}^{-\frac{n-2}{2}} \tu_{1}\left( \frac{\T^{-1}(x)}{k_{1,\eps} }  \right)\right|^{1-\frac{\pe}{\crits-2}}=0.$\\
\item[{\bf(A8)}]
For any $\delta>0$ and any $1 \leq  i \leq N-1$, we have 
$$\lim \limits_{R\to +\infty}\lim \limits_{\eps\to 0}\sup \limits_{\delta k_{i+1, \eps}\geq |x|\geq R \kei}|x|^{\frac{n-2}{2}}\left|\ue(x)-\mu_{i+1,\eps}^{-\frac{n-2}{2}} \tu_{i+1}\left( \frac{\T^{-1}(x)}{k_{i+1,\eps} }  \right)\right|^{1-\frac{\pe}{\crits-2}}=0.$$\\
\item[{\bf(A9)}]
For any $i\in\{1,...,N\}$, there exists $t_i\in (0,1]$ such that
$\lim_{\eps\to 0}\mei^{\pe}=t_i.$
\end{enumerate}
\end{proposition}
\medskip

\noindent The proof of this proposition is inspired by \cite{dhr} and proceeds in five steps.

\medskip

\noindent
Since  $s>0$, the subcriticality $\crits<\crit$  of equations $(E_{\eps})$ along with \eqref{weak:lim:ue} yields that  $\ue \to u_{0}$ in $C^2_{loc}(\Omegabar\setminus\{0\})$. So the only blow-up point is the origin.

\medskip\noindent{\bf Step \ref{sec:exh}.1:}
 The construction of the $\mei$'s proceeds by induction. This step is the initiation.
\medskip

\noindent
By the regularity Theorem \ref{th:hopf} and the definition of $\tau$ in  \eqref{tau} it follows that for any $\eps >0$ there exists $\xeun \in \Omegabar\setminus\{0\}$ such that 
\begin{align}\label{def: xe}
\sup \limits_{ x \in \Ono} |x|^{\tau} |\ue (x)|= |\xeun|^{\tau} |\ue (\xeun)|. 
\end{align}
\noindent
We define  $\meun$ and $\keun>0$ as follows 
\begin{align}\label{def:me:ke}
\meun^{-\frac{n-2}{2}}:=|\ue(\xeun)| ~ \hbox{ and } ~\keun:=\meun^{1-\frac{\pe}{\crit-2}}.
\end{align}
Since blow-up occurs, that is  (\ref{hyp:blowup}) holds and since   $\ue \to u_{0}$ in $C^2_{loc}(\Omegabar\setminus\{0\})$, we have that 
$$\lim \limits_{\eps \to 0} \xeun=0 \in \bdry \qquad \hbox{ and  }\qquad  \lim \limits_{\eps \to 0} \meun =0.$$
It follows that $\ue$ satisfies the hypothesis \eqref{scale lem:hyp 1 on u} of  Lemma \ref{scaling lemma 1} with $\ye=\xeun$, $\nu_{\eps}=\meun$. Therefore 
$$|\xeun| = O\left(\keun\right)~ \hbox{ as } \eps \to 0. $$
\noindent
In fact, we claim that  there exists $c_{1}>0$ such that
\begin{align}\label{def:c1}
\lim \limits_{\eps \to 0} \frac{|\xeun|}{\keun}= c_{1}.
\end{align}
We argue by contradiction and we assume that $|\xeun|=o(\keun)$ as $\eps \to 0$. Let $\displaystyle \tilde{x}_{1,\eps}:= \T^{-1}(\xeun)  \in \rnm$. Since $|\xeun|=o(\keun)$ as $\eps \to 0$,  so also $| \tilde{x}_{1,\eps}|=o(k_{1,\eps})$ as $\eps \to 0$.
\medskip

\noindent
We define  for $\eps >0$
 $$\tilde{v}_{\eps}(x):=\meun^{\frac{n-2}{2}}\ue( \T ( |\tilde{x}_{1,\eps}| ~x)) \qquad \hbox{ for } x \in  \frac{U}{|\tilde{x}_{1,\eps}|} \cap \rnmpbar  $$
Using  $(E_{\eps})$ we obtain  that $\tilde{v}_{\eps}$ satisfies the equation 
$$-\Delta\tilde{v}_{\eps}- \frac{\gamma}{\left| \frac{\T \left(|\tilde{x}_{1,\eps}| x\right)}{|\tilde{x}_{1,\eps}|}\right|^2} \tilde{v}_{\eps}+|\xeun|^2\he \circ \T(|\tilde{x}_{1,\eps}|x) ~\tilde{v}_{\eps}=\left(\frac{|\tilde{x}_{1,\eps}|}{\keun}\right)^{2-s-\pe}\frac{|\tilde{v}_{\eps}|^{\crits-2-\pe}\tilde{v}_{\eps}}{\left| \frac{\T \left(|\tilde{x}_{1,\eps}| x\right)}{|\tilde{x}_{1,\eps}|}\right|^s}$$
The definition \eqref{def: xe} yields  as $\eps \to 0$, $\left| x \right|^\tau |\tilde{v}_{\eps}(x)|\leq 2$ for all $x\in \rnm. $\\
Standard elliptic theory then yields the existence of $\tilde{v}\in C^2(\rnmpbar)$ such that $\tv_{\eps} \to \tv$ in $C^2_{loc}(\rnmpbar)$ where
 $$\left\{ \begin{array}{llll}
-\Delta \tv -\frac{\gamma}{|x|^2}\tv &=&0  &\hbox{ in } \rnm \\
\hfill \tv&=&0 & \hbox{ on } \partial \rnmp .
\end{array}\right.$$ 
In addition, we have that  $ \left|\tilde{v}_{\eps} \left(|\tilde{x}_{1,\eps} |^{-1} \tilde{x}_{1,\eps} \right) \right|=1$ and   so $\tv \not \equiv 0$. Also  since $|x|^\tau |\tilde{v}(x)|\leq 2$  in $\rnmpbar$,  we have the bound that
\begin{align}\label{sing sol:grad control}
|x|^{\tau+1} |\tilde{v}(x)|\leq 2 |x_{1}|  \qquad \hbox{ for  all } ~x =(x_{1}, \tilde{x}) \hbox{ in } \rnm, 
\end{align}
which implies that 
$$|\tv(x)| < 4 \frac{|x_{1}|}{|x|^{\bp}}+ 4 \frac{|x_{1}|}{|x|^{\bm}}  \qquad \hbox{ for  all } ~x =(x_{1}, \tilde{x}) \hbox{ in } \rnm.$$
Therefore $x\mapsto \tilde{V}(x):=4 \frac{|x_{1}|}{|x|^{\bp}}+ 4 \frac{|x_{1}|}{|x|^{\bm}} -\tv(x)$ is a positive solution to  $-\Delta \tilde{V}-\frac{\gamma}{|x|^2} \tilde{V}=0$ in $\rnm$. Proposition \ref{prop:liouville} yields the existence of  $A,B \in \R$ such that 
 $$\tv(x) =A \frac{|x_{1}|}{|x|^{\bp}}+ B \frac{|x_{1}|}{|x|^{\bm}}  \qquad \hbox{ for  all } ~x  \hbox{ in } \rnm.$$
But the pointwise control \eqref{sing sol:grad control}  then implies  $A=B=0$ by letting $|x|\to 0$ and $\to \infty$. This contradicts $\tv \not \equiv 0$. This proves Claim \eqref{def:c1}.
\medskip

\noindent
We rescale and define for all $\eps >0$
$$\tu_{1,\eps}(x):=\meun^{\frac{n-2}{2}} \ue( \T (k_{1,\eps} ~x)) \qquad \hbox{ for } x \in k_{1,\epsilon}^{-1} U \cap \rnmpbar $$ 
It follows from  \eqref{def: xe} and  \eqref{def:c1}  that $\tueun$ satisfies the hypothesis \eqref{scale lem:hyp 2 on u} of  Lemma \ref{scaling lemma 2} with $\ye=x_{1,\eps}$, $\nu_{\eps}=\mu_{1,\eps}$. Then using Lemma  \ref{scaling lemma 2} we get that  there exists $\tu_{1} \in H_{1,0}^2(\rnm) \cap C^1(\rnmpbar)$ weakly satisfying the equation:
$$\left\{ \begin{array}{llll}
-\Delta \tu_{1}-\frac{\gamma}{|x|^2}\tu_{1} &=& \frac{|\tu_{1}|^{\crits-2} \tui}{|x|^s}  &\hbox{ in } \rnm \\
\hfill \tu_{1}&=&0 & \hbox{ on } \partial \rnmp .
\end{array}\right.$$ 
and 
\begin{align*}
\tu_{1,\eps} \longrightarrow &~ \tu_{1} \qquad \hbox{in } C^{1}_{loc}(\rnmpbar)  \qquad \hbox{ as } \eps \to0, \notag\\
\tu_{1,\eps} \rightharpoonup  &~ \tu_{1} \qquad \hbox{ weakly in } H_{1,0}^2(\rnm)  \quad \hbox{ as } \eps \to0.
\end{align*}\\
It follows from the definition that $\left| \tu_{i_{o},\eps}\left( \frac{\tilde{x}_{1,\eps}}{\keun}\right) \right|=1$. From  \eqref{def:c1} we therefore have that $\tu_{1} \not \equiv 0$.  And hence  again from  Lemma \ref{scaling lemma 2} we get that 
$$\int \limits_{\rnm} \frac{|\tu_{1} |^{\crits}}{|x|^{s}} \geq \mu_{\gamma, s}(\rnm)^{\frac{\crits}{\crits-2}}.$$
Moreover, there exists $t_{1}\in (0,1]$ such that $\lim \limits_{\eps\to 0}\meun^{\pe}=t_{1}$.
Since $\frac{|x|^{\bm}}{|x_{1}|}\tu_{1}\in C^0(\R^{n})$, we  get  as $\eps\to 0$
$$|\ye|^{\frac{n-2}{2}}\left|\mu_{1, \eps}^{-\frac{n-2}{2}}\tu_{1}\left(\frac{\T^{-1}(\ye)}{k_{1,\eps}}\right)\right|^{1-\frac{\pe}{\crits-2}}=O\left(|\tilde{y}_{\eps}|\right)^{\frac{n}{2}-\bm}=o(1),$$
and 
$$\lim_{\delta \to 0}  \lim_{\eps\to 0}~\sup_{  B_{\delta k_{1,\eps}}(0) \cap \Omega } |x|^{\frac{n-2}{2}}\left|\ue(x)-\mu_{1,\eps}^{-\frac{n-2}{2}} \tu_{1}\left( \frac{\T^{-1}(x)}{k_{1,\eps} }  \right)\right|^{1-\frac{\pe}{\crits-2}}=0.$$  \qed
\bigskip

\noindent{\bf Step \ref{sec:exh}.2:} We claim that there exists $C>0$ such that
\begin{align}\label{ineq:est:1}
|x|^{\frac{n-2}{2}}|\ue(x)|^{1-\frac{\pe}{\crits-2}} \leq C \quad \hbox{for all $\eps>0$ and all $x\in\Ono$.}
\end{align}
We argue by contradiction and let $(\ye)_{\eps>0} \in\Ono$ be such that
\begin{align}
\label{hyp:step32}
\sup_{x\in\Ono}|x|^{\frac{n-2}{2}}|\ue(x)|^{1-\frac{\pe}{\crits-2}}=|\ye|^{\frac{n-2}{2}}|\ue(\ye)|^{1-\frac{\pe}{\crits-2}}\to +\infty ~ \hbox{ as } \eps\to 0.
\end{align}
By the regularity Theorem \ref{th:hopf},  it follows that the sequence $(\ye)_{\eps>0}$ is well-defined and moreover $\lim \limits_{\eps \to 0} y_{\eps} =0$, since $\ue \to u_{0}$ in $C^2_{loc}(\Omegabar\setminus\{0\})$. 
For $\eps >0$ we let
$$\nu_{\eps}:=|\ue(\ye)|^{-\frac{2}{n-2}},~ \elle:=\nu_{\eps}^{1-\frac{\pe}{\crits-2}} \hbox{ and } \kappa_{\eps}:= \left| y_{\epsilon} \right|^{s/2} \elle^{\frac{2-s}{2}} .$$
Then it follows from (\ref{hyp:step32}) that
\begin{align}\label{lim:infty:34}
\lim_{\eps\to 0}\nu_{\eps}=0,  ~\lim_{\eps\to 0}\frac{|\ye|}{\elle}=+\infty  \hbox{ and } \lim \limits_{\eps\to 0}\frac{\kappa_{\eps}}{|\ye|}=0.
\end{align}

\noindent
Let $R>0$ and let $x\in B_R(0)$ be  such that $\ye+\kappa_{\eps} x\in\Ono$. It follows from the definition (\ref{hyp:step32}) of $\ye$ that for all $\eps>0$
$$|\ye+\kappa_{\eps} x|^{\frac{n-2}{2}}|\ue(\ye+\kappa_{\eps} x)|^{1-\frac{\pe}{\crits-2}}\leq |\ye|^{\frac{n-2}{2}}|\ue(\ye)|^{1-\frac{\pe}{\crits-2}}$$
and then, for all $\eps>0$ 
$$\left(\frac{|\ue(\ye+\kappa_{\eps} x)|}{|\ue(\ye)|}\right)^{1-\frac{\pe}{\crits-2}}\leq\left(\frac{1}{1-\frac{\kappa_{\eps}}{|\ye|}R}\right)^{\frac{n-2}{2}}$$
for all  $x\in B_R(0)$ such that $\ye+\kappa_{\eps} x\in\Ono$. Using  \eqref{lim:infty:34}, we get that there exists $C(R)>0$ such that the hypothesis \eqref{scale lem:hyp 1 on u} of  Lemma \ref{scaling lemma 1}  is satisfied and therefore  one has  $|\ye|=O(\elle)$ when $\eps\to 0$, contradiction to (\ref{lim:infty:34}). This proves (\ref{ineq:est:1}).
\qed
\bigskip

\noindent 
Let $\I \in\mathbb{N}^\star$. We consider the following assertions:
\begin{enumerate}
\item[{\bf(B1)}]
$0< \meun< ...<\mu_{\I,\eps}.$\\
\item[{\bf(B2)}]
 $\lim_{\eps\to 0}\mu_{\eps,\I}=0\hbox{ and }~\lim_{\eps\to 0}\frac{\mu_{i+1,\eps}}{\mei}=+\infty\hbox{ for all } 1 \leq i \leq \I-1$\\
\item[{\bf(B3)}]
 For all $ 1 \leq i \leq \I$, there exists $\tui \in H_{1,0}^2(\rnm) \cap C^2(\rnmpbar)$  such that $\tui $ weakly solves the equation 
$$\left\{ \begin{array}{llll}
-\Delta \tu_{i}-\frac{\gamma}{|x|^2}\tu_{i} &=& \frac{|\tu_{i}|^{\crits-2} \tui}{|x|^s}  &\hbox{ in } \rnm \\
\hfill \tu_{i}&=&0 & \hbox{ on } \partial \rnmp, 
\end{array}\right.$$ 
with 
$$\int \limits_{\rnm} \frac{|\tui|^{\crits}}{|x|^{s}} \geq \mu_{\gamma, s}(\rnm)^{\frac{\crits}{\crits-2}},$$ 
and
\begin{align*}
\tuei \longrightarrow &~ \tui \qquad \hbox{in } C^{1}_{loc}(\rnmpbar)  \qquad \hbox{ as } \eps \to0, \notag\\
\tuei \rightharpoonup  &~ \tui \qquad \hbox{ weakly in } H_{1,0}^2(\rnm)  \quad \hbox{ as } \eps \to0, 
\end{align*}
where for $\eps >0$, we have set $\kei=\mei^{1-\frac{\pe}{\crit-2}}$ and 
$$\tu_{i,\eps}(x):=\meun^{\frac{n-2}{2}} \ue( \T (k_{i,\eps} ~x)) \qquad \hbox{ for } x \in k_{i,\epsilon}^{-1} U \cap \rnmpbar.$$ 
\item[{\bf(B4)}]
For all $1 \leq i \leq \I $, there exists $t_i\in (0,1]$ such that
$\lim_{\eps\to 0}\mei^{\pe}=t_i.$
\end{enumerate}
\medskip

\noindent
We shall then say that $({\mathcal H}_{\I})$ holds if there exists $\I$ sequences  $(\mei)_{\eps>0}$, $i=1,...,\I$ such that items (B1), (B2) (B3) and (B4) holds. Note that it follows from Step \ref{sec:exh}.1 that $({\mathcal H}_{1})$ holds.
Next we  show the following:
\medskip

\noindent
{\bf Step \ref{sec:exh}.3} \label{prop:HN}
Let $I\geq 1$. We assume that $({\mathcal H}_{\I})$ holds.
Then, either
$$\lim_{R\to +\infty}\lim_{\eps\to 0}\sup_{ \Omega \setminus  B_{R k_{\I,\eps}}(0)} |x|^{\frac{n-2}{2}}|\ue(x)-u_0(x)|^{1-\frac{\pe}{\crits-2}}=0,$$
or ${\mathcal H}_{\I+1}$ holds.
\medskip

\noindent{\it Proof of Step \ref{sec:exh}.3:}
Suppose
$\lim\limits_{R\to +\infty}\lim\limits_{\eps\to 0}\sup_{\Omega \setminus  B_{R k_{\I,\eps}}(0)} |x|^{\frac{n-2}{2}}|\ue(x)-u_0(x)|^{1-\frac{\pe}{\crit-2}}\neq 0.$
Then, there exists a sequence of points $(\ye)_{\eps>0}\in\Ono$ such that

\bequa\label{hyp:lim:ye:HN}
\lim_{\eps\to 0}\frac{|\ye|}{k_{\I,\eps}}=+\infty\hbox{ and }\lim_{\eps\to 0}|\ye|^{\frac{n-2}{2}}|\ue(\ye)-u_0(\ye)|^{1-\frac{\pe}{\crits-2}}=a>0.
\eequa
Since $\ue \to u_{0}$ in $C^2_{loc}(\Omegabar\setminus\{0\})$ it follows that  $\lim \limits_{\eps \to 0} y_{\eps} =0$. Then by the regularity Theorem \ref{th:hopf} and since $\bm < \frac{n-2}{2}$, we get
\begin{align}\label{ref:a}
\lim_{\eps\to 
0}|\ye|^{\frac{n-2}{2}}|\ue(\ye)|^{1-\frac{\pe}{\crits-2}}=a>0
\end{align}
for some positive constant $a$.
In particular, $\lim \limits_{\eps\to 0}|\ue(\ye)|=+\infty$. Let
$$\mu_{\I+1, \eps}:=|\ue(\ye)|^{-\frac{2}{n-2}}\hbox{ and } ~ k_{\I+1,\eps}:=\mu_{\I+1,\eps}^{1-\frac{\pe}{\crit-2}}.$$
As a consequence we have 
\begin{align} \label{hyp:lim:ye:1:HN}
\lim \limits_{\eps\to 0}\mu_{\I +1,\eps}=0 \qquad \hbox{ and } \qquad  \lim \limits_{\eps\to 0}\frac{|\ye|}{k_{\I+1,\eps}}=a>0.
\end{align}
 
\medskip

\noindent
We rescale and define
$$\tu_{\I+1,\eps}(x):=\mu_{\I+1,\eps}^{\frac{n-2}{2}} \ue( \T (k_{I+1,\eps} ~x)) \qquad \hbox{ for } x \in k_{\I+1,\epsilon}^{-1} U \cap \rnmpbar. $$
It follows from \eqref{ineq:est:1} that for all $\eps >0$
$$\left| \frac{ \T (k_{I+1,\eps} ~x)}{k_{\I+1,\eps}}\right|^{\frac{n-2}{2}}|\tu_{\I+1,\eps}(x)|^{1-\frac{\pe}{\crit-2}}\leq C \qquad \hbox{ for } x \in k_{\I+1,\epsilon}^{-1}\Omega\setminus\{0\},$$
and so hypothesis \eqref{scale lem:hyp 2 on u} of Lemma \ref{scaling lemma 2} is satisfied.  Using Lemma \ref{scaling lemma 2}, we then get that  there exists $\tu_{\I+1}\in H_{1,0}^2(\rnm) \cap C^1(\rnmpbar)$ that satisfies weakly  the equation:
$$-\Delta \tu_{\I+1}-\frac{\gamma}{|x|^2}\tu_{\I+1}= \frac{|\tu_{\I+1}|^{\crits-2} \tu_{\I+1}}{|x|^s}\hbox{ in }\rnm.$$
while 
\begin{align*}
 \tu_{\I+1,\eps} \rightharpoonup   ~\tu_{\I+1} \quad \hbox{weakly in $H_{1,0}^2(\rnm)$ \quad 
 and \quad 
$ \tu_{\I+1,\eps} \rightarrow   ~\tu_{\I+1} \quad \hbox{in } C^{1}_{loc}(\rnmpbar)$,} 
\end{align*}
 as  $\eps \rightarrow 0.$
\medskip

\noindent
We denote $\displaystyle \tye:=\frac{ \T^{-1}(\ye) }{k_{\I+1,\eps}} \in \rnm$. From \eqref{hyp:lim:ye:1:HN}  it follows that that $\lim \limits_{\eps\to0} |\tye|:= |\tilde{y}_0| >a/2\neq 0$. 
Therefore 
$$|\tu_{\I+1}(\tilde{y}_0)|=\lim_{\eps\to 0}|\tu_{\I+1,\eps}(\tye)|=1.$$
Since $\tu_{\I+1}\equiv 0$ on $ \partial \rnmp$ so $\tye \notin \partial \rnm$ and hence $\tu_{\I+1}\not\equiv 0$. Hence  again from  Lemma \ref{scaling lemma 2}, we get
$$\int \limits_{\rnm} \frac{|\tu_{\I+1}|^{\crits}}{|x|^{s}} \geq \mu_{\gamma, s}(\rnm)^{\frac{\crits}{\crits-2}}$$
and there exists $t_{\I+1}\in (0,1]$ such that 
$\lim \limits_{\eps\to0}\mu_{\I+1,\eps}^{\pe}=t_{\I+1}$. Moreover, it follows from \eqref{hyp:lim:ye:HN} and  \eqref{hyp:lim:ye:1:HN}  that
$$\lim \limits_{\eps\to 0}\frac{\mu_{\I+1,\eps}}{\mu_{\I,\eps}}=+\infty \hbox{ and } \lim_{\eps\to 0}\mu_{\I+1,\eps}=0.$$
Hence the families $(\mei)_{\eps>0}$, $1\leq i \leq \I+1$ satisfy ${\mathcal H}_{\I+1}$.\qed
\bigskip

\noindent
The next step is  equivalent to step \ref{sec:exh}.3  at intermediate scales.
\medskip

\noindent
{\bf Step \ref{sec:exh}.4} \label{prop:HN+}
Let $I\geq 1$. We assume that $({\mathcal H}_{\I})$ holds. Then, for any $1\leq  i \leq \I-1$ and for any $\delta>0$,  either
$$\lim \limits_{R\to +\infty} \lim \limits_{\eps\to 0} \sup_{\Omega \cap B_{\delta k_{i+1,\eps}}(0)\setminus \overline{B}_{R k_{i,\eps}(0)}}|x|^{\frac{n-2}{2}}\left|\ue(x)-\mu_{i+1,\eps}^{-\frac{n-2}{2}}\tu_{i+1}\left(\frac{ \T^{-1}(x)}{k_{i+1,\eps}}\right)\right|^{1-\frac{\pe}{\crit-2}}=0$$
or $({\mathcal H}_{\I+1})$ holds.
\medskip

\noindent{\it Proof of Step \ref{sec:exh}.4:}
We assume that there exist an $i\leq \I-1$ and   $\delta>0$ such that $$\lim \limits_{R\to +\infty} \lim \limits_{\eps\to 0} \sup_{\Omega \cap B_{\delta k_{i+1,\eps}}(0)\setminus \overline{B}_{R k_{i,\eps}}(0)}|x|^{\frac{n-2}{2}}\left|\ue(x)-\mu_{i+1,\eps}^{-\frac{n-2}{2}}\tu_{i+1}\left(\frac{ \T^{-1}(x)}{k_{i+1,\eps}}\right)\right|^{1-\frac{\pe}{\crits-2}}>0.$$
It then follows that there exists a sequence  $(\ye)_{\eps>0}\in\Omega$ such that
\beqn
&&\lim_{\eps\to 0}\frac{|\ye|}{k_{i, \eps}}=+\infty,\qquad |\ye|\leq \delta
k_{i+1,\eps}\hbox{ for all }\eps>0\label{hyp:lim:ye:1:bis}\\
&&  |\ye|^{\frac{n-2}{2}}\left|\ue(\ye)-\mu_{i+1,\eps}^{-\frac{n-2}{2}}\tu_{i+1}\left(\frac{\T^{-1}(\ye)}{k_{i+1,\eps}}\right)\right|^{1-\frac{\pe}{\crits-2}}=a>0,\label{hyp:lim:ye:1:ter}
\eeqn
for some positive constant $a$.  Note that $ a< +\infty$ since $$|x|^{\frac{n-2}{2}}\left|\ue(x)-\mu_{i+1,\eps}^{-\frac{n-2}{2}}\tu_{i+1}\left(\frac{ \T^{-1}(x)}{k_{i+1,\eps}}\right)   \right|^{1-\frac{\pe}{\crits-2}}$$
is uniformly bounded for all $x \in \Omega \cap B_{\delta k_{i+1,\eps}}(0)\setminus \overline{B}_{R k_{i,\eps}}(0)$.
\medskip

\noindent
 Let  $\tye^{*} \in\rnm$ be such that $ \T^{-1}(\ye)= k_{i+1,\eps}~\tye^{*}$. It follows  that $|\tye^{*}|\leq \delta$ for all $\eps>0$. We rewrite (\ref{hyp:lim:ye:1:ter}) as
$$\lim_{\eps\to 0} |\tye^{*}|^{\frac{n-2}{2}}\left|\tu_{i+1, \eps}(\tye^{*})-\tu_{i+1}(\tye^{*})\right|^{1-\frac{\pe}{\crits-2}}=a>0.$$
Then  from point (B3) of ${\mathcal H}_{\I}$ it follows that $\tye^{*}\to 0$ as $\eps\to 0$. Since $\frac{|x|^{\bm}}{|x_{1}|} \tu_{i+1}\in C^0(\R^{n})$, we  get  as $\eps\to 0$
$$|\ye|^{\frac{n-2}{2}}\left|\mu_{i+1, \eps}^{-\frac{n-2}{2}}\tu_{i+1}\left(\frac{\ye}{k_{i+1,\eps}}\right)\right|^{1-\frac{\pe}{\crits-2}}=O\left(\frac{|\ye|}{k_{i+1,\eps}}\right)^{\frac{n}{2}-\bm}=o(1)$$
Then  (\ref{hyp:lim:ye:1:ter}) becomes 
\begin{align}\label{hyp:lim:ye:bis}
\lim \limits_{\eps\to0}|\ye|^{\frac{n-2}{2}}|\ue(\ye)|^{1-\frac{\pe}{\crits-2}}=a>0.
\end{align}
In particular, $\lim \limits_{\eps\to 0}|\ue(\ye)|=+\infty$. We let
$$\nu_{\eps}:=|\ue(\ye)|^{-\frac{2}{n-2}}\hbox{ and }\elle:=\nu_{\eps}^{1-\frac{\pe}{\crits-2}}.$$
Then we have  
\begin{align} \label{hyp:lim:ye:1:HN+}
\lim \limits_{\eps\to 0}\nu_{\eps}=0 \qquad \hbox{ and } \qquad  \lim \limits_{\eps\to 0}\frac{|\ye|}{\elle}=a>0.
\end{align}
\medskip

\noindent
We rescale and define
$$\tu_{\eps}(x):=\nu_{\eps}^{\frac{n-2}{2}}\ue(\T(\elle~x)) \qquad \hbox{ for } x \in \elle^{-1} U \cap \rnmpbar . $$
It follows from \eqref{ineq:est:1} that for all $\eps >0$
$$|x|^{\frac{n-2}{2}}|\tu_{\eps}(x)|^{1-\frac{\pe}{\crits-2}}\leq C \qquad \hbox{ for } x \in  \elle^{-1} U \cap \rnmpbar, $$
so that hypothesis \eqref{scale lem:hyp 2 on u} of Lemma \ref{scaling lemma 2} is satisfied. We can then use it to get that  there exists $\tu \in D^{1.2}(\rnm) \cap C^{1}(\rnmpbar)$ that satisfies weakly  the equation:
$$-\Delta \tu-\frac{\gamma}{|x|^2} \tu= \frac{|\tu|^{\crits-2} \tu}{|x|^s}\hbox{ in }\rnm,$$
while
\begin{align*}
 \tu_{ \eps} \rightharpoonup  & ~\tu \qquad \hbox{ weakly in } H_{1,0}^2(\rnm) \quad \text{ as } \eps \rightarrow 0 \notag\\
  \tu_{ \eps} \rightarrow  & ~\tu \qquad \hbox{in } C^{1}_{loc}(\rnmpbar)  \qquad \text{ as } \eps \rightarrow 0. 
\end{align*}
\noindent
We denote $\displaystyle \tye:=\frac{ \T^{-1}(\ye) }{\elle} \in \rnm$. From \eqref{hyp:lim:ye:bis}  it follows that that $\lim \limits_{\eps\to0} |\tye|:= |\tilde{y}_0| >a/2\neq 0$.
Therefore 
$$|\tu(\tilde{y}_0)|=\lim_{\eps\to 0}|\tu_{\eps}(\tye)|=1.$$
Since $\tu \equiv 0$ on $ \partial \rnmp$ so $\tye \notin \partial \rnm$ and hence $\tu\not\equiv 0$. Hence  again from  Lemma \ref{scaling lemma 2} we get
$$\int \limits_{\rnm} \frac{|\tu|^{\crits}}{|x|^{s}} \geq \mu_{\gamma, s}(\rnm)^{\frac{\crits}{\crits-2}},$$
and there exists $t\in (0,1]$ such that 
$\lim \limits_{\eps\to0}\nu_{\eps}^{\pe}=t$.  Moreover, from (\ref{hyp:lim:ye:bis}), (\ref{hyp:lim:ye:1:bis}), and since $\lim \limits_{\eps \to 0} \frac{|\ye|}{k_{i+1,\eps}}=0$, it follows that
$$\lim \limits_{\eps\to 0}\frac{\nu_{\eps}}{\mu_{i,\eps}}=+\infty~\hbox{ and }~\lim \limits_{\eps\to 0}\frac{\mu_{i+1, \eps}}{\nu_{\eps}}=+\infty .$$
Hence  the families $(\meun)$,..., $(\mei)$, $(\nu_{\eps})$, $(\mu_{i+1, \eps})$,..., $(\mu_{\I,\eps})$ satisfy $({\mathcal H}_{\I+1})$.\qed
\medskip

\noindent
The last step tells us that the process of constructing $\{ {\mathcal H}_{\I} \}$ stops after a  finite number of steps. \\

\noindent{\bf Step \ref{sec:exh}.5:} Let $N_0=\max\{\I : ({\mathcal H}_\I) \hbox{ holds } \}$. Then $N_0<+\infty$ and the conclusion of Proposition \ref{prop:exhaust} holds with $N=N_0$.
\medskip

\noindent{\it Proof of Step \ref{sec:exh}.5:}
Indeed, assume that $({\mathcal H}_{\I})$ holds. Since $\mei=o(\mu_{i+1, \eps})$ for all $1 \leq i \leq N-1$,  we get with a change of variable and the
definition of $\tuei$  that for any $R > \delta>0$
\begin{align*}
\int \limits_{\Omega}  \frac{|\ue|^{\crits -\pe}}{|x|^{s}} dx & \geq \sum_{i=1}^{\I} \int \limits_{ \T \left( B_{R \kei}(0)\setminus \overline{B}_{\delta\kei}(0) \cap \rnm \right)}  \frac{|\ue|^{\crits-\pe}}{|x|^{s}} dx \notag \\
&\geq  \sum_{i=1}^{\I}  \int \limits_{  B_{R \kei}(0)\setminus \overline{B}_{\delta\kei}(0) \cap \rnm } \frac{|\tu_{i,\eps}|^{\crits-\pe}}{|x|^{s}} dv_{g_{i,\eps}}. 
\end{align*}
Here  $g_{i,\eps}$ is the metric such that $(g_{\eps,i})_{qr}=(\partial_q\T(\kei x),\partial_r\T(\kei x))$ for all $q,r\in\{1,...,n\}$. Then  from \eqref{bnd:ue} we have 
\begin{align}
~\Lambda \geq  \sum_{i=1}^{\I}  \int \limits_{  B_{R \kei}(0)\setminus \overline{B}_{\delta\kei}(0) \cap \rnm } \frac{|\tu_{i,\eps}|^{\crits-\pe}}{|x|^{s}} dv_{g_{i,\eps}}. 
\end{align}
Passing to the limit $\eps\to 0$ and then $\delta \to0$, $R \to +\infty$ we obtain using point (B3) of ${\mathcal H}_{\I}$, that
$$\Lambda \geq \I \mu_{\gamma, s}(\rnm)^{\frac{\crits}{\crits-2}},$$
from which it follows that $N_0<+\infty$. \hfill $\Box$\\

\noindent To complete the proof, we let families $(\meun)_{\eps>0}$,..., $(\mu_{N_0, \eps})_{\eps>0}$ be such that ${\mathcal H}_{N_0}$ holds. We argue
by contradiction and assume that the conclusion of Proposition \ref{prop:exhaust} does not hold with $N=N_0$. Assertions (A1), (A2), (A3),(A4), (A5) , (A7) and (A9) hold. Assume that (A6) or (A8) does not hold. It then follows from Steps 4.3, 4.4 and 4.5 that ${\mathcal H}_{N+1}$ holds. A contradiction with the choice of $N=N_0$. Hence the proposition is proved. \qed

\section{\, Strong pointwise estimates}\label{sec:se:1}

\noindent
The objective of this section is to obtain pointwise controls on $\ue$ and $\nabla \ue $. The core is the proof of the following proposition in the spirit of Druet-Hebey-Robert \cite{dhr}:

\begin{proposition}\label{prop:fund:est}    
Let $\Omega$ be a smooth bounded domain of $\rn$, $n\geq 3$,  such that $0\in \bdry $ and  assume that  $0< s < 2$,   $\gamma<\frac{n^2}{4}$.  Let $(\ue)$, $(\he)$ and $(\pe)$ be such that $(E_\eps)$, \eqref{hyp:he},  \eqref{lim:pe} and \eqref{bnd:ue} holds. Assume that blow-up occurs, that is
$$\lim_{\eps\to 0}\Vert |x|^{\tau} \ue\Vert_{L^\infty(\Omega)}=+\infty ~\hbox{ where }  ~ \bm-1<\tau<\frac{n-2}{2} .$$ 
Consider $\mu_{1,\eps},...,\mu_{N,\eps}$  from Proposition  \ref{prop:exhaust}.  Then, there exists $C>0$ such that for all $\eps>0$
\begin{align}\label{eq:est:global}
|\ue(x)| \leq~  C& \left(~\sum_{i=1}^N\frac{\mei^{\frac{\bp-\bm}{2}} |x| }{ \mei^{\bp-\bm}|x|^{\bm}+|x|^{\bp}}+ \frac{  \Vert  |x|^{\bm-1}u_0 \Vert_{L^{\infty}(\Omega)} }{|x|^{\bm}} |x| \right) .
\end{align}
for all $x\in \Omega $.
\end{proposition}
\noindent The proof of this estimate,  inspired by the methodology of \cite{dhr},  proceeds in  seven steps. 

\medskip

\noindent{\bf Step \ref{sec:se:1}.1:} We claim that for any  $\alpha >0$  small and any $R>0$, there exists $C(\alpha, R)>0$ such that for all $\eps>0$ sufficiently small, we have  for all $x\in\Omega\setminus \overline{B}_{ Rk_{N,\eps}}(0)$, 
\begin{align}\label{estim:alpha:1}
|\ue(x)|\leq  C(\alpha,R)~\left( \frac{\mu_{N,\eps}^{\frac{\bp-\bm}{2}-\alpha } |x|}{|x|^{\bp-\alpha}}+ \frac{\Vert  |x|^{\bm-1}u_0 \Vert_{L^{\infty}(\Omega)}}{|x|^{\bm +\alpha}}|x| \right). 
\end{align}

\smallskip\noindent{\it Proof of  Step \ref{sec:se:1}.1:} We fix $\gamma'$ such that $\gamma<\gamma'<\frac{n^2}{4}$. Since the  operator $-\Delta - \frac{\gamma}{|x|^2}-h_{0}(x)$ is coercive, taking $\gamma'$ close to $\gamma$ it follows that the operator $-\Delta-\frac{\gamma'}{|x|^2} -h_{0}$ is  also coercive in  $\Omega$. From Theorem \ref{th:sing}, there exists   $H\in C^2(\overline{\Omega}\setminus\{0\})$ such that 
\begin{align} \label{eqn:sup-sol}
\left\{\begin{array}{ll}
-\Delta H-\frac{\gamma'}{|x|^2}H-h_{0}(x) H=0 &\hbox{ in } \Omega \\
\hfill H>0&\hbox{ in } \Omega\\
\hfill H=0&\hbox{ on }\partial\Omega \setminus \{ 0\}.
\end{array}\right.
\end{align}
And we have the following bound on $H$, that there exists  $C_{1}>0$ such that
\begin{align}\label{esti:sup-sol}
\frac{1}{C_{1}}\frac{d(x,\bdry)}{|x|^{\beta_{+}(\gamma') }} \leq H(x) \leq C_{1}\frac{d(x,\bdry)}{|x|^{\beta_{+}(\gamma') }} \qquad \hbox{ for all } x  \in \Omega.
\end{align}
\medskip

\noindent We let $\lambda^{\gamma'}_{1}>0$ be  the first eigenvalue of the coercive operator $-\Delta-\frac{\gamma'}{|x|^2} -h_{0}$  on $\Omega$ and we let
$\varphi\in C^2(\overline{\Omega}\setminus\{0\})\cap \huno$ be the unique eigenfunction such that
\begin{align}\label{eqn:sub-sol}
\left\{\begin{array}{ll}
-\Delta\varphi-\frac{\gamma'}{|x|^2}\varphi-h_{0}(x) \varphi~=\lambda^{\gamma'}_{1}\varphi ~  &\hbox{ in } \Omega\\
\hfill \varphi~ >0   &\hbox{ in }\Omega\\
\hfill \varphi~=0  &\hbox{ on }\partial\Omega\ \setminus \{ 0\}.\
\end{array}\right.
\end{align}
It follows from  the regularity result, Theorem \ref{th:hopf} that  there exists $C_2>0$ such that
\begin{align}\label{esti:sub-sol}
\frac{1}{C_{2}}\frac{d(x,\bdry)}{|x|^{\beta_{-}(\gamma') }} \leq \varphi(x) \leq C_{2}\frac{d(x,\bdry)}{|x|^{\beta_{-}(\gamma') }} \qquad \hbox{ for all } x \in \Omega.
\end{align}
\noindent
We define the operator
\begin{equation}\label{def:L}
\mathcal{L}_\eps:=-\Delta-\left(\frac{\gamma}{|x|^2}+h_{\eps} \right)-\frac{|\ue|^{\crits-2-p_{\eps}}}{|x|^s}.\end{equation}
\noindent
{\it Step \ref{sec:se:1}.1.1:} We claim that given any  $\gamma<\gamma'<\frac{n^2}{4}$  there exist $\delta_0>0$ and $R_0>0$ such that for any  $0<\delta < \delta_0$ and $R>R_0$, we have  for  $\eps>0$ sufficiently small
\begin{align}\label{ineq:LG}
\mathcal{L}_\eps H(x) >0\hbox{ and } \mathcal{L}_\eps \varphi(x) &~>0 \qquad \hbox{for all }x\in B_\delta(0)\setminus \overline{B}_{R\keN}(0) \cap \Omega. \notag\\
 \mathcal{L}_\eps H(x) &~>0 \qquad \hbox{for all }x\in \Omega\setminus \overline{B}_{R\keN}(0),  \hbox{ if } u_{0} \equiv0. 
\end{align}
\smallskip\noindent We prove the claim. As one checks for all $\eps >0$ and $ x \in \Omega$
\begin{align*}
\frac{\mathcal{L}_\eps H(x)}{ H(x)}= \frac{\gamma' -\gamma}{|x|^2}+(h_{0} -h_{\eps}) -\frac{|\ue|^{\crits-2-p_{\eps}}}{|x|^s}.
\end{align*}
and 
\begin{align*}
\frac{\mathcal{L}_\eps \varphi (x)}{ \varphi(x)}= \frac{\gamma' -\gamma}{|x|^2} +(h_{0} -h_{\eps}) -\frac{|\ue|^{\crits-2-p_{\eps}}}{|x|^s}+ \lambda^{\gamma'}_{1}.
\end{align*}
\noindent
One has for $\eps >0$ sufficiently small $\Vert h_{0}- h_{\eps}\Vert_\infty \leq  \frac{\gamma'-\gamma}{4(1+\sup \limits_{\Omega}|x|^{2} )}$ and we choose  $0 <\delta_0 <1$  such that 
\begin{align}\label{def:delta0}
 ~ \delta_{0}^{\left(\crits-2 \right)\left(\frac{n}{2}-\bm \right)} \Vert  |x|^{\bm-1}u_0 \vert|_{L^{\infty}(\Omega)}^{\crits-2} \leq\frac{\gamma'-\gamma}{ 2^{\crits+3}}.
\end{align}
This choice is possible thanks to (\ref{hyp:he}) and the regularity Theorem \ref{th:hopf} respectively. It follows from point (A6) of Proposition \ref{prop:exhaust} that, there exists $R_0>0$ such that for any $R>R_0$, we have for all $\eps>0$ sufficiently small
$$|x|^{\frac{n-2}{2}}|\ue(x)-u_0(x)|^{1-\frac{\pe}{\crits-2}}\leq \left(\frac{\gamma'-\gamma}{2^{\crits+2}}\right)^{\frac{1}{\crits-2}} \qquad \hbox{for all }  x\in \Omega\setminus \overline{B}_{R\keN}(0) $$
\noindent
With this choice of $\delta_0$  and $R_{0}$ we get that for any $0<\delta < \delta_0$ and $R>R_0$, we have  for $\eps>0$ small enough
\begin{align*}
|x|^{2-s}|\ue(x)|^{\crits-2-\pe}&\leq 2^{\crits-1-\pe}|x|^{2-s}|\ue(x)-u_0(x)|^{\crits-2-\pe} \\
&\quad +  2^{\crits-1-\pe}|x|^{2-s}|u_0(x)|^{\crits-2-\pe} \notag \\
&\leq  2^{-\pe}\frac{\gamma' -\gamma}{4} \leq   \frac{\gamma' -\gamma}{4} 
\end{align*}
for all  $x\in  B_\delta(0) \setminus \overline{B}_{R\keN}(0) \cap \Omega$,  if  $u_{0} \not\equiv0$, and  
\[
 ~|x|^{2-s}|\ue(x)|^{\crits-2-\pe}  \leq \frac{\gamma' -\gamma}{4}    \qquad  \hbox{for all }x\in \Omega\setminus \overline{B}_{R\keN}(0),  \hbox{ if } u_{0} \equiv0. 
\]
Hence  we obtain that for $\eps>0$ small enough
\begin{align}
\frac{\mathcal{L}_\eps H(x)}{ H(x)} = &~ \frac{\gamma' -\gamma}{|x|^2} +h_{0}-h_{\eps} -\frac{|\ue|^{\crits-2-p_{\eps}}}{|x|^s}\notag\\ 
\geq& ~\frac{\gamma' -\gamma}{|x|^2} +h_{0}-h_{\eps} - \frac{\gamma' -\gamma}{4 |x|^2} \notag\\
\geq& ~\frac{\gamma' -\gamma}{|x|^2} -\frac{\gamma' -\gamma}{4 |x|^2}  - \frac{\gamma' -\gamma}{4 |x|^2} = \frac{\gamma' -\gamma}{2 |x|^2} \notag\\ 
>&~0 \qquad \hbox{ for all } x\in   B_\delta (0)\setminus \overline{B}_{R\keN}(0) \cap \Omega  \hbox{ if } u_{0} \not\equiv0  \notag\\
\hbox{and }~~\frac{\mathcal{L}_\eps H(x)}{ H(x)}  >&~0  \qquad  \hbox{for all }x\in \Omega\setminus \overline{B}_{R\keN}(0),  \hbox{ if } u_{0} \equiv0.\notag
\end{align} 
Similarly we have 
\begin{align*}
\frac{\mathcal{L}_\eps \varphi (x)}{ \varphi(x)} >~0 \qquad \hbox{ for all } x\in    B_\delta (0)\setminus \overline{B}_{R\keN}(0) \cap \Omega.
\end{align*}
\qed
\medskip

\noindent{\it Step \ref{sec:se:1}.1.2:} It follows from point (A4) of Proposition \ref{prop:exhaust} that there exists $C'_1(R)>0$ such that for  all $\eps>0$  small
\begin{align*}
|\ue(x)|\leq C'_1(R)  \frac{\meN^{\frac{\beta_{+}(\gamma')-\beta_{-}(\gamma')}{2}} d(x, \bdry) }{|x|^{\beta_{+}(\gamma')}}\qquad \hbox{ for all } x\in \Omega \cap \partial B_{R\keN}(0).
\end{align*}
By estimate \eqref{esti:sup-sol} on $H$, we then have for some constant $C_1(R)>0$ 
\begin{align}\label{sup:ue:boundary:2}
|\ue(x)|\leq C_1(R) \meN^{\frac{\beta_{+}(\gamma')-\beta_{-}(\gamma')}{2}} H(x)\qquad \hbox{ for all } x\in \Omega \cap \partial B_{R\keN}(0).
\end{align}
\medskip

\noindent
It follows from point (A1) of Proposition \ref{prop:exhaust}  and the regularity Theorem \ref{th:hopf}, that there exists $C'_2(\delta)>0$ such that for all $\eps>0$  small
\begin{align}\label{sup:ue:boundary:1}
|\ue(x)|\leq  C'_2(\delta) \Vert { |x|^{\bm-1}u_0} \Vert_{L^{\infty}(\Omega)} \frac{ d(x, \bdry) }{|x|^{\beta_{-}(\gamma') }} \qquad \hbox{ for all } x\in \Omega \cap \partial B_{\delta}(0),  \hbox{ if } u_{0}\not \equiv0.
\end{align}
And then by the estimate \eqref{esti:sub-sol} on $\varphi$ we have for some constant $C_2(\delta)>0$ 
\begin{align}\label{sup:ue:boundary:1:bis}
|\ue(x)|\leq C_2(\delta) \Vert { |x|^{\bm-1}u_0} \Vert_{L^{\infty}(\Omega)}~ \varphi(x) \qquad \hbox{ for all } x\in \Omega \cap \partial B_{\delta}(0),  \hbox{ if } u_{0}\not \equiv0.
\end{align}
\medskip

\noindent
We now let for all $\eps>0$ 
$$\Psi_{\eps}(x):=  C_1(R)   \meN^{\frac{\beta_{+}(\gamma')-\beta_{-}(\gamma')}{2}} H(x) +C_2(\delta) \Vert { |x|^{\bm-1}u_0} \Vert_{L^{\infty}(\Omega)}~ \varphi(x)~~ \hbox{ for } x \in \Omega.$$
\noindent
Then  \eqref{sup:ue:boundary:1} and \eqref{sup:ue:boundary:2} imply that for all $\eps>0$  small
\begin{align}\label{control:ue:bord}
|u_{\eps}(x)| \leq \Psi_{\eps}(x) \qquad  \hbox{ for all } x \in \partial \left(B_\delta(0)\setminus \overline{B}_{R\keN}(0) \cap \Omega \right), \hbox{ if }  u_{0} \not \equiv 0
\end{align}
and 
\begin{align}\label{control:ue:bord:u0}
|u_{\eps}(x)| \leq \Psi_{\eps}(x) \qquad  \hbox{ for all } x \in \partial(\Omega\setminus \overline{B}_{R\keN}(0)), \hbox{ if } u_{0}\equiv0.
\end{align}
\noindent
Therefore when $ u_{0}\not \equiv0$ it follows from \eqref{ineq:LG} and \eqref{control:ue:bord} that for all  $\eps>0$ sufficiently small
$$\left\{\begin{array}{ll}
\mathcal{L}_\eps \Psi_\eps \geq 0=\mathcal{L}_\eps\ue & \hbox{ in } B_\delta(0)\setminus \overline{B}_{R\keN}(0) \cap \Omega\\
\Psi_\eps \geq \ue& \hbox{ on } \partial\left(B_\delta(0)\setminus \overline{B}_{R\keN}(0) \cap \Omega \right) \\
\mathcal{L}_\eps \Psi_\eps \geq 0=-\mathcal{L}_\eps\ue & \hbox{ in } B_\delta(0)\setminus \overline{B}_{R\keN}(0) \cap \Omega\\
\Psi_\eps \geq -\ue& \hbox{ on } \partial\left(B_\delta(0)\setminus \overline{B}_{R\keN}(0) \cap \Omega \right).
\end{array}\right.$$
and  from \eqref{ineq:LG} and \eqref{control:ue:bord:u0}, in case $u_{0} \equiv 0$,  we have for  $\eps>0$ sufficiently small
$$\left\{\begin{array}{ll}
\mathcal{L}_\eps \Psi_\eps \geq 0=\mathcal{L}_\eps \ue & \hbox{ in } \Omega \setminus \overline{B}_{R\keN}(0) \\
\Psi_\eps \geq \ue& \hbox{ on } \partial(\Omega\setminus \overline{B}_{R\keN}(0))\\
\mathcal{L}_\eps \Psi_\eps \geq 0=-\mathcal{L}_\eps\ue & \hbox{ in } \Omega \setminus \overline{B}_{R\keN}(0)\\
\Psi_\eps \geq -\ue& \hbox{ on } \partial(\Omega \setminus \overline{B}_{R\keN}(0)).
\end{array}\right.$$
Since $\Psi_\eps>0$  and $\mathcal{L}_\eps \Psi_\eps>0$, it follows from the comparison principle of Berestycki-Nirenberg-Varadhan \cite{bnv} that  the operator $\mathcal{L}_\eps$ satisfies the comparison principle on $ B_\delta(0)\setminus \overline{B}_{R\keN}(0) \cap \Omega$. Therefore
\begin{align*}
|\ue(x)|\leq&~ \Psi_\eps(x) \qquad  \hbox{ for all } x \in B_\delta(0)\setminus \overline{B}_{R\keN}(0) \cap \Omega, \notag\\
\hbox{ and }~ |\ue(x)|\leq &~ \Psi_\eps(x) \qquad  \hbox{ for all } x \in \Omega \setminus \overline{B}_{R\keN}(0) ~ \hbox{ if } u_{0} \equiv 0.
\end{align*}
Therefore   when $ u_{0}\not \equiv0$ we have for for all $\eps>0$  small
\begin{align*}
|u_{\eps}(x)| \leq C_1(R)   \meN^{\frac{\beta_{+}(\gamma')-\beta_{-}(\gamma')}{2}} H(x) +C_2(\delta) \Vert {|x|^{\bm-1}u_0} \Vert_{L^{\infty}(\Omega)}~ \varphi(x) 
\end{align*}
for all $ x \in B_\delta(0)\setminus \overline{B}_{R\keN}(0) \cap \Omega$, for $R$ large and $\delta$ small. 
\medskip

\noindent
Then, when $ u_{0}\not \equiv0$,  using the estimates \eqref{esti:sup-sol} and \eqref{esti:sub-sol}, we have or all $\eps>0$  small
\begin{align*}
|u_{\eps}(x)| &~\leq  C_1(R) \frac{ \meN^{\frac{\beta_{+}(\gamma')-\beta_{-}(\gamma')}{2}}  d\left( x, \bdry\right) }{|x|^{\beta_{+}(\gamma') }}+
C_2(\delta) \Vert {|x|^{\bm-1}u_0} \Vert_{L^{\infty}(\Omega)} \frac{d\left( x, \bdry\right)}{|x|^{\beta_{-}(\gamma') }}  \notag\\ 
 &~\leq  C_1(R) \frac{ \meN^{\frac{\beta_{+}(\gamma')-\beta_{-}(\gamma')}{2}}  |x| }{|x|^{\beta_{+}(\gamma') }}+C_2(\delta) \frac{\Vert {|x|^{\bm-1}u_0} \Vert_{L^{\infty}(\Omega)} }{|x|^{\beta_{-}(\gamma') }}|x|  .
\end{align*}
for all $ x \in B_\delta(0)\setminus \overline{B}_{R\keN}(0) \cap \Omega$, for $R$ large and $\delta$ small. 
\medskip

\noindent
Similarly if $u_{0} \equiv 0$, then all $\eps>0$  small and $R >0 $ large
$$|u_{\eps}(x)|\leq   C_1(R)  \frac{\meN^{\frac{\beta_{+}(\gamma')-\beta_{-}(\gamma')}{2}} |x|}{|x|^{\beta_{+}(\gamma')}} ~ \hbox{ for all } x \in \Omega \setminus \overline{B}_{R\keN}(0).$$
Taking $\gamma'$ close to $\gamma$, along with  points (A1) and (A4) of Proposition \ref{prop:exhaust}, it then follows that estimate \eqref{estim:alpha:1} holds on  $\Omega\setminus \overline{B}_{ Rk_{\eps,N}}(0)$ for all $R>0$.\qed 
\bigskip

\noindent{\bf Step \ref{sec:se:1}.2:} Let $1 \leq i \leq N-1$. We\ claim that for any $\alpha >0$ small  and any $R,\rho>0$, there exists  $C(\alpha,R,\rho)>0$ such that  all $\eps>0$.\par
\begin{align}\label{estim:alpha:1:bis}
|\ue(x)|\leq C(\alpha,R,\rho)\left(\frac{\mu_{i,\eps}^{\frac{\bp-\bm}{2}-\alpha}|x|}{|x|^{\bp-\alpha}}+\frac{|x|}{\mu_{i+1,\eps}^{\frac{\bp-\bm}{2}-\alpha} |x|^{\bm +\alpha}}  \right)   
\end{align}
for all $x\in B_{\rho k_{i+1,\eps}}(0)\setminus \overline{B}_{R k_{i,\eps}}(0) \cap \Omega$. 

\smallskip\noindent{\it Proof of  Step \ref{sec:se:1}.2:} We let $i\in \{1,...,N-1\}$. We emulate  the proof of Step \ref{sec:se:1}.1. Fix $\gamma'$ such that $\gamma<\gamma'<\frac{n^2}{4}$. Consider the functions $H$ and $\varphi$ defined in  Step \ref{sec:se:1}.1 satisfying \eqref{eqn:sup-sol} and \eqref{eqn:sup-sol} respectively.

\medskip

\noindent{\it Step \ref{sec:se:1}.2.1:} 
We claim that given any  $\gamma<\gamma'<\frac{n^2}{4}$  there exist  $\rho_0>0$ and $R_0>0$ such that for any  $0<\rho < \rho_0$ and $R>R_0$, we have  for  $\eps>0$ sufficiently small
\begin{align}
\label{ineq:LG2}
\mathcal{L}_\eps H(x) >0\hbox{ and } \mathcal{L}_\eps \varphi(x) >0 \qquad \hbox{for all }x\in B_{\rho k_{i+1,\eps}}(0)  \setminus \overline{B}_{R k_{i, \eps}}(0)  \cap \Omega 
\end{align}
where $\mathcal{L}\eps$ is as in \eqref{def:L}.

\smallskip\noindent We prove the claim. As one checks for all $\eps >0$ and $ x \in \Omega$
\begin{align*}
\frac{\mathcal{L}_\eps H(x)}{ H(x)}&~= \frac{\gamma' -\gamma}{|x|^2} + h_{0} -h_{\eps} -\frac{|\ue|^{\crits-2-p_{\eps}}}{|x|^s}, \notag \\
\frac{\mathcal{L}_\eps \varphi (x)}{ \varphi(x)}&~\geq \frac{\gamma' -\gamma}{|x|^2} + h_{0}-h_{\eps} -\frac{|\ue|^{\crits-2-p_{\eps}}}{|x|^s}.
\end{align*}

\noindent
We choose  $0<\rho_0<1$ such that
\begin{align}\label{def:rho0}
& \rho_0^2\sup \limits_{\Omega}| h_{0}-h_{\eps}| \leq  \frac{\gamma'-\gamma}{4} \quad \hbox{ for all } \eps>0 \hbox{ small and } ~ \notag \\
&\rho_{0}^{\left(\crits-2 \right)\left(\frac{n}{2}-\bm \right)} \Vert  |x|^{\bm-1} \tu_{i+1} \vert|_{L^\infty(B_2(0)\cap \rnm)}^{\crits-2} \leq\frac{\gamma'-\gamma}{ 2^{\crits+3}}
\end{align}
It follows from point (A8) of Proposition \ref{prop:exhaust} that there exists $R_0>0$ such that for any $R>R_0$ and any $0<\rho < \rho_0$, we have for all $\eps>0$ sufficiently small
$$ |x|^{\frac{n-2}{2}}\left|\ue(x)-\mu_{i+1,\eps}^{-\frac{n-2}{2}}\tu_{i+1}\left( \frac{\T^{-1}(x)}{k_{i+1,\eps} }  \right)\right|^{1-\frac{\pe}{\crits-2}} \leq  \left(\frac{\gamma'-\gamma}{2^{\crits+2}}\right)^{\frac{1}{\crits-2}} $$
for all $x\in B_{\rho k_{i+1, \eps}}(0)\setminus \overline{B}_{R\kei}(0) \cap \Omega$. 
\medskip

\noindent
With this choice of $\rho_0$ and $R_{0}$ we get that for any $0<\rho < \rho_0$ and $R>R_0$, we have  for $\eps>0$ small enough
\begin{align}
|x|^{2-s}|\ue(x)|^{\crits-2-\pe} \leq &~ 2^{\crits-1-\pe}|x|^{2-s}\left|\ue(x)-\mu_{i+1,\eps}^{-\frac{n-2}{2}}\tu_{i+1}\left(\frac{\T^{-1}(x)}{k_{i+1, \eps}} \right)\right|^{\crits-2-\pe} \notag\\
&+2^{\crits-1-\pe}\left(\frac{|x|}{k_{i+1,\eps}}\right)^{2-s} \left|\tu_{i+1} \left(\frac{ \T^{-1}(x)}{k_{i+1, \eps}}\right) \right|^{\crits-2-\pe} \notag\\
& \leq   \frac{\gamma' -\gamma}{4} \qquad \hbox{ for all } x\in B_{\rho k_{i+1, \eps}}(0)\setminus \overline{B}_{R\kei}(0).\notag
\end{align}
Hence  as in Step \ref{sec:se:1}.1 we have that for $\eps>0$ small enough
\begin{align*}
\frac{\mathcal{L}_\eps H(x)}{ H(x)} >0 ~ \hbox{ and } \frac{\mathcal{L}_\eps \varphi (x)}{ \varphi(x)}>0 ~ \hbox{ for all }x\in B_{\rho k_{i+1,\eps}}(0)  \setminus \overline{B}_{R k_{i, \eps}}(0)  \cap \Omega.  
\end{align*}
\medskip

\noindent{\it Step \ref{sec:se:1}.2.2:}  Let $i\in \{1,...,N-1\}$. It follows from point (A4) of Proposition \ref{prop:exhaust} that there exists $C'_1(R)>0$ such that for  all $\eps>0$  small
\begin{align*}
|\ue(x)|\leq C'_1(R)  \frac{\mei^{\frac{\beta_{+}(\gamma')-\beta_{-}(\gamma')}{2}} d(x, \bdry) }{|x|^{\beta_{+}(\gamma')}}\qquad \hbox{ for all } x\in \Omega \cap \partial B_{R\kei}(0),
\end{align*}
And then by the estimate \eqref{esti:sup-sol} on $H$ we have for some constant $C_1(R)>0$ 
\begin{align}\label{sup:ue:boundary:2:bis}
|\ue(x)|\leq C_1(R) \mei^{\frac{\beta_{+}(\gamma')-\beta_{-}(\gamma')}{2}} H(x)\qquad \hbox{ for all } x\in \Omega \cap \partial B_{R\kei}(0).
\end{align}
\medskip

\noindent
Again from point (A4) of Proposition \ref{prop:exhaust} it follows  that there exists $C'_2(\rho)>0$ such that for all $\eps>0$  small
\begin{align*}
|\ue(x)|\leq C'_2(\rho)  \frac{d(x,\bdry)}{\mu_{i+1,\eps}^{\frac{\beta_{+}(\gamma')-\beta_{-}(\gamma')}{2}} |x|^{\beta_{-}(\gamma')}} \qquad \hbox{ for all } x\in \Omega\cap \partial B_{\rho k_{i+1, \eps}}(0),
\end{align*}
and then by the estimate \eqref{esti:sub-sol} on $\varphi$ we have for some constant $C_2(\delta)>0$ 
\begin{align}\label{sup:ue:boundary:3:bis}
|\ue(x)|\leq C_2(\rho)  \frac{\varphi(x) }{\mu_{i+1,\eps}^{\frac{\beta_{+}(\gamma')-\beta_{-}(\gamma')}{2}} }  \qquad \hbox{ for all } x\in \Omega\cap \partial B_{\rho k_{i+1, \eps}}(0).
\end{align}
\medskip

\noindent
We let for all $\eps>0$ 
$$\tilde{\Psi}_{\eps}(x):=  C_1(R)   \mei^{\frac{\beta_{+}(\gamma')-\beta_{-}(\gamma')}{2}} H(x) +C_2(\rho)  \frac{\varphi(x) }{\mu_{i+1,\eps}^{\frac{\beta_{+}(\gamma')-\beta_{-}(\gamma')}{2}} }  \qquad \hbox { for } x \in \Omega. $$
Then  \eqref{sup:ue:boundary:2:bis} and \eqref{sup:ue:boundary:3:bis} implies that for all $\eps>0$  small
\begin{align}\label{control:ue:bord:bis}
|u_{\eps}(x)| \leq \tilde{\Psi}_{\eps}(x) \qquad  \hbox{ for all } x \in \partial \left(B_{\rho k_{i+1, \eps}}(0)\setminus \overline{B}_{R\kei}(0 ) \cap \Omega \right).
\end{align}

\noindent
Therefore it follows from \eqref{ineq:LG2} and \eqref{control:ue:bord:bis} that $\eps>0$ sufficiently small
$$\left\{\begin{array}{ll}
\mathcal{L}_\eps \tilde{\Psi}_\eps \geq 0=\mathcal{L}_\eps\ue & \hbox{ in } B_{\rho k_{i+1, \eps}}(0)\setminus \overline{B}_{R\kei}(0) \cap \Omega \\
\tilde{\Psi}_\eps \geq \ue& \hbox{ on } \partial \left(B_{\rho k_{i+1, \eps}}(0)\setminus \overline{B}_{R\kei}(0) \cap \Omega \right) \\
\mathcal{L}_\eps \tilde{\Psi}_\eps \geq 0=-\mathcal{L}_\eps\ue & \hbox{ in } B_{\rho k_{i+1, \eps}}(0)\setminus \overline{B}_{R\kei}(0) \cap \Omega \\
\tilde{\Psi}_\eps \geq -\ue& \hbox{ on } \partial \left(B_{\rho k_{i+1, \eps}}(0)\setminus \overline{B}_{R\kei}(0 ) \cap \Omega\right).
\end{array}\right.$$
Since $\tilde{\Psi}_\eps>0$  and $\mathcal{L}_\eps \tilde{\Psi}_\eps>0$, it follows from the comparison principle of Berestycki-Nirenberg-Varadhan \cite{bnv} that  the operator $\mathcal{L}_\eps$ satisfies the comparison principle on $ B_{\rho k_{i+1, \eps}}(0)\setminus \overline{B}_{R\kei}(0) \cap \Omega$. Therefore
$$|\ue(x)|\leq \tilde{\Psi}_\eps(x) \qquad  \hbox{ for all } x \in B_{\rho k_{i+1, \eps}}(0)\setminus \overline{B}_{R\kei}(0) \cap \Omega.$$
So for all $\eps>0$  small
$$|u_{\eps}(x)|\leq    C_1(R)   \mei^{\frac{\beta_{+}(\gamma')-\beta_{-}(\gamma')}{2}} H(x) +C_2(\rho)  \frac{\varphi(x) }{\mu_{i+1,\eps}^{\frac{\beta_{+}(\gamma')-\beta_{-}(\gamma')}{2}} } $$
 for all $ x \in B_{\rho k_{i+1, \eps}}(0)\setminus \overline{B}_{R\kei}(0) \cap \Omega$, for $R$ large and $\rho$ small.  Then using the estimates \eqref{esti:sup-sol} and \eqref{esti:sub-sol} we have or all $\eps>0$  small
\begin{align*}
|u_{\eps}(x)| &~\leq  C_1(R) \frac{ \mei^{\frac{\beta_{+}(\gamma')-\beta_{-}(\gamma')}{2}}  d\left( x, \bdry\right) }{|x|^{\beta_{+}(\gamma') }}+
C_2(\rho)  \frac{d(x,\bdry)}{\mu_{i+1,\eps}^{\frac{\beta_{+}(\gamma')-\beta_{-}(\gamma')}{2}} |x|^{\beta_{-}(\gamma')}}   \notag\\ 
 &~\leq  C_1(R) \frac{ \mei^{\frac{\beta_{+}(\gamma')-\beta_{-}(\gamma')}{2}}  |x| }{|x|^{\beta_{+}(\gamma') }}+C_2(\rho)  \frac{|x|}{\mu_{i+1,\eps}^{\frac{\beta_{+}(\gamma')-\beta_{-}(\gamma')}{2}} |x|^{\beta_{-}(\gamma')}}  .
\end{align*}
 for all $ x \in B_{\rho k_{i+1, \eps}}(0)\setminus \overline{B}_{R\kei}(0) \cap \Omega$, for $R$ large and $\rho$ small.

Taking $\gamma'$ close to $\gamma$, along with  point (A4) of Proposition \ref{prop:exhaust} it then follows that estimate \eqref{estim:alpha:1:bis} holds on  $B_{\rho k_{i+1, \eps}}(0)\setminus \overline{B}_{R\kei}(0) \cap \Omega$, for all $R \rho>0$.\qed 
\bigskip

\noindent{\bf Step \ref{sec:se:1}.3:}  We\ claim that for any $\alpha >0$ small  and any $\rho>0$, there exists  $C(\alpha,\rho)>0$ such that  all $\eps>0$.\par
\begin{align}\label{estim:alpha:1:tis}
|\ue(x)|\leq C(\alpha,\rho)\frac{|x|}{\mu_{1,\eps}^{\frac{\bp-\bm}{2}-\alpha} |x|^{\bm +\alpha}}   \qquad \hbox{ for all } x\in B_{\rho k_{1,\eps}}(0) \cap \Omega.
\end{align}

\smallskip\noindent{\it Proof of  Step \ref{sec:se:1}.3:}  Fix $\gamma'$ such that $\gamma<\gamma'<\frac{n^2}{4}$. Consider the function $\varphi$ defined in  Step \ref{sec:se:1}.1 satisfying  \eqref{eqn:sup-sol}. 

\smallskip\noindent{\it Step \ref{sec:se:1}.3.1:} 
We claim that given any  $\gamma<\gamma'<\frac{n^2}{4}$  there exist  $\rho_0>0$  such that for any  $0<\rho < \rho_0$ we have  for  $\eps>0$ sufficiently small
\begin{align}
\label{ineq:LG3}
\mathcal{L}_\eps \varphi(x) >0 \qquad \hbox{for all }x\in B_{\rho k_{1,\eps}}(0) \cap \Omega, 
\end{align}
where $\mathcal{L}\eps$ is as in \eqref{def:L}.

\smallskip\noindent Indeed, for all $\eps >0$ and $ x \in \Omega$
\begin{align*}
\frac{\mathcal{L}_\eps \varphi (x)}{ \varphi(x)}&~\geq \frac{\gamma' -\gamma}{|x|^2} -h_{\eps} -\frac{|\ue|^{\crits-2-p_{\eps}}}{|x|^s}.
\end{align*}
We choose  $0<\rho_0<1$ such that
\begin{align}
& \rho_0^2\sup \limits_{\Omega}| h_{\eps}| \leq  \frac{\gamma'-\gamma}{4} \quad \hbox{ for all } \eps>0 \hbox{ small and } ~ \notag \\
&\rho_{0}^{\left(\crits-2 \right)\left(\frac{n}{2}-\bm \right)} \Vert  |x|^{\bm-1} \tu_{1} \vert|_{L^\infty(B_2(0)\cap \rnm)}^{\crits-2} \leq\frac{\gamma'-\gamma}{ 2^{\crits+3}}\notag
\end{align}
It follows from point (A7) of Proposition \ref{prop:exhaust} that for any $0<\rho < \rho_0$, we have for all $\eps>0$ sufficiently small
$$ |x|^{\frac{n-2}{2}}\left|\ue(x)-\mu_{1,\eps}^{-\frac{n-2}{2}}\tu_{1}\left( \frac{\T^{-1}(x)}{k_{1,\eps} }  \right)\right|^{1-\frac{\pe}{\crits-2}} \leq  \left(\frac{\gamma'-\gamma}{2^{\crits+2}}\right)^{\frac{1}{\crits-2}} $$
for all $x\in B_{\rho k_{1,\eps}}(0) \cap \Omega$. 
\medskip

\noindent
With this choice of $\rho_0$  we get that for any $0<\rho < \rho_0$  we have  for $\eps>0$ small enough
\begin{align}
|x|^{2-s}|\ue(x)|^{\crits-2-\pe} \leq &~ 2^{\crits-1-\pe}|x|^{2-s}\left|\ue(x)-\mu_{1,\eps}^{-\frac{n-2}{2}}\tu_{1}\left(\frac{\T^{-1}(x)}{k_{1, \eps}} \right)\right|^{\crits-2-\pe} \notag\\
&+2^{\crits-1-\pe}\left(\frac{|x|}{k_{1,\eps}}\right)^{2-s} \left|\tu_{1} \left(\frac{ \T^{-1}(x)}{k_{1, \eps}}\right) \right|^{\crits-2-\pe} \notag\\
& \leq   \frac{\gamma' -\gamma}{4} \qquad \hbox{ for all } x\in B_{\rho k_{1,\eps}}(0) \cap \Omega.\notag
\end{align}
Hence  as in Step \ref{sec:se:1}.1 we have that for $\eps>0$ small enough
\begin{align*}
\frac{\mathcal{L}_\eps \varphi (x)}{ \varphi(x)}>0 ~ \hbox{ for all } x\in B_{\rho k_{1,\eps}}(0) \cap \Omega. 
\end{align*}
\hfill $\Box$

\noindent{\it Step \ref{sec:se:1}.3.2:}   It follows from point (A4) of Proposition \ref{prop:exhaust}  that there exists $C'_2(\rho)>0$ such that for all $\eps>0$  small
\begin{align*}
|\ue(x)|\leq C'_2(\rho)  \frac{d(x,\bdry)}{\mu_{1,\eps}^{\frac{\beta_{+}(\gamma')-\beta_{-}(\gamma')}{2}} |x|^{\beta_{-}(\gamma')}} \qquad \hbox{ for all } x\in \partial B_{\rho k_{1, \eps}}(0) \cap \Omega
\end{align*}
and then by the estimate \eqref{esti:sub-sol} on $\varphi$ we have for some constant $C_2(\delta)>0$ 
\begin{align}\label{sup:ue:boundary:3:tis}
|\ue(x)|\leq C_2(\rho)  \frac{\varphi(x) }{\mu_{1,\eps}^{\frac{\beta_{+}(\gamma')-\beta_{-}(\gamma')}{2}} }  \qquad \hbox{ for all } x\in  \partial B_{\rho k_{1, \eps}}(0) \cap \Omega.
\end{align}

\noindent
We let for all $\eps>0$ 
$$\Psi^{0}_{\eps}(x):=  C_2(\rho)  \frac{\varphi(x) }{\mu_{1,\eps}^{\frac{\beta_{+}(\gamma')-\beta_{-}(\gamma')}{2}} }\qquad \hbox{ for } x  \in \Omega. $$
Then   \eqref{sup:ue:boundary:3:tis} implies that for all $\eps>0$  small
\begin{align}\label{control:ue:bord:tis}
|u_{\eps}(x)| \leq \Psi^{0}_{\eps}(x) \qquad  \hbox{ for all } x \in \partial\left(  B_{\rho k_{1, \eps}}(0) \cap \Omega \setminus \{ 0\}\right) .
\end{align}
 
\noindent
Therefore it follows from \eqref{ineq:LG3} and \eqref{control:ue:bord:tis} that $\eps>0$ sufficiently small
$$\left\{\begin{array}{ll}
\mathcal{L}_\eps\Psi^{0}_\eps \geq 0=\mathcal{L}_\eps\ue & \hbox{ in }  B_{\rho k_{1,\eps}}(0) \cap \Omega \\
\Psi^{0}_\eps \geq \ue& \hbox{ on } \partial\left(  B_{\rho k_{1, \eps}}(0) \cap \Omega \setminus \{ 0\}\right) \\
\mathcal{L}_\eps \Psi^{0}_\eps \geq 0=-\mathcal{L}_\eps\ue & \hbox{ in }  B_{\rho k_{1,\eps}}(0) \cap \Omega\\
\Psi^{0}_\eps \geq -\ue& \hbox{ on } \partial\left(  B_{\rho k_{1, \eps}}(0) \cap \Omega \setminus \{ 0\}\right).
\end{array}\right.$$
Since  the operator $\mathcal{L}_\eps$ satisfies the comparison principle on $ B_{\rho k_{1,\eps}}(0)$. Therefore
$$|\ue(x)|\leq\Psi^{0}_\eps(x) \qquad  \hbox{ for all } x  \in B_{\rho k_{1,\eps}}(0) \cap \Omega.$$
And so  for all $\eps >0$ small
$$|u_{\eps}(x)|\leq C_2(\rho) \frac{ \varphi(x)}{\mu_{1,\eps}^{\frac{\beta_{+}(\gamma')-\beta_{-}(\gamma')}{2}} }  ~  \hbox{ for all } x  \in B_{\rho k_{1,\eps}}(0)\cap \Omega.$$
for  $\rho$ small. Using the estimate \eqref{esti:sub-sol} we have or all $\eps>0$  small
\begin{align*}
|u_{\eps}(x)| &~\leq  C_2(\rho) \frac{d(x,\bdry)}{\mu_{1,\eps}^{\frac{\beta_{+}(\gamma')-\beta_{-}(\gamma')}{2}} |x|^{\beta_{-}(\gamma')}}  \notag\\ 
 &~\leq C_2(\rho) \frac{|x|}{\mu_{1,\eps}^{\frac{\beta_{+}(\gamma')-\beta_{-}(\gamma')}{2}} |x|^{\beta_{-}(\gamma')}}  .
\end{align*}
for  $\rho$ small. It then follows from point (A4) of Proposition \ref{prop:exhaust} that estimate \eqref{estim:alpha:1:tis} holds on  $x  \in B_{\rho k_{1,\eps}}(0) \cap \Omega$ for all $\rho>0$.\qed \\

\noindent{\bf Step \ref{sec:se:1}.4:}  Combining the previous three steps,  it follows from \eqref{estim:alpha:1}, \eqref{estim:alpha:1:bis}, \eqref{estim:alpha:1:tis} and Proposition \ref{prop:exhaust} that for any $\alpha>0$ small, there exists $C(\alpha)>0$ such that for all $\eps >0$  we have for all $x \in \Omega$, 
\begin{eqnarray}\label{eq:est:alpha:global}
|\ue(x)|&\leq&
C(\alpha) \left(~\sum_{i=1}^N\frac{\mei^{\frac{\bp-\bm}{2}-\alpha} ~ |x|}{ \mei^{\left(\bp-\bm\right)-2\alpha}|x|^{\bm+\alpha}+|x|^{\bp -\alpha}}\right.\nonumber\\
&&\left.+ \frac{\Vert |x|^{\bm-1}u_0 \vert|_{L^{\infty}(\Omega)}}{|x|^{\bm+\alpha}} |x| \right). 
\end{eqnarray}

\medskip

\noindent
Next we improve the above estimate and show that one can take  $\alpha =0$ in  \eqref{eq:est:alpha:global}.\\

\noindent
We let $G_{\epsilon}$ be the Green's function for the coercive operator  $-\Delta-\frac{\gamma}{|x|^{2}}-h_{\epsilon}$ on $\Omega$ with Dirichlet boundary condition. Green's representation formula, the pointwise bounds on the Green's function \eqref{est:G:up} and the regularity Theorem \ref{th:hopf}, yields for any $ z \in \Omega$, 
$$\ue(z) = \int \limits_\Omega G_{\epsilon}(z,x)\left(  \frac{|\ue (x)|^{\crits-2-\pe} ~\ue(x)}{|x|^s}\right) dx,$$
and therefore,
\begin{align}\label{ineq:1} 
|\ue(z)| \leq & \int \limits_\Omega G_{\epsilon}(z,x)  \frac{|\ue (x)|^{\crits-1-\pe} }{|x|^s} ~ dx \notag \\
\leq&~ C  \int \limits_\Omega  \left(\frac{\max\{|z|,|x|\}}{\min\{|z|,|x|\}}\right)^{\bm} \frac{d(x,\bdry) d(z,\bdry)}{|x-z|^{n}} \frac{|\ue(x)|^{\crits-1-\pe}}{|x|^s}~ dx.
\end{align}
\noindent
To estimate the above integral  we break it into three parts. \\

\noindent{\bf Step \ref{sec:se:1}.5:} There exist $C>0$ such that for any sequence $(\ze)$ with $\ze\in \Omega\setminus \overline{B}_{k_{N,\eps}}(0)$,   we have
\begin{align}\label{estim:1:ze}
\int \limits_\Omega G_{\epsilon}(\ze,x)  \frac{|\ue (x)|^{\crits-1-\pe} }{|x|^s}~dx \leq  C~\left(\frac{\mu_{N,\eps}^{\frac{\bp-\bm}{2}} |\ze|}{|\ze|^{\bp}}+\frac{\Vert |x|^{\bm-1}u_0 \vert|_{L^{\infty}(\Omega)}}{|\ze|^{\bm }}|\ze| \right). 
\end{align}
 \noindent{\it Proof of  Step \ref{sec:se:1}.5:} To estimate the right-hand-side of  \eqref{ineq:1} in this case, we split  $\Omega$ into four subdomains as:  $\Omega = \bigcup \limits_{i=1}^{4} D^{N}_{i,\eps}$ where
\begin{itemize}
\item
$D^{N}_{1,\eps}:= B_{k_{N,\eps}}(0) \cap \Omega,\; D^{N}_{2,\eps} :=\{ k_{N,\eps}<|x|< \frac{1}{2}|\ze|\} \cap \Omega,$
\item
$D^{N}_{3,\eps}:=\{\frac{1}{2}|\ze|<|x|< 2|\ze|\} \cap \Omega,\; D^{N}_{4,\eps}:=\{ 2|\ze|<|x|\} \cap \Omega. $
\end{itemize}
Note that one has $\frac{1}{2}|\ze|< |x-\ze|$ in $D^{N}_{2,\eps}$ and $\frac{1}{2}|x|< |x-\ze|$ in $D^{N}_{4,\eps}$.

\noindent
Using point (A5) of Proposition \ref{prop:exhaust} and a change of variable, we get
\begin{align}\label{estim:1:I1:N}
I^{N}_{1}: =&~C~ \int \limits_{D^{N}_{1,\eps}}  \left(\frac{\max\{|\ze|,|x|\}}{\min\{|\ze|,|x|\}}\right)^{\bm} \frac{d(x,\bdry) d(\ze,\bdry)}{|x-\ze|^{n}} \frac{|\ue(x)|^{\crits-1-\pe}}{|x|^s}~ dx \notag \\
\leq &~C~d(\ze,\bdry) \int \limits_{D^{N}_{1,\eps}} \frac{|\ze|^{\bm}}{|x|^{\bm-1}} |x-\ze|^{-n}\frac{|\ue(x)|^{\crits-1-\pe}}{|x|^s}~ dx  \notag \\
\leq & ~C~ \frac{d(\ze,\bdry)}{|\ze|^{\bp}} \int \limits_{D^{N}_{1,\eps}} \frac{|\ue(x)|^{\crits-1-\pe}}{|x|^{\bm-1 +s}} ~ dx  \notag \\
\leq & ~C~ \frac{d(\ze,\bdry)}{|\ze|^{\bp}} \int \limits_{B_{k_{N,\eps}}(0)} \frac{1}{|x|^{\bm -1+s+ (\crits-1-\pe)\frac{n-2}{2}}} ~ dx  \notag \\
\leq & ~C ~\frac{\mu_{N,\eps}^{\frac{\bp-\bm}{2}} d(\ze,\bdry)}{  |\ze|^{\bp} } \int \limits_{B_{1}(0)} \frac{1}{|x|^{ n- \frac{\bp-\bm}{2} -\pe \frac{n-2}{2} }} ~ dx \notag \\
\leq & ~C~ \frac{\mu_{N,\eps}^{\frac{\bp-\bm}{2}}|\ze|}{  |\ze|^{\bp} }  . 
\end{align}
\noindent
Now we estimate
\begin{align*}
I^{N}_{2}: =&~C \int \limits_{D^{N}_{2,\eps}}  \left(\frac{\max\{|\ze|,|x|\}}{\min\{|\ze|,|x|\}}\right)^{\bm} \frac{d(x,\bdry) d(\ze,\bdry)}{|x-\ze|^{n}} \frac{|\ue(x)|^{\crits-1-\pe}}{|x|^s}~ dx \notag \\
\leq &~C~d(\ze,\bdry)\frac{|\ze|^{\bm}}{|\ze|^{n}} \int \limits_{D^{N}_{2,\eps}}  \frac{|\ue(x)|^{\crits-1-\pe}}{|x|^{\bm-1 +s}}~ dx.  \notag \\
\end{align*}
Using \eqref{estim:alpha:1} we get for $0<\alpha < \frac{\crits-2}{\crits-1}\left( \frac{\bp-\bm}{2}\right)$
\begin{align}\label{estim:1:I2:N}
I^{N}_{2}\leq &~C~ \frac{d(\ze,\bdry)}{|\ze|^{\bp}}  \int \limits_{D^{N}_{2,\eps}}  \frac{|\ue(x)|^{\crits-1-\pe}}{|x|^{\bm -1+s}}~ dx  \notag \\
\leq &~C~ \frac{d(\ze,\bdry)}{|\ze|^{\bp}}   \mu_{N,\eps}^{\left( \frac{\bp-\bm}{2}-\alpha\right)(\crits-1-\pe)}  \int \limits_{D^{N}_{2,\eps}}  \frac{|x|^{-\bm+1 -s} }{|x|^{(\crits-1-\pe)(\bp-1-\alpha)}}~ dx  \notag \\
& + C~ \frac{d(\ze,\bdry)}{|\ze|^{\bp}} \int \limits_{D^{N}_{2,\eps}}  \frac{\Vert |x|^{\bm-1}u_0 \vert|^{\crits-1-\pe}_{L^{\infty}(\Omega)} }{|x|^{(\crits-1-\pe)(\bm-1+\alpha)+\bm -1+s}}~ dx   \notag \\
\leq &~C~  \frac{d(\ze,\bdry)}{|\ze|^{\bp}} \mu_{N,\eps}^{\frac{\bp-\bm}{2}}  \int \limits_{1 \leq |x| }  \frac{1}{|x|^{n+ (\crits-2-\pe)\left( \frac{\bp-\bm}{2}\right)-\alpha  (\crits-1-\pe)}}~ dx  \notag \\
& + C~ \frac{d(\ze,\bdry)}{|\ze|^{\bp}}  \int \limits_{|x| \leq \frac{1}{2}|\ze|}  \frac{\Vert |x|^{\bm-1}u_0 \vert|^{\crits-1-\pe}_{L^{\infty}(\Omega)} }{|x|^{ (\crits-\pe) (\bm-1) +s +\alpha(\crits-1-\pe)}}~ dx   \notag \\
\leq &~C~  \ \frac{\mu_{N,\eps}^{\frac{\bp-\bm}{2}} |\ze|}{|\ze|^{\bp}}  \int \limits_{1 \leq |x| }  \frac{1}{|x|^{n+ (\crits-2-\pe)\left( \frac{\bp-\bm}{2}\right)-\alpha  (\crits-1-\pe)}}~ dx  \notag \\
& + C~ \frac{|\ze|^{(\crits-2-\pe)\left( \frac{\bp-\bm}{2}\right)-\alpha  (\crits-1-\pe)} \Vert |x|^{\bm-1}u_0 \vert|^{\crits-1-\pe}_{L^{\infty}(\Omega)} }{|\ze|^{\bm}} |\ze|  \notag \\
\leq &~ C \left(\ \frac{\mu_{N,\eps}^{\frac{\bp-\bm}{2}} |\ze|}{|\ze|^{\bp}}  + \frac{ \Vert |x|^{\bm-1}u_0 \vert|^{\crits-1-\pe}_{L^{\infty}(\Omega)} }{|\ze|^{\bm}} |\ze|  \right).
\end{align}
\noindent
For the next integral
\begin{align*}
I^{N}_{3} :=&~C ~\int \limits_{D^{N}_{3,\eps}}  \left(\frac{\max\{|\ze|,|x|\}}{\min\{|\ze|,|x|\}}\right)^{\bm} \frac{d(x,\bdry) d(\ze,\bdry)}{|x-\ze|^{n}} \frac{|\ue(x)|^{\crits-1-\pe}}{|x|^s}~ dx\notag \\
\leq &~C~d(\ze,\bdry) \int \limits_{D^{N}_{3,\eps}} \frac{|x| }{|x-\ze|^{n}}\frac{|\ue(x)|^{\crits-1-\pe}}{|x|^s}~ dx. 
\end{align*}
From \eqref{estim:alpha:1} we get for $0<\alpha < \frac{\crits-2}{\crits-1}\left( \frac{\bp-\bm}{2}\right)$
\begin{align}\label{estim:1:I3:N}
I^{N}_{3} \leq &~C ~d(\ze,\bdry) \mu_{N,\eps}^{\left( \frac{\bp-\bm}{2}-\alpha\right)(\crits-1-\pe)}  \int \limits_{D^{N}_{3,\eps}}  \frac{|x|^{1-s} |x-\ze|^{-n}}{|x|^{(\bp-1 -\alpha)(\crits -1-\pe)}}~ dx  \notag \\
& + Cd(\ze,\bdry) \int \limits_{D^{N}_{3,\eps}}  \frac{ |x|^{1-s} |x-\ze|^{-n}}{|x|^{(\bm-1 +\alpha)(\crits -1-\pe)}}\Vert |x|^{\bm-1}u_0 \vert|^{\crits-1-\pe}_{L^{\infty}(\Omega)} ~dx  \notag \\
\leq &~C ~\frac{ \mu_{N,\eps}^{\left( \frac{\bp-\bm}{2}-\alpha\right)(\crits-1-\pe)} d(\ze,\bdry)}{|\ze|^{(\bp -1-\alpha)(\crits -1-\pe)+s-1}}  \int \limits_{D^{N}_{3,\eps}}   |x-\ze|^{-n} ~ dx  \notag \\
& + C \frac{ \Vert |x|^{\bm-1}u_0 \vert|^{\crits-1-\pe}_{L^{\infty}(\Omega)}d(\ze,\bdry) }{|\ze|^{(\bm -1+\alpha)(\crits -1-\pe)+s-1}} \int \limits_{D^{N}_{3,\eps}} |x-\ze|^{-n}~ dx  \notag \\
\leq &~C ~ \frac{\mu_{N,\eps}^{\frac{\bp-\bm}{2}} d(\ze,\bdry) }{|\ze|^{\bp}}  \left( \frac{\mu_{N,\eps}}{|\ze|}\right)^{(\crits-2-\pe)\left( \frac{\bp-\bm}{2}\right)-\alpha  (\crits-1-\pe)} \notag \\
& + C~\frac{ \Vert |x|^{\bm-1}u_0 \vert|^{\crits-1-\pe}_{L^{\infty}(\Omega)}d(\ze,\bdry)}{|\ze|^{\bm}} |\ze|^{(\crits-2-\pe)\left( \frac{\bp-\bm}{2}\right)-\alpha  (\crits-1-\pe)}  \notag \\
\leq &~ C~\left(\frac{\mu_{N,\eps}^{\frac{\bp-\bm}{2}}|\ze|}{|\ze|^{\bp}}+\frac{\Vert |x|^{\bm-1}u_0 \vert|^{\crits-1-\pe}_{L^{\infty}(\Omega)}}{|\ze|^{\bm -1}} |\ze|\right). 
\end{align}
\noindent
Finally we estimate 
\begin{align*}
I^{N}_{4} :=&~C~ \int \limits_{D^{N}_{4,\eps}}\left(\frac{\max\{|\ze|,|x|\}}{\min\{|\ze|,|x|\}}\right)^{\bm} \frac{d(x,\bdry) d(\ze,\bdry)}{|x-\ze|^{n}} \frac{|\ue(x)|^{\crits-1-\pe}}{|x|^s}~ dx \notag \\
\leq &~C~\frac{d(\ze,\bdry)}{|\ze|^{\bm}}  \int \limits_{2 |\ze| \leq |x|} |x|^{\bm+1-n} \frac{|\ue(x)|^{\crits-1-\pe}}{|x|^s}~ dx \notag \\
\leq &~C~ \frac{d(\ze,\bdry)}{|\ze|^{\bm}} \int \limits_{2 |\ze| \leq |x|}  \frac{|\ue(x)|^{\crits-1-\pe}}{|x|^{\bp+s-1}}~ dx. \notag \\
\end{align*}
Then from \eqref{estim:alpha:1} we get for $0<\alpha < \frac{\crits-2}{\crits-1}\left( \frac{\bp-\bm}{2}\right)$
\begin{align}\label{estim:1:I4:N}
I^{N}_{4}\leq &~C~\frac{\mu_{N,\eps}^{\left( \frac{\bp-\bm}{2}-\alpha\right)(\crits-1-\pe)}d(\ze,\bdry)}{|\ze|^{\bm}}  \int \limits_{2 |\ze| \leq |x|}  \frac{|x|^{\alpha(\crits-1-\pe)}}{|x|^{(\crits-\pe )(\bp-1)+s}}~ dx  \notag \\
& +C~ \frac{d(\ze,\bdry)}{|\ze|^{\bm}}  \int \limits_{2 |\ze| \leq |x|}  \frac{\Vert |x|^{\bm-1}u_0 \vert|^{\crits-1-\pe}_{L^{\infty}(\Omega)}}{|x|^{ \bp+s+\bm(\crits-1-\pe)+\alpha(\crits-1-\pe)}}~ dx  \notag \\
\leq &~C~\frac{\mu_{N,\eps}^{\left( \frac{\bp-\bm}{2}-\alpha\right)(\crits-1-\pe)}d(\ze,\bdry)}{|\ze|^{\bm}}  \int \limits_{2 |\ze| \leq |x|}  \frac{|x|^{\alpha(\crits-1-\pe)}}{|x|^{n+( \crits-\pe)\left( \frac{\bp-\bm}{2}\right)}}~ dx  \notag \\
& +C~ \frac{d(\ze,\bdry)}{|\ze|^{\bm}}  \int \limits_{2 |\ze| \leq |x|}  \frac{\Vert |x|^{\bm-1}u_0 \vert|^{\crits-1-\pe}_{L^{\infty}(\Omega)}}{|x|^{ n- \left[(\crits-2-\pe)\left( \frac{\bp-\bm}{2}\right)-\alpha  (\crits-1-\pe)\right]}}~ dx  \notag \\
\leq &~C~\frac{\mu_{N,\eps}^{\left( \frac{\bp-\bm}{2}-\alpha\right)(\crits-1-\pe)}}{|\ze|^{\bm}}    \frac{d(\ze,\bdry)}{|\ze|^{ (\crits-\pe)\left( \frac{\bp-\bm}{2}\right)-\alpha  (\crits-1-\pe)}}  \notag \\
& + C~\frac{\Vert |x|^{\bm-1}u_0 \vert|^{\crits-1-\pe}_{L^{\infty}(\Omega)}d(\ze,\bdry)}{|\ze|^{\bm}}   \notag \\
\leq &~C~\frac{\mu_{N,\eps}^{ \frac{\bp-\bm}{2}}d(\ze,\bdry)}{|\ze|^{\bp}}     \left( \frac{\mu_{N,\eps}}{|\ze|}\right)^{(\crits-2-\pe)\left( \frac{\bp-\bm}{2}\right)-\alpha  (\crits-1-\pe)}   \notag \\
&+ C~\frac{\Vert |x|^{\bm-1}u_0 \vert|^{\crits-1-\pe}_{L^{\infty}(\Omega)}d(\ze,\bdry)}{|\ze|^{\bm}}   \notag \\
\leq &~ C~\left(\frac{\mu_{N,\eps}^{\frac{\bp-\bm}{2}}|\ze|}{|\ze|^{\bp}}+\frac{\Vert |x|^{\bm-1}u_0 \vert|^{\crits-1-\pe}_{L^{\infty}(\Omega)}}{|\ze|^{\bm }} |\ze|\right). 
\end{align} 
\noindent
Combining \eqref{estim:1:I1:N}, \eqref{estim:1:I2:N}, \eqref{estim:1:I3:N} and \eqref{estim:1:I4:N},   we then obtain  for some constant $C>0$
\begin{eqnarray*}
\int \limits_\Omega G_{\epsilon}(\ze,x)  \frac{|\ue (x)|^{\crits-1-\pe} }{|x|^s}~dx  &\leq& ~ C~\left(\frac{\mu_{N,\eps}^{\frac{\bp-\bm}{2}}|\ze|}{|\ze|^{\bp}}\right.\\
&&\left.+\frac{\Vert |x|^{\bm-1}u_0 \vert|^{\crits-1-\pe}_{L^{\infty}(\Omega)}}{|\ze|^{\bm }} |\ze| \right),
\end{eqnarray*}
which we write as
\begin{eqnarray*}
&&\int \limits_\Omega G_{\epsilon}(\ze,x)  \frac{|\ue (x)|^{\crits-1-\pe} }{|x|^s}~dx \\
&&\leq ~ C~\left(\frac{\mu_{N,\eps}^{\frac{\bp-\bm}{2}}|\ze|}{|\ze|^{\bp}}+\frac{\Vert |x|^{\bm-1}u_0 \vert|_{L^{\infty}(\Omega)}}{|\ze|^{\bm }}|\ze| \right)
\end{eqnarray*}
for some $C>0$. This proves \eqref{estim:1:ze}. \qed \\

\noindent{\bf Step \ref{sec:se:1}.6:} There exists $C>0$ such that for sequence of points $(\ze)$ in $ B_{k_{1,\eps}}(0) \cap \Omega$  we have 
\begin{align}\label{estim:1:tis:ze}
\int \limits_\Omega G_{\epsilon}(\ze,x)  \frac{|\ue (x)|^{\crits-1-\pe} }{|x|^s}~dx  \leq C~ \frac{|\ze|}{\mu_{1,\eps}^{\frac{\bp-\bm}{2}} |\ze|^{\bm }}    .
\end{align}

\smallskip\noindent{\it Proof of  Step \ref{sec:se:1}.6:}  Here again, to estimate the right-hand-side of \eqref{ineq:1} in this case, we split  $\Omega$ into four subdomains as: $\Omega = \bigcup \limits_{i=1}^{4} D^{1}_{i,\eps}(R)$ where

\begin{itemize}
\item
$D^{1}_{1,\eps} :=\{|x|< \frac{1}{2}|\ze|\} \cap \Omega,$
\item
$D^{1}_{2,\eps}:=\{\frac{1}{2}|\ze|<|x|< 2|\ze|\} \cap \Omega,$
\item
$D^{1}_{3,\eps}:=\{ 2|\ze|<|x| \leq  k_{1,\eps} \}  \cap \Omega, $
\item
$D^{1}_{4,\eps}:=\{  k_{1,\eps} <|x|\}\cap \Omega . $
\end{itemize}
Note that one has $\frac{1}{2}|\ze|< |x-\ze|$ in $D^{1}_{1,\eps}$ and $\frac{1}{2}|x|< |x-\ze|$ in $D^{1}_{3,\eps}$.
We then have 
\begin{align*}
I^{1}_{1} :=&~C \int \limits_{D^{1}_{1,\eps}} \left(\frac{\max\{|\ze|,|x|\}}{\min\{|\ze|,|x|\}}\right)^{\bm} \frac{d(x,\bdry) d(\ze,\bdry) }{|x-\ze|^{n}} \frac{|\ue(x)|^{\crits-1-\pe}}{|x|^s}~ dx \notag \\
\leq &~C d(\ze,\bdry)\frac{|\ze|^{\bm}}{|\ze|^{n-2}} \int \limits_{D^{1}_{1,\eps}}  \frac{|\ue(x)|^{\crits-1-\pe}}{|x|^{\bm +s-1}}~ dx.  \notag \\
\end{align*}
Using \eqref{estim:alpha:1:tis} we get for $0<\alpha < \frac{\crits-2}{\crits-1}\left( \frac{\bp-\bm}{2}\right)$
\begin{align}\label{estim:1:I1:tis}
I^{1}_{1}\leq &~C~ \frac{d(\ze,\bdry)}{|\ze|^{\bp}} \int \limits_{D^{1}_{1,\eps}}  \frac{|\ue(x)|^{\crits-1-\pe}}{|x|^{\bm +s-1}}~ dx  \notag \\
\leq &~   C~ \frac{  \mu_{1,\eps}^{-\left( \frac{\bp-\bm}{2}-\alpha\right)(\crits-1-\pe)} d(\ze,\bdry)}{|\ze|^{\bp}}  \int \limits_{D^{1}_{1,\eps}}  \frac{|x|^{-\bm -s+1}}{|x|^{(\crits-1-\pe)(\bm-1+\alpha)}}~ dx  \notag \\
\leq &~ C~\frac{  \mu_{1,\eps}^{-\left( \frac{\bp-\bm}{2}-\alpha\right)(\crits-1-\pe)} d(\ze,\bdry)}{|\ze|^{\bp}} \int \limits_{|x| \leq \frac{1}{2}|\ze|}  \frac{|x|^{-\alpha(\crits-1-\pe)}}{|x|^{( \crits -\pe)(\bm-1) +s }}~ dx  \notag \\
\leq &~  C~ \left(\frac{|\ze| }{\mu_{1,\eps}} \right)^{(\crits-2-\pe)\left( \frac{\bp-\bm}{2}\right)-\alpha  (\crits-1-\pe)} \frac{d(\ze,\bdry)}{\mu_{1,\eps}^{\frac{\bp-\bm}{2}}|\ze|^{\bm}} \notag \\
\leq &~   C~\frac{|\ze|}{\mu_{1,\eps}^{\frac{\bp-\bm}{2}}|\ze|^{\bm}}.
\end{align}
\noindent
Next we have
\begin{align*}
I^{1}_{2} :=&~C \int \limits_{D^{1}_{2,\eps}} \left(\frac{\max\{|\ze|,|x|\}}{\min\{|\ze|,|x|\}}\right)^{\bm} \frac{d(x,\bdry) d(\ze,\bdry) }{|x-\ze|^{n}} \frac{|\ue(x)|^{\crits-1-\pe}}{|x|^s}~ dx  \notag \\
\leq &~C ~d(\ze,\bdry) \int \limits_{D^{1}_{2,\eps}}  \frac{d(x,\bdry) }{|x-\ze|^{n}} \frac{|\ue(x)|^{\crits-1-\pe}}{|x|^s}~ dx. 
\end{align*}
From \eqref{estim:alpha:1:tis} we get for $0<\alpha < \frac{\crits-2}{\crits-1}\left( \frac{\bp-\bm}{2}\right)$
\begin{align}\label{estim:1:I2:tis}
I^{1}_{2} \leq &~C~  \mu_{1,\eps}^{-\left( \frac{\bp-\bm}{2}-\alpha\right)(\crits-1-\pe)} d(\ze,\bdry) \int \limits_{D^{1}_{2,\eps}}  \frac{|x|^{-s+1} |x-\ze|^{-n}}{|x|^{(\bm-1 +\alpha)(\crits -1-\pe)}}~ dx  \notag \\
\leq &~ C~  \frac{\mu_{1,\eps}^{-\left( \frac{\bp-\bm}{2}-\alpha\right)(\crits-1-\pe)}d(\ze,\bdry)}{|\ze|^{(\bm -1+\alpha)(\crits -1-\pe)+s-1}} \int \limits_{D^{1}_{2,\eps}} |x|^{-n}~ dx  \notag \\
\leq &~  C~ \left(\frac{|\ze| }{\mu_{1,\eps}} \right)^{(\crits-2-\pe)\left( \frac{\bp-\bm}{2}\right)-\alpha  (\crits-1-\pe)} \frac{d(\ze,\bdry) }{\mu_{1,\eps}^{\frac{\bp-\bm}{2}}|\ze|^{\bm}} \notag \\
\leq &~   C~\frac{|\ze|}{\mu_{1,\eps}^{\frac{\bp-\bm}{2}}|\ze|^{\bm}}.
\end{align}
\noindent
For 
\begin{align*}
I^{1}_{3}: =&~C \int \limits_{D^{1}_{3,\eps}}  \left(\frac{\max\{|\ze|,|x|\}}{\min\{|\ze|,|x|\}}\right)^{\bm} \frac{d(x,\bdry) d(\ze,\bdry) }{|x-\ze|^{n}} \frac{|\ue(x)|^{\crits-1-\pe}}{|x|^s}~ dx  \notag \\
\leq &~C ~\frac{d(\ze,\bdry)}{|\ze|^{\bm}}  \int \limits_{2 |\ze| \leq |x|\leq  k_{1,\eps}} |x|^{\bm+1-n} \frac{|\ue(x)|^{\crits-1-\pe}}{|x|^s}~ dx \notag \\
\leq &~C~\frac{d(\ze,\bdry)}{|\ze|^{\bm}} \int \limits_{2 |\ze| \leq |x| \leq  k_{1,\eps}}  \frac{|\ue(x)|^{\crits-1-\pe}}{|x|^{\bp+s-1}}~ dx. \notag \\
\end{align*}
Then from \eqref{estim:alpha:1:tis} we get for $0<\alpha < \frac{\crits-2}{\crits-1}\left( \frac{\bp-\bm}{2}\right)$
\begin{align*}
I^{1}_{3}\leq &~C~ \frac{\mu_{1,\eps}^{-\left( \frac{\bp-\bm}{2}-\alpha\right)(\crits-1-\pe)}d(\ze,\bdry)}{|\ze|^{\bm}}  A_\eps 
\end{align*} 

$$A_\eps:=\int \limits_{2 |\ze| \leq |x| \leq  k_{1,\eps}}  \frac{1}{|x|^{ \bp+s-1+(\bm-1-\pe)(\crits-1)+\alpha(\crits-1-\pe)}}~ dx$$
With a change of variable, we get that
\begin{eqnarray*}
A_\eps &\leq &\int \limits_{2 |\ze| \leq |x| \leq  k_{1,\eps}}  \frac{dx}{|x|^{ n- \left[(\crits-2-\pe)\left( \frac{\bp-\bm}{2}\right)-\alpha  (\crits-1-\pe)\right]}}\\
&\leq & \mu_{1,\eps}^{\left( \frac{\bp-\bm}{2}\right)(\crits-2-\pe)-\alpha(\crits-1-\pe)}\int \limits_{B_{1}(0)}  \frac{|x|^{ -\alpha  (\crits-1-\pe)}dx}{|x|^{ n- (\crits-2)\left( \frac{\bp-\bm}{2}\right)}}\\
&\leq & C \mu_{1,\eps}^{\left( \frac{\bp-\bm}{2}\right)(\crits-2-\pe)-\alpha(\crits-1-\pe)}
\end{eqnarray*}
and then
\begin{align}\label{estim:1:I3:tis}
I^{1}_{3}\leq &~  C~  \frac{d(\ze,\bdry)}{\mu_{1,\eps}^{\frac{\bp-\bm}{2}}|\ze|^{\bm}}   \leq  C~ \frac{|\ze|}{\mu_{1,\eps}^{\frac{\bp-\bm}{2}}|\ze|^{\bm}}.
\end{align}

\noindent
For the last integral, we use point (A5) of Proposition \ref{prop:exhaust} and a change of variable to get
\begin{align}
I^{1}_{4} :=&~C \int \limits_{D^{1}_{4,\eps}} \left(\frac{\max\{|\ze|,|x|\}}{\min\{|\ze|,|x|\}}\right)^{\bm} \frac{d(x,\bdry) d(\ze,\bdry) }{|x-\ze|^{n}} \frac{|\ue(x)|^{\crits-1-\pe}}{|x|^s}~ dx \notag \\
\leq &~C~\frac{d(\ze,\bdry)}{|\ze|^{\bm}} \int \limits_{|x| \geq  k_{1,\eps}}  \frac{|\ue(x)|^{\crits-1-\pe}}{|x|^{\bp+s-1}}~ dx \notag \\
\leq &~C ~\frac{d(\ze,\bdry)}{|\ze|^{\bm}} \int \limits_{|x| \geq  k_{1,\eps}}  \frac{1}{|x|^{\bp+s -1+\frac{n-2}{2}(\crits-1-\pe) }}~ dx \label{estim:1:I4:tis}\\
\leq &~ C~  \frac{d(\ze,\bdry)}{\mu_{1,\eps}^{\frac{\bp-\bm}{2}}|\ze|^{\bm}} \int \limits_{|x| \geq 1 }\frac{d}{|x|^{ n+\frac{\bp-\bm}{2}} }\leq  C~ \frac{|\ze|}{\mu_{1,\eps}^{\frac{\bp-\bm}{2}}|\ze|^{\bm}}.\notag
\end{align}
\medskip

\noindent
Combining \eqref{estim:1:I1:tis}, \eqref{estim:1:I2:tis}, \eqref{estim:1:I3:tis} and \eqref{estim:1:I4:tis},   we then obtain \eqref{estim:1:tis:ze}.  \hfill $\Box$
\bigskip

\noindent{\bf Step \ref{sec:se:1}.7:} Let $1 \leq i \leq N-1$. There exists $C>0$ such that for sequence of points $(\ze)$ in $B_{ k_{i+1,\eps}}(0)\setminus \overline{B}_{ k_{i,\eps}}(0) \cap \Omega$  we have 
\begin{align}\label{estim:1:bis:ze}
\int \limits_\Omega G_{\epsilon}(\ze,x)  &\frac{|\ue (x)|^{\crits-1-\pe} }{|x|^s}~dx  \leq C~\left(\frac{\mu_{i,\eps}^{\frac{\bp-\bm}{2}}|\ze|}{|\ze|^{\bp}}+\frac{|\ze|}{\mu_{i+1,\eps}^{\frac{\bp-\bm}{2}} |\ze|^{\bm }}  \right).
\end{align}

\smallskip\noindent{\it Proof of  Step \ref{sec:se:1}.7:}  To estimate the right-hand-side of \eqref{ineq:1} in this case, we split  $\Omega$ into  five subdomains as: $\Omega = \bigcup \limits_{j=1}^{5} D^{i}_{j,\eps}$ where

\begin{itemize}
\item
$D^{i}_{1,\eps}:= B_{k_{i,\eps}}(0) \cap \Omega,$\\
\item
$D^{i}_{2,\eps} :=\{ k_{i,\eps}<|x|< \frac{1}{2}|\ze|\} \cap \Omega,$\\
\item
$D^{i}_{3,\eps}:=\{\frac{1}{2}|\ze|<|x|< 2|\ze|\} \cap \Omega,$\\
\item
$D^{i}_{4,\eps}:=\{ 2|\ze|<|x| < k_{i+1,\eps} \} \cap \Omega, $\\
\item
$D^{i}_{5,\eps}:=\{  k_{i+1,\eps} <|x| \} \cap \Omega. $
\end{itemize}
Note that one has $\frac{1}{2}|\ze|< |x-\ze|$ in $D^{i}_{2,\eps}$ and $\frac{1}{2}|x|< |x-\ze|$ in $D^{i}_{4,\eps}$.
\medskip

\noindent
First we have using point (A5) of Proposition \ref{prop:exhaust} and a change of variable
\begin{align}\label{estim:1:I1:bis}
I^{i}_{1} :=&~C \int \limits_{D^{i}_{1,\eps}} \left(\frac{\max\{|\ze|,|x|\}}{\min\{|\ze|,|x|\}}\right)^{\bm} \frac{d(x,\bdry) d(\ze,\bdry) }{|x-\ze|^{n}} \frac{|\ue(x)|^{\crits-1-\pe}}{|x|^s}~ dx\notag \\
\leq &~C~d(\ze,\bdry) \int \limits_{D^{i}_{1,\eps}} \frac{|\ze|^{\bm}}{|x|^{\bm-1}} |x-\ze|^{-n}\frac{|\ue(x)|^{\crits-1-\pe}}{|x|^s}~ dx  \notag \\
\leq & ~C~ \frac{d(\ze,\bdry)}{|\ze|^{\bp}} \int \limits_{D^{i}_{1,\eps}} \frac{|\ue(x)|^{\crits-1-\pe}}{|x|^{\bm-1 +s}} ~ dx  \notag \\
\leq & ~C~ \frac{d(\ze,\bdry)}{|\ze|^{\bp}} \int \limits_{B_{k_{i,\eps}}(0)} \frac{1}{|x|^{\bm -1+s+ (\crits-1-\pe)\frac{n-2}{2}}} ~ dx  \notag \\
\leq & ~C ~\frac{\mu_{i,\eps}^{\frac{\bp-\bm}{2}} d(\ze,\bdry)}{  |\ze|^{\bp} } \int \limits_{B_{1}(0)} \frac{1}{|x|^{ n- \frac{\bp-\bm}{2}} } ~ dx\leq C~ \frac{\mu_{i,\eps}^{\frac{\bp-\bm}{2}}|\ze|}{  |\ze|^{\bp} }  .  
\end{align}

\noindent
Now we estimate
\begin{align*}
I^{i}_{2} :=&~C \int \limits_{D^{i}_{2,\eps}} \left(\frac{\max\{|\ze|,|x|\}}{\min\{|\ze|,|x|\}}\right)^{\bm} \frac{d(x,\bdry) d(\ze,\bdry) }{|x-\ze|^{n}} \frac{|\ue(x)|^{\crits-1-\pe}}{|x|^s}~ dx\notag \\
\leq &~C~d(\ze,\bdry)\frac{|\ze|^{\bm}}{|\ze|^{n}} \int \limits_{D^{i}_{2,\eps}}  \frac{|\ue(x)|^{\crits-1-\pe}}{|x|^{\bm-1 +s}}~ dx.  \notag \\ \notag \\
\end{align*}
Using \eqref{estim:alpha:1:bis} we get for $0<\alpha < \frac{\crits-2}{\crits-1}\left( \frac{\bp-\bm}{2}\right)$
\begin{align}\label{estim:1:I2:bis}
I^{i}_{2}\leq &~C~ \frac{d(\ze,\bdry)}{|\ze|^{\bp}} \int \limits_{D^{i}_{2,\eps}}  \frac{|\ue(x)|^{\crits-1-\pe}}{|x|^{\bm +s}}~ dx  \notag \\
\leq &~C~ \frac{d(\ze,\bdry)}{|\ze|^{\bp}}   \mu_{i,\eps}^{\left( \frac{\bp-\bm}{2}-\alpha\right)(\crits-1-\pe)}  \int \limits_{D^{i}_{2,\eps}}  \frac{|x|^{-\bm+1 -s} \, dx}{|x|^{(\crits-1-\pe)(\bp-1-\alpha)}}  \notag \\
& + C~ \frac{  \mu_{i+1,\eps}^{-\left( \frac{\bp-\bm}{2}-\alpha\right)(\crits-1-\pe)}d(\ze,\bdry) }{|\ze|^{\bp}}   \int \limits_{D^{i}_{2,\eps}}  \frac{|x|^{-\bm+1 -s} \, dx}{|x|^{(\crits-1-\pe)(\bm-1+\alpha)}}  \notag \\
\leq &~C~  \frac{d(\ze,\bdry)}{|\ze|^{\bp}} \mu_{i,\eps}^{\frac{\bp-\bm}{2}}  \int \limits_{1 \leq |x| }  \frac{dx}{|x|^{n+ (\crits-2-\pe)\left( \frac{\bp-\bm}{2}\right)-\alpha  (\crits-1-\pe)}}  \notag \\
& + C~ \frac{  \mu_{i+1,\eps}^{-\left( \frac{\bp-\bm}{2}-\alpha\right)(\crits-1-\pe)}d(\ze,\bdry) }{|\ze|^{\bp}}  \int \limits_{|x| \leq \frac{1}{2}|\ze|}  \frac{|x|^{-\alpha(\crits-1-\pe)}\, dx}{|x|^{ \crits (\bm-1) +s }}  \notag \\
\leq &~C~  \frac{\mu_{i,\eps}^{\frac{\bp-\bm}{2}}d(\ze,\bdry)}{|\ze|^{\bp}}  \int \limits_{1 \leq |x| }  \frac{dx}{|x|^{n+ (\crits-2-\pe)\left( \frac{\bp-\bm}{2}\right)-\alpha  (\crits-1-\pe)}}\notag \\
& + C~ \left(\frac{|\ze| }{\mu_{i+1,\eps}} \right)^{(\crits-2-\pe)\left( \frac{\bp-\bm}{2}\right)-\alpha  (\crits-1-\pe)}  \frac{d(\ze,\bdry)}{ \mu^{\frac{\bp-\bm}{2}}_{i+1,\eps} |\ze|^{\bm}}  \notag \\
\leq &~  C ~\left (\frac{\mu_{i,\eps}^{\frac{\bp-\bm}{2}}|\ze|}{|\ze|^{\bp}}     +  \frac{|\ze|}{ \mu^{\frac{\bp-\bm}{2}}_{i+1,\eps} |\ze|^{\bm}}   \right).
\end{align}
And next
\begin{align*}
I^{i}_{3} :=&~C \int \limits_{D^{i}_{3,\eps}} \left(\frac{\max\{|\ze|,|x|\}}{\min\{|\ze|,|x|\}}\right)^{\bm} \frac{d(x,\bdry) d(\ze,\bdry) }{|x-\ze|^{n}} \frac{|\ue(x)|^{\crits-1-\pe}}{|x|^s}~ dx\notag \\
\leq &~C~d(\ze,\bdry) \int \limits_{D^{i}_{3,\eps}} \frac{|x| }{|x-\ze|^{n}}\frac{|\ue(x)|^{\crits-1-\pe}}{|x|^s}~ dx. 
\end{align*}
From \eqref{estim:alpha:1:bis} we get for $0<\alpha < \frac{\crits-2}{\crits-1}\left( \frac{\bp-\bm}{2}\right)$
\begin{align}\label{estim:1:I3:bis}
I^{i}_{3} \leq &~C ~d(\ze,\bdry) \mu_{i,\eps}^{\left( \frac{\bp-\bm}{2}-\alpha\right)(\crits-1-\pe)}  \int \limits_{D^{i}_{3,\eps}}  \frac{|x|^{1-s} |x-\ze|^{-n}\, dx}{|x|^{(\bp-1 -\alpha)(\crits -1-\pe)}}  \notag \\
& + C~d(\ze,\bdry) \mu_{i+1,\eps}^{-\left( \frac{\bp-\bm}{2}-\alpha\right)(\crits-1-\pe)} \int \limits_{D^{i}_{3,\eps}}  \frac{|x|^{1-s} |x-\ze|^{-n}\, dx}{|x|^{(\bm-1 +\alpha)(\crits -1-\pe)}}  \notag \\
\leq &~C ~\frac{ \mu_{i,\eps}^{\left( \frac{\bp-\bm}{2}-\alpha\right)(\crits-1-\pe)}d(\ze,\bdry)}{|\ze|^{(\bp -1-\alpha)(\crits -1-\pe)+s-1}}  \int \limits_{D^{i}_{3,\eps}}   |x-\ze|^{-n} ~ dx  \notag \\
& + C~ \frac{  \mu_{i+1,\eps}^{-\left( \frac{\bp-\bm}{2}-\alpha\right)(\crits-1-\pe)}d(\ze,\bdry) }{|\ze|^{(\bm -1+\alpha)(\crits -1-\pe)+s-1}} \int \limits_{D^{i}_{3,\eps}} |x-\ze|^{-n}~ dx  \notag \\
\leq &~C ~ \frac{\mu_{i,\eps}^{\frac{\bp-\bm}{2}} d(\ze,\bdry) }{|\ze|^{\bp}}  \left( \frac{\mu_{i,\eps}}{|\ze|}\right)^{(\crits-2-\pe)\left( \frac{\bp-\bm}{2}\right)-\alpha  (\crits-1-\pe)} \notag \\
& + C~\left(\frac{|\ze| }{\mu_{i+1,\eps}} \right)^{(\crits-2-\pe)\left( \frac{\bp-\bm}{2}\right)-\alpha  (\crits-1-\pe)}  \frac{ d(\ze,\bdry)}{|\ze|^{\bm}} \notag \\
\leq &~ C~\left(\frac{\mu_{i,\eps}^{\frac{\bp-\bm}{2}}|\ze|}{|\ze|^{\bp}}+\frac{|\ze| }{ \mu^{\frac{\bp-\bm}{2}}_{i+1,\eps} |\ze|^{\bm}}  \right). 
\end{align}
The next integral becomes
\begin{align*}
I^{i}_{4} :=&~C \int \limits_{D^{i}_{4,\eps}} \left(\frac{\max\{|\ze|,|x|\}}{\min\{|\ze|,|x|\}}\right)^{\bm} \frac{d(x,\bdry) d(\ze,\bdry)}{|x-\ze|^{n}} \frac{|\ue(x)|^{\crits-1-\pe}}{|x|^s}~ dx \notag \\
\leq &~C~\frac{d(\ze,\bdry)}{|\ze|^{\bm}} \int \limits_{2 |\ze| \leq |x| <  k_{i+1,\eps}}  |x|^{\bm+1-n} \frac{|\ue(x)|^{\crits-1-\pe}}{|x|^s}~ dx \notag \\
\leq &~C~ \frac{d(\ze,\bdry)}{|\ze|^{\bm}} \int \limits_{2 |\ze| \leq |x| <  k_{i+1,\eps}} \frac{|\ue(x)|^{\crits-1-\pe}}{|x|^{\bp+s-1}}~ dx. \notag \\\notag \\
\end{align*}
Then from \eqref{estim:alpha:1:bis} we get for $0<\alpha < \frac{\crits-2}{\crits-1}\left( \frac{\bp-\bm}{2}\right)$
\begin{eqnarray*}
I^{i}_{4}&\leq& C\frac{d(\ze,\bdry)}{|\ze|^{\bm}}\left(\mu_{i,\eps}^{\left( \frac{\bp-\bm}{2}-\alpha\right)(\crits-1-\pe)} A_\eps \right.\\
&&\left.+ \mu_{i+1,\eps}^{-\left( \frac{\bp-\bm}{2}-\alpha\right)(\crits-1-\pe)}  B_\eps\right)
\end{eqnarray*} 
with
\begin{eqnarray*}
A_\eps&:=& \int \limits_{2 |\ze| \leq |x| <  k_{i+1,\eps}}  \frac{|x|^{\alpha(\crits-1-\pe)}\, dx}{|x|^{(\crits-\pe) (\bp-1)+s}}\\
&\leq &  \int \limits_{2 |\ze| \leq |x|<  k_{i+1,\eps}}  \frac{|x|^{\alpha(\crits-1-\pe)}\, dx}{|x|^{n+ (\crits-\pe)\left( \frac{\bp-\bm}{2}\right)}}~ dx
\end{eqnarray*}

\begin{eqnarray*}
B_\eps&:=&\int \limits_{2 |\ze| \leq |x|<  k_{i+1,\eps}}  \frac{dx}{|x|^{ \bp+s+\bm(\crits-1-\pe)+\alpha(\crits-1-\pe)}}\\
&\leq & \int \limits_{2 |\ze| \leq |x| <  k_{i+1,\eps}}  
\frac{ dx}{|x|^{ n- \left[(\crits-2-\pe)\left( \frac{\bp-\bm}{2}\right)-\alpha  (\crits-1-\pe)\right]}}
\end{eqnarray*}
Arguing as in the case $i=1$, with a change of variable we get that
\begin{eqnarray}\label{estim:1:I4:bis}
I^{i}_{4}&\leq& C\frac{\mu_{i,\eps}^{ \frac{\bp-\bm}{2}}d(\ze,\bdry)}{|\ze|^{\bp}}     \left( \frac{\mu_{i,\eps}}{|\ze|}\right)^{(\crits-2-\pe)\left( \frac{\bp-\bm}{2}\right)-\alpha  (\crits-1-\pe)}   \notag \\
&& + C\left(\frac{|\ze| }{\mu_{i+1,\eps}} \right)^{(\crits-2-\pe)\left( \frac{\bp-\bm}{2}\right)-\alpha  (\crits-1-\pe)}  \frac{ d(\ze,\bdry)}{|\ze|^{\bm}} \notag \\
&\leq& C\left(\frac{\mu_{i,\eps}^{\frac{\bp-\bm}{2}}|\ze|}{|\ze|^{\bp}}+\frac{|\ze| }{ \mu^{\frac{\bp-\bm}{2}}_{i+1,\eps} |\ze|^{\bm}}  \right). 
\end{eqnarray} 
Finally we get for the last integral from point (A5) of Proposition \ref{prop:exhaust} and a change of variable
\begin{align}\label{estim:1:I5:bis}
I^{i}_{5} :=&~C \int \limits_{D^{i}_{5,\eps}}\left(\frac{\max\{|\ze|,|x|\}}{\min\{|\ze|,|x|\}}\right)^{\bm} \frac{d(x,\bdry) d(\ze,\bdry)}{|x-\ze|^{n}} \frac{|\ue(x)|^{\crits-1-\pe}}{|x|^s}~ dx  \notag \\
\leq &~C~\frac{d(\ze,\bdry)}{|\ze|^{\bm}} \int \limits_{|x| \geq  k_{i+1,\eps}}  \frac{|\ue(x)|^{\crits-1-\pe}}{|x|^{\bp+s-1}}~ dx \notag \\
\leq &~C ~\frac{d(\ze,\bdry)}{|\ze|^{\bm}} \int \limits_{|x| \geq  k_{i+1,\eps}}  \frac{1}{|x|^{\bp+s -1+\frac{n-2}{2}(\crits-1-\pe) }}~ dx \notag\\
\leq &~ C~  \frac{d(\ze,\bdry)}{\mu_{i+1,\eps}^{\frac{\bp-\bm}{2}}|\ze|^{\bm}} \int \limits_{|x| \geq 1 }\frac{1}{|x|^{ n+\frac{\bp-\bm}{2}} } \notag\\
\leq &~  C~ \frac{|\ze|}{\mu_{i+1,\eps}^{\frac{\bp-\bm}{2}}|\ze|^{\bm}}.
\end{align}
\noindent
Then from \eqref{estim:1:I1:bis}, \eqref{estim:1:I2:bis}, \eqref{estim:1:I3:bis}, \eqref{estim:1:I4:bis} and \eqref{estim:1:I5:bis} we get the estimate \eqref{estim:1:bis:ze}. \hfill $\Box$
\bigskip

\noindent
Combining the estimates  \eqref{ineq:1}, \eqref{estim:1:ze},  \eqref{estim:1:tis:ze} and  \eqref{estim:1:bis:ze} we get that, there exists a constant $C>0$ such that  for any sequence of points $(\ze)$ in $\Omega$ we have 
\begin{align*}
|\ue(\ze)| \leq &~ C \left( \sum_{i=1}^N\frac{\mei^{\frac{\bp-\bm}{2}}|\ze|}{ \mei^{\bp-\bm}|\ze|^{\bm}+|\ze|^{\bp}}+ \frac{\Vert |x|^{\bm}u_0 \vert|_{L^{\infty}(\Omega)}}{|\ze|^{\bm}} |\ze| \right). \notag
\end{align*}
This completes the proof of  Proposition \ref{prop:fund:est}. \hfill $\Box$\\

\noindent
In our next  result we  obtain a pointwise control on the gradient. 
\begin{proposition}\label{prop:fund:est:grad}    
Let $\Omega$ be a smooth bounded domain of $\rn$, $n\geq 3$,  such that $0\in \bdry $ and  assume that  $0< s < 2$,   $\gamma<\frac{n^2}{4}$.  Let $(\ue)$, $(\he)$ and $(\pe)$ be such that $(E_\eps)$, \eqref{hyp:he},  \eqref{lim:pe} and \eqref{bnd:ue} holds. Assume that blow-up occurs, that is
$$\lim_{\eps\to 0}\Vert |x|^{\tau} \ue\Vert_{L^\infty(\Omega)}=+\infty ~\hbox{ where }  ~ \bm-1<\tau<\frac{n-2}{2} .$$ 
Consider the $\mu_{1,\eps},...,\mu_{N,\eps}$  from Proposition  \ref{prop:exhaust}.  Then there exists $C>0$ such that for all $\eps>0$ 
\begin{align}\label{eq:est:grad:global}
|\nabla \ue(x)|\leq
C \left(~\sum_{i=1}^N\frac{\mei^{\frac{\bp-\bm}{2}}  }{ \mei^{\bp-\bm}|x|^{\bm}+|x|^{\bp}}+\frac{ \Vert  |x|^{\bm-1}u_0 \Vert_{L^{\infty}(\Omega)} }{|x|^{\bm}}  \right)
\end{align}
for all $x\in \overline{\Omega} \setminus \{ 0\}$.
\end{proposition}
\noindent{\it Proof of Proposition \ref{prop:fund:est:grad}:} We let $G_{\epsilon}$ be the Green's function of the coercive operator  $-\Delta-\frac{\gamma}{|x|^{2}}-h_{\epsilon}$ on $\Omega$ with Dirichlet boundary condition. Differentiating the  Green's representation formula, and then using the pointwise bounds on the gradient  Green's function \eqref{ineq:grad:G} and the regularity result Theorem \ref{th:hopf} yields for any $ z \in \Omega$
\begin{align}\label{estim:grad:global:z}
\ue(z) =& \int \limits_\Omega G_{\epsilon}(z,x) \frac{|\ue (x)|^{\crits-2-\pe} ~\ue(x)}{|x|^s}dx \notag \\
|\nabla \ue(z)| \leq & ~C~\int \limits_\Omega \left|\nabla_{z}  G_{\epsilon}(z,x) \right|   \frac{|\ue (x)|^{\crits-1-\pe} }{|x|^s} dx \notag \\
\leq & ~C\int \limits_\Omega \left|\nabla_{z}  G_{\epsilon}(z,x) \right|  \frac{|\ue (x)|^{\crits-1-\pe} }{|x|^s}~dx\notag \\
\leq&~  \int \limits_\Omega  \left(\frac{\max\{|z|,|x|\}}{\min\{|z|,|x|\}}\right)^{\bm} \frac{d(x,\bdry) }{|x-z|^{n}} \frac{|\ue(x)|^{\crits-1-\pe}}{|x|^s}~ dx \notag.
\end{align}

\noindent
Then  using  the pointwise estimates \eqref{eq:est:global}  the proof goes exactly as in Proposition \ref{prop:fund:est}. 
\hfill $\Box$

\section{\, Sharp blow-up rates and the proof of Compactness}\label{sec:cpct}

The proof of compactness rely on the following two key propositions.
\begin{proposition}\label{prop:rate:sc}
Let $\Omega$ be a smooth bounded domain of $\rn$, $n\geq 3$,  such that $0\in \bdry $ and  assume that  $0< s < 2$,   $\gamma<\frac{n^2}{4}$.  Let $(\ue)$, $(\he)$ and $(\pe)$ be such that $(E_\eps)$, \eqref{hyp:he},  \eqref{lim:pe} and \eqref{bnd:ue} holds. Assume that blow-up occurs, that is
\begin{equation}\label{blowup}
\lim_{\eps\to 0}\Vert |x|^{\tau} \ue\Vert_{L^\infty(\Omega)}=+\infty ~\hbox{ for some }  ~ \bm-1<\tau<\frac{n-2}{2} .
\end{equation} 
Consider the $\mu_{1,\eps},...,\mu_{N,\eps}$  and $t_{1},...,t_{N}$  from Proposition \ref{prop:exhaust}. Suppose  that
\begin{equation}
\hbox{either }\{\,\bp-\bm>2\,\}\hbox{ or }\{\bp-\bm>1\hbox{ and }u_0\equiv 0\}.
\end{equation}
Then, we have following blow-up rates:
\begin{align}\label{rate:1}
\lim \limits_{\eps \to 0} \frac{\pe}{\mu_{N,\eps}}&~=c_{n,s,t_N}\frac{ \int \limits_{\partial\rnm} II_{0}(x,x)  |\nabla \tu_{N}|^2  ~d\sigma}{  \sum \limits_{i=1}^{N}  \frac{1}{t_{i}^{ \frac{n-2}{\crits-2}}} \int \limits_{\rnm }   \frac{|\tu_{i}|^{\crits }}{|x|^{s}} ~ dx }. \end{align}
Here $II_{0}$ denotes the second fundamental form of $\bdry$ at $0 \in \bdry$ and
\begin{equation*}
c_{n,s,t_N}:=\frac{n-s}{(n-2)^{2}}  \frac{1}{ t_{N}^{\frac{n-1}{\crits-2}}}.
\end{equation*}
\end{proposition}

\begin{proposition}[\textit{The positive case}]\label{prop:rate:sc:2}
Let $\Omega$ be a smooth bounded domain of $\rn$, $n\geq 3$,  such that $0\in \bdry $ and  assume that  $0< s < 2$,   $\gamma<\frac{n^2}{4}$.  
Let $(\ue)$, $(\he)$ and $(\pe)$ be as in Proposition \ref{prop:rate:sc} and let $H(0)$ denote the mean curvature of $\partial\Omega$ at $0$. Assume that blow-up occurs as in \eqref{blowup}. 
Consider the $\mu_{1,\eps},...,\mu_{N,\eps}$  and $t_{1},...,t_{N}$  from Proposition \ref{prop:exhaust}. Suppose in addition that
\begin{equation}
\ue>0 \qquad \hbox{for all } \eps >0.
\end{equation}
We define
\begin{equation}
C_{n,s,(t_i), (\tu)_i}:=\frac{c_{n,s,t_N} \int \limits_{\partial\rnm} |x|^2  |\nabla \tu_{N}|^2  ~d\sigma}{  (n-1)\sum \limits_{i=1}^{N}  \frac{1}{t_{i}^{ \frac{n-2}{\crits-2}}} \int \limits_{\rnm }   \frac{|\tu_{i}|^{\crits }}{|x|^{s}} ~ dx }
\end{equation}
Then,  we have the following blow-up rates:

\medskip\noindent 1) When $\bp-\bm\geq 2$, then 
\begin{align*}
\lim \limits_{\eps \to 0} \frac{\pe}{\mu_{N,\eps}}&~=C_{n,s,(t_i), (\tu)_i}\cdot H(0)\,\, \hbox{ if }\left\{\begin{array}{l}
\bp-\bm>2\\
\hbox{ or }\bp-\bm=2\hbox{ and }u_0\equiv 0
\end{array}\right\}. 
\end{align*}

\begin{align*}
\lim \limits_{\eps \to 0} \frac{\pe}{\mu_{N,\eps}}&~=C_{n,s,(t_i), (\tu)_i}\cdot H(0)-\tilde{K} \,\, \hbox{ if }\bp-\bm=2\hbox{ and }u_0>0. \end{align*}
for some $\tilde{K}>0$.  

\medskip\noindent 2) When $\bp-\bm<2$, then $u_0\equiv 0$ and

\begin{align*}
\lim \limits_{\eps \to 0} \frac{\pe}{\mu_{N,\eps}}&~=C_{n,s,(t_i), (\tu)_i}\cdot H(0)\quad \hbox{ if }\bp-\bm>1.\end{align*}

\begin{equation}\label{asymp:crit:log}
\lim \limits_{\eps \to 0} \frac{\pe}{\mu_{N,\eps}\ln\frac{1}{\mu_{\eps,N}}}=C'_{n,s,(t_i), (\tu)_i}\cdot H(0)\quad \hbox{ if }\bp-\bm=1
\end{equation}
\begin{equation}
\lim_{\eps\to 0}\frac{p_\eps}{\mu_{N,\eps}^{\bp-\bm}}=-\chi\cdot m_{\gamma,h}(\Omega)\quad \hbox{ if }\bp-\bm<1\label{est:mass}
\end{equation}
where
$$C'_{n,s,(t_i), (\tu)_i}:=\frac{n-s}{(n-2)^{2}}\frac{K^2\omega_{n-2}}{(n-1) \sum \limits_{i=1}^{N}  \frac{1}{t_{i}^{ \frac{n-2}{\crits-2}}} \int \limits_{\rnm }   \frac{|\tu_{i}|^{\crits }}{|x|^{s}} ~ dx}, $$
the constant $K$ is as in \eqref{def:K}, $\chi>0$ is a constant and $m_{\gamma,h}(\Omega)$ is the boundary mass defined in Theorem \ref{thdefi:mass}.
\end{proposition}

 \medskip\noindent {\it Proof of Theorems \ref{th:cpct:sc}, \ref{th:cpct:sc:2} and \ref{th:cpct:sc:3}:} We argue by contradiction and assume that the family is not pre-compact. Then, up to a subsequence, it blows up. We then apply Propositions \ref{prop:rate:sc} and \ref{prop:rate:sc:2} to get the blow-up rate (that is nonegative). However, the hypothesis of Theorems \ref{th:cpct:sc}, \ref{th:cpct:sc:2} and \ref{th:cpct:sc:3} yield exactly negative blow-up rates. This is a contradiction, and therefore the family is pre-compact. This proves the Theorems.\qed

\medskip\noindent We now establish Propositions \ref{prop:rate:sc} and \ref{prop:rate:sc:2}. The proof is divided in 13 steps in Sections \ref{loc pohoaev id} to \ref{pf blow-up rates}. These steps are numbered Steps P1, P2, etc. 

\section{\, Estimates on the localized Pohozaev identity}\label{loc pohoaev id}

\smallskip\noindent In the sequel, we let $(\ue)$, $(\he)$ and $(\pe)$ be such that $(E_\eps)$, \eqref{hyp:he},  \eqref{lim:pe} and \eqref{bnd:ue} hold. We assume that blow-up occurs. Note that 
$$\gamma<\frac{n^2}{4}-1\; \Leftrightarrow\;\bp-\bm>2,$$
and 
$$\gamma<\frac{n^2-1}{4}\; \Leftrightarrow\;\bp-\bm>1.$$
In the sequel, we will permanently use the following consequence of (A9) of Proposition \ref{prop:exhaust}: for all $i=1,...,N$, there exists $c_i>1$ such that
\begin{equation}\label{bnd:k:mu}
c_i^{-1} \mu_{\eps,i}\leq k_{\eps,i}\leq c_i \mu_{\eps,i}. 
\end{equation}

\begin{step}[Pohozaev identity]\label{step:p1} We let $(\ue)$, $(\he)$ and $(\pe)$ be such that $(E_\eps)$, \eqref{hyp:he},  \eqref{lim:pe} and \eqref{bnd:ue} hold. We assume that blow-up occurs. We define 
\begin{align}\label{def:Fe}
F_{\eps}(x):=&~( x, \nu) \left( \frac{|\nabla \ue|^2}{2}  -\frac{\gamma}{2}\frac{\ue^{2}}{|x|^{2}} - \frac{h_{\eps}(x)}{2} \ue^{2}  - \frac{1}{\crits -p_{\eps}} \frac{|\ue|^{\crits-p_{\eps} }}{|x|^{s}} \right)  \notag\\
& ~-  \left(x^i\partial_i \ue+\frac{n-2}{2} \ue \right)\partial_\nu \ue 
\end{align}
We let $\T$ be a chart at $0$ as in \eqref{def:T:bdry}. We define $\re:=\sqrt{\mu_{N,\eps}}$. Then
\begin{align}\label{PohoId3}
& \int \limits_{\T\left( \rnm \cap B_{ r_{\eps} }(0) \setminus B_{ k_{1,\eps}^{3} }(0) \right)} \left( h_{\eps}(x) +\frac{ \left( \nabla h_{\eps}, x \right)}{2} \right)\ue^{2}~ dx   \notag\\
 &  +\frac{p_{\eps}}{\crits} \left( \frac{n-s}{\crits -p_{\eps}}\right) \int \limits_{\T\left( \rnm \cap B_{ r_{\eps} }(0) \setminus B_{ k_{1,\eps}^{3} }(0) \right) }   \frac{|\ue|^{\crits-p_{\eps} }}{|x|^{s}} ~ dx \notag \\
 =~& - \int \limits_{ \T \left( \rnm \cap \partial  B_{ r_{\eps} }(0)  \right) } F_{\eps}(x)~d\sigma +\int \limits_{ \T \left( \rnm \cap \partial  B_{ k_{1,\eps}^{3} }(0)  \right) } F_{\eps}(x)~d\sigma \notag\\
& + \int \limits_{   \T\left( \partial \rnm \cap B_{ r_{\eps} }(0) \setminus B_{ k_{1,\eps}^{3} }(0) \right) }    ( x, \nu)  \frac{|\nabla \ue|^2}{2}   ~d\sigma
\end{align}
and, for $\delta_0>0$ small enough,
\begin{align}\label{PohoId3:bis}
& \int \limits_{\T\left( \rnm \cap B_{\delta_0 }(0) \setminus B_{ k_{1,\eps}^{3} }(0) \right)} \left( h_{\eps}(x) +\frac{ \left( \nabla h_{\eps}, x \right)}{2} \right)\ue^{2}~ dx   \notag\\
 & \qquad  \qquad +\frac{p_{\eps}}{\crits} \left( \frac{n-s}{\crits -p_{\eps}}\right) \int \limits_{\T\left( \rnm \cap B_{ \delta_0 }(0) \setminus B_{ k_{1,\eps}^{3} }(0) \right) }   \frac{|\ue|^{\crits-p_{\eps} }}{|x|^{s}} ~ dx \notag \\
 =~& -\int \limits_{ \T \left( \rnm \cap \partial  B_{ \delta_0 }(0)  \right) } F_{\eps}(x)~d\sigma + \int \limits_{ \T \left( \rnm \cap \partial  B_{ k_{1,\eps}^{3} }(0)  \right) } F_{\eps}(x)~d\sigma \notag\\
& + \int \limits_{   \T\left( \partial \rnm \cap B_{ \delta_0 }(0) \setminus B_{ k_{1,\eps}^{3} }(0) \right) }    ( x, \nu)  \frac{|\nabla \ue|^2}{2}   ~d\sigma
\end{align}
\end{step}

\medskip\noindent{\it Proof of Step P1:} We apply the Pohozaev identity \eqref{PohoId} with $y_{0}=0$ and
 $$U_{\eps}= \T\left( \rnm \cap B_{ r_{\eps} }(0) \setminus B_{ k_{1,\eps}^{3} }(0) \right) \subset \Omega.$$ 
This yields
\begin{align}\label{PohoId2}
  - \int \limits_{U_{\eps}} \left( h_{\eps}(x) +\frac{ \left( \nabla h_{\eps}, x \right)}{2} \right)\ue^{2}~ dx  
  &~ - \frac{p_{\eps}}{\crits} \left( \frac{n-s}{\crits -p_{\eps}}\right) \int \limits_{U_{\eps} }   \frac{|\ue|^{\crits-p_{\eps} }}{|x|^{s}} ~ dx \notag \\
 &~ =  \int \limits_{\partial U_{\eps}}  F_{\eps}(x)~d\sigma. 
\end{align}
It follows from the properties of the boundary map that 
\begin{align*}
 \partial U_{\eps} =&~\partial \left( \T \left( \rnm \cap B_{ r_{\eps} }(0) \setminus B_{ k_{1,\eps}^{3} }(0) \right) \right)&~  \notag\\
  = &~\T \left( \rnm \cap \partial B_{ r_{\eps} }(0)  \right) \cup   \T \left( \rnm \cap \partial B_{ k_{1,\eps}^{3} }(0)  \right) \cup  \T\left( \partial \rnm \cap B_{ r_{\eps} }(0) \setminus B_{ k_{1,\eps}^{3} }(0) \right) 
\end{align*}
Since for all $\eps >0$, $\ue \equiv 0$ on $\bdry$, identity \eqref{PohoId2} yields \eqref{PohoId3}. Concerning \eqref{PohoId3:bis}, we apply the Pohozaev identity \eqref{PohoId} with $y_{0}=0$ and
 $$V_{\eps}= \T\left( \rnm \cap B_{ \delta_0 }(0) \setminus B_{ k_{1,\eps}^{3} }(0) \right) \subset \Omega.$$
 The argument is similar. This ends the proof of Step \ref{step:p1}.\qed

\medskip\noindent We will estimate each of the terms in the above integral identities and calculate the limit as $\eps \to 0$.  \medskip

\subsection{Estimates of the $L^{\crits}$ and $L^2-$terms in the localized Pohozaev identity}

\begin{step}\label{step:p2} We let $(\ue)$, $(\he)$ and $(\pe)$ be such that $(E_\eps)$, \eqref{hyp:he},  \eqref{lim:pe} and \eqref{bnd:ue} hold. We assume that blow-up occurs. We claim that, as $\eps \to 0$

\begin{equation}\label{eq:mystery}
 \int \limits_{\T\left( \rnm  \cap B_{r_\eps }(0) \setminus B_{ k_{1,\eps}^3 }(0)\right)}
 \frac{|\ue|^{\crits-p_{\eps} }}{|x|^{s}} ~ dx =  \sum \limits_{i=1}^{N}  \frac{1}{t_{i}^{ \frac{n-2}{\crits-2}}} \int \limits_{\rnm }   \frac{|\tu_{i}|^{\crits }}{|x|^{s}} ~ dx +o(1).
\end{equation}
and

\begin{align}\label{cpct:2*term:2}
 \int \limits_{\T \left( \rnm  \cap B_{ \delta_0 }(0) \setminus B_{ k_{1,\eps}^3 }(0)  \right)}   \frac{|\ue|^{\crits-p_\eps }}{|x|^s} ~ dx=  \sum \limits_{i=1}^{N}  \frac{1}{t_i^{\frac{n-2}{\crits-2}}} \int \limits_{\rnm }   \frac{|\tu_{i}|^{\crits }}{|x|^{s}} ~ dx +o(1)\hbox{ if }u_0\equiv 0.
\end{align}
\end{step}

\noindent {\it Proof of Step \ref{step:p2}:} For any $R,\rho >0$ we decompose  the above integral as
\begin{align*}
 \int \limits_{   \T \left( \rnm \cap B_{r_{\eps} }(0) \setminus B_{k_{1,\eps}^{3} }(0) \right)}     \frac{|\ue|^{\crits-p_{\eps} }}{|x|^{s}} ~ dx&= 
 \int \limits_{   \T \left( \rnm \cap B_{ r_{\eps} }(0) \setminus \overline{B}_{R k_{N,\eps}}(0)  \right)}   \frac{|\ue|^{\crits-p_{\eps} }}{|x|^{s}} ~ dx \notag\\
 &+ \sum  \limits_{i=1}^{N}   \int \limits_{  \T \left( \rnm \cap B_{R k_{i,\eps}}(0)  \setminus \overline{B}_{\rho k_{i,\eps}}(0)  \right)}     \frac{|\ue|^{\crits-p_{\eps} }}{|x|^{s}} ~ dx  \notag\\
&+ \sum  \limits_{i=1}^{N-1}   \int \limits_{  \T \left( \rnm \cap B_{\rho k_{i+1,\eps}}(0)  \setminus \overline{B}_{R k_{i,\eps}}(0)  \right)}     \frac{|\ue|^{\crits-p_{\eps} }}{|x|^{s}} ~ dx  \notag\\
&  +  
 \int \limits_{   \T \left( \rnm \cap B_{\rho k_{1,\eps}}(0) \setminus B_{ k_{1,\eps}^{3} }(0)  \right) }     \frac{|\ue|^{\crits-p_{\eps} }}{|x|^{s}} ~ dx. 
\end{align*}
We will evaluate each of the above terms and calculate the limit $\lim\limits_{R \to + \infty} \lim\limits_{\rho \to 0} \lim\limits_{\eps \to 0}$. 
\medskip

\noindent
From the estimate \eqref{eq:est:global}, we  get as $\eps \to 0$
\begin{align*}
 &\int \limits_{   \T \left( \rnm \cap B_{ r_{\eps} }(0) \setminus \overline{B}_{R k_{N,\eps}}(0)  \right)}   \frac{|\ue|^{\crits-p_{\eps} }}{|x|^{s}} ~ dx  \notag\\
\leq &~  C~\int \limits_{   \T \left( \rnm \cap B_{ r_{\eps} }(0) \setminus \overline{B}_{R k_{N,\eps}}(0)  \right)} \left[  \frac{ \mu_{N,\eps}^{\frac{\bp-\bm}{2}(\crits-\pe)}}{|x|^{(\bp-1)(\crits-\pe)+s}}+ 
\frac{1}{ |x|^{(\bm-1)(\crits-\pe)+s}} \right] ~dx\notag\\
\leq &~  C~\int \limits_{ \rnm \cap  B_{ r_{\eps}}(0)  \setminus \overline{B}_{{R k_{N,\eps}}}(0)  }   \frac{ \mu_{N,\eps}^{\frac{\bp-\bm}{2}(\crits-\pe)}}{|x|^{(\bp-1)(\crits-\pe)+s}} \left| \hbox{Jac } \T(x) \right|~dx \notag \\
+ &~  C~\int \limits_{ \rnm \cap  B_{ r_{\eps}}(0)  \setminus \overline{B}_{{R k_{N,\eps}}}(0)  }    \frac{1}{ |x|^{(\bm-1)(\crits-\pe)+s}}\left| \hbox{Jac } \T(x) \right|~dx \notag \\
\leq &~  C~ \int \limits_{   \rnm \cap  B_{ \frac{ r_{\eps}}{k_{N,\eps}}}(0)  \setminus \overline{B}_{R }(0)  }    \frac{ 1}{|x|^{n+ \crits \left( \frac{\bp-\bm}{2}\right)-\pe(\bp-1)}}\left| \hbox{Jac } \T(k_{N,\eps}x) \right|~dx \notag\\
+ &~  C~\int \limits_{ \rnm \cap  B_{ 1}(0)  \setminus \overline{B}_{\frac{R k_{N,\eps}}{r_{\eps}}}(0)  }    \frac{ 1 }{ |x|^{n-\crits \left( \frac{\bp-\bm}{2}\right)-\pe(\bm-1)}} \left| \hbox{Jac } \T( r_{\eps} x) \right|~dx \notag \\
 \leq &C\left( R^{-\crits \left( \frac{\bp-\bm}{2}\right)-\pe(\bp-1)}  + r_{\eps}^{\crits \left( \frac{\bp-\bm}{2}\right)+\pe(\bm-1)}\right)  .
\end{align*}
Therefore 
\begin{align}{\label{cpct:2*term:1}}
\lim\limits_{R \to + \infty}  \lim\limits_{\eps \to 0}\int \limits_{   \T \left( \rnm \cap B_{ r_{\eps} }(0) \setminus \overline{B}_{R k_{N,\eps}}(0)  \right)}   \frac{|\ue|^{\crits-p_{\eps} }}{|x|^{s}} ~ dx =0.
\end{align}
\medskip

\noindent
It follows from Proposition \ref{prop:exhaust} that for any $1 \leq i\leq N$
\begin{align}{\label{cpct:2*term:2:lala}}
\lim\limits_{R \to + \infty} \lim\limits_{\rho \to 0} \lim\limits_{\eps \to 0}  \int \limits_{  \T \left( \rnm \cap B_{R k_{i,\eps}}(0)  \setminus \overline{B}_{\rho k_{i,\eps}}(0)  \right)}     \frac{|\ue|^{\crits-p_{\eps} }}{|x|^{s}} ~ dx 
= \frac{1}{t_{i}^{ \frac{n-2}{\crits-2}}} \int \limits_{\rnm }   \frac{|\tu_{i}|^{\crits }}{|x|^{s}} ~ dx.
\end{align}
\medskip

\noindent
Let $1 \leq i\leq N-1$. In Proposition \ref{prop:fund:est}, we had obtained the following pointwise estimates: For any $R, \rho>0$ and all $\eps>0$ we have 
\begin{align*}
|\ue(x)| \leq  C~ \frac{\mei^{\frac{\bp-\bm}{2}} |x|}{|x|^{\bp}}+C~ \frac{ |x| }{ \mu_{ i+1,\eps}^{\frac{\bp-\bm}{2}}|x|^{\bm}}
\end{align*}
for all $x \in B_{\rho k_{i+1,\eps}}(0)  \setminus \overline{B}_{R k_{i,\eps}}(0) $.
\medskip

\noindent
Then we have  as $\eps \to 0$
\begin{align*}
 &\int \limits_{  \T \left( \rnm \cap B_{\rho k_{i+1,\eps}}(0)  \setminus \overline{B}_{R k_{i,\eps}}(0)  \right)}     \frac{|\ue|^{\crits-p_{\eps} }}{|x|^{s}} ~ dx  \notag\\
 \leq &~  C~\int \limits_{  \T \left( \rnm \cap B_{\rho k_{i+1,\eps}}(0)  \setminus \overline{B}_{R k_{i,\eps}}(0)  \right)}   \left[   \frac{ \mei^{\frac{\bp-\bm}{2}(\crits-\pe)}}{|x|^{(\bp-1)(\crits-\pe)+s}}+ \frac{ \mu_{i+1,\eps}^{-\frac{\bp-\bm}{2}(\crits-\pe)} }{ |x|^{(\bm-1)(\crits-\pe)+s}}  \right]~dx \notag\\
  \leq &~  C~\int \limits_{ \rnm \cap  B_{ \rho k_{i+1,\eps}}(0)  \setminus \overline{B}_{R k_{i,\eps}}(0)  }  \left[   \frac{ \mei^{\frac{\bp-\bm}{2}(\crits-\pe)}}{|x|^{(\bp-1)(\crits-\pe)+s}}+ \frac{ \mu_{i+1,\eps}^{-\frac{\bp-\bm}{2}(\crits-\pe)} }{ |x|^{(\bm-1)(\crits-\pe)+s}}  \right] \,dx \notag \\
\leq &~  C~ \int \limits_{   \rnm \cap  B_{ \frac{\rho k_{i+1,\eps}}{k_{i,\eps}}}(0)  \setminus \overline{B}_{R }(0)  }    \frac{ 1}{|x|^{n+ \crits \left( \frac{\bp-\bm}{2}\right)-\pe(\bp-1)}}\,dx \notag\\
 &+ ~  C~\int \limits_{  \rnm \cap  B_{2\rho }(0)  \setminus \overline{B}_{ \frac{R k_{i,\eps}}{k_{i+1,\eps}}}(0)  }    \frac{ 1 }{ |x|^{n-\crits \left( \frac{\bp-\bm}{2}\right)-\pe(\bm-1)}} \, dx \notag\\
 \leq &C\left( R^{-\crits \left( \frac{\bp-\bm}{2}\right)-\pe(\bp-1)}  + \rho^{\crits \left( \frac{\bp-\bm}{2}\right)+\pe(\bm-1)}\right).
\end{align*}
And so 
\begin{align}{\label{cpct:2*term:3}}
\lim\limits_{R \to + \infty} \lim\limits_{\rho \to 0} \lim\limits_{\eps \to 0} \int \limits_{  \T \left( \rnm \cap B_{\rho k_{i+1,\eps}}(0)  \setminus \overline{B}_{R k_{i,\eps}}(0)  \right)}     \frac{|\ue|^{\crits-p_{\eps} }}{|x|^{s}} ~ dx =0.
\end{align}
\medskip

\noindent
Again, from the pointwise  estimates of  Proposition \ref{prop:fund:est}, we have as $\eps \to 0$
\begin{align*}
&\int \limits_{   \T \left( \rnm \cap B_{\rho k_{1,\eps}}(0) \setminus B_{ k_{1,\eps}^{3} }(0)  \right) }     \frac{|\ue|^{\crits-p_{\eps} }}{|x|^{s}} ~ dx \notag\\
\leq &~  C~\int \limits_{   \T \left( \rnm \cap B_{\rho k_{1,\eps}}(0) \setminus B_{ k_{1,\eps}^{3} }(0)  \right) }   \frac{ \mu_{1,\eps}^{-\frac{\bp-\bm}{2}(\crits-\pe)} }{ |x|^{(\bm-1)(\crits-\pe)+s}}  ~dx \notag\\
 \leq &~  C~\int \limits_{ \rnm \cap  B_{\rho k_{1,\eps}}(0)  \setminus \overline{B}_{ k_{1,\eps}^{3}}(0)  }  \frac{ \mu_{1,\eps}^{-\frac{\bp-\bm}{2}(\crits-\pe)} }{ |x|^{(\bm-1)(\crits-\pe)+s}}  \left| \hbox{Jac } \T(x) \right|~dx \notag \\
\leq& ~  C~\int \limits_{  \rnm \cap  B_{\rho }(0)  \setminus \overline{B}_{ k_{1,\eps}^{2}}(0)  }    \frac{ 1 }{ |x|^{n-\crits \left( \frac{\bp-\bm}{2}\right)-\pe(\bm-1)}} \left| \hbox{Jac } \T(k_{1,\eps}x) \right|~dx \notag\\
 \leq &C   ~\rho^{\crits \left( \frac{\bp-\bm}{2}\right)+\pe(\bm-1)}.
\end{align*}
Therefore
\begin{align}{\label{cpct:2*term:4}}
 \lim\limits_{\rho \to 0} \lim\limits_{\eps \to 0}   \int \limits_{   \T \left( \rnm \cap B_{\rho k_{1,\eps}}(0) \setminus B_{ k_{1,\eps}^{3} }(0)  \right) }     \frac{|\ue|^{\crits-p_{\eps} }}{|x|^{s}}~dx=0.
\end{align}
\medskip

\noindent
Combining  \eqref{cpct:2*term:1}, \eqref{cpct:2*term:2:lala}, \eqref{cpct:2*term:3} and \eqref{cpct:2*term:4} we obtain \eqref{eq:mystery}.

\medskip\noindent We now prove \eqref{cpct:2*term:2} under the assumption that $u_0\equiv 0$. We decompose the integral as
\begin{align*}
 \int \limits_{   \T \left( \rnm \cap B_{\delta_0 }(0) \setminus B_{k_{1,\eps}^{3} }(0) \right)}     \frac{|\ue|^{\crits-p_{\eps} }}{|x|^{s}} ~ dx&=
 \int \limits_{   \T \left( \rnm \cap B_{ \delta_0 }(0) \setminus \overline{B}_{r_\eps}(0)  \right)}   \frac{|\ue|^{\crits-p_{\eps} }}{|x|^{s}} ~ dx\notag\\
 & \qquad +  \int \limits_{   \T \left( \rnm \cap B_{r_\eps }(0) \setminus B_{k_{1,\eps}^{3} }(0) \right)}     \frac{|\ue|^{\crits-p_{\eps} }}{|x|^{s}} ~ dx,  
\end{align*}
with $r_\eps:=\sqrt{\mu_{N,\eps}}$. From the estimate \eqref{eq:est:global} and $u_0\equiv 0$, we  get as $\eps \to 0$
\begin{equation*}
 \int \limits_{   \T \left( \rnm \cap B_{ \delta_0 }(0) \setminus \overline{B}_{r_\eps}(0)  \right)}   \frac{|\ue|^{\crits-p_{\eps} }}{|x|^{s}} ~ dx 
\leq  C~\int \limits_{   \rnm \cap B_{ \delta_0 }(0) \setminus \overline{B}_{r_{\eps}}(0) } \left[  \frac{ \mu_{N,\eps}^{\frac{\bp-\bm}{2}(\crits-\pe)}}{|x|^{(\bp-1)(\crits-\pe)+s}} \right] ~dx
\end{equation*}
Since $(\bp-1)\crits+s>n$, we then get that
$$ \int \limits_{   \T \left( \rnm \cap B_{ \delta_0 }(0) \setminus \overline{B}_{r_\eps}(0)  \right)}   \frac{|\ue|^{\crits-p_{\eps} }}{|x|^{s}} ~ dx \leq C\left(\frac{\mu_{N,\eps}}{r_\eps}\right)^{\frac{\crits}{2}(\bp-\bm)}=o(1)$$
as $\eps\to 0$. Therefore, with \eqref{cpct:2*term:1}, we get \eqref{cpct:2*term:2}. This proves \eqref{cpct:2*term:2}.

\smallskip\noindent This ends the proof of Step \ref{step:p2}.\qed

\begin{step}\label{step:p3} We let $(\ue)$, $(\he)$ and $(\pe)$ be such that $(E_\eps)$, \eqref{hyp:he},  \eqref{lim:pe} and \eqref{bnd:ue} hold. We assume that blow-up occurs. We claim that

\begin{eqnarray}
 &&\int \limits_{\T \left( \rnm  \cap B_{ r_{\eps} }(0) \setminus B_{ k_{1,\eps}^{3} }(0)  \right)}  \left( h_{\eps}(x) +\frac{ \left( \nabla h_{\eps}, x \right)}{2} \right)\ue^{2}~ dx \nonumber\\
 &&=
 \left\{ \begin{array}{llll}
   O(\mu_{N,\eps}^2) ~ & \hbox{ if } \bp-\bm > 2\\ 
  O(\mu_{N,\eps}^2\ln \frac{1}{\mu_{N,\eps}}) ~ & \hbox{ if } \bp-\bm = 2\\ 
    O(\mu_{N,\eps}^{1+\frac{\bp-\bm}{2}}) ~ &\hbox{ if } \bp-\bm < 2.
\end{array} \right. \label{cpct:L2term:sc} 
\end{eqnarray}
And if $u_{0}\equiv 0$
\begin{eqnarray}
&& \int \limits_{\T \left( \rnm  \cap B_{ \delta_{0} }(0) \setminus B_{ k_{1,\eps}^{3} }(0)  \right)}  \left( h_{\eps}(x) +\frac{ \left( \nabla h_{\eps}, x \right)}{2} \right)\ue^{2}~ dx \nonumber\\
 &&=
 \left\{ \begin{array}{llll}
   O(\mu_{N,\eps}^2) ~ & \hbox{ if } \bp-\bm > 2\\ 
  O(\mu_{N,\eps}^2\ln \frac{1}{\mu_{N,\eps}}) ~ & \hbox{ if } \bp-\bm = 2\\ 
    O(\mu_{N,\eps}^{\bp-\bm}) ~ &\hbox{ if } \bp-\bm < 2.
\end{array} \right. \label{cpct:L2term:u=0} 
\end{eqnarray}

\end{step}
\noindent{\it Proof of Step \ref{step:p3}:}
From estimate \eqref{eq:est:global} and after a change of variables, we  get as $\eps \to 0$,
\begin{eqnarray}
 && \int \limits_{\T \left( \rnm  \cap B_{ r_{\eps} }(0) \setminus B_{ k_{1,\eps}^{3} }(0)  \right)}  \left( h_{\eps}(x) +\frac{ \left( \nabla h_{\eps}, x \right)}{2} \right)\ue^{2}~ dx \hfil \notag \\
 && \leq  C \int \limits_{\T \left( \rnm  \cap B_{ r_{\eps} }(0) \setminus B_{ k_{1,\eps}^{3} }(0)  \right)}  \ue^{2}~ dx   \notag\\
&& \hfil \leq ~  C~\int \limits_{\T \left( \rnm  \cap B_{ r_{\eps} }(0) \setminus B_{ k_{1,\eps}^{3} }(0)  \right)}   \left[ \frac{ \mu_{N,\eps}^{\bp-\bm}}{|x|^{2(\bp-1)}}~dx + 
   \frac{ 1}{|x|^{2(\bm-1)}}~dx \right] \label{ineq:2018}\\\
 && \leq ~  C~ \int \limits_{ \rnm \cap  B_{ r_{\eps}}(0)  \setminus \overline{B}_{{R k_{1,\eps}^{3}}}(0)  }  \left(\sum_{i=1}^N\frac{\mu_{i,\eps}^{\bp-\bm}|x|^2}{\mu_{i,\eps}^{2(\bp-\bm)}|x|^{2\bm}+|x|^{2\bp}}\right.\notag\\
 &&\hskip4cm \left.+\frac{ 1}{|x|^{2(\bm-1)}}\right)\, dx. \notag
 \end{eqnarray}
\noindent{\it Case 1:} Assuming that $\bp-\bm<2$, we then have the following rough bound from \eqref{ineq:2018},
 \begin{eqnarray}
&& \int_{\T \left( \rnm  \cap B_{ r_{\eps} }(0) \setminus B_{ k_{1,\eps}^{3} }(0)  \right)}  \ue^{2}\, dx\nonumber\\
&&\leq C\int \limits_{ \rnm \cap  B_{ r_{\eps}}(0)  \setminus \overline{B}_{{R k_{1,\eps}^{3}}}(0)  } \left(\frac{ \mu_{N,\eps}^{\bp-\bm}}{|x|^{2(\bp-1)}} + \frac{ 1}{|x|^{2(\bm-1)}}\right)\, dx\nonumber\\
 &&\leq C \mu_{N,\eps}^{1+\frac{\bp-\bm}{2}}\hbox{ if }\bp-\bm<2.\label{ineq:111}
 \end{eqnarray}
 
 \medskip\noindent{\it Case 2:} Assuming $\bp-\bm\geq 2$, then via a change of variable in \eqref{ineq:2018}, we get
 \begin{eqnarray*}
 \int \limits_{\T \left( \rnm  \cap B_{ r_{\eps} }(0) \setminus B_{ k_{1,\eps}^{3} }(0)  \right)}  \ue^{2}~ dx 
&&\leq C\sum_{i=1}^N\mu_{i,\eps}^2  \int_{  B_{ \frac{r_{\eps}}{\mu_{i,\eps}}}(0)  \setminus \overline{B}_{ \frac{k_{1,\eps}^3}{\mu_{i,\eps}}}(0) }\frac{|x|^2\, dx}{|x|^{2\bm}+|x|^{2\bp}}\notag \\
&&+ C\int_{B_{ r_{\eps} }(0) \setminus B_{ k_{1,\eps}^{3} }(0) }|x|^{2-2\bm}\, dx.
 \end{eqnarray*}
 Therefore, if $\bp-\bm>2$, then
  \begin{eqnarray}
&& \int \limits_{\T \left( \rnm  \cap B_{ r_{\eps} }(0) \setminus B_{ k_{1,\eps}^{3} }(0)  \right)}  \ue^{2}~ dx \leq C\sum_{i=1}^N\mu_{i,\eps}^2  + Cr_\eps^{n+2-2\bm}\leq C\mu_{N,\eps}^2.\label{ineq:222}
 \end{eqnarray}
When $\bp-\bm=2$, we get that
 \begin{eqnarray*}
 \int \limits_{\T \left( \rnm  \cap B_{ r_{\eps} }(0) \setminus B_{ k_{1,\eps}^{3} }(0)  \right)}  \ue^{2}~ dx  &&\leq C\sum_{i=1}^N\mu_{i,\eps}^2  \left(1+\int_{  B_{ \frac{r_{\eps}}{\mu_{i,\eps}}}(0)  \setminus \overline{B}_1(0) }|x|^{2-\bp}\, dx\right)\\
 &&+ Cr_\eps^{2+\bp-\bm}\\
&&\leq C\mu_{N,\eps}^2\ln\frac{1}{\mu_{N,\eps}}+C\sum_{i=1}^{N-1}\mu_{i,\eps}^2\ln\frac{1}{\mu_{i,\eps}}.
 \end{eqnarray*}
Since $\mu_{N,\eps}\to 0$ and $\lim_{\eps\to 0}\mu_{i,\eps}/\mu_{N,\eps}$ is finite for all $i=1,...,N-1$, we get that
\begin{equation}
\int \limits_{\T \left( \rnm  \cap B_{ r_{\eps} }(0) \setminus B_{ k_{1,\eps}^{3} }(0)  \right)}  \ue^{2}\, dx=O\left(\mu_{N,\eps}^2\ln\frac{1}{\mu_{N,\eps}}\right),\label{ineq:333}
\end{equation}
since $\bp-\bm=2$.
Inequality \eqref{ineq:2018} put together with \eqref{ineq:111}, \eqref{ineq:222} and \eqref{ineq:333} yield \eqref{cpct:L2term:sc}. 

\medskip\noindent When $u_0\equiv 0$ we decompose the integral and proceed as in the proof of \eqref{cpct:2*term:2} to obtain  \eqref{cpct:L2term:u=0}. This ends Step \ref{step:p3}.\qed \bigskip

\subsection{Estimate of the curvature term in the Pohozaev identity when $\bp-\bm>1$}

\begin{step}\label{step:p4} We let $(\ue)$, $(\he)$ and $(\pe)$ be such that $(E_\eps)$, \eqref{hyp:he},  \eqref{lim:pe} and \eqref{bnd:ue} hold. We assume that blow-up occurs and that $\bp -\bm >1$. We claim that, as $\eps \to 0$
\begin{eqnarray}
&&\int \limits_{   \T\left( \partial \rnm \cap B_{ r_{\eps} }(0) \setminus B_{ k_{1,\eps}^{3} }(0) \right) }    ( x, \nu)  \frac{|\nabla \ue|^2}{2}   ~d\sigma\nonumber\\
&&=  \frac{\mu_{N,\eps}}{2} \left(\frac{1}{ t_{N}^{\frac{n-1}{\crits-2}}} \int \limits_{\partial\rnm} II_{0}(x,x)  \frac{|\nabla \tu_{N}|^2}{2}  ~d\sigma +o(1) \right).\label{cpct:bdry:term:C}
\end{eqnarray}
Here, see Proposition \ref{prop:rate:sc}, $II_0$ denotes the second fundamental form. Moreover, when $u_0\equiv 0$, we claim that as $\eps\to 0$, 
\begin{eqnarray}
&&\int \limits_{   \T\left( \partial \rnm \cap B_{ \delta_0 }(0) \setminus B_{ k_{1,\eps}^{3} }(0) \right) }    ( x, \nu)  \frac{|\nabla \ue|^2}{2}   ~d\sigma\nonumber\\
&&=  \frac{\mu_{N,\eps}}{2} \left(\frac{1}{ t_{N}^{\frac{n-1}{\crits-2}}} \int \limits_{\partial\rnm} II_{0}(x,x)  \frac{|\nabla \tu_{N}|^2}{2}  ~d\sigma +o(1) \right).\label{cpct:bdry:term:C:2}
\end{eqnarray}

\end{step}
\noindent {\it Proof of Step \ref{step:p4}:} We have  for any  $R, \rho >0$,
\begin{eqnarray}\label{decomp:bdry:term:C}
&&\int \limits_{   \T\left( \partial \rnm \cap B_{ r_{\eps} }(0) \setminus B_{ k_{1,\eps}^{3} }(0) \right) }     ( x, \nu)  \frac{|\nabla \ue|^2}{2}   ~d\sigma\\
&&= 
  \int \limits_{   \T \left( \partial \rnm \cap B_{ r_{\eps} }(0) \setminus \overline{B}_{R k_{N,\eps}}(0)  \right)}  ( x, \nu)  \frac{|\nabla \ue|^2}{2}   ~d\sigma  \notag\\
  &&+  \sum  \limits_{i=1}^{N} \int \limits_{ \T \left(  \partial \rnm \cap B_{R k_{i,\eps}}(0)  \setminus \overline{B}_{\rho k_{i,\eps}}(0)  \right)}  ( x, \nu)  \frac{|\nabla \ue|^2}{2}   ~d\sigma\notag\\
  &&   + \sum  \limits_{i=1}^{N-1}   \int \limits_{ \T \left(  \partial \rnm \cap B_{\rho k_{i+1,\eps}}(0)  \setminus \overline{B}_{R k_{i,\eps}}(0)  \right)}   ( x, \nu)  \frac{|\nabla \ue|^2}{2}   ~d\sigma  \notag\\
  &  &+  
 \int \limits_{ \T \left( \partial \rnm \cap B_{\rho k_{1,\eps}}(0) \setminus B_{ k_{1,\eps}^{3} }(0)  \right)}  ( x, \nu)  \frac{|\nabla \ue|^2}{2}   ~d\sigma .
\end{eqnarray}
\noindent
We consider the second fundamental form associated to $\partial\Omega$, 
$II_0(x,y)=(d\nu_px,y)$
for  $0\in\partial\Omega$ and all $x,y\in T_{0}\partial\Omega$ ($\nu$ is the outward normal vector at the hypersurface $\partial\Omega$). In the canonical basis of $\partial\rnm=T_0\partial\Omega$, the matrix of the bilinear form $II_{0}$ is $-D^2_0\T_0$, where $D^2_0\T_0$ is the Hessian matrix of $\T_0$ at $0$. Using the expression of $\T$ (see \eqref{def:T:bdry}), we can write  for all $x \in U \cap \partial \rnm$
$$\nu(\T(x))=\frac{(1,-\partial_2\T_0(x), ...,-\partial_n\T_0(x))}{\sqrt{1+\sum_{i=2}^n(\partial_i\T_0(x))^2}}.$$
 With the expression of $\T$, we then
get that
\begin{align*}
\left(\nu\circ \T(x),\T(x) \right)= \frac{\T_0(x)-\sum_{p=2}^n x^p\partial_p\T_0(x)}{\sqrt{1+\sum_{p=2}^n(\partial_p\T_0(x))^2}}
\end{align*}
And so for  all $x\in U\cap\partial \rnm$.
\begin{align}\label{est:T:2}
|(\T (x),\nu\circ\T(x))|\leq C|x|^2
\end{align}
Since $\T_0(0)=0$ and $\nabla\T_0(0)=0$ (see \eqref{def:T:bdry}), we then get as $|x| \to 0$
\begin{align}\label{est:T:3}
\left(\nu\circ \T(x),\T(x) \right)=-\frac{1}{2}\sum \limits_{p,q=2}^n x^px^q \partial_{pq} \T_0(0)+O(|x|^{3})  
\end{align}
and therefore for all $\eps>0$ and all $x\in B_R(0) \cap \partial \rnm $ 
\begin{align}\label{est:T:4}
\left(\T(k_{N,\eps}x ), \nu \ \circ  T(k_{N,\eps} x)\right)&~=-\frac{1}{2}\keN^2 \sum_{p,q=2}^n x^px^q\partial_{pq}\T_0(0)+\theta_{\eps,R}(x) \keN^2 \notag\\
&~= \frac{1}{2}\keN^2 II_{0}(x,x)+\theta_{\eps,R}(x) \keN^2
\end{align}
where $\lim \limits_{\eps\to 0}\sup \limits_{B_R(0)\cap\{x_1=0\}}|\theta_{\eps,R}|=0$ for any $R>0$.

\medskip\noindent{\it Step \ref{step:p4}.1:}  
Let $1 \leq i\leq N-1$. In Proposition \ref{prop:fund:est:grad} we have obtained the pointwise estimates, that for any $R, \rho>0$ and all $\eps>0$ we have for all $x \in B_{\rho k_{i+1,\eps}}(0)  \setminus \overline{B}_{R k_{i,\eps}}(0) $, 
\begin{align*}
|\nabla \ue(x)| \leq  C~ \frac{\mei^{\frac{\bp-\bm}{2}} }{|x|^{\bp}}+C~ \frac{ 1 }{ \mu_{ i+1,\eps}^{\frac{\bp-\bm}{2}}|x|^{\bm}}.
\end{align*}
 \noindent
For clearness, we write in this step
$$D_\eps:= \T \left(  \partial \rnm \cap B_{\rho k_{i+1,\eps}}(0)  \setminus \overline{B}_{R k_{i,\eps}}(0)  \right).$$
As $\eps \to 0$ we have that
\begin{eqnarray*}
 \left|\int \limits_{D_\eps}     ( x, \nu)  \frac{|\nabla \ue|^2}{2}  ~d\sigma\right|& \leq &~  C~\int \limits_{D_\eps}    (x, \nu)  \left[   \frac{ \mei^{\bp-\bm}}{|x|^{2\bp}}+ \frac{ 1 }{ \mu_{ i+1,\eps}^{\bp-\bm}|x|^{2\bm}}  \right]~d\sigma \notag\\
 & \leq &~  C~\int \limits_{ D_\eps}    |x|^{2} \left[   \frac{ \mei^{\bp-\bm}}{|x|^{2\bp}}+ \frac{ 1 }{ \mu_{ i+1,\eps}^{\bp-\bm}|x|^{2\bm}}  \right]~d\sigma \notag\\
& \leq &~  C~\mu_{i,\eps} \int \limits_{  \partial \rnm \cap  B_{ \frac{\rho k_{i+1,\eps}}{k_{i,\eps}}}(0)  \setminus \overline{B}_{R }(0)  }    \frac{ d\sigma}{|x|^{(n-1)+(\bp-\bm-1)}} \notag\\
 &&+ ~  C~\mu_{i+1,\eps}\int \limits_{  \partial \rnm \cap  B_{\rho }(0)  \setminus \overline{B}_{ \frac{R k_{i,\eps}}{k_{i+1,\eps}}}(0)  }    \frac{ d\sigma }{ |x|^{(n-1)-(\bp-\bm+1)}}  \notag\\
 & \leq& C~\left( \mu_{i,\eps} R^{-(\bp-\bm-1)}+ \mu_{i+1,\eps} \rho^{\bp-\bm+1} \right).
\end{eqnarray*}
So then  for all  $1 \leq i\leq N-1$
\begin{align}\label{cpct:bdry:term:C:3}
\lim\limits_{R \to + \infty} \lim\limits_{\rho \to 0} \lim\limits_{\eps \to 0} \left(  \mu_{N,\eps}^{-1} \int \limits_{ \T \left(  \partial \rnm \cap B_{\rho k_{i+1,\eps}}(0)  \setminus \overline{B}_{R k_{i,\eps}}(0)  \right)}     ( x, \nu)  \frac{|\nabla \ue|^2}{2}  ~d\sigma \right)=0.
\end{align}
This ends Step \ref{step:p4}.1.

\medskip\noindent{\it Step \ref{step:p4}.2:} Again from the estimates of Proposition \ref{prop:fund:est:grad}, we have as $\eps \to 0$
\begin{eqnarray*}
&& \left|\int \limits_{ \T \left( \partial \rnm \cap B_{\rho k_{1,\eps}}(0) \setminus B_{ k_{1,\eps}^{3} }(0)  \right)}    ( x, \nu)  \frac{|\nabla \ue|^2}{2}   ~d\sigma \right|
 \notag \\
 &&  \leq C \int \limits_{ \T \left( \partial \rnm \cap B_{\rho k_{1,\eps}}(0) \setminus B_{ k_{1,\eps}^{3} }(0)  \right)}      \frac{  |( x, \nu)|\,d\sigma}{ \mu_{ 1,\eps}^{\bp-\bm}|x|^{2\bm}}     \notag\\
&&  \leq ~C~   \int \limits_{  \partial \rnm \cap B_{\rho k_{1,\eps}}(0) \setminus B_{ k_{1,\eps}^{3} }(0) }     \frac{  |x|^{2}\,d\sigma}{ \mu_{ 1,\eps}^{\bp-\bm}|x|^{2\bm}}  \notag\\ 
&&   \leq ~  C~k_{1,\eps}\int \limits_{  \partial \rnm \cap  B_{\rho}(0) \setminus B_{ k^{2}_{1,\eps} }(0) }    \frac{ d\sigma }{ |x|^{2\bm-2}}   \leq ~C~\mu_{1,\eps}~ \rho^{\bp-\bm+1}. \notag
\end{eqnarray*}
Then, using again \eqref{bnd:k:mu}, we get that
\begin{align}\label{cpct:bdry:term:C:4}
\lim\limits_{R \to + \infty} \lim\limits_{\rho \to 0} \lim\limits_{\eps \to 0} \left( \mu_{N,\eps}^{-1} \int \limits_{ \T \left( \partial \rnm \cap B_{\rho k_{1,\eps}}(0) \setminus B_{ k_{1,\eps}^{3} }(0)  \right)}    ( x, \nu)  \frac{|\nabla \ue|^2}{2}   ~d\sigma \right)=0.
\end{align}
This ends Step \ref{step:p4}.2.

\medskip\noindent{\it Step \ref{step:p4}.3:}  With the pointwise estimates of Proposition \ref{prop:fund:est:grad}, we obtain as $\eps \to 0$
\begin{eqnarray*}
&&\int \limits_{   \T \left( \partial \rnm \cap B_{ r_{\eps} }(0) \setminus \overline{B}_{R k_{N,\eps}}(0)  \right)}    ( x, \nu)  \frac{|\nabla \ue|^2}{2}   ~d\sigma  \notag\\
&& \leq~  C~\int \limits_{   \T \left( \partial \rnm \cap B_{ r_{\eps} }(0) \setminus \overline{B}_{R k_{N,\eps}}(0)  \right)}   |x|^{2} \left[   \frac{ \mu_{N,\eps}^{\bp-\bm}}{|x|^{2\bp}}+ \frac{ 1 }{|x|^{2\bm}}  \right]~d\sigma \notag\\
&&\leq ~  C~k_{N,\eps} \int \limits_{  \partial \rnm \cap  B_{ \frac{ r_{\eps}}{k_{N,\eps}}}(0)  \setminus \overline{B}_{R }(0)  }    \frac{ 1}{|x|^{2\bp-2}}~d\sigma \notag\\
& &+ ~  C~r_{\eps}^{\bp-\bm+1}\int \limits_{  \partial \rnm \cap  B_{1}(0)  \setminus \overline{B}_{ \frac{R k_{N,\eps}}{ r_{\eps}}}(0)  }    \frac{ 1 }{ |x|^{2\bm-2}} ~d\sigma \notag\\
 &&\leq~  C~k_{N,\eps} \int \limits_{  \partial \rnm \cap  B_{ \frac{ k_{N,\eps}}{k_{N-1,\eps}}}(0)  \setminus \overline{B}_{R /2}(0)  }    \frac{ 1}{|x|^{(n-1)+(\bp-\bm-1)}}~d\sigma, \notag
\end{eqnarray*}
and then
\begin{eqnarray*}
&&\int \limits_{   \T \left( \partial \rnm \cap B_{ r_{\eps} }(0) \setminus \overline{B}_{R k_{N,\eps}}(0)  \right)}    ( x, \nu)  \frac{|\nabla \ue|^2}{2}   ~d\sigma  \notag\\
& &+ ~  C~r_{\eps}^{\bp-\bm+1} \int \limits_{  \partial \rnm \cap  B_{1 }(0)  \setminus \overline{B}_{ \frac{R k_{N,\eps}}{2r_{\eps}}}(0)  }    \frac{ 1 }{ |x|^{(n-1)-(\bp-\bm+1)}} ~d\sigma \notag\\
&& \leq C~k_{N,\eps}~\left(  R^{-(\bp-\bm-1)}+  r_{\eps}^{\bp-\bm-1} \right)~d\sigma. \notag\\
\end{eqnarray*}
So if $\bp-\bm>1$
\begin{align}\label{cpct:bdry:term:C:1}
\lim\limits_{R \to + \infty} \lim\limits_{\rho \to 0} \lim\limits_{\eps \to 0} \left( \mu_{N,\eps}^{-1} \int \limits_{ \T \left( \partial \rnm \cap B_{ r_{\eps} }(0) \setminus \overline{B}_{R k_{N,\eps}}(0)  \right)}    ( x, \nu)  \frac{|\nabla \ue|^2}{2}   ~d\sigma  \right)=0.
\end{align}
This ends Step \ref{step:p4}.3.

\medskip
\noindent{\it Step \ref{step:p4}.4:}   Let $1 \leq i\leq N$.  When $\bp -\bm >1$, we have 
\begin{align}\label{cpct:bdry:term:C:i}
\lim\limits_{R \to + \infty} \lim\limits_{\rho \to 0} \lim\limits_{\eps \to 0} &\left( \mu_{i,\eps}^{-1}  \int \limits_{ \T \left(  \partial \rnm \cap B_{R k_{i,\eps}}(0)  \setminus \overline{B}_{\rho k_{i,\eps}}(0)  \right)}    ( x, \nu)  \frac{|\nabla \ue|^2}{2}   ~d\sigma   \right) \notag\\
&\qquad \qquad \qquad =\frac{1}{2}\frac{1}{ t_{i}^{\frac{n-1}{\crits-2}}} \int \limits_{\partial\rnm} II_{0}(x,x)  \frac{|\nabla \tui|^2}{2}  ~d\sigma,
\end{align}
where $II_{0}(x,x)$ is the second fundamental form of the boundary $\bdry$ at $0$.
\medskip

\noindent{\it Proof of Step \ref{step:p4}.4:}   Consider $\tu_{i}$ obtained in  Proposition \ref{prop:exhaust}.  It follows  that  for some constant $C>0$, 
\begin{align*}
|\nabla \tu_{i}(x)| \leq \frac{C}{|x|^{\bm}+|x|^{\bp}} \qquad \hbox{ for all } x \in \rnmpbar.
\end{align*}
So when $\bp-\bm >1$,  the function  $|x|^{2}|\nabla \tu_{i}| \in L^{2}(\R^{n-1}) $.
\medskip

\noindent
With a change of variable and the definition of $\tu_{i,\eps}$ we then obtain 
\begin{eqnarray}
&&\mu_{i,\eps}^{-1}  \int \limits_{ \T \left(  \partial \rnm \cap B_{R k_{i,\eps}}(0)  \setminus \overline{B}_{\rho k_{i,\eps}}(0)  \right)}    ( x, \nu)  \frac{|\nabla \ue|^2}{2}   ~d\sigma  \notag\\
&& = \frac{k_{i,\eps}^{n-3}}{ \mu_{i,\eps}^{n-1}}  ~ \int \limits_{   \partial \rnm \cap B_{R }(0)  \setminus \overline{B}_{\rho }(0)  }   \left(\T(k_{N,\eps}x ), \nu \ \circ  T(k_{N,\eps} x)\right) \frac{|\nabla \tu_{i,\eps}|^2}{2}   ~d\sigma   \notag\\
&& = -\frac{k_{i,\eps}^{n-3}}{ \mu_{i,\eps}^{n-1}}   \left(~ \int \limits_{   \partial \rnm \cap B_{R }(0)  \setminus \overline{B}_{\rho }(0)  }   \frac{1}{2} \keN^2 \sum_{p,q=2}^n\partial_{pq}\T_0(0)x^px^q  \frac{|\nabla \tu_{i,\eps}|^2}{2}   ~d\sigma  +\theta_{\eps,R}(x) \keN^2  \right)  \notag\\
&& =  -\left( \frac{k_{i,\eps}}{ \mu_{i,\eps}} \right)^{n-1}    \left(~ \int \limits_{   \partial \rnm \cap B_{R }(0)  \setminus \overline{B}_{\rho }(0)  }   \frac{1}{2}  \sum_{p,q=2}^n\partial_{pq}\T_0(0)x^px^q  \frac{|\nabla \tu_{i,\eps}|^2}{2}   ~d\sigma  +\theta_{\eps,R}(x)   \right).  \notag
\end{eqnarray}
Since   $|x|^{2}|\nabla \tu_{i}| \in L^{2}(\R^{n-1}) $, passing to the limits it follows from the expression of the second fundamental form in \eqref{est:T:4}, that 
\begin{eqnarray*}
&&\lim\limits_{R \to + \infty} \lim\limits_{\rho \to 0} \lim\limits_{\eps \to 0} \left( \mu_{i,\eps}^{-1}  \int \limits_{ \T \left(  \partial \rnm \cap B_{R k_{i,\eps}}(0)  \setminus \overline{B}_{\rho k_{i,\eps}}(0)  \right)}    ( x, \nu)  \frac{|\nabla \ue|^2}{2}   ~d\sigma   \right)  \notag\\
&& = - \frac{1}{2}\frac{1}{ t^{\frac{n-1}{\crits-2}}}  \int \limits_{   \partial \rnm \cap B_{R }(0)  \setminus \overline{B}_{\rho }(0)  }   \sum_{p,q=2}^n\partial_{pq}\T_0(0)x^px^q  \frac{|\nabla \tui|^2}{2}   ~d\sigma  \notag\\
&& = \frac{1}{2}\frac{1}{ t_{i}^{\frac{n-1}{\crits-2}}} \int \limits_{\partial\rnm} II_{0}(x,x)  \frac{|\nabla \tui|^2}{2}  ~d\sigma.
\end{eqnarray*}
This ends Step \ref{step:p4}.4.

\smallskip\noindent Plugging \eqref{cpct:bdry:term:C:1}, \eqref{cpct:bdry:term:C:i}, \eqref{cpct:bdry:term:C:3} and  \eqref{cpct:bdry:term:C:4} in the integral \eqref{decomp:bdry:term:C}, we get \eqref{cpct:bdry:term:C}. This proves the first identity of Step \ref{step:p4}.

\medskip\noindent{\it Step \ref{step:p4}.5:} We now assume that $u_0\equiv 0$ and $\bp-\bm>1$. We prove \eqref{cpct:bdry:term:C:2}. We write
\begin{align}\label{decomp:bdry:term:C:2}
\int \limits_{   \T\left( \partial \rnm \cap B_{ \delta_0 }(0) \setminus B_{ k_{1,\eps}^{3} }(0) \right) }     ( x, \nu)  \frac{|\nabla \ue|^2}{2}   ~d\sigma=& 
  \int \limits_{   \T \left( \partial \rnm \cap B_{ \delta }(0) \setminus \overline{B}_{r_\eps}(0)  \right)}    ( x, \nu)  \frac{|\nabla \ue|^2}{2}   ~d\sigma \notag\\
  &+  \int \limits_{   \T \left( \partial \rnm \cap  \overline{B}_{r_\eps}(0) \setminus B_{ k_{1,\eps}^{3} }(0) \right)}    ( x, \nu)  \frac{|\nabla \ue|^2}{2}   ~d\sigma 
 \end{align}

\medskip\noindent With the pointwise estimates of Proposition \ref{prop:fund:est:grad} with $u_0\equiv 0$, and using that $\bp-\bm>1$, we obtain as $\eps \to 0$
\begin{eqnarray*}
&&\left|\int \limits_{   \T \left( \partial \rnm \cap B_{\delta_0 }(0) \setminus \overline{B}_{r_{\eps}}(0)  \right)}    ( x, \nu)  \frac{|\nabla \ue|^2}{2}   ~d\sigma\right|  \\&&\leq   C~\int \limits_{   \partial \rnm \cap B_{ \delta_0 }(0) \setminus \overline{B}_{r_\eps}(0) }   |x|^{2} \left[   \frac{ \mu_{N,\eps}^{\bp-\bm}}{|x|^{2\bp}}  \right]~d\sigma \notag\\
&&\leq  C \frac{ \mu_{N,\eps}^{\bp-\bm}}{\re^{2\bp-2-n+1}}\leq C \mu_{N,\eps}^{1+\frac{\bp-\bm-1}{2}}=o(\mu_{N,\eps}),
\end{eqnarray*}
since $\bp-\bm>1$. Then, with \eqref{cpct:bdry:term:C}, we get \eqref{cpct:bdry:term:C:2}. This ends Step \ref{step:p4}.5.

\smallskip\noindent These five substeps prove Step \ref{step:p4}.\qed
\bigskip

\subsection{Estimates of the boundary terms}

\begin{step}\label{step:p5} We let $(\ue)$, $(\he)$ and $(\pe)$ be such that $(E_\eps)$, \eqref{hyp:he},  \eqref{lim:pe} and \eqref{bnd:ue} hold. We assume that blow-up occurs. We fix a chart $\T$ as in \eqref{def:T:bdry} and, for any $\eps>0$, we define
\begin{align*}
\tve(x):= r_{\eps}^{\bm-1}  \ue( \T (r_{\eps} x)) \qquad \hbox{for } x \in r_{\eps}^{-1} U \cap \rnmpbar,
\end{align*}
where $r_{\eps}:= \sqrt{\mu_{N,\eps}}$. We claim that there exists  $\tv \in C^{1} ( \rnmpbar) $ such that
$$\lim \limits_{\eps \to 0}  \tve(x)= \tv~\hbox{ in } C^{1}_{loc} ( \rnmpbar) $$
where $\tv$ is a solution of 
\begin{equation}\label{blowup:eqn:cpct}
\left\{ \begin{array}{llll}
 -\Delta \tv -\frac{\gamma}{|x|^{2}} \tv&=&0  & \hbox{ in } \rnm \\
\hfill \tv&=&0 & \hbox{ on } \partial \rnm \setminus \{ 0\}.
\end{array}\right.
\end{equation}
\end{step}
\noindent {\it Proof of Step \ref{step:p5}:} For any $i,j=1,...,n$, we let $(\tge)_{ij}=(\T^*Eucl)(r_\eps x)_{ij}=(\partial_i\T(r_{\eps} x),\partial_j\T( r_{\eps} x))$,
where $(\cdot,\cdot)$ denotes the Euclidean scalar product on $\rn$. We consider $\tge$ as a metric on $\rn$. In the sequel, we let $\Delta_g=div_g(\nabla)$ be the Laplace-Beltrami
operator with respect to a metric $g$. From $(E_{\eps})$ it follows that for all $\eps >0$, the rescaled functions $\tve $  weakly satisfies the equation 
\begin{align}{\label{eqn:cpct:1}}
-\Delta_{\tge}  \tve - \frac{\gamma}{\left| \frac{\T ( r_{\eps} x)}{r_{\eps}}\right|^{2}}  \tve -r_{\eps}^2~\he \circ \T(r_{\eps} x)~ \tve =r_{\eps}^{\theta +\pe \bm}\frac{| \tve |^{\crits-2-\pe} \tve}{\left| \frac{\T (r_{\eps} x)}{r_{\eps}}\right|^{s}}.
\end{align}
with  $\theta:=(\crits-2)\frac{\bp-\bm}{2}>0$ and $  \tve \equiv 0$ on $  \partial \rnm \setminus \{ 0\}$.
\medskip

\noindent
Using the pointwise estimates  \eqref{eq:est:global}  we obtain the bound, that as $\eps \to 0$ we have for $ x \in \rnm$
\begin{align*}
|\tv_{\eps}(x)|  \leq & ~C~r_{\eps}^{\bm-1} \sum_{i=1}^N\frac{\mei^{\frac{\bp-\bm}{2}} | \T (r_{\eps} x)|}{ \mei^{\bp-\bm}| \T (r_{\eps} x)|^{\bm}+| \T (r_{\eps} x)|^{\bp}}  \notag\\
&~+~C~ r_{\eps}^{\bm-1}\frac{\Vert |x|^{\bm-1}u_0 \vert|_{L^{\infty}(\Omega)}}{| \T (r_{\eps} x)|^{\bm}}  |\T (r_{\eps} x)| \notag\\
 \leq & ~C~\sum_{i=1}^N \frac{\left( \frac{\mei}{\mu_{N,\eps}} \right)^{\frac{\bp-\bm}{2}}  \left| \frac{\T (r_{\eps} x)}{r_{\eps}} \right|}{ \left( \frac{\mei}{ \sqrt{\mu_{N,\eps}}} \right)^{\bp-\bm} \left| \frac{\T (r_{\eps} x)}{r_{\eps}} \right|^{\bm}+ \left| \frac{\T (r_{\eps} x)}{r_{\eps}} \right|^{\bp}}  \notag\\
&~+~C~  \frac{\Vert |x|^{\bm-1}u_0 \vert|_{L^{\infty}(\Omega)}   }{ \left| \frac{\T (r_{\eps} x)}{r_{\eps}} \right|^{\bm}}  \left| \frac{\T (r_{\eps} x)}{r_{\eps}} \right| \notag\\
 \leq & ~C \left(~\sum_{i=1}^N\frac{ \left( \frac{\mei}{\mu_{N,\eps}} \right)^{\frac{\bp-\bm}{2}} |x| }{ \left( \frac{\mei}{ \sqrt{\mu_{N,\eps}}} \right)^{\bp-\bm}| x|^{\bm}+|x|^{\bp}}\right.\\
 &\left.+ \frac{\Vert |x|^{\bm-1}u_0 \vert|_{L^{\infty}(\Omega)}+ }{|x|^{\bm}} |x| \right) \notag\\
  \leq &~C ~ \left(\frac{1}{ |x|^{\bp-1}} + \frac{\Vert |x|^{\bm-1}u_0 \vert|_{L^{\infty}(\Omega)}}{ |x|^{\bm-1} }\right).
\end{align*} 
Then  passing to limits in the  equation \eqref{eqn:cpct:1}, standard elliptic theory  yields the existence of $\tilde{v}\in C^2(\rnmpbar)$ such that $\tv_{\eps} \to \tv$ in $C^2_{loc}(\rnmpbar)$ and  $\tv$ satisfies the equation: 
\begin{align*}
\left\{ \begin{array}{llll}
-\Delta \tv -\frac{\gamma}{|x|^{2}} \tv&=&0  \ \ & \hbox{ in } \rnm \\
\hfill \tv&=&0 & \hbox{ on } \partial \rnm \setminus \{ 0\}.
\end{array}\right.
\end{align*}
and we have the following bound on $\tv$
$$|\tv(x)|  \leq C ~ \left(\frac{|x_{1}|}{ |x|^{\bp}} + \frac{\Vert |x|^{\bm-1}u_0 \vert|_{L^{\infty}(\Omega)}}{ |x|^{\bm} }|x_{1}| \right) \qquad \hbox{ for all  } x=(x_{1},\tilde{x}) \hbox{ in } \rnm.$$ 
This ends the proof of Step \ref{step:p5}.\qed

\begin{step}\label{step:p6} We let $(\ue)$, $(\he)$ and $(\pe)$ be such that $(E_\eps)$, \eqref{hyp:he},  \eqref{lim:pe} and \eqref{bnd:ue} hold. We assume that blow-up occurs. We claim that, as $\eps \to 0$,
\begin{align}\label{cpct:bdry:term1and23}
\int \limits_{ \T \left( \rnm  \cap \partial B_{ r_{\eps} }(0) \right)}  F_{\eps}(x)~d\sigma =~
\mu_{N,\eps}^{\frac{\bp-\bm}{2}} \left(  \mathcal{F}_0 +o(1) \right)
\end{align}
with
\begin{align}\label{def:F1}
\mathcal{F}_0 := \int \limits_{\rnm \cap \partial  B_{1 }(0) } ( x, \nu) \left( \frac{|\nabla \tv|^2}{2}  -\frac{\gamma}{2}\frac{\tv^{2}}{|x|^{2}} \right) -  \left(x^i\partial_i \tv+\frac{n-2}{2} \tv \right)\partial_\nu \tv ~d \sigma 
\end{align}
and 
\begin{align}\label{cpct:bdry:term1and2}
 ~ \int \limits_{ \T \left( \rnm  \cap \partial B_{ k_{1,\eps}^{3} }(0) \right)} F_{\eps}(x)~d\sigma&~=  o\left( \mu_{N,\eps}^{\bp-\bm}\right).
\end{align}
\end{step}

\noindent{\it Proof of Step \ref{step:p6}:} We keep the notations of Step \ref{step:p5}.
With a change of variable and the definition of $\tve$, and $\theta:=(\crits-2)\frac{\bp-\bm}{2}>0$, we get
\begin{align*}
& \int \limits_{ \T \left( \rnm  \cap \partial B_{ r_{\eps} }(0) \right)} F_{\eps}(x)~d\sigma= \notag\\
& ~ r_{\eps}^{\bp-\bm} \int \limits_{\rnm \cap \partial  B_{1  }(0) } ( x, \nu)_{\tge} \left( \frac{|\nabla_{\tge} \tve|^2}{2}  -\frac{\gamma}{2}\frac{\tve^{2}}{|x|_{\tge}^{2}} \right) -  \left(x^i\partial_i \tve+\frac{n-2}{2} \tve \right)\partial_{\nu} \tve ~d \sigma_{\tge}  \notag\\
& ~-r_{\eps}^{\bp-\bm} \int \limits_{\rnm \cap \partial  B_{1 }(0) } \left(r_{\eps}^{2}\frac{h_{\eps}( \T (r_{\eps}x))}{2} \tve^{2}  - \frac{ r_{\eps}^{\theta+(\bm -1)\pe} }{\crits -p_{\eps} } 
\frac{|\tve|^{\crits-p_{\eps} }}{|x|_{\tge}^{s}}\right) d\sigma_{\tge}.
\end{align*}
From the convergence result of Step \ref{step:p5}, we  then get \eqref{cpct:bdry:term1and23}.
\medskip

\noindent
For the next boundary term,  from  the estimates \eqref{eq:est:global} and  \eqref{eq:est:grad:global} we obtain
\begin{eqnarray*}
&&\left| \int_{ \T \left( \rnm  \cap \partial B_{ k_{1,\eps}^{3} }(0) \right)} F_{\eps}(x)~d\sigma \right|\\
&& \leq ~\frac{C}{\mu_{1,\eps}^{\bp-\bm}}   \int \limits_{ \T \left( \rnm  \cap \partial B_{ k_{1,\eps}^{3} }(0) \right)} |x| \left( \frac{1}{|x|^{2\bm}}+\frac{|x|^{2}}{|x|^{2\bm}}  \right)~dx \notag\\
 &&+C   \int \limits_{ \T \left( \rnm  \cap \partial B_{ k_{1,\eps}^{3} }(0) \right)} |x| \frac{ \mu_{1,\eps}^{-\crits \left(\frac{\bp -\bm}{2} \right)+ \pe \left(\frac{\bp -\bm}{2} \right)} }{|x|^{(\bm-1)(\crits-\pe)+s}}~dx \notag\\
&& \leq ~\frac{C}{\mu_{1,\eps}^{\bp-\bm}}   \int \limits_{  \rnm  \cap \partial B_{ k_{1,\eps}^{3} }(0) } |x| \left( \frac{1}{|x|^{2\bm}}+\frac{|x|^{2}}{|x|^{2\bm}}  \right)~dx \notag\\
 && +C~   \int \limits_{  \rnm  \cap \partial B_{ k_{1,\eps}^{3} }(0) } |x| \frac{ \mu_{1,\eps}^{-\crits \left(\frac{\bp -\bm}{2} \right)+ \pe \left(\frac{\bp -\bm}{2} \right)} }{|x|^{(\bm-1)(\crits-\pe)+s}}~dx \notag\\
  & &\leq  C \mu_{1,\eps}^{\bp-\bm} \left( \mu_{1,\eps}^{\bp-\bm}+\mu_{1,\eps}^{(\bp-\bm)\left( \frac{2-s}{n-2}\right) +\pe\left( \frac{n-2}{2}\right)}  \right).
\end{eqnarray*}
And so 
\begin{align}\label{cpct:bdry:term2}
\int \limits_{ \T \left( \rnm  \cap \partial B_{ k_{1,\eps}^{2} }(0) \right)} F_{\eps}(x)~d\sigma=  o\left( \mu_{N,\eps}^{\bp-\bm}\right).
\end{align}
This ends Step \ref{step:p6}.\qed

\begin{step}\label{step:p7} We let $(\ue)$, $(\he)$ and $(\pe)$ be such that $(E_\eps)$, \eqref{hyp:he},  \eqref{lim:pe} and \eqref{bnd:ue} hold. We assume that blow-up occurs. We assume that $u_0\equiv 0$. We define 
\begin{equation}\label{def:bue}
\bar{u}_\eps:=\frac{\ue}{\mu_{N,\eps}^{\frac{\bp-\bm}{2}}}.
\end{equation}
We claim that 
there exists $\bar{u}\in C^2(\overline{\Omega}\setminus\{0\})$ such that
\begin{equation}\label{cv:bue}
\lim_{\eps\to 0}\bar{u}_\eps=\bar{u}\hbox{ in }C^2_{loc}(\overline{\Omega}\setminus\{0\})\hbox{ with }\left\{\begin{array}{ll}
-\Delta \bar{u}-\left(\frac{\gamma}{|x|^2}+h_0\right)\bar{u}=0 &\hbox{ in }\Omega\\
\bar{u}=0&\hbox{ in }\partial\Omega\setminus\{0\}
\end{array}\right.
\end{equation}
\end{step}
\noindent{\it Proof of Step \ref{step:p7}:} Since $u_0\equiv 0$, it follows from \eqref{eq:est:global} that there exists $C>0$ such that 
\begin{equation}\label{control:bue}
|\bar{u}_\eps(x)|\leq C |x|^{1-\bp}\hbox{ for all }x\in\Omega\hbox{ and }\eps>0.
\end{equation}
Moreover, equation  $(E_\eps)$ rewrites
$$-\Delta \bar{u}_\eps-\left(\frac{\gamma}{|x|^2}+h_\eps\right)\bar{u}_\eps=\mu_{N,\eps}^{\frac{\bp-\bm}{2}(\crits-2-p_\eps)}\frac{|\bar{u}_\eps|^{\crits-2-p_\eps}\bar{u}_\eps}{|x|^s}\;\hbox{ in }\Omega,
$$
and $\bar{u}_\eps=0$ on $\partial\Omega$. It then follows from standard elliptic theory that the claim holds. This ends Step \ref{step:p7}.\qed

\begin{step}\label{step:p8} We let $(\ue)$, $(\he)$ and $(\pe)$ be such that $(E_\eps)$, \eqref{hyp:he},  \eqref{lim:pe} and \eqref{bnd:ue} hold. We assume that blow-up occurs. We assume that $u_0\equiv 0$. We claim that
\begin{align}\label{cpct:bdry:term1and2:lala}
\int \limits_{ \T \left( \rnm  \cap \partial B_{ \delta_0 }(0) \right)}  F_{\eps}(x)~d\sigma&~=\left( \mathcal{F}_{\delta_{0}}+o(1)\right)\mu_{N,\eps}^{\bp-\bm}, 
\end{align}
and 
\begin{align}\label{id:tralala}
 ~ \int \limits_{ \T \left( \rnm  \cap \partial B_{ k_{1,\eps}^{3} }(0) \right)} F_{\eps}(x)~d\sigma&~=  o\left( \mu_{N,\eps}^{\bp-\bm}\right), 
\end{align}
where
\begin{equation}\label{def:Fd}
 \mathcal{F}_{\delta_{0}} := \int \limits_{ \T \left( \rnm  \cap \partial B_{ \delta_0 }(0) \right)} ( x, \nu) \left( \frac{|\nabla \bar{u}|^2}{2}  -\left(\frac{\gamma}{|x|^2}+h_0\right)\frac{\bar{u}^{2}}{2} \right) -  \left(x^i\partial_i \bar{u}+\frac{n-2}{2} \bar{u} \right)\partial_{\nu} \bar{u} ~d \sigma.
\end{equation}
\end{step}
\noindent{\it Proof of Step \ref{step:p8}:} The second term has already been estimated in \eqref{cpct:bdry:term1and2}. We are left with the first term. 
With a change of variable, the definition of $\bar{u}_\eps$ and the convergence \eqref{cv:bue}, we get
\begin{eqnarray}
&& \int \limits_{ \T \left( \rnm  \cap \partial B_{ \delta_0 }(0) \right)} F_{\eps}(x)~d\sigma\\
&&=  \mu_{N,\eps}^{\bp-\bm} \int \limits_{ \T \left( \rnm  \cap \partial B_{ \delta_0 }(0) \right)} ( x, \nu) \left( \frac{|\nabla \bar{u}_\eps|^2}{2}  -\left(\frac{\gamma}{|x|^{2}}+h_\eps\right)\frac{\bar{u}_\eps^{2}}{2}\right)~d \sigma  \notag \\
 && - \mu_{N,\eps}^{\frac{\bp-\bm}{2}(\crit-2-\eps)}\int \limits_{ \T \left( \rnm  \cap \partial B_{ \delta_0 }(0) \right)}\frac{|\bar{u}_\eps|^{\crit-2-\eps}\bar{u}_\eps}{|x|^2}~d \sigma  \notag \\
 &&- \int \limits_{ \T \left( \rnm  \cap \partial B_{ \delta_0 }(0) \right)}  \left(x^i\partial_i \bar{u}_\eps+\frac{n-2}{2} \bar{u}_\eps \right)\partial_{\nu} \bar{u}_\eps ~d \sigma  \notag\\
&&=\mu_{N,\eps}^{\bp-\bm} \left(   \mathcal{F}_{\delta_{0}}+o(1) \right).\label{est:Fe}
\end{eqnarray}
where $\mathcal{F}_{\delta_{0}}$ is as above.
Arguing as in the proof of \eqref{cpct:bdry:term2}, we get that
\begin{align}\label{cpct:bdry:term2:bis}
\int \limits_{ \T \left( \rnm  \cap \partial B_{ k_{1,\eps}^{3} }(0) \right)} F_{\eps}(x)~d\sigma=  o\left( \mu_{N,\eps}^{\bp-\bm}\right)\hbox{ as }\eps\to 0.
\end{align}
This ends Step \ref{step:p8}.\qed

\begin{step}\label{step:p9} We let $(\ue)$, $(\he)$ and $(\pe)$ be such that $(E_\eps)$, \eqref{hyp:he},  \eqref{lim:pe} and \eqref{bnd:ue} hold. We assume that blow-up occurs. We assume that $\ue>0$ for all $\eps>0$. Then $\mathcal{F}_0\geq 0$ and
$${\mathcal F}_0>0\;\Leftrightarrow\; u_0>0.$$
where $\mathcal{F}_0$ is as in \eqref{def:F1}.
\end{step}
\noindent{\it Proof of Step \ref{step:p9}:} We let $\tv$ be defined as in Step \ref{step:p5}. It follows from Step \ref{step:p5} that $\tv$ satisfies \eqref{blowup:eqn:cpct} and we have the following bound on $\tv$
\begin{equation}
|\tv(x)|  \leq C ~ \left(\frac{|x_{1}|}{ |x|^{\bp}} + \frac{\Vert |x|^{\bm-1}u_0 \vert|_{L^{\infty}(\Omega)}}{ |x|^{\bm} }|x_{1}| \right) \qquad \hbox{ for all  } x=(x_{1},\tilde{x}) \hbox{ in } \rnm.\label{est:tv}
\end{equation} 
Given $\alpha\in\rr$, we define $v_\alpha(x):=x_1|x|^{-\alpha}$ for all $x\in \rnm$. Since $\tv\geq 0$, it follows from Proposition 6.4 in Ghoussoub-Robert \cite{gr4} that  there exists $A,B\geq 0$ such that
\begin{equation}\label{eq:tv:A:B}
\tv:=A v_{\bp}+B v_{\bm}.
\end{equation}

\noindent{\it Step \ref{step:p9}.1:} We claim that $B=0$ when $u_0\equiv 0$.

\smallskip\noindent This is a direct consequence of controling \eqref{eq:tv:A:B} with \eqref{est:tv} when $u_0\equiv 0$ and letting $|x|\to\infty$.

\medskip\noindent{\it Step \ref{step:p9}.2:} We claim that $B>0$ when $u_0>0$.

\smallskip\noindent We prove the claim. We fix $x\in \rnm$. Green's representation formula yields
\begin{eqnarray*}
\tv_\eps(x)&=& \int_\Omega r_{\eps}^{\bm-1} G_\eps(\T (r_{\eps} x), y)\frac{\ue^{\crits-1}(y)}{|y|^s}\, dy.
\end{eqnarray*}
We fix $\omega\subset\subset \Omega$. Then there exists $c(\omega)>0$ such that $|y|\geq d(y,\partial\Omega)\geq c(\omega)$ for all $y\in \omega$. Moreover, the control \eqref{est:G:up} of the Green's function yields
\begin{eqnarray*}
\tv_\eps(x)&\geq & c\int_\omega r_{\eps}^{\bm-1} \frac{r_\eps x_1}{r_\eps^{\bm}|x|^{\bm}}|c(\omega)-r_\eps |x||^{-n}\frac{\ue^{\crits-1}(y)}{|y|^s}\, dy,
\end{eqnarray*}
and then, passing to the limit $\eps\to 0$,  we get that
\begin{eqnarray*}
\tv(x)&\geq &  \frac{c x_1}{|x|^{\bm}}\int_\omega \frac{u_0^{\crits-1}(y)}{|y|^s}\, dy,
\end{eqnarray*}
for all $x\in \rnm$. As one checks, this yields $B\geq c\int_\omega \frac{u_0^{\crits-1}(y)}{|y|^s}\, dy >0$ when $u_0>0$. This ends Step \ref{step:p9}.2.

\medskip\noindent{\it Step \ref{step:p9}.3:} We claim that $A>0$.

\smallskip\noindent The proof is similar to Step \ref{step:p9}.2. We fix $x\in \rnm $ and $\omega\subset\subset \rnm$. Green's representation formula and the pointwise control \eqref{est:G:up} yield
\begin{eqnarray*}
\tv_\eps(x)&\geq & \int_{\T (\mu_{N,\eps}\omega)} r_{\eps}^{\bm-1} G_\eps(\T (r_{\eps} x), y)\frac{\ue^{\crits-1}(y)}{|y|^s}\, dy\\
&\geq & \int_{\omega} r_{\eps}^{\bm-1} G_\eps(\T (r_{\eps} x), \T(\mu_{N,\eps} y))\mu_{N,\eps}^n \frac{\ue(\T(\mu_{N,\eps} y))^{\crits-1}}{|\mu_{N,\eps} y|^s}\, dy\\
&\geq & \int_{\omega} r_{\eps}^{\bm-1} \left(\frac{r_{\eps}|x|}{\mu_{N,\eps}|y|}\right)^{\bm}K_\eps(x,y) \mu_{N,\eps}^{\frac{n-2}{2}}\frac{\tuei(y)^{\crits-1}}{|y|^s}\, dy\\
&\geq & \int_{\omega} r_{\eps}^{2\bm-n} |x|^{\bm} \left|x-\frac{\mu_{N,\eps}}{r_\eps} y \right|^{-n}x_1 y_1 \mu_{N,\eps}^{\frac{n}{2}-\bm}\frac{\tuei(y)^{\crits-1}}{|y|^s}\, dy
\end{eqnarray*}
with
$$K_\eps(x,y)=|r_\eps x-\mu_{N,\eps} y|^{2-n}\min\left\{1,\frac{\mu_{N,\eps}}{r_\eps}\frac{x_1 y_1}{|x-\frac{\mu_{N,\eps}}{r_\eps} y|^2}\right\} $$
Since $r_\eps:=\sqrt{\mu_{N,\eps}}$, letting $\eps\to 0$, we get with the convergence (A4) of Proposition \ref{prop:exhaust} that
\begin{eqnarray*}
\tv_\eps(x)&\geq & \int_{\omega} r_{\eps}^{2\bm-n} |x|^{\bm} \left|x-\frac{\mu_{N,\eps}}{r_\eps} y \right|^{-n}x_1 y_1 \mu_{N,\eps}^{\frac{n}{2}-\bm}\frac{\tuei(y)^{\crits-1}}{|y|^s}\, dy\\
&\geq & \frac{x_1}{|x|^{\bp}}\int_{\omega} \frac{\tilde{u}_i(y)^{\crits-1}}{|y|^s}\, dy
\end{eqnarray*}
for all $x\in\rnm$. Therefore, as one checks, $A\geq \int_{\omega} \frac{\tilde{u}_i(y)^{\crits-1}}{|y|^s}\, dy>0$. This ends Step \ref{step:p9}.3.

\medskip\noindent{\it Step \ref{step:p9}.4:} We claim that
\begin{equation}\label{id:F:A:B}
\mathcal{F}_0=\frac{\omega_{n-1}}{n}\left(\frac{n^2}{4}-\gamma\right)\cdot AB.
\end{equation}
We prove the claim. The definition \eqref{def:F1} reads \begin{align}\label{def:F0:2}
\mathcal{F}_0 := \int \limits_{\rnm \cap \partial  B_{1 }(0) } ( x, \nu) \left( \frac{|\nabla \tv|^2}{2}  -\frac{\gamma}{2}\frac{\tv^{2}}{|x|^{2}} \right) -  \left(x^i\partial_i \tv+\frac{n-2}{2} \tv \right)\partial_\nu \tv ~d \sigma 
\end{align}
For simplicity, we define the bilinear form
\begin{eqnarray*}
{\mathcal H}_\delta(u,v)&=&\int \limits_{\rnm \cap \partial  B_{\delta}(0) } \left[( x, \nu) \left( (\nabla u,\nabla v)  -\gamma \frac{uv}{|x|^{2}} \right) -  \left(x^i\partial_i u+\frac{n-2}{2} u \right)\partial_\nu v\right.\\
&&\left.-\left(x^i\partial_i v+\frac{n-2}{2} v \right)\partial_\nu u\right] ~d \sigma
\end{eqnarray*}
As one checks,
\begin{eqnarray}
\mathcal{F}_0&=&\frac{1}{2}{\mathcal H}_1(A v_{\bp}+B v_{\bm},A v_{\bp}+B v_{\bm})\nonumber\\
&=& \frac{A^2}{2}{\mathcal H}_1(v_{\bp},v_{\bp})+AB{\mathcal H}_1(v_{\bp},v_{\bm})\nonumber\\
&&+\frac{B^2}{2}{\mathcal H}_1(v_{\bm},v_{\bm})\nonumber
\end{eqnarray}
In full generality, we compute ${\mathcal H}_\delta(v_\alpha,v_\beta)$ for all $\alpha,\beta\in\rr$ and all $\delta>0$. As one checks, for any $i=1,...,n$, we have that $\partial_i v_\alpha=\left(\delta_{i,1}-\alpha\frac{x_1 x_i}{|x|^2}\right)|x|^{-\alpha}$ for all $x\in\rnm$. Moreover, for $x\in \partial B_\delta(0)$, we have that $\partial_\nu v_\alpha=\frac{x^i}{|x|}\partial_i v_\alpha$. Consequently, straightforward computations yield
$$  \left(x^i\partial_i v_\alpha+\frac{n-2}{2} v_\alpha \right)\partial_\nu v_\beta=-(\beta-1)\left(\frac{n}{2}-\alpha\right)\frac{v_\alpha v_\beta}{|x|}$$
and
$$(x,\nu)\left((\nabla v_\alpha,\nabla v_\beta)-\frac{\gamma}{|x|^2}v_\alpha v_\beta\right)=|x|^{1-\alpha-\beta}+(\alpha\beta-\alpha-\beta-\gamma)\frac{v_\alpha v_\beta}{|x|}$$
and then
\begin{eqnarray*}
{\mathcal H}_\delta(v_\alpha,v_\beta) &=& \int_{\rnm\cap \partial B_\delta(0)}\left(|x|^{1-\alpha-\beta}+\left(\frac{n}{2}(\alpha+\beta)-n-\alpha\beta-\gamma\right)\frac{v_\alpha v_\beta}{|x|}\right)\, d\sigma
\end{eqnarray*}
We have that
$$\int_{\rnm\cap \partial B_\delta(0)}|x|^{1-\alpha-\beta}\, d\sigma=\frac{1}{2}\int_{B_\delta(0)}|x|^{1-\alpha-\beta}\, d\sigma=\frac{\omega_{n-1}}{2}\delta^{n-\alpha-\beta}$$
and
\begin{eqnarray*}
\int_{\rnm\cap \partial B_\delta(0)}\frac{v_\alpha v_\beta}{|x|}\, d\sigma&=&\frac{1}{2}\int_{B_\delta(0)}x_1^2|x|^{-\alpha-\beta-1}\, d\sigma\\
&=&\frac{1}{2n}\int_{B_\delta(0)}|x|^{-\alpha-\beta+1}\, d\sigma=\frac{\omega_{n-1}}{2n}\delta^{n-\alpha-\beta}
\end{eqnarray*}
Plugging all these identities together yields
$${\mathcal H}_\delta(v_\alpha,v_\beta) =\frac{\omega_{n-1}}{2n}\delta^{n-\alpha-\beta}\left(\frac{n}{2}(\alpha+\beta)-\alpha\beta-\gamma\right).$$
Since $\bp,\bm$ are solutions to $X^2-nX+\gamma=0$, we get that 
$${\mathcal H}_\delta(v_{\bm},v_{\bm})= {\mathcal H}_\delta(v_{\bp},v_{\bp})= 0.$$
 Since $\bp+\bm=n$ and $\bp\bm=\gamma$, we get that
$${\mathcal H}_\delta(v_{\bm},v_{\bp})=\frac{\omega_{n-1}}{n}\left(\frac{n^2}{4}-\gamma\right).$$
Plugging all these results together yields \eqref{id:F:A:B}. This ends Step \ref{step:p9}.4.

\smallskip\noindent These substeps end the proof of Step \ref{step:p9}.\qed

\begin{step}\label{step:p10} We let $(\ue)$, $(\he)$ and $(\pe)$ be such that $(E_\eps)$, \eqref{hyp:he},  \eqref{lim:pe} and \eqref{bnd:ue} hold. We assume that blow-up occurs. We assume that $\bp-\bm<2$ and $\ue>0$ for all $\eps>0$. Then $u_0\equiv 0$.
\end{step}
\noindent{\it Proof of Step \ref{step:p10}:} We claim that, as $\eps\to 0$,
\begin{equation}\label{est:2:1}
\int \limits_{   \T\left( \partial \rnm \cap B_{ r_{\eps} }(0) \setminus B_{ k_{1,\eps}^{3} }(0) \right) }    ( x, \nu)  \frac{|\nabla \ue|^2}{2}   ~d\sigma=o\left(\mu_{N,\eps}^{\frac{\bp-\bm}{2}}\right)\hbox{ when }\bp-\bm<2
\end{equation}
Indeed, if $\bp-\bm>1$, the claim follows from \eqref{cpct:bdry:term:C} and $1>\frac{\bp-\bm}{2}$. If now $\bp-\bm<1$, then \eqref{est:T:2} and the control \eqref{eq:est:grad:global} yield that
\begin{eqnarray*}
&&\left|\int \limits_{   \T\left( \partial \rnm \cap B_{ r_{\eps} }(0) \setminus B_{ k_{1,\eps}^{2} }(0) \right) }    ( x, \nu)  \frac{|\nabla \ue|^2}{2}   ~d\sigma\right| \\
&&\leq  C\int_{   \partial \rnm \cap B_{ r_{\eps} }(0)  }|x|^2\left(\sum_{i=1}^N\frac{\mu_{i,\eps}^{\bp-\bm}}{|x|^{2\bp}}\, dx+ \frac{dx}{|x|^{2\bm}}\right)   ~d\sigma\\
&&\leq C\sum_{i=1}^N\mu_{i,\eps}^{\bp-\bm}r_\eps^{n-1-2(\bp-1)}+Cr_\eps^{\bp-\bm+1}= o(\mu_{N,\eps}^{\frac{\bp-\bm}{2}})
\end{eqnarray*}
as $\eps\to 0$. The limit case $\bp-\bm=1$ is similar. This proves the claim.

\medskip\noindent Plugging \eqref{eq:mystery}, \eqref{cpct:L2term:sc}, \eqref{cpct:bdry:term1and23}, \eqref{cpct:bdry:term1and2} and \eqref{est:2:1} into the Pohozaev identity \eqref{PohoId3}, we get

\begin{align}\label{id:poho:36}
  & \frac{p_{\eps}}{\crits} \left( \frac{n-s}{\crits}\right)  \left(\sum \limits_{i=1}^{N}  \frac{1}{t_{i}^{ \frac{n-2}{\crits-2}}} \int \limits_{\rnm }   \frac{|\tu_{i}|^{\crits }}{|x|^{s}} ~ dx +o(1)\right)
 =-\left(\mathcal{F}_0+o(1)\right)\mu_{N,\eps}^{\frac{\bp-\bm}{2}}
\end{align}
as $\eps\to 0$, where $\mathcal{F}_0$ is as in \eqref{def:F0:2}. Therefore $\mathcal{F}_0\leq 0$. Since $\ue>0$, it then follows from \eqref{id:F:A:B} of Step \ref{step:p9} that $u_0\equiv 0$. This proves Step \ref{step:p10}.\qed \bigskip

\section{\, Proof of the sharp blow-up rates}\label{pf blow-up rates}
We now prove the sharp blow-up rates claimed in  Propositions \ref{prop:rate:sc} and \ref{prop:rate:sc:2}. We start with the case when $\bp-\bm\neq 1$. As a preliminary estimate, we claim that
\begin{eqnarray}
&  & \frac{p_{\eps}}{\crits} \left( \frac{n-s}{\crits -p_{\eps}}\right)  \left(\sum \limits_{i=1}^{N}  \frac{1}{t_{i}^{ \frac{n-2}{\crits-2}}} \int \limits_{\rnm }   \frac{|\tu_{i}|^{\crits }}{|x|^{s}} ~ dx +o(1)\right)\notag\\
 &&=\int \limits_{   \T\left( \partial \rnm \cap B_{ r_{\eps} }(0) \setminus B_{ k_{1,\eps}^{3} }(0) \right) }    ( x, \nu)  \frac{|\nabla \ue|^2}{2}   ~d\sigma
   -\left( \mathcal{F}_0+o(1)\right)\mu_{N,\eps}^{\frac{\bp-\bm}{2}} \label{est:1}
\end{eqnarray}
as $\eps\to 0$, where $\mathcal{F}_0 $ is as in \eqref{def:F1}; and, when $u_0\equiv 0$, we claim that
\begin{eqnarray}\label{est:4:1}
 && \frac{p_{\eps}}{\crits} \left( \frac{n-s}{\crits -p_{\eps}}\right)  \left(\sum \limits_{i=1}^{N}  \frac{1}{t_{i}^{ \frac{n-2}{\crits-2}}} \int \limits_{\rnm }   \frac{|\tu_{i}|^{\crits }}{|x|^{s}} ~ dx +o(1)\right)   \notag\\
&& = \int \limits_{   \T\left( \partial \rnm \cap B_{ \delta_0 }(0) \setminus B_{ k_{1,\eps}^{3} }(0) \right) }    ( x, \nu)  \frac{|\nabla \ue|^2}{2}   ~d\sigma - \left(\mathcal{F}_{\delta_0} +o(1)\right)\mu_{N,\eps}^{\bp-\bm}  \notag\\
&& {+ \underbrace{o(\mu_{N,\eps})}_{\text{when  $\bp-\bm \geq 2$}} +    \underbrace{O(\mu_{N,\eps}^{\bp-\bm})}_{\text{when  $\bp-\bm < 2$}}},
 \end{eqnarray}
where $\mathcal{F}_{\delta_{0}}$ is as in \eqref{def:Fd}.

\smallskip\noindent We prove the claim. Collecting  the first estimate of Step P2, \eqref{cpct:L2term:sc}, \eqref{cpct:bdry:term1and23} and \eqref{cpct:bdry:term1and2} of the terms of the  Pohozaev identity \eqref{PohoId3} gives \eqref{est:1}. Similarly, the second estimate of Step P2, \eqref{cpct:L2term:u=0}, \eqref{cpct:bdry:term1and2:lala} and \eqref{id:tralala} of the terms of the  Pohozaev identity \eqref{PohoId3:bis} gives \eqref{est:4:1}.

\subsection{Proof of the sharp blow-up rates when $\bp-\bm\neq 1$}

We first assume $\ue>0$ and $\bp-\bm<1$.

\begin{step}\label{step:p11} We let $(\ue)$, $(\he)$ and $(\pe)$ be such that $(E_\eps)$, \eqref{hyp:he},  \eqref{lim:pe} and \eqref{bnd:ue} hold. We assume that blow-up occurs. We assume that $\ue>0$ and $\bp-\bm<1$. Then \eqref{est:mass} holds, that is
\begin{equation}
\lim_{\eps\to 0}\frac{p_\eps}{\mu_{N,\eps}^{\bp-\bm}}=-\frac{\frac{\omega_{n-1}\crits^2}{n}\left(\frac{n^2}{4}-\gamma\right)A^2}{(n-s)\sum \limits_{i=1}^{N}  \frac{1}{t_{i}^{ \frac{n-2}{\crits-2}}} \int \limits_{\rnm }   \frac{|\tu_{i}|^{\crits }}{|x|^{s}} ~ dx}\cdot m_{\gamma,h}(\Omega)
\label{est:mass:bis}
\end{equation}
for some $A>0$, where $m_{\gamma,h}(\Omega)$ is the boundary mass.
\end{step}

\noindent{\it Proof of Step \ref{step:p11}:} It follows from Step \ref{step:p10} that $u_0\equiv 0$. 

\smallskip\noindent{\it Step \ref{step:p11}:1:} We now claim that
\begin{equation}
\frac{p_{\eps}}{\crits} \left( \frac{n-s}{\crits }\right)  \left(\sum \limits_{i=1}^{N}  \frac{1}{t_{i}^{ \frac{n-2}{\crits-2}}} \int \limits_{\rnm }   \frac{|\tu_{i}|^{\crits }}{|x|^{s}} ~ dx +o(1)\right) = \mu_{N,\eps}^{\bp-\bm}\left(M_{\delta_0}+o(1)\right)\nonumber\end{equation}
where
\begin{eqnarray}\label{est:inf:1}
M_{\delta_0}&:=&-\int_{\T\left( \rnm \cap B_{\delta_0 }(0)  \right)} \left( h_0(x) +\frac{ \left( \nabla h_0, x \right)}{2} \right)\bar{u}^{2}~ dx- \mathcal{F}_{\delta_0}\nonumber\\
&&+\int_{   \T\left( \partial \rnm \cap B_{ \delta_0}(0)  \right) }    ( x, \nu)  \frac{|\nabla \bar{u}|^2}{2}   ~d\sigma,
\end{eqnarray}
and $\mathcal{F}_{\delta_{0}}$ is as in \eqref{def:Fd} and $\bar{u}$ is as in \eqref{cv:bue}.

\smallskip\noindent Indeed, the Pohozaev identity \eqref{PohoId3}, the convergence \eqref{def:bue}, \eqref{control:bue}, \eqref{cv:bue} and $\bp-\bm<1$ yield
\begin{eqnarray}\label{est:l2}
&&\int \limits_{\T\left( \rnm \cap B_{\delta_0 }(0) \setminus B_{ k_{1,\eps}^{3} }(0) \right)} \left( h_{\eps}(x) +\frac{ \left( \nabla h_{\eps}, x \right)}{2} \right)\ue^{2}~ dx\\
&&=\mu_{N,\eps}^{\bp-\bm}\left(\int_{\T\left( \rnm \cap B_{\delta_0 }(0)  \right)} \left( h_0(x) +\frac{ \left( \nabla h_0, x \right)}{2} \right)\bar{u}^{2}~ dx+o(1)\right)\nonumber
\end{eqnarray}
With $u_0\equiv 0$ and the control \eqref{eq:est:grad:global}, we get that $|\nabla \ue(x)|\leq C\mu_{N,\eps}^{\frac{\bp-\bm}{2}}|x|^{-\bp}$ for all $\eps>0$ and $x\in \Omega$. Therefore, with \eqref{def:bue} and \eqref{cv:bue}, we get that
\begin{eqnarray}\label{est:grad}
&&\int \limits_{   \T\left( \partial \rnm \cap B_{ \delta_0}(0) \setminus B_{ k_{1,\eps}^{3} }(0) \right) }    ( x, \nu)  \frac{|\nabla \ue|^2}{2}   ~d\sigma\nonumber\\
&&=\mu_{N,\eps}^{\bp-\bm}\left(\int_{   \T\left( \partial \rnm \cap B_{ \delta_0}(0)  \right) }    ( x, \nu)  \frac{|\nabla \bar{u}|^2}{2}   ~d\sigma+o(1)\right)
\end{eqnarray}
as $\eps\to 0$. Plugging \eqref{cpct:bdry:term1and2}, \eqref{est:l2} and \eqref{est:grad} into \eqref{PohoId3:bis}, we get \eqref{est:inf:1}. This proves the claim and ends Step \ref{step:p11}.1.

\smallskip\noindent We fix $\delta<\delta'$. Taking $U:=\T(\rnm\cap B_{\delta'}(0)\setminus B_\delta(0))$, $K=0$ and $u=\bar{u}$ in \eqref{PohoId}, and using \eqref{cv:bue}, we get that $M_\delta$ is independent of the choice of $\delta>0$ small enough.

\medskip\noindent{\it Step \ref{step:p11}.2:} We claim that $\bar{u}>0$.

\smallskip\noindent We prove the claim. Since $\bar{u}\geq 0$ is a solution to \eqref{cv:bue}, it is enough to prove that $\bar{u}\not\equiv 0$. We argue as in the proof of Step \ref{step:p9}. We fix $x\in \Omega$. Green's identity and$\ue>0$  yield
\begin{eqnarray*}
&&\quad \bue(x)=\mu_{N,\eps}^{-(\bp-\bm)/2}\int_\Omega G_\eps(x,y)\frac{\ue(y)^{\crits-1-p_\eps}}{|y|^s}\, dy\\
&&\geq  \mu_{N,\eps}^{-(\bp-\bm)/2}\int_{A_\eps} G_\eps(x,y)\frac{\ue(y)^{\crits-1-p_\eps}}{|y|^s}\, dy\\
&&\geq  C \mu_{N,\eps}^{n-s-(\bp-\bm)/2}\int_{A} G_\eps(x,\T(\mu_{N,\eps}y))\frac{\ue(\T(\mu_{N,\eps}y))^{\crits-1-p_\eps}}{|y|^s}\, dy,
\end{eqnarray*}
where $A_\epsilon:=\T(\rnm\cap B_{2\mu_{N,\eps}}(0)\setminus B_{\mu_{N,\eps}}(0))$, $A:=\rnm\cap B_2(0)\setminus B_1(0)$. With the pointwise control \eqref{est:G:up}, we get
\begin{eqnarray*}
&&\quad \bue(x)\geq\\
&&  C \int_A \left(\frac{|x|}{|y|}\right)^{\bm}|x-\T(\mu_{N,\eps}y)|^{2-n}\left(\frac{d(x,\partial\Omega) |y_1|}{|x-\T(\mu_{N,\eps}y)|^2}\right)
\frac{u_{\eps,i}(y)^{\crits-1-p_\eps}}{|y|^s}\, dy
\end{eqnarray*}
where $u_{\eps,i}$ is as in Proposition \ref{prop:exhaust}. Letting $\eps\to 0$ and using the convergence (A4) of Proposition \ref{prop:exhaust}, we get that
$$\bar{u}(x)\geq C\frac{d(x,\partial\Omega)}{|x|^{\bp}} \hbox{ for all }x\in \Omega.$$
And then $\bar{u}>0$ in $\Omega$. This proves the claim and Step \ref{step:p11}.2.

\medskip\noindent We fix $r_0>0$ and $\eta\in C^\infty(\rn)$ such that $\eta(x)=1$ in $B_{r_0}(0)$ and $\eta(x)=0$ in $\rn\setminus B_{2r_0}(0)$. It then follows from \cites{gr4,gr5} that, for $r_0>0$ small enough, there exists $A>0$ and  $\beta\in H_{0}^1(\Omega)$ such that
$$\bar{u}(x)=A\left(\frac{\eta(x) d(x,\partial\Omega)}{|x|^{\bp}}+\beta(x)\right)\hbox{ for all }x\in \Omega$$
with
$$\beta(x)=m_{\gamma,h}(\Omega)\frac{\eta(x) d(x,\partial\Omega)}{|x|^{\bm}}+o\left(\frac{\eta(x) d(x,\partial\Omega)}{|x|^{\bm}}\right)$$
as $\eps\to 0$. Here, $m_{\gamma,h}(\Omega)$ is the boundary mass. 

\medskip\noindent{\it Step \ref{step:p11}.3:}  We claim that
\begin{equation}\label{est:Md}
\lim_{\delta\to 0}M_\delta= - \frac{\omega_{n-1}}{n}\left(\frac{n^2}{4}-\gamma\right)A^2\cdot m_{\gamma,h}(\Omega)
\end{equation}
We prove the claim. Since $\bar{u}$ is a solution to \eqref{cv:bue}, it follows from standard elliptic theory that there exists $C>0$ such that $\bar{u}(x)+|x||\nabla\bar{u}(x)|\leq C|x|^{1-\bp}$ for all $x\in \Omega$. Therefore, since $\bp-\bm<1$, we get that 
$$\lim_{\delta\to 0}\int_{\T(\rnm\cap B_\delta(0))}\bar{u}^2\, dx+\int_{\T(\rnm\cap \partial B_\delta(0))}\bar{u}^2\, d\sigma+\int_{\T(\partial\rnm\cap B_\delta(0))}|x|^2|\nabla\bar{u}|^2\, d\sigma=0.$$
Therefore,
\begin{equation*}
M_\delta=-\frac{A^2}{2}\bar{\mathcal H}_\delta(\bar{v}_{\bp}+\bar{v}_{\bm},\bar{v}_{\bp}+\bar{v}_{\bm} )+o(1)
\end{equation*}
as $\delta\to 0$, where
\begin{eqnarray*}
\bar{\mathcal H}_\delta(u,v)&:=&\int \limits_{ \T \left( \rnm  \cap \partial B_{ \delta_0 }(0) \right)} \left[( x, \nu) \left( (\nabla u,\nabla v)  -\frac{\gamma}{|x|^2}uv \right) -  \left(x^i\partial_i u+\frac{n-2}{2} u \right)\partial_{\nu} v\right.\\
&&\left.- \left(x^i\partial_i v+\frac{n-2}{2} v \right)\partial_{\nu} u\right]\,d \sigma
\end{eqnarray*}
and
$$\bar{v}_{\bp}(x):=\frac{\eta(x) d(x,\partial\Omega)}{|x|^{\bp}}\hbox{ and }\bar{v}_{\bm}(x)=\beta(x)\hbox{ for all }x\in \Omega.$$
We then get that
\begin{eqnarray*}
M_\delta&=&-\frac{A^2}{2}\bar{\mathcal H}_\delta(\bar{v}_{\bp},\bar{v}_{\bp} )-A^2\bar{\mathcal H}_\delta(\bar{v}_{\bp},\bar{v}_{\bm} )\\
&&-\frac{A^2}{2}\bar{\mathcal H}_\delta(\bar{v}_{\bm},\bar{v}_{\bm} )+o(1)
\end{eqnarray*}
as $\delta\to 0$. For any $x\in \rnm\cap B_\delta(0)$, with the chart $\T$ and the definition of $\beta$, we get
$$\bar{v}_{\bp}(\T(x)):=\frac{|x_1|}{|x|^{\bp}}+O(|x|^{2-\bp})=v_{\bp}+O(|x|^{2-\bp})$$
$$\hbox{ and }\bar{v}_{\bm}(\T(x))=m_{\gamma,h}(\Omega)\frac{|x_1|}{|x|^{\bm}}+O(|x|^{2-\bm})=m\cdot v_{\bm}+O(|x|^{2-\bm}).$$
Moreover, elliptic theory yields
$$\nabla (\bar{v}_{\bp}\circ\T(x)):=\nabla v_{\bp}+O(|x|^{1-\bp}).$$
$$\hbox{ and }\nabla (\bar{v}_{\bm}\circ\T(x))=m_{\gamma,h}(\Omega)\cdot \nabla v_{\bm}+O(|x|^{1-\bm})\hbox{ for all }x\in \rnm\cap B_\delta(0), $$
where $v_\beta$ is defined in the proof of Step \ref{step:p9}. Since $\bp-\bm<1$ and $\bp+\bm=n$, we get with a change of variable that as $\delta\to 0$, 
\begin{eqnarray*}
\bar{\mathcal H}_\delta(\bar{v}_{\bp},\bar{v}_{\bp} )&=&{\mathcal H}_\delta(v_{\bp},v_{\bp} )+O(\delta^{1-(\bp-\bm)})\\
\bar{\mathcal H}_\delta(\bar{v}_{\bp},\bar{v}_{\bm} )&=&m_{\gamma,h}(\Omega)\cdot {\mathcal H}_\delta(v_{\bp},v_{\bm} )+O(\delta^{1-(\bp-\bm)})\\
\bar{\mathcal H}_\delta(\bar{v}_{\bm},\bar{v}_{\bm} )&=&O(\delta^{n-2\bm}).
\end{eqnarray*}
Using the computations performed in the proof of Step \ref{step:p9}, we then get \eqref{est:Md}. This proves the claim and ends Step \ref{step:p11}.3.

\medskip\noindent{\it End of the proof of Step \ref{step:p11}:} Since $M_\delta$ is independent of $\delta$ small, we then get that $M_{\delta_0}=- \frac{\omega_{n-1}}{n}\left(\frac{n^2}{4}-\gamma\right)A^2m_{\gamma,h}(\Omega)$. Putting this estimate in \eqref{est:inf:1}, we then get \eqref{est:mass:bis}. This end Step \ref{step:p11}.\qed

\medskip\noindent {\it Proof of Proposition \ref{prop:rate:sc} when $\bp-\bm>2$:} Plugging \eqref{cpct:bdry:term:C} into \eqref{est:1} and using that $ \bp-\bm>2$, we obtain 
\begin{align*}
\lim \limits_{\eps \to 0} \frac{\pe}{\mu_{N,\eps}}&~=\frac{n-s}{(n-2)^{2}}  \frac{1}{ t_{N}^{\frac{n-1}{\crits-2}}}\frac{ \int \limits_{\partial\rnm} II_{0}(x,x)  |\nabla \tu_{N}|^2  ~d\sigma}{  \sum \limits_{i=1}^{N}  \frac{1}{t_{i}^{ \frac{n-2}{\crits-2}}} \int \limits_{\rnm }   \frac{|\tu_{i}|^{\crits }}{|x|^{s}} ~ dx }. 
\end{align*}
This yields \eqref{rate:1} when $\bp-\bm>2$. \bigskip 

\medskip\noindent{\it Proof of Proposition \ref{prop:rate:sc} when $\bp-\bm>1$ and $u_0\equiv 0$.} Plugging \eqref{cpct:bdry:term:C:2} into \eqref{est:4:1} and using that $ \bp-\bm>1$, we obtain also \eqref{rate:1}.\\

\medskip\noindent{\it Proof of Proposition \ref{prop:rate:sc:2}  when $\bp-\bm>1$.} Since $\ue>0$, we get that $\tilde{u}_N>0$. Therefore, it follows from Ghoussoub-Robert \cite{gr4} that $\bar{u}_N(x_1,x')=\bar{U}_N(x_1,|x'|)$ for all $(x_1,x')\in(0,+\infty)\times \rr^{n-1}$. Due to this symmetry, when $\bp-\bm>1$, we get that 
 
\begin{eqnarray}
&&\int_{\partial\rnm} II_{0}(x,x)  |\nabla \tu_{N}|^2  ~d\sigma=\sum_{i,j=1}^{n-1}\int_{\partial\rnm} II_{0,ij}x^ix^j |\nabla \tu_{N}|^2  ~d\sigma\label{est:3:1}\\
&&=\frac{\sum_{i=1}^{n-1}II_{0,ii}}{n-1}\int_{\partial\rnm} |x|^2 |\nabla \tu_{N}|^2  ~d\sigma=\frac{\int_{\partial\rnm} |x|^2 |\nabla \tu_{N}|^2  ~d\sigma}{n-1}H(0).\notag
\end{eqnarray}
When $\bp-\bm>2$ or $\{\bp-\bm=2\hbox{ and }u_0\equiv 0\}$, Proposition \ref{prop:rate:sc:2} follows from \eqref{rate:1} and \eqref{est:3:1}. When $\{\bp-\bm=2\hbox{ and }u_0>0\}$, Proposition \ref{prop:rate:sc:2} follows from \eqref{est:1}, \eqref{id:F:A:B} of Step \ref{step:p9}, \eqref{cpct:bdry:term:C} and \eqref{est:3:1}. When $1<\bp-\bm<2$, Proposition \ref{prop:rate:sc:2} follows from Step \ref{step:p10}, \eqref{est:4:1}, \eqref{cpct:bdry:term:C:2} and \eqref{est:3:1}.

\medskip\noindent{\it Proof of Proposition \ref{prop:rate:sc:2}  when $\bp-\bm<1$:} This is a direct consequence of Steps \ref{step:p10} and \ref{step:p11}.

 \subsection{Proof of the sharp blow-up rates when $\bp-\bm = 1$}
 
We start with the following refined asymptotics when $\ue>0$, $\bp-\bm=1$ and $u_0\equiv 0$.

\begin{step}\label{prop:11} We let $(\ue)$, $(\he)$ and $(\pe)$ be such that $(E_\eps)$, \eqref{hyp:he},  \eqref{lim:pe} and \eqref{bnd:ue} hold. We assume that blow-up occurs. We assume that $\ue>0$ and $u_0\equiv 0$. We fix a family of parameters $(\lambda_\eps)_{\eps>0}\in (0,+\infty)$ such that
\begin{equation}\label{hyp:lambda}
\lim_{\eps\to 0}\lambda_\eps=0\hbox{ and }\lim_{\eps\to 0}\frac{\mu_{N,\eps}}{\lambda_\eps}=0.
\end{equation}
Then, for all $x\in\overline{\rnm}$, $x\neq 0$, we have that
\begin{equation*}
\lim_{\eps\to 0}\frac{\lambda_\eps^{\bp-1}}{\mu_{N,\eps}^{\frac{\bp-\bm}{2}}}\ue(\T(\lambda_\eps x))=K\cdot\frac{|x_1|}{|x|^{\bp}},
\end{equation*}
where $\T$ is as in \eqref{def:T:bdry},
\begin{equation}\label{def:K}
K:=t_N^{-\frac{\bp-1}{\crits-2}}L_{\gamma,\Omega}\int_{\rnm}\frac{|y_1|}{|y|^{\bm}}\frac{\tilde{u}_N^{\crits-1}(y)}{|y|^s}\, dy>0
\end{equation}
and $L_{\gamma,\Omega}>0$ is given by \eqref{est:G:4}. Moreover, this limit holds in $C^2_{loc}(\overline{\rnm}\setminus\{0\})$.
\end{step}
\noindent{\it Proof of Step \ref{prop:11}:} We define
$$w_\eps(x):=\frac{\lambda_\eps^{\bp-1}}{\mu_{N,\eps}^{\frac{\bp-\bm}{2}}}\ue(\T(\lambda_\eps x))$$
for all $x\in \rnm\cap \lambda_\eps^{-1}U$. As in the proof of \eqref{cpct:bdry:term1and23}, for any $i,j=1,...,n$, we let $(\tge)_{ij}=(\partial_i\T(r_{\eps} x),\partial_j\T( r_{\eps} x))$,
where $(\cdot,\cdot)$ denotes the Euclidean scalar product on $\rn$. We consider $\tge$ as a metric on $\rn$. We let $\Delta_g=div_g(\nabla)$, the Laplace-Beltrami
operator with respect to the metric $g$.  From   $(E_{\eps})$ it follows that for all $\eps >0$, we have that
\begin{equation*}
\left\{\begin{array}{cl}
-\Delta_{\tge}  w_\eps - \frac{\gamma}{\left| \frac{\T ( \lambda_{\eps} x)}{\lambda_{\eps}}\right|^{2}}  w_\eps -\lambda_{\eps}^2~\he \circ \T(\lambda_{\eps} x) w_\eps =s_\eps\frac{w_\eps^{\crits-1-\pe}}{\left| \frac{\T (\lambda_{\eps} x)}{\lambda_{\eps}}\right|^{s}}& \hbox{ in }\rnm\cap \lambda_\eps^{-1}U\\
 w_\eps>0 & \hbox{ in }\rnm\cap \lambda_\eps^{-1}U\\
 w_\eps=0 & \hbox{ on }(\partial\rnm\setminus\{0\})\cap \lambda_\eps^{-1}U.
\end{array}\right.
\end{equation*}
With
$$s_\eps:=\left(\frac{\mu_{N,\eps}^{\frac{\bp-\bm}{2}}}{\lambda_\eps^{\bp-1}}\right)^{\crit-2-p_\eps}\lambda_\eps^{2-s} .$$
Since $\mu_{N,\eps}^{p_\eps}\to t_N>0$ (see (A9) of Proposition \ref{prop:exhaust}) and $$(\bp-1)(\crits-2)-(2-s)=(\crits-2)\frac{\bp-\bm}{2},$$ then using the hypothesis \eqref{hyp:lambda}, we get that
$$\left(\frac{\mu_{N,\eps}^{\frac{\bp-\bm}{2}}}{\lambda_\eps^{\bp-1}}\right)^{\crit-2-p_\eps}\lambda_\eps^{2-s}\leq C
\left(\frac{\mu_{N,\eps}}{\lambda_\eps}\right)^{(\bp-1)(\crits-2-p_\eps)-(2-s)}=o(1)$$
as $\eps\to 0$. Since $u_0\equiv 0$, it follows from the pointwise control \eqref{eq:est:global} that there exists $C>0$ such that $0<w_\eps(x)\leq C|x_1|\cdot |x|^{-\bp}$ for all $x\in \rnm\cap \lambda_\eps^{-1}U$. It then follows from standard elliptic theory that there exists $w\in C^2(\overline{\rnm}\setminus\{0\})$ such that
\begin{equation}\label{lim:w:c2}
\lim_{\eps\to 0}w_\eps=w\hbox{ in }C^2_{loc}(\overline{\rnm}\setminus\{0\})
\end{equation}
with
\begin{equation*}
\left\{\begin{array}{ll}
-\Delta w - \frac{\gamma}{|x|^2}  w=0&\hbox{ in }\rnm\\
0\leq w(x)\leq C|x_1|\cdot |x|^{-\bp}&\hbox{ in }\rnm\\
w=0&\hbox{ on }\partial\rnm\setminus\{0\}.
\end{array}\right.
\end{equation*}
It follows from Lemma 4.2 in Ghoussoub-Robert \cite{gr4} (see also Pinchover-Tintarev  \cite{PT}) that there exists $\Lambda\geq 0$ such that $w(x)= \Lambda |x_1|\cdot |x|^{-\bp}$ for all $x\in\rnm$. We are left with proving that $\Lambda=K$ defined in \eqref{def:K}. We fix $x\in \rnm$. Green's representation formula yields
\begin{eqnarray}
w_\eps(x)&=&\int_\Omega\frac{\lambda_\eps^{\bp-1}}{\mu_{N,\eps}^{\frac{\bp-\bm}{2}}}G_\eps(\T(\lambda_\eps x),y)\frac{\ue(y)^{\crits-1-p_\eps}}{|y|^s}\,dy\nonumber\\
&=& \int_{\T(\rnm\cap (B_{R k_{N,\eps}}(0)\setminus B_{\delta k_{N,\eps}}(0))}+\int_{\Omega\setminus\T(\rnm\cap (B_{R k_{N,\eps}}(0)\setminus B_{\delta k_{N,\eps}}(0))}\label{est:w:0}
\end{eqnarray}

\medskip\noindent{\it Step \ref{prop:11}.1:} We estimate the first term of the right-hand-side. Since $D_{0} \mathcal{T} = \mathbb{I}_{\R^{n}}$, a change of variable yields
\begin{eqnarray*}
&&\int_{\T(\rnm\cap (B_{Rk_{N,\eps}}(0)\setminus B_{\delta k_{N,\eps}}(0))}\frac{\lambda_\eps^{\bp-1}}{\mu_{N,\eps}^{\frac{\bp-\bm}{2}}}G_\eps(\T(\lambda_\eps x),y)\frac{\ue(y)^{\crits-1-p_\eps}}{|y|^s}\,dy\\
&&=s_\eps^{(1)}\int_{\rnm\cap (B_{R}(0)\setminus B_{\delta}(0))}G_\eps(\T(\lambda_\eps x),\T(k_{N,\eps}z))\frac{\tilde{u}_{N,\eps}(z)^{\crits-1-p_\eps}}{|z|^s}(1+o(1))\,dz
\end{eqnarray*}
with
$$s_\eps^{(1)}:=\frac{\lambda_\eps^{\bp-1}}{\mu_{N,\eps}^{\frac{\bp-\bm}{2}}}k_{N,\eps}^{n-s}\mu_{N,\eps}^{-\frac{n-2}{2}(\crits-1-p_\eps)}$$
It follows from \eqref{est:G:4} that for any $z\in \rnm$, we have that
$$G_\eps(\T(\lambda_\eps x),\T(k_{N,\eps}z))=(L_{\gamma,\Omega}+o(1))\frac{\lambda_\eps|x_1|}{\lambda_\eps^{\bp}|x|^{\bp}}\cdot \frac{k_{N,\eps}|y_1|}{k_{N,\eps}^{\bm}|z|^{\bm}},$$
and that the convergence is uniform with repect to $z\in \rnm\cap (B_{R}(0)\setminus B_{\delta}(0))$. Plugging this estimate in the above equality, using that $k_{N,\eps}=\mu_{N,\eps}^{1-p_\eps/(\crits-2)}$, $\mu_{N,\eps}^{p_\eps}\to t_N>0$ and the convergence of $\tilde{u}_{N,\eps}$ to $\tilde{u}_N$ (see Proposition \ref{prop:exhaust}), we get that
\begin{eqnarray*}
&&\int_{\T(\rnm\cap (B_{R k_{N,\eps}}(0)\setminus B_{\delta k_{N,\eps}}(0))}\frac{\lambda_\eps^{\bp-1}}{\mu_{N,\eps}^{\frac{\bp-\bm}{2}}}G_\eps(\T(\lambda_\eps x),y)\frac{\ue(y)^{\crits-1-p_\eps}}{|y|^s}\,dy\\
&&=L_{\gamma,\Omega}\frac{|x_1|}{|x|^{\bp}}t_N^{-\frac{\bp-1}{\crits-2}}\int_{\rnm\cap (B_{R}(0)\setminus B_{\delta}(0))}\frac{|y_1|}{|y|^{\bm}}\frac{\tilde{u}_{N}(z)^{\crits-1}}{|z|^s}\,dz+o(1)
\end{eqnarray*}
as $\eps\to 0$. Therefore,
\begin{eqnarray}
&&\lim_{R\to +\infty, \delta\to 0}\lim_{\eps\to 0}\int_{\T(\rnm\cap (B_{R k_{N,\eps}}(0)\setminus B_{\delta k_{N,\eps}}(0))}\frac{\lambda_\eps^{\bp-1}}{\mu_{N,\eps}^{\frac{\bp-\bm}{2}}}G_\eps(\T(\lambda_\eps x),y)\frac{\ue(y)^{\crits-1}}{|y|^s}\,dy\nonumber\\
&&=K\frac{|x_1|}{|x|^{\bp}}\label{est:w:1}
\end{eqnarray}
where $K$ is as in \eqref{def:K}. 

\medskip\noindent{\it Step  \ref{prop:11}.2:} With the control \eqref{est:G:up} on the Green's function and the pointwise control \eqref{eq:est:global} on $\ue$, we get that
\begin{eqnarray}
&&\int_{\Omega\setminus\T(\rnm\cap (B_{R k_{N,\eps}}(0)\setminus B_{\delta k_{N,\eps}}(0))}\frac{\lambda_\eps^{\bp-1}}{\mu_{N,\eps}^{\frac{\bp-\bm}{2}}}G_\eps(\T(\lambda_\eps x),y)\frac{\ue(y)^{\crits-1-\pe}}{|y|^s}\,dy\nonumber\\
&&\leq \sum_{i=1}^{N-1}A_{i,\eps} +B_\eps(R)+C_\eps(\delta)\label{est:w:0:bis}
\end{eqnarray}
where
\begin{equation*}
A_{i,\eps}:=C\frac{\lambda_\eps^{\bp-1}}{\mu_{N,\eps}^{\frac{\bp-\bm}{2}}}\int_{B_{R_0}(0)}\frac{\ell_\epsilon (x, y)^{\bm}r_\epsilon (x , y)}{|\T(\lambda_\eps x)-y|^{n-2}|y|^s}\left(\frac{\mu_{i,\eps}^{\frac{\bp-\bm}{2}}|y|}{\mu_{i,\eps}^{\bp-\bm}|y|^{\bm}+|y|^{\bp}}\right)^{\crits-1}\, dy
\end{equation*}
\begin{equation*}
B_\eps(R):=C\frac{\lambda_\eps^{\bp-1}}{\mu_{N,\eps}^{\frac{\bp-\bm}{2}}}\int_{B_{R_0}(0)\setminus B_{R k_{N,\eps}}(0)}\frac{\ell_\epsilon (x, y)^{\bm}r_\epsilon (x, y)}{|\T(\lambda_\eps x)-y|^{n-2}}\frac{\mu_{N,\eps}^{\frac{\bp-\bm}{2}(\crits-1)}}{|y|^{(\bp-1)(\crits-1)+s}}\, dy
\end{equation*}

\begin{equation*}
C_\eps(\delta):=C(x)\frac{\lambda_\eps^{\bp-1+2-n+\bm-1}}{\mu_{N,\eps}^{\frac{\bp-\bm}{2}\cdot\crits}}\int_{ B_{\delta k_{N,\eps}}(0)}\frac{dy}{ |y|^{(\bm-1)(\crits-1)+s+\bm-1}}
\end{equation*}
where $\ell_\epsilon (x, y):=\frac{\max\{\lambda_\eps|x|,|y|\}}{\min\{\lambda_\eps|x|,|y|\}}$, and $r_\epsilon (x, y)=\min\left\{1,\frac{\lambda_\eps|x_1|\cdot|y|}{|\T(\lambda_\eps x)-y|^2}\right\}$.

\medskip\noindent{\it Step \ref{prop:11}.3.} We first estimate $C_\eps(\delta)$. Since $n>s+\crits(\bm-1)$ (this is a consequence of $\bm<n/2$), straightforward computations yield
$$C_\eps(\delta)\leq C(x)\delta^{\frac{\crits}{2}(\bp-\bm)},$$
and therefore
\begin{equation}\label{est:w:2}
\lim_{\delta\to 0}\lim_{\eps\to 0}C_\eps(\delta)=0.
\end{equation}

\medskip\noindent{\it Step \ref{prop:11}.4.} We estimate $B_\eps(R)$. We split the integral as
$$B_\eps(R)=\int_{R k_{\eps,N}<|y|<\frac{\lambda_\eps|x|}{2}}I_\eps(y)\, dy+\int_{\frac{\lambda_\eps|x|}{2}<|y|<2\lambda_\eps|x|}I_\eps(y)\, dy+\int_{|y|>2\lambda_\eps|x|}I_\eps(y)\, dy$$
where $I_\eps(y)$ is the integrand. Since 
$$n-(s+(\bp-1)(\crits-1)+\bm-1)=-\frac{\crits-2}{2}(\bp-\bm)<0,$$ straightforward computations yield
\begin{eqnarray*}
&&\int_{R k_{N,\eps}<|y|<\frac{\lambda_\eps|x|}{2}}I_\eps(y)\, dy\\
&&\leq C(x)\frac{\lambda_\eps^{\bp-1+\bm-1+2-n}}{\mu_{N,\eps}^{\frac{\bp-\bm}{2}}} \int_{R k_{N,\eps}<|y|<\frac{\lambda_\eps|x|}{2}}\frac{\mu_{N,\eps}^{\frac{\bp-\bm}{2}(\crits-1)}}{|y|^{(\bp-1)(\crits-1)+s+\bm-1}}\, dy\\
&&\leq C(x)R^{-\frac{\crits-2}{2}(\bp-\bm)},
\end{eqnarray*}
For the next term, a change of variable yields
\begin{eqnarray*}
&&\int_{\frac{\lambda_\eps|x|}{2}<|y|<2\lambda_\eps|x|}I_\eps(y)\, dy\\
&&\leq C(x)\frac{\lambda_\eps^{\bp-1}}{\mu_{N,\eps}^{\frac{\bp-\bm}{2}}} \int_{\frac{\lambda_\eps|x|}{2}<|y|<2\lambda_\eps|x|}|\T(\lambda_\eps x)-y|^{2-n}\frac{\mu_{N,\eps}^{\frac{\bp-\bm}{2}(\crits-1)}}{|y|^{(\bp-1)(\crits-1)+s}}\, dy\\
&&\leq C(x) \left(\frac{\mu_{N,\eps}}{\lambda_\eps}\right)^{\frac{\crits-2}{2}(\bp-\bm)}\int_{\frac{|x|}{2}<|z|<2|x|}|x-z|^{2-n}\, dz=o(1)
\end{eqnarray*}
as $\eps\to 0$. Finally, since $\bp+\bm=n$ and $n-s-(\bp-1)\crits=\frac{\crits}{2}(\bp-\bm)$, we estimate the last term
\begin{eqnarray*}
&&\int_{|y|>2\lambda_\eps|x|}I_\eps(y)\, dy\\
&&\leq C(x) \mu_{N,\eps}^{\frac{\crits-2}{2}(\bp-\bm)}\lambda_\eps^{\bp-\bm}\int_{|y|>2\lambda_\eps|x|}\frac{|y|^{\bm+1-n-s}\, dy}{|y|^{(\bp-1)(\crits-1)}}\\
&&\leq C(x) \left(\frac{\mu_{N,\eps}}{\lambda_\eps}\right)^{\frac{\crits-2}{2}(\bp-\bm)}=o(1)
\end{eqnarray*}
as $\eps\to 0$. All these inequalities yield
\begin{equation}\label{est:w:3}
\lim_{R\to +\infty}\lim_{\eps\to 0}B_\eps(R)=0.
\end{equation}

\medskip\noindent{\it Step \ref{prop:11}.5.} We fix $i\in \{1,...,N-1\}$ and estimate $A_{i,\eps}$. As above, we split the integral as
$$A_{i,\eps}=\int_{|y|<\frac{\lambda_\eps|x|}{2}}J_{i,\eps}(y)\, dy+\int_{\frac{\lambda_\eps|x|}{2}<|y|<2\lambda_\eps|x|}J_{i,\eps}(y)\, dy+\int_{|y|>2\lambda_\eps|x|}J_{i,\eps}(y)\, dy,$$
where $J_{i,\eps}$ is the integrand. Since $\mu_{i,\eps}\leq \mu_{N,\eps}$, as one checks, the second and the third integral of the right-hand-side are controled from above respectively by $\int_{\frac{\lambda_\eps|x|}{2}<|y|<2\lambda_\eps|x|}I_\eps(y)\, dy$ and $\int_{|y|>2\lambda_\eps|x|}I_\eps(y)\, dy$ that have been computed just above and go to $0$ as $\eps\to 0$. We are then left with the first term. With a change of variables, we have that
\begin{eqnarray*}
&&\int_{|y|<\frac{\lambda_\eps|x|}{2}}J_{i,\eps}(y)\, dy\\
&&\leq  C(x)\frac{\lambda_\eps^{\bp-1+\bm+2-n-1}}{\mu_{N,\eps}^{\frac{\bp-\bm}{2}}} \\
&&\times\int_{|y|<\frac{\lambda_\eps|x|}{2}}
\left(\frac{\mu_{i,\eps}^{\frac{\bp-\bm}{2}}|y|}{\mu_{i,\eps}^{\bp-\bm}|y|^{\bm}+|y|^{\bp}}\right)^{\crits-1}
\frac{dy}{|y|^{s+\bm-1}}\\
&&\leq C(x) \frac{\mu_{i,\eps}^{1+n-s-\bm-\frac{n-2}{2}(\crits-1)}}{\mu_{N,\eps}^{\frac{\bp-\bm}{2}}}\\
&&\times\int_{|z|<\frac{\lambda_\eps|x|}{2\mu_{i,\eps}}}\frac{1}{|z|^{(\bm-1)+s}}\left(\frac{|z|}{|z|^{\bm}+|z|^{\bp}}\right)^{\crits-1}\, dz\\
&&\leq C(x)\left(\frac{\mu_{i,\eps}}{\mu_{N,\eps}}\right)^{\frac{\bp-\bm}{2}}
\end{eqnarray*}
since $n>s+(\crits(\bm-1))$ and $n<(\bm-1)+s+(\crits-1)(\bp-1)$. Since $\mu_{i,\eps}=o(\mu_{N,\eps})$ as $\eps\to 0$, we get that
\begin{equation}\label{est:w:4}
\lim_{\eps\to 0}A_{i,\eps}=0.
\end{equation}

\medskip\noindent{\it Step \ref{prop:11}.6:} Plugging \eqref{est:w:1}, \eqref{est:w:2}, \eqref{est:w:3} and \eqref{est:w:4} into \eqref{est:w:0} and \eqref{est:w:0:bis} yields $\lim_{\eps\to 0}w_\eps(x)=K\frac{|x_1|}{|x|^{\bp}}$ for all $x\in \rnm$. With \eqref{lim:w:c2}, we then get that $\Lambda=K$. This proves Step \ref{prop:11}.\\

Now we can prove Proposition \ref{prop:rate:sc:2} when $\bp-\bm=1$ in the case when $\ue>0.$

\begin{step}\label{step:p13} We let $(\ue)$, $(\he)$ and $(\pe)$ be such that $(E_\eps)$, \eqref{hyp:he},  \eqref{lim:pe} and \eqref{bnd:ue} hold. We assume that blow-up occurs. We assume that $\ue>0$ and $\bp-\bm=1$. Then $u_0\equiv 0$ and 

\begin{eqnarray}
&&   \frac{p_{\eps}}{\crits} \left( \frac{n-s}{\crits}\right)  \left(\sum \limits_{i=1}^{N}  \frac{1}{t_{i}^{ \frac{n-2}{\crits-2}}} \int \limits_{\rnm }   \frac{|\tu_{i}|^{\crits }}{|x|^{s}} ~ dx +o(1)\right)\nonumber\\
&& = \frac{K^2\omega_{n-2}H(0)}{4(n-1)}\mu_{N,\eps}\ln\frac{1}{\mu_{N,\eps}}+o\left(\mu_{N,\eps}\ln\frac{1}{\mu_{N,\eps}}\right).\label{est:inter}
\end{eqnarray}
\end{step}

The case $\bp-\bm=1$ of Proposition \ref{prop:rate:sc:2} is a consequence of Step \ref{step:p13}.\\

\smallskip\noindent{\it Proof of Step \ref{step:p13}:} First remark that since $\bp+\bm=n$, we then have that
$$\bp=\frac{n+1}{2}\hbox{ and }\bm=\frac{n-1}{2}.$$
It follows from Step \ref{step:p10} that $u_0\equiv 0$. We use \eqref{est:4:1} that writes
{\small 
\begin{equation}
   \frac{p_{\eps}}{\crits} \left( \frac{n-s}{\crits -p_{\eps}}\right)  \left(\sum \limits_{i=1}^{N}  \frac{1}{t_{i}^{ \frac{n-2}{\crits-2}}} \int \limits_{\rnm }   \frac{|\tu_{i}|^{\crits }}{|x|^{s}} ~ dx +o(1)\right)
 = \int \limits_{ {\mathcal T}_\epsilon }    ( x, \nu)  \frac{|\nabla \ue|^2}{2}   ~d\sigma +O\left(\mu_{N,\eps}\right)\label{est:5}.
\end{equation}
}
where ${\mathcal T}_\epsilon:= \T\left( \partial \rnm \cap B_{ \delta_0 }(0) \setminus B_{ k_{1,\eps}^{3} }(0) \right)$. 
It follows from \eqref{est:T:3} that
\begin{eqnarray}
&&\ds \int \limits_{  {\mathcal T}_\epsilon}    ( x, \nu)  \frac{|\nabla \ue|^2}{2}   ~d\sigma=\ds-\frac{1}{4}\int \limits_{  {\mathcal T}_\epsilon} \sum_{p,q=2}^{n}  x^px^q \partial_{pq} \T_{0}(0)     |\nabla (\ue\circ\T)|_{\T^\star\eucl}^2   (1+O(|x|)~d\sigma\nonumber\\
&&\ds +O\left(\int \limits_{   \partial \rnm \cap B_{ \delta_0 }(0)  }|x|^3 |\nabla (\ue\circ\T)|_{\T^\star\eucl}^2  \,d\sigma\right)\nonumber\\
&&=\ds-\frac{1}{4}\int \limits_{   \partial \rnm \cap B_{ \delta_0 }(0) \setminus B_{ k_{1,\eps}^{3} }(0) }  \sum_{p,q=2}^{n}  x^px^q \partial_{pq} \T_{0}(0)       |\nabla (\ue\circ\T)|^2  \,d\sigma\nonumber\\
&&\ds +O\left(\int \limits_{   \partial \rnm \cap B_{ \delta_0 }(0)  }|x|^3 |\nabla (\ue\circ\T)|^2  \,d\sigma\right).\label{est:23}
\end{eqnarray}
With the control \eqref{eq:est:grad:global} and $\bp-\bm=1$, we get that
\begin{eqnarray}
\int_{   \partial \rnm \cap B_{ \delta_0 }(0)  }|x|^3 |\nabla (\ue\circ\T)|^2  \,d\sigma&\leq& C\sum_{i=1}^N\int_{   \partial \rnm \cap B_{ \delta_0 }(0)  }|x|^3 \frac{\mu_{i,\eps}^{\bp-\bm}}{|x|^{2\bp}}\,d\sigma\nonumber\\
&\leq& C\mu_{N,\eps}^{\bp-\bm}= C\mu_{N,\eps}\label{est:24}
\end{eqnarray}

\medskip\noindent We need an intermediate result. We let $(s_\eps)_\eps,(t_\eps)_\eps\in [0,+\infty)$ such that $0\leq s_\eps\leq t_\eps$, and $\mu_{\eps,N}=o(t_\eps)$ as $\eps\to 0$. We claim that
\begin{eqnarray}\label{est:square}
\int_{\partial \rnm \cap \left(B_{ t_\eps }(0)\setminus B_{s_\eps}(0)  \right)}|x|^2 |\nabla (\ue\circ\T)|^2  \,d\sigma&\leq& C\sum_{i}\mu_{i,\eps}\ln\left(\frac{t_\eps}{\max\{s_\eps,\mu_{i,\eps}\}}\right)
\end{eqnarray}
Indeed, with the pointwise control \eqref{eq:est:grad:global}, $u_0\equiv 0$ and $2\bp=n+1$, we get that
 
\begin{eqnarray*}
&&\ds \int \limits_{   \partial \rnm \cap \left(B_{ t_\eps }(0)\setminus B_{s_\eps}(0)  \right)}|x|^2 |\nabla (\ue\circ\T)|^2  \,d\sigma\\
&&\leq C\sum_{i=1,...,N}\mu_{i,\eps}^{\bp-\bm}\int_{s_\eps}^{t_\eps}\frac{r^{2+(n-1)-1}\, dr}{\mu_{i,\eps}^{2(\bp-\bm)}r^{2\bm}+r^{2\bp}}\\
&&\leq  C\sum_{i=1,...,N}\mu_{i,\eps}\int_{\frac{s_\eps}{\mu_{i,\eps}}}^{\frac{t_\eps}{\mu_{i,\eps}}}\frac{r^{2\bp-1}\, dr}{r^{2\bm}+r^{2\bp}}
\end{eqnarray*}
Distinguishing the cases $s_\eps\leq \mu_{i,\eps}$ and  $s_\eps\geq \mu_{i,\eps}$, we get \eqref{est:square}. This proves the claim.

\medskip\noindent We define $\theta_\eps:=\frac{1}{\sqrt{|\ln\mu_{N,\eps}|}}$, $\alpha_\eps:=\mu_{N,\eps}^{\theta_\eps}$ and $\beta_\eps:=\mu_{N,\eps}^{1-\theta_\eps}$. As one checks, we have that
\begin{equation}
\left\{\begin{array}{ccc}
\mu_{\eps,N}=o(\beta_\eps) & \beta_\eps=o(\alpha_\eps) & \alpha_\eps=o(1) \\
\ln \frac{\alpha_\eps}{\beta_\eps}\simeq \ln\frac{1}{\mu_{N,\eps}} & \ln \frac{\beta_\eps}{\mu_{N,\eps}}=o\left(\ln\frac{1}{\mu_{N,\eps}} \right) & \ln\alpha_\eps=o(\ln\mu_{N,\eps})
\end{array}\right\}\label{ppty:a:b}
\end{equation}
as $\eps\to 0$. It then follows from \eqref{est:square} and the properties \eqref{ppty:a:b} that
\begin{equation}
\left\{\begin{array}{l}
\ds \int\limits_{   \partial \rnm \cap \left(B_{ \delta_0}(0)\setminus B_{\alpha_\eps}(0)  \right)}|x|^2 |\nabla (\ue\circ\T)|^2  =o\left(\mu_{N,\eps} \ln\frac{1}{\mu_{N,\eps}}\right);\\
\ds \int \limits_{   \partial \rnm \cap B_{ \beta_\eps}(0)}|x|^2 |\nabla (\ue\circ\T)|^2  =o\left(\mu_{N,\eps} \ln\frac{1}{\mu_{N,\eps}}\right)
\end{array}\right\}\label{est:25}
\end{equation}
Since $\mu_{N,\eps}=o(\beta_\eps)$ and $\alpha_\eps=o(1)$ as $\eps\to 0$, it follows from Proposition \ref{prop:11} that
\begin{equation}\label{asymp:ue}
\lim_{\eps\to 0}\sup_{x\in \partial\rnm\cap B_{\alpha_\eps}(0)\setminus B_{\beta_\eps}(0)}\left|\frac{|x|^{2\bp}|\nabla (\ue\circ\T)|^2(x)}{\mu_{N,\eps}^{\bp-\bm}}-K^2\right|=0
\end{equation}
We fix $i,j\in \{2,...,n\}$. It follows from \eqref{asymp:ue} and $\bp-\bm=1$ that
\begin{eqnarray}
&&\ds \int \limits_{\partial\rnm\cap B_{\alpha_\eps}(0)\setminus B_{\beta_\eps}(0)}x^ix^j  \partial_{ij} \T_{0}(0)|\nabla(\ue\circ\T)|^2\, dx\nonumber\\
&&= \ds \int \limits_{\partial\rnm\cap B_{\alpha_\eps}(0)\setminus B_{\beta_\eps}(0)}\mu_{N,\eps}\frac{x^ix^j   \partial_{ij} \T_{0}(0)}{|x|^{2\bp}}K^2\, dx\nonumber\\
&&\ds +\int \limits_{\partial\rnm\cap B_{\alpha_\eps}(0)\setminus B_{\beta_\eps}(0)}\mu_{N,\eps}\frac{x^ix^j  \partial_{ij} \T_{0}(0) }{|x|^{2\bp}}\left(\frac{|x|^{2\bp}|\nabla(\ue\circ\T)|^2}{\mu_{N,\eps}}-K^2\right)\, dx\nonumber\\
&&\ds =\int \limits_{\partial\rnm\cap B_{\alpha_\eps}(0)\setminus B_{\beta_\eps}(0)}\mu_{N,\eps}\frac{x^ix^j   \partial_{ij} \T_{0}(0) }{|x|^{2\bp}}K^2\, dx\nonumber\\
&&+o\left(\int_{\partial\rnm\cap B_{\alpha_\eps}(0)\setminus B_{\beta_\eps}(0)}\mu_{N,\eps}\frac{|x|^2}{|x|^{2\bp}}\, dx\right)\label{est:234}
\end{eqnarray}
Independently, with a change of variable and $2\bp=n+1$, we get that
\begin{align*}
\ds \int \limits_{\partial\rnm\cap B_{\alpha_\eps}(0)\setminus B_{\beta_\eps}(0)}\frac{x^ix^j   \partial_{ij} \T_{0}(0) }{|x|^{2\bp}}\, dx&= \partial_{ij} \T_{0}(0) \left(\int_{\beta_\eps}^{\alpha_\eps}\frac{dr}{r}\right)\left(\int_{\mathbb{S}^{n-2}}\sigma^i\sigma^j\, d\sigma\right) \\
 &=\delta_{ij}   \partial_{ij} \T_{0}(0)\frac{ \omega_{n-2}}{n-1}\ln\frac{\alpha_\eps}{\beta_\eps}, 
\end{align*}
where $\omega_{n-2}$ is the volume of the round $(n-2)-$unit sphere. This equality, \eqref{est:234} and the properties \eqref{ppty:a:b} yield
\begin{eqnarray}\label{est:26}
&&\ds \int \limits_{\partial\rnm\cap B_{\alpha_\eps}(0)\setminus B_{\beta_\eps}(0)}x^ix^j \delta_{ij}   \partial_{ij} \T_{0}(0)|\nabla(\ue\circ\T)|^2\, dx\nonumber\\
&&=\delta_{ij}  \partial_{ij} \T_{0}(0) \frac{K^2\omega_{n-2}}{n-1}\mu_{N,\eps}\ln\frac{1}{\mu_{N,\eps}}+o\left(\mu_{N,\eps}\ln\frac{1}{\mu_{N,\eps}}\right).
\end{eqnarray}
Therefore, plugging \eqref{est:24}, \eqref{est:25} and \eqref{est:26} into \eqref{est:23} yields
\begin{eqnarray*}
&& \int \limits_{   \T\left( \partial \rnm \cap B_{ \delta_0 }(0) \setminus B_{ k_{1,\eps}^{3} }(0) \right) }    ( x, \nu)  \frac{|\nabla \ue|^2}{2} \,d\sigma\\
&&= -\frac{K^2\omega_{n-2}  \sum_{i=2}^{n} \partial_{ii} \T_{0}(0)}{4(n-1)}\mu_{N,\eps}\ln\frac{1}{\mu_{N,\eps}}+o\left(\mu_{N,\eps}\ln\frac{1}{\mu_{N,\eps}}\right)\\
&&=\frac{K^2\omega_{n-2}  \sum_{i=2}^{n} II_{0,ii}}{4(n-1)}\mu_{N,\eps}\ln\frac{1}{\mu_{N,\eps}}+o\left(\mu_{N,\eps}\ln\frac{1}{\mu_{N,\eps}}\right)\\
&&=  \frac{K^2\omega_{n-2} H(0)}{4(n-1)}\mu_{N, \eps}\ln\frac{1}{\mu_{N,\eps}}+o\left(\mu_{N,\eps}\ln\frac{1}{\mu_{N,\eps}}\right).
\end{eqnarray*}
Plugging this latest estimate into \eqref{est:5} yields \eqref{est:inter}. This ends the proof of Step \ref{step:p13}.\qed

\section{\, Proof of multiplicity}\label{sec:proof:th}
\noindent {\bf Proof of Theorem \ref{th:cpct:sc}:} We fix $\gamma<n^2/4$ and $h\in C^1(\overline{\Omega})$ such that $-\Delta-\gamma|x|^{-2}-h$ is coercive. For each  $2<p\leq \crit$, we consider the $C^2$-functional
\begin{equation*}
I_{p, \gamma}(u)=\frac{1}{2}\int_\Omega\left(\vert \nabla u\vert^2\, dx-\frac{\gamma}{2} {|u|^2}{|x|^2}-hu^2\right)\, dx-
\frac{1}{p}\int_{\Omega} \frac {|u|^p}{|x|^s}\, dx
\end{equation*}
on $\huno$, whose critical points are the weak solutions of
\begin{equation}
\label{main}
\left\{ \begin{array}{lll}
-\Delta u-\frac{\gamma}{|x|^2}u-h u&= \frac{|u|^{p-2}u}{|x|^s} &{\rm on} \ \Omega \\
  \hfill   u  &=   0  &{\rm on} \ \partial\Omega.
\end{array} \right.
\end{equation}
For a fixed $u\in \huno$, $u\not\equiv 0$, we have that
$$I_{p, \gamma}(\lambda u)=\frac{\lambda^2}{2}\int_\Omega\vert \nabla u\vert^2\, dx-\frac{\gamma \lambda^2}{2}\int_{\Omega} \frac {|u|^2}{|x|^2}\, dx-\lambda^2\int_\Omega hu^2\, dx-
\frac{\lambda^p}{p}\int_{\Omega} \frac {|u|^p}{|x|^s}\, dx$$
Then, since coercivity holds, we have that that $\lim_{\lambda \to \infty}I_{p, \gamma}(\lambda u)=-\infty$, which means 
that for each finite dimensional subspace $E_k \subset E:=\huno$, there 
exists $R_k>0$ such that
\begin{equation}
\label{neg}
\sup \{I_{p, \gamma}(u); u\in E_k, \Vert u\Vert_{H_1^2} >R_k\} <0
\end{equation}
when $p\to\crit (s)$. Let $(E_k)^\infty_{k=1}$ be an increasing sequence of subspaces of
$\huno$ such that
$\dim E_k=k $ and $\overline{\cup^\infty_{k=1} E_k}=E:=\huno$ and define
the min-max values:
$$c_{p,k} = {\ds \inf_{g \in {\bf H}_k} \sup_{x\in E_k} I_{p, \gamma}(g(x))},$$
where
$$ {\bf H}_k=\{g \in C(E,E); \hbox{ $g$ is odd and $g(v)=v$
   for $\|v\| > R_k$ for some $R_k>0 \}$}.$$
\begin{proposition} \label{critical.values} With the above notation and assuming  $n\ge 3$, we have:
  \begin{enumerate}
  \item For each $k\in \nn$, $c_{p,k}>0$ and  $\lim\limits_{p\to \crits} c_{p,k}=c_{\crits,k}:=c_k.$
   \item If $2<p<{\crits}$, there exists for each $k$, functions  $u_{p,k}\in
\huno$ such that $I'_{p, \gamma}(u_{p,k})=0$, and $I_{p, \gamma}(u_{p,k})=c_{p,k}$.
   \item  For each  $2<p < {\crits}$, we have  
$ c_{p,k}\ge D_{n,p} k^{\frac{p+1}{p-1}\frac{2}{n}}$ where $D_{n,p}>0$ is 
such that
$\lim \limits_{p\to  {\crits}}D_{n,p}=0$.
\item $\lim\limits_{k\to \infty} c_k=\lim\limits_{k\to \infty} 
c_{\crits,k}=+\infty$.
\end{enumerate}
\end{proposition}
\noindent {\bf Proof:} (1) Coercivity yields the existence of $a_0>0$ such that
\begin{equation}\label{coer}
\int_\Omega\left(|\nabla u|^2-\frac{\gamma}{|x|^2}u^2-hu^2\right)\, dx\geq a_0\int_\Omega|\nabla u|^2\, dx\hbox{ for all }u\in \huno.
\end{equation}
With \eqref{coer}, the Hardy and the Hardy-Sobolev inequality \eqref{HS-ineq}, there exists $C>0$ and $\alpha>0$ such that
\[
I_{p, \gamma}(u)\geq \frac{a_0}{2}\|\nabla u\|_2^2-C\|\nabla u\|_2^p =\|\nabla 
u\|_2^2\left(\frac{a_0}{2}-C\|\nabla u\|_2^{p-2}\right)\geq \alpha >0
\]
for all $u\in\huno$ such that provided   $\Vert \nabla u\Vert_2= \rho$ for some $\rho>0$ small enough.  
Then the sphere $S_\rho=\{u\in E; \|u\|_{\huno}= \rho\}$ intersects every  image $g(E_k)$  by an odd continuous function $g$. It follows that
\[
c_{p,k}\geq \inf \{I_{p, \gamma}(u); u\in S_\rho \} \geq \alpha >0.
\]
In view of (\ref{neg}), it follows that for each $g\in {\bf H}_k$, we have 
that
\[
\sup\limits_{x\in
E_k}I_{p_i, \gamma}(g(x))=\sup\limits_{x\in D_k}I_{p, \gamma}(g(x))
\]
  where $D_k$ denotes the
ball in $E_k$ of radius $R_k$.
  Consider now a sequence $p_i \to \crits$ and note
first that for each $u\in E$,  we have that $I_{p_i, \gamma}(u) \to I_{\crits, \gamma}
(u)$.   Since $g(D_k)$ is compact and the family of
functionals $(I_{p, \gamma})_p$ is equicontinuous, it follows that
   $\sup\limits_{x\in E_k}I_{p, \gamma}(g(x))\to \sup\limits_{x\in E_k}I_{\crits, \gamma
}(g(x))$, from which follows that
   $\limsup\limits_{i\in \nn}c_{p_i,k}\leq \sup\limits_{x\in E_k}I_{\crits, \gamma}
(g(x))$. Since this holds for any $g\in {\bf H}_k$, it follows that
   \[
    \limsup\limits_{i\in \nn}c_{p_i,k}\leq c_{\crits, k}=c_k.
   \]
On the other hand, the function $f(r)=\frac{1}{p}r^p-\frac{1}{\crits}r^{\crits}$ 
attains its maximum on $[0, +\infty)$ at $r=1$ and therefore
$f(r) \leq \frac{1}{p}-\frac{1}{\crits}$ for all $r>0$. It follows
  \begin{eqnarray*}   I_{\crits, \gamma}(u) &= &I_{p, \gamma}(u)+\int_\Omega \frac{1}{|x|^s}
\left(\frac{1}{p}|u(x)|^p-\frac{1}{\crits}|u(x)|^{\crits }\right)\, dx\\
&\leq&
I_{p, \gamma}(u)+\int_\Omega \frac{1}{|x|^s} \left(\frac{1}{p}-\frac{1}{\crits}\right)\, dx
   \end{eqnarray*}
from which follows that $c_k\leq  \liminf\limits_{i\in \nn}c_{p_i,k}$, and 
claim (1) is proved. \\
If now $p< {\crits}$, we are in the subcritical case, that is we have
compactness in the Sobolev embedding $\huno \to L^p(\Omega; |x|^{-s}dx)$ 
and
therefore $I_{p, \gamma}$ has the Palais-Smale condition. It is then standard to
find critical points $u_{p,k}$ for $I_{p, \gamma}$ at each level $c_{p,k}$  (see for
example the book \cite{Gh}). Consider now the functional 
\begin{equation*}
I_{p, 0}(u)=\frac{1}{2}\int_\Omega\vert \nabla u\vert^2\, dx- 
\frac{1}{p}\int_{\Omega} \frac {|u|^p}{|x|^s}\, dx
\end{equation*}
and its critical values 
$$c^0_{p,k} = {\ds \inf_{g \in {\bf H}_k} \sup_{x\in E_k} I_{p, 0}(g(x))}.$$
It has been shown in \cite{gr2} that (1), (2) and (3) of Proposition \ref{critical.values} hold, with  
$c^0_{p,k}$ and $c^0_{k}$ replacing $c_{p,k}$ and  $c_{k}$ respectively.  In particular, $\lim\limits_{k\to \infty} c^0_k=\lim\limits_{k\to \infty} 
c^0_{\crits,k}=+\infty$.

\noindent On the other hand, with the coercivity \eqref{coer}, we have that 
\[
I_{p, \gamma}(u) \geq a_0^{\frac{p}{p-2}} I_{p,0}(v) \quad \hbox{for every $u\in \huno$,}
\]
where  $v=a_0^{-\frac{1}{p-2}}u$.  It then follows that $\lim\limits_{k\to \infty} c_k=\lim\limits_{k\to \infty} 
c_{\crits,k}=+\infty$.

\smallskip\noindent To complete the proof of Theorem  \ref{th:cpct:sc},   notice that since for each $k$, we
have $$ \lim\limits_{p_i\to \crits}I_{p_i, \gamma}(u_{p_i,k})=\lim\limits_{p_i\to\crits}c_{p_i,k}=c_k,$$
it follows that the sequence  $(u_{p_i,k})_i $ is
uniformly bounded in $\huno$. Moreover, since $I_{p_i}'(u_{p_i,k})=0$, it 
follows from the compactness result that by letting $p_i\to \crits$, we
get a solution $u_k$ of (\ref{main}) in such a way that
$I_{\crit(s), \gamma}(u_k)=\lim\limits_{p\to \crits}I_{p, \gamma}(u_{p,k})=\lim\limits_{p\to
\crits}c_{p,k}=c_k$. Since the latter sequence goes to infinity, it follows 
that   (\ref{main}) has  an infinite number of critical levels. 

\section{Proof of the non-existence result}\label{sec:nonex}  

\medskip\noindent{\bf Proof of Theorem \ref{thm:non:ter}:} We argue by contradiction. We fix $\gamma<\gamma_H(\Omega)\leq\frac{n^2}{4}$ and $\Lambda>0$. We assume that there is a family $(\ue)_{\eps>0}\in \huno$ of solutions to 
\begin{equation}\label{eq:non}\left\{ \begin{array}{cl}
-\Delta \ue-\gamma \frac{\ue}{|x|^2}-h_\eps \ue
= \frac{\ue^{\crits-1}}{|x|^s}  &\text{in } \Omega,\\
\ue>0 &\hbox{ in }\Omega\\
\ue=0 &\hbox{ on }\bono
\end{array}\right.\end{equation}
with $\Vert \nabla \ue\Vert_2\leq\Lambda$ and $\lim_{\eps\to 0}h_\eps=h_0$ in $C^1(\overline{\Omega})$.

\medskip\noindent We claim that $(\ue)_{\eps>0}$ is not pre-compact in $\huno$. Otherwise, up to extraction, there would be $u_0\in \huno$, $u_0\geq 0$, such that $\ue\to u_0$ in $\huno$ as $\eps\to 0$. Passing to the limit in the equation, we get that $u_0\geq 0$ and
\begin{equation}\label{eq:non:2}\left\{ \begin{array}{cl}
-\Delta u_0-\gamma \frac{u_0}{|x|^2}-h_0 u_0= \frac{u_0^{\crits-1}}{|x|^s}  &\text{in } \Omega,\\
u_0\geq 0 &\hbox{ in }\Omega\\
u_0=0 &\hbox{ on }\bono.
\end{array}\right.\end{equation}
The coercivity of $-\Delta u_0-\gamma |x|^{-2}-h_0$ and the convergence of $(h_\eps)_\eps$ yield
\begin{eqnarray*}
&&C\left(\int_{\Omega}\frac{\ue^{\crits}}{|x|^s}dx\right)^{{2}/{\crits}}\leq \int_{\Omega} |\nabla \ue|^2\,dx-\int_\Omega\left(\frac{\gamma}{|x|^2}+\he\right)\ue^2\, dx\leq\int_{\Omega}\frac{\ue^{\crits}}{|x|^s}dx,
\end{eqnarray*}
for small $\epsilon>0$, and then, since $\ue>0$, there exists $c_0>0$ such that 
$$\int_{\Omega}\frac{\ue^{\crits}}{|x|^s}dx\geq c_0$$
for all $\eps>0$. Passing to the limit yields $u_0\not\equiv 0$. Therefore, $u_0>0$ is a solution to \eqref{eq:non} with $\epsilon=0$. 
This is not possible simply by the hypothesis. 

The family $(\ue)_\eps$ is not pre-compact and it therefore blows-up with bounded energy. Let $u_0\in \huno$ be its weak limit, which is necessarily a solution to \eqref{eq:non:2}, and hence must be the trivial solution $u_0\equiv 0$.  Proposition \ref{prop:rate:sc:2} then yields that 
either 
\begin{equation}
\hbox{$\bp-\bm\geq 1$ and therefore $H(0)=0$,}
\end{equation}
or 
\begin{equation}\label{res:blowup}
\hbox{$\bp-\bm<1$ and therefore  $m_{\gamma,h_0}(\Omega)=0$.}
\end{equation}
It now suffices to note that when $\gamma\leq (n^2-1)/4$ then $\bp-\bm\geq 1$ and the above contradicts our assumption that $H(0)\neq 0$. Similarly, if $\gamma>(n^2-1)/4$, then $\bp-\bm<1$ and the above contradicts our assumption that the mass is non-zero. In either case,  this means that no such a family of positive solutions $(\ue)_{\eps>0}$ exist. \qed 

\medskip\noindent{\bf Proof of Corollary \ref{thm:non}:} First note that if $h_0$ satisfies 
\begin{equation}\label{hyp:h00}
h_0(x)+\frac{1}{2}(\nabla h_0(x), x)\leq 0\hbox{ for all }x\in \Omega, 
\end{equation}
then by differentiating for any $x\in \Omega$, the function $t\mapsto t^2 h_0(tx)$ (which is well defined for $t\in [0,1]$ since $\Omega$ is starshaped), we get that $h_0\leq 0$. Therefore $-\Delta-\gamma|x|^{-2}-h_0$ is coercive.\\
Assume now there is a positive variational solution $u_0$ corresponding to $h_0$, the Pohozaev identity \eqref{PohoId} then gives
$$\int_{\partial\Omega}(x,\nu)\frac{(\partial_\nu u_0)^2}{2}\, d\sigma-\int_\Omega \left(h_0+\frac{1}{2}(\nabla h_0, x)\right)u_0^2\, dx=0.$$
Hopf's strong comparison principle yields $\partial_\nu u_0 <0$. Since $\Omega$ is starshaped with respect to $0$, we get that $(x,\nu)\geq 0$ on $\partial\Omega$. Therefore, with \eqref{hyp:h00}, we get that $(x,\nu)=0$ for all $x\in \Omega$, which is a contradiction since $\Omega$ is smooth and bounded. \\
If now $\gamma\leq (n^2-1)/4$, the result follows from Theorem \ref{thm:non:ter} since we have assumed that $H(0)\neq 0$. \\
If  $\gamma > (n^2-1)/4$, we use 
Theorem 7.1 in Ghoussoub-Robert \cite{gr4} to find $\mathcal{K}\in C^2(\overline{\Omega}\setminus\{0\})$ and $A>0$ such that 
$$\left\{\begin{array}{ll}
-\Delta\mathcal{K}-\frac{\gamma}{|x|^2}\mathcal{K}-h_0 \mathcal{K}=0 &\hbox{ in }\Omega\\
\mathcal{K}>0 &\hbox{ in }\Omega\\
\mathcal{K}=0 &\hbox{ on }\partial\Omega\setminus\{0\}.
\end{array}\right.$$
and  such that 
$$\mathcal{K}(x)=A\left(\frac{\eta(x)d(x,\partial\Omega)}{|x|^{\bp}}+\beta(x)\right)\hbox{ for all }x\in\Omega,$$
where $\eta\in C^\infty_c(\rn)$ and $\beta\in \huno$ are as in Step \ref{step:p11}. We now apply the Pohozaev identity \eqref{PohoId} to $\mathcal{K}$ on the domain $U:=\Omega\setminus \T(B_\delta(0))$ for $\T$ as in \eqref{def:T:bdry}: using that $\mathcal{K}^2\in L^1(\Omega)$ and $(\cdot,\nu)(\partial_\nu \mathcal{K})^2\in L^1(\partial\Omega)$ when $\bp-\bm<1$, we get that
$$\int_{\partial\Omega} (x,\nu)\frac{(\partial_\nu \mathcal{K})^2}{2}\, d\sigma-\int_\Omega \left(h_0+\frac{1}{2}(\nabla h_0, x)\right)\mathcal{K}^2\, dx=M_{\delta}$$
where $M_\delta$ is defined in \eqref{est:inf:1}. With \eqref{est:Md}, we then get
$$\int_{\partial\Omega} (x,\nu)\frac{(\partial_\nu \mathcal{K})^2}{2}\, d\sigma-\int_\Omega \left(h_0+\frac{1}{2}(\nabla h_0, x)\right)\mathcal{K}^2\, dx=- \frac{\omega_{n-1}}{n}\left(\frac{n^2}{4}-\gamma\right)A^2\cdot m_{\gamma,h_0}(\Omega).$$
Since $\Omega$ is star-shaped and $h_0$ satisfies (\ref{hyp:h00}), it follows that $m_{\gamma,h_0}(\Omega)<0$ and Theorem \ref{thm:non:ter} then applies to complete our corollary.

\section{Appendix A: The Pohozaev identity}\label{sec:app:poho}
\begin{proposition} Let $U\subset\rn$ be a smooth bounded domain and let $u\in C^2(\overline{U})$ be a solution of 
\begin{equation}\label{poh}
- \Delta u -\gamma \frac{u}{|x|^2}-hu =K\frac{|u|^{\crits -2-p}}{|x|^{s}}u \quad \hbox{on $U$}.
\end{equation}
 Then,  we have 
\begin{eqnarray*}
  &&- \int \limits_{U} \left( h(x) +\frac{ \left( \nabla h, x \right)}{2} \right)u^{2}~ dx  
  ~ - \frac{p}{\crits} \left(\frac{n-s}{\crits -p}\right) \int \limits_{U } K  \frac{|u|^{\crits-p }}{|x|^{s}} ~ dx \\
 &&=  \int \limits_{\partial U}  F(x)~d\sigma, 
\end{eqnarray*}
where
\begin{eqnarray*}
F(x)&:=&( x, \nu) \left( \frac{|\nabla u|^2}{2}  -\frac{\gamma}{2}\frac{u^{2}}{|x|^{2}} - \frac{h(x)}{2} u^{2}  - \frac{K}{\crits -p} \frac{|u|^{\crits-p }}{|x|^{s}} \right)  \\
&& -  \left(x^i\partial_i u+\frac{n-2}{2} u \right)\partial_\nu u. 
\end{eqnarray*}
\end{proposition}
\smallskip\noindent{\it Proof:}
For any $y_{0}\in \rn$, the classical Pohozaev identity yields
\begin{eqnarray*}
&&-\int \limits_{U}\left((x-y_{0})^i\partial_i u+\frac{n-2}{2}u\right)\Delta u~ dx\\
&& = \int \limits_{\partial U}\left[(x-y_{0},\nu)\frac{|\nabla u|^2}{2}-\left((x-y_{0})^i\partial_i u+\frac{n-2}{2}u\right)\partial_\nu u\right]d\sigma,
\end{eqnarray*}
where $\nu$ is the outer normal to the boundary $\partial U$. 
\medskip

\noindent
One has  for $1 \leq j \leq n$
\begin{align*}
\partial_{j} \left( \frac{|u|^{\crits-p }}{|x|^{s}} \right) =  -s \frac{x^{j}}{|x|^{s+2}} |u|^{\crits -p}+ (\crits -p) \frac{|u|^{\crits -2-p}}{|x|^{s}} u \partial_{j}u
\end{align*}
So
\begin{eqnarray*}
&& \left( x-y_{0}, \nabla u \right)\frac{|u|^{\crits -2-p}}{|x|^{s}}u = \frac{1}{\crits -p}(x-y_{0})^{j}\partial_{j} \left( \frac{|u|^{\crits-p }}{|x|^{s}} \right) \\
&&+ \frac{s}{\crits -p} \frac{|u|^{\crits -p}}{|x|^{s}}  -  \frac{s}{\crits -p} \frac{(x,y_{0})}{|x|^{s+2}}|u|^{\crits -p}.
\end{eqnarray*}
Then integration by parts yields
\begin{align*}
\int \limits_{U}  \left( x-y_{0}, \nabla u \right)\frac{|u|^{\crits -2-p}}{|x|^{s}}u ~ dx= \frac{1}{\crits -p} \int \limits_{U} (x-y_{0})^{j}\partial_{j} \left( \frac{|u|^{\crits-p }}{|x|^{s}} \right) dx \\
+ \frac{s}{\crits -p} \int \limits_{U}\frac{|u|^{\crits -p}}{|x|^{s}} dx  - \frac{s}{\crits -p} \int \limits_{U}\frac{(x,y_{0})}{|x|^{s+2}}|u|^{\crits -p} dx
\end{align*}
\begin{align*}
=& -\frac{n-s}{\crits -p} \int \limits_{U}   \frac{|u|^{\crits-p }}{|x|^{s}} ~ dx - \frac{s}{\crits -p} \int \limits_{U}\frac{(x,y_{0})}{|x|^{s+2}}|u|^{\crits -p} dx 
  \\ &  +  \frac{1}{\crits -p} \int \limits_{\partial U} ( x-y_{0}, \nu) \frac{|u|^{\crits-p }}{|x|^{s}}~ d\sigma.
\end{align*}
\noindent
Similarly, 
\begin{align*}
 \left( x-y_{0}, \nabla u \right)\frac{u}{|x|^{2}} = \frac{1}{2}(x-y_{0})^{j}\partial_{j} \left( \frac{u^{2}}{|x|^{2}} \right) +\frac{u^{2}}{|x|^{2}}  -  \frac{(x,y_{0})}{|x|^{4}}u^{2}
\end{align*}
\begin{align*}
\int \limits_{U}   \left( x-y_{0}, \nabla u \right)\frac{u}{|x|^{2}}  ~ dx =& -\frac{n-2}{2} \int \limits_{U}  \frac{u^{2}}{|x|^{2}} ~ dx  - \int \limits_{U}\frac{(x,y_{0})}{|x|^{4}}u^{2} dx  \notag \\ 
&+ \frac{1}{2}\int \limits_{\partial U} ( x-y_{0}, \nu) \frac{u^{2}}{|x|^{2}}~ d\sigma
\end{align*}
and 
\begin{align*}
\int \limits_{U}   \left( x-y_{0}, \nabla u \right) h(x) u ~ dx = & -\frac{n}{2} \int \limits_{U}  h(x) u^{2}~ dx 
-\frac{1}{2}\int \limits_{U}    \left( \nabla h, x-y_{0}\right) u^{2} ~ dx \\ &+ \frac{1}{2}\int \limits_{\partial U} ( x-y_{0}, \nu) h(x)u^{2}~ d\sigma
\end{align*}
\noindent
Combining the above, we obtain for any $K$ and any $y_0\in \rn$, 
\begin{align}\label{PohoId}
&\int \limits_{U}\left((x-y_{0})^i\partial_i u+\frac{n-2}{2}u\right)\left(- \Delta u -\gamma \frac{u}{|x|^2}-hu -K\frac{|u|^{\crits -2-p}}{|x|^{s}}u  \right)~ dx \notag \\
&\quad  - \int \limits_{U}  h(x) u^{2}~ dx  -  \frac{1}{2}\int \limits_{U}    \left( \nabla h, x-y_{0}\right) u^{2} ~ dx \notag\\
&\quad - \frac{p}{\crits} \left( \frac{n-s}{\crits -p}\right) \int \limits_{U}   K\frac{|u|^{\crits-p }}{|x|^{s}} ~ dx \notag \\
&\quad - \gamma \int \limits_{U}\frac{(x,y_{0})}{|x|^{4}}u^{2} dx - \frac{s}{\crits -p} \int \limits_{U}\frac{(x,y_{0})}{|x|^{s+2}}K|u|^{\crits -p} dx \notag \\
&  =  \int \limits_{\partial U}\left[( x-y_{0}, \nu) \left( \frac{|\nabla u|^2}{2}  -\frac{\gamma}{2}\frac{u^{2}}{|x|^{2}} - \frac{h(x)}{2} u^{2}  - \frac{K}{\crits -p} \frac{|u|^{\crits-p }}{|x|^{s}} \right) \right]d\sigma \notag \\
&\quad - \int \limits_{\partial U} \left[ \left((x-y_{0})^i\partial_i u+\frac{n-2}{2}u\right)\partial_\nu u \right]d\sigma. 
\end{align}
We conclude by taking $y_0=0$ and using that $u$ satisfies (\ref{poh}) on $U$. 
\section[Appendix B: A continuity property for the first eigenvalue]{Appendix B: A continuity property of the first eigenvalue of Schr\"odinger operators}\label{sec:app:lemma}
 
\begin{lemma}\label{lem:vp} Let $\Omega\subset\rn$, $n\geq 3$, be a smooth bounded domain. Let $(V_k)_k:\Omega\to \rr$ and $V_\infty:\Omega\to\rr$ be measurable functions and let {$(x_k)_k\in \overline{\Omega}$} be a sequence of points. We assume that
\begin{enumerate}[i)]
\item $\lim \limits_{k\to +\infty}V_k(x)=V_\infty(x)$ for a.e. $x\in\Omega$,
\item There exists $C>0$ such that $|V_k(x)|\leq C|x-x_k|^{-2}$ for all $k\in\nn$ and $x\in\Omega$.
\item $\lim \limits_{k\to +\infty}x_k=0\in\partial\Omega.$
\item For some $\gamma_0<n^2/4$, there exists $\delta>0$ such that $|V_k(x)|\leq \gamma_0|x-x_k|^{-2}$ for all $k\in\nn$ and {$x\in B_\delta(0)\cap\Omega$}.  
\item The first eigenvalue $\lambda_1(-\Delta+V_k)$ is achieved for all $k\in\nn$.
\end{enumerate}
Then,
\begin{equation}
\lim_{k\to +\infty}\lambda_1(-\Delta+V_k)=\lambda_1(-\Delta+V_\infty).
\end{equation}
\end{lemma}
\smallskip\noindent{\it Proof:} 
 We first claim that $(\lambda_1(-\Delta+V_k))_k$ is bounded. Indeed, fix $\varphi\in \huno\setminus\{0\}$ and use the Hardy inequality to write for all $k\in\nn$, 
$$\lambda_1(-\Delta+V_k)\leq \frac{\int_\Omega (|\nabla \varphi|^2+V_k\varphi^2)\, dx}{\int_\Omega \varphi^2\, dx}\leq \frac{\int_\Omega (|\nabla \varphi|^2+C|x-x_k|^{-2}\varphi^2)\, dx}{\int_\Omega \varphi^2\, dx}:=M<+\infty$$
For the lower bound, we have for any $\varphi\in \huno$,
\begin{eqnarray}
\int_\Omega (|\nabla \varphi|^2+V_k\varphi^2)\, dx&= & \int_\Omega |\nabla \varphi|^2\, dx+\int_{B_\delta(0)}V_k\varphi^2\, dx+\int_{\Omega\setminus B_\delta(0)}V_k\varphi^2\, dx\nonumber \\
&\geq & \int_\Omega |\nabla \varphi|^2\, dx-\gamma_0\int_{B_\delta(0)}|x-x_k|^{-2}\varphi^2\, dx\nonumber\\
&&-4C\delta^{-2}\int_{\Omega\setminus B_\delta(0)}\varphi^2\, dx\nonumber \\
&\geq & \left(1-4\gamma_0/n^2\right)\int_\Omega |\nabla \varphi|^2\, dx-4C\delta^{-2}\int_{\Omega}\varphi^2\, dx.\label{ineq:123}
\end{eqnarray}
Since $\gamma_0<n^2/4$, we then get that $\lambda_1(-\Delta+V_k)\geq -4C\delta^{-2}$ for large $k$, which proves the lower bound.  

\smallskip\noindent Up to a subsequence, we can now assume that $(\lambda_1(-\Delta+V_k))_k$ converges as $k\to +\infty$. We now show that 
\begin{equation}\label{ineq:vp:100}
\liminf_{k\to +\infty}\lambda_1(-\Delta+V_k)\geq \lambda_1(-\Delta+V_\infty).
\end{equation}
For $k\in\nn$, we let $\varphi_k\in \huno$ be a minimizer of $\lambda_1(-\Delta+V_k)$ such that $\int_\Omega\varphi_k^2\, dx=1$. In particular,
\begin{equation}\label{eq:phi:vp}
-\Delta\varphi_k+V_k\varphi_k=\lambda_1(-\Delta+V_k)\varphi_k\hbox{ weakly in }\huno.
\end{equation}
Inequality \eqref{ineq:123} above yields the boundedness of $(\varphi_k)_k$ in $\huno$. Up to a subsequence, we let $\varphi\in \huno$ such that, as $k\to +\infty$, $\varphi_k\rightharpoonup \varphi$ weakly in $\huno$, $\varphi_k\to \varphi$ strongly in $L^2(\Omega)$ (then $\int_\Omega\varphi^2\, dx=1$) and $\varphi_k(x)\to\varphi(x)$ for a.e. $x\in \Omega$. Letting
 $k\to +\infty$ in \eqref{eq:phi:vp}, the hypothesis on $(V_k)$ allow us to conclude that 
$$-\Delta\varphi+V_\infty\varphi=\lim_{k\to +\infty}\lambda_1(-\Delta+V_k)\varphi\hbox{ weakly in }\huno.$$
Since $\int_\Omega\varphi^2\, dx=1$ and we have extracted subsequences, we then get \eqref{ineq:vp:100}.

\medskip\noindent Finally, we prove the reverse inequality. For $\eps>0$, let $\varphi\in \huno$ be such that 
$$ \frac{\int_\Omega (|\nabla \varphi|^2+V_\infty\varphi^2)\, dx}{\int_\Omega \varphi^2\, dx}\leq \lambda_1(-\Delta+V_\infty)+\eps.$$
We have 
$$\lambda_1(-\Delta+V_k)\leq \lambda_1(-\Delta+V_\infty)+\eps +\frac{\int_\Omega |V_k-V_\infty|\varphi^2\, dx}{\int_\Omega \varphi^2\, dx}.$$
The hypothesis of Lemma \ref{lem:vp} allow us to conclude that $\int_\Omega |V_k-V_\infty|\varphi^2\, dx\to 0$ as $k\to +\infty$. Therefore $\limsup_{k\to +\infty}\lambda_1(-\Delta+V_k)\leq \lambda_1(-\Delta+V_\infty)+\eps$ for all $\eps>0$. Letting $\eps\to 0$, we get the reverse inequality and 
 the conclusion of Lemma \ref{lem:vp}.\qed

\section[Appendix C: Regularity and the Hardy-Schr\"odinger operator]{Appendix C: Regularity and the Hardy-Schr\"odinger operator on $\R^{n}_{-}$}\label{sec:app:regul}

In this section, we  collect some important results from the paper \cite{gr4}  used in the proof of the compactness theorems. First we state the following regularity result: 
\begin{theorem}[\cite{gr4}, see also \cite{ff} ]\label{th:hopf} 
Let $\Omega$ be a smooth bounded domain of $\mathbb{R}^{n}$ $(n \geq 3)$  such that $0 \in \partial \Omega$.  We fix $\gamma<\frac{n^2}{4}$ and $f: \Omega\times \rr\to \rr$ is a Caratheodory function such that
$$|f(x, v)|\leq C|v| \left(1+\frac{|v|^{\crits-2}}{|x|^s}\right)\hbox{\rm  for all }x\in \Omega\hbox{ \rm and }v\in \rr.$$
Let $u\in \huno$ be a weak solution of 
\begin{align}\label{regul:eq}
-\Delta u-\frac{\gamma+O(|x|^\theta)}{|x|^2}u=f(x,u) \qquad  \hbox{ in } \left(\huno \right)'
\end{align}
for some $\theta>0$. Then there exists $K\in\rr$ such that
\begin{align}
\label{eq:hopf}
\lim_{x\to 0}\frac{u(x)}{d(x,\bdry) |x|^{-\bm}}=K.
\end{align}
Moreover, if $u\geq 0$ and $u\not\equiv 0$, we have that $K>0$.
\end{theorem}
The following result characterizes the positive solution to the singular global equation 
\begin{proposition}[\cite{gr4}]\label{prop:liouville} 
Let $\gamma<\frac{n^2}{4}$ and let $u\in C^2(\rn\setminus \{0\})$ be a nonnegative function such that
\begin{align*}
 \left\{ \begin{array}{llll}
-\Delta u -\frac{\gamma}{|x|^2} u &=&0  &\hbox{ in } \rnm \\
\hfill u&=&0 & \hbox{ on } \partial \rnmp .
\end{array}\right.
\end{align*}
Then there exist $C_-,C_+\geq 0$ such that
$$u(x)=C_{-}\frac{|x_{1}|}{|x|^{\bp}}+ C_{+} \frac{|x_{1}|}{|x|^{\bm}}  \qquad \hbox{ for  all } x\in \rnm.$$
\end{proposition}
\medskip

\noindent Next, we recall the existence and behaviour of the singular solution to the homogeneous equation.  
\begin{theorem}[\cite{gr4}] \label{th:sing} 
Let $\Omega$ be a smooth bounded domain of $\mathbb{R}^{n}$ $(n \geq 3)$  such that $0 \in \partial \Omega$.  Fix $\gamma<\frac{n^2}{4}$ and $h\in C^{1}(\overline\Omega)$ be such that the operator  $\Delta-\gamma|x|^{-2}-h$ is coercive.  There exists then  $\mathcal{H}\in C^2(\overline{\Omega}\setminus\{0\})$ such that 
$$\qquad\qquad \left\{\begin{array}{ll}
-\Delta \mathcal{H}-\frac{\gamma}{|x|^2}\mathcal{H}+h(x) \mathcal{H}=0 &\hbox{ in } \Omega\\
\hfill \mathcal{H}>0&\hbox{ in } \Omega\\
\hfill \mathcal{H}=0&\hbox{ on }\partial\Omega \setminus\{0\}.
\end{array}\right.$$
These solutions are unique up to a positive multiplicative constant, and there exists $c>0$ such that 
$\mathcal{H}(x)\simeq_{x\to 0}c\frac{d(x,\bdry)}{|x|^{\beta_+(\gamma)}}.$
\end{theorem}

\section[Appendix D: Green's function on a bounded domain]{Appendix D: Green's function for the Hardy-Schr\"odinger operator with boundary singularity on a bounded domain}\label{sec:app:c}

\begin{defi} Let $\Omega$ be a smooth bounded domain of $\rn$, $n\geq 3$, such that $0\in\partial\Omega$. We fix $\gamma<n^2/4$ and $h\in C^{0,\theta}(\overline{\Omega})$, $\theta\in (0,1)$ such that  $-\Delta-(\gamma|x|^{-2}+h)$ is coercive. We say that $G:\Omega\times\Omega\setminus\{(x,x)/\, x\in \Omega\}$ is a Green's function for  $-\Delta-\gamma|x|^{-2}-h$ if

\smallskip\noindent{\bf $\bullet$} For any $p\in\Omega$, $G_p:=G(p,\cdot)\in L^1(\Omega)$.

\smallskip\noindent{\bf $\bullet$} For all $f\in C^\infty_c(\Omega)$ and all $p\in \Omega$, then
\begin{equation*}
\varphi(p)=\int_{\Omega}G_p(x)f(x)\, dx.
\end{equation*}
where $\varphi\in \huno\cap C^0(\Omega)$ is the unique solution to
\begin{equation*}
-\Delta\varphi-\left(\frac{\gamma}{|x|^2}+h(x)\right)\varphi=f\hbox{ in }\Omega\; ; \; \varphi_{|\partial\Omega}=0.
\end{equation*}
\end{defi}
\noindent This appendix is devoted to the proof of the following result.  

\begin{theorem}\label{th:green:gamma:domain} Let $\Omega$ be a smooth bounded domain of $\rn$ such that $0\in\partial\Omega$. We fix $\gamma<\frac{n^2}{4}$. We let $h\in C^{0,\theta}(\Omegabar)$ be such that $-\Delta-\gamma|x|^{-2}-h$ is coercive. Then there exists a unique Green's function $G$ for $-\Delta-\gamma|x|^{-2}-h$. Moreover:

 \noindent{\bf I. Properties of $G$.} The Green's function $G$ is such that 

\smallskip\noindent{\bf (a)} $G_p\in C^{2,\theta}(\overline{\Omega}\setminus\{0,p\})$ and $G_p>0$ for all $p\in\Omega$. 

\smallskip\noindent{\bf (b)} For all $p\in\Omega$ and all $\eta\in C^\infty_c(\rn\setminus\{p\})$, we have that $\eta G_p\in \huno$.

\smallskip\noindent{\bf (c)} For all $f\in L^{\frac{2n}{n+2}}(\Omega)\cap L^q(\Omega\setminus B_\delta(0))$, for all $\delta>0$ and some $q>n/2$,  we have for any $p\in\Omega$
\begin{equation}\label{id:172}
\varphi(p)=\int_{\Omega}G_p(x)f(x)\, dx.
\end{equation}
where $\varphi\in \huno\cap C^0(\Omega)$ is the unique solution to
\begin{equation}\label{eq:G:dist}
-\Delta\varphi-\left(\frac{\gamma}{|x|^2}+h(x)\right)\varphi=f\hbox{ in }\Omega\; ; \; \varphi_{|\partial\Omega}=0,
\end{equation}
In particular,
\begin{equation}\label{eq:G:c2}
\left\{\begin{array}{ll}
-\Delta G_p-\left(\frac{\gamma}{|x|^2}+h(x)\right)G_p=0 &\hbox{ in }\Omega\setminus\{p\},\\
G_p>0  &\hbox{ in }\Omega\setminus\{p\},\\
G_p=0  &\hbox{ in }\partial\Omega\setminus\{0\}.
\end{array}\right.
\end{equation}

\noindent{\bf II. Asymptotics.}\label{th:green:gamma:asymp} $G$ satisfies the following properties:\\

\smallskip\noindent{\bf (d)} For all $p\in\Omega\setminus\{0\}$, there exists $c_0(p)>0$ such that
\begin{equation}\label{asymp:G}
G_p(x)\sim_{x\to 0} c_0(p)\frac{d(x,\partial\Omega)}{|x|^{\bm}}\hbox{ and }G_p(x)\sim_{x\to p}\frac{1}{(n-2)\omega_{n-1}|x-p|^{n-2}}
\end{equation}
where
$$\bm:=\frac{n}{2}-\sqrt{\frac{n^2}{4}-\gamma}\hbox{ and }\bp:=\frac{n}{2}+\sqrt{\frac{n^2}{4}-\gamma}.$$

\smallskip\noindent{\bf (e)}  There exists $c>0$ depending only on $\gamma$, the coercivity constant and an upper-bound for $\Vert h\Vert_{C^{0,\theta}}$ such that
\begin{equation}\label{est:G:up}
c^{-1}H_p(x)< G_p(x)< c H_p(x)\hbox{ for }x\in\Omega-\{0,p\},
\end{equation}
where
\begin{equation}\label{def:Hp:1}
H_p(x):=\left(\frac{\max\{|p|,|x|\}}{\min\{|p|,|x|\}}\right)^{\bm}|x-p|^{2-n}\min\left\{1,\frac{d(x,\partial\Omega)d(p,\partial\Omega)}{|x-p|^2}\right\}.
\end{equation}
And
\begin{equation}\label{ineq:grad:G}
|\nabla G_p(x)|\leq c \left(\frac{\max\{|p|,|x|\}}{\min\{|p|,|x|\}}\right)^{\bm}|x-p|^{1-n}\min\left\{1,\frac{d(p,\partial\Omega)}{|x-p|}\right\}\hbox{ for }x\in\Omega-\{0,p\}.
\end{equation}
\smallskip\noindent{\bf (f)} There exists $L_{\gamma,\Omega}>0$ such that for any $(h_i)_i\in C^{0,\theta}(\Omega)$ such that $\lim \limits_{i\to+\infty}h_i=h$ in $C^{0,\theta}$, then for any sequences $(x_i)_i,(y_i)_i\in \Omega$ such that 
$$y_i=o(|x_i|)\hbox{ and }x_i=o(1)\hbox{ as }i\to +\infty,$$
then, as $i\to +\infty$ we have that
\begin{equation}\label{est:G:4}
G_{h_i}(x_i,y_i)=(L_{\gamma,\Omega}+o(1))\frac{d(x_i,\partial\Omega)}{|x_i|^{\bp}}\frac{d(y_i,\partial\Omega)}{|y_i|^{\bm}}
\end{equation}
\end{theorem}

\medskip\noindent{\bf Notations:} In order to simplify notations, we will often drop the dependence in the domain $\Omega$ and the dimension $n\geq 3$. If $F: A\times B\to\rr$ is a function, then for any $x\in A$, we define $F_x: B\to\rr$ by $F_x(y):=F(x,y)$ for all $y\in B$. Finally, we will write $\hbox{Diag}(A):=\{(x,x)/\, x\in A\}$ for any set $A$. \\

We split the proof into several parts.

\subsection{Proof of existence and uniqueness of the Green function} 
We let $\eta_\eps(x):=\tilde{\eta}(\eps^{-1}|x|)$ for all $x\in\rn$ and $\eps>0$, where $\tilde{\eta}\in C^\infty(\rr)$ is nondecreasing and such that $\tilde{\eta}(t)=0$ for $t<1$  and $\tilde{\eta}(t)=1$ for $t>1$. It follows from Lemma \ref{lem:vp} (see Appendix B) and the coercivity of $-\Delta-\left(\gamma|x|^{-2}+h\right)$ that there exists $\eps_0>0$ and $c>0$ such that such that for all $\varphi\in \huno$ and $\eps\in (0,\eps_0)$, 
$$\int_\Omega\left(|\nabla \varphi|^2-\left(\frac{\gamma\eta_\eps}{|x|^2}+h(x)\right)\varphi^2\right)\, dx\geq c\int_\Omega\varphi^2\, dx.$$
As a consequence, there exists $c>0$ such that for all $\varphi\in \huno$ and $\eps\in (0,\eps_0)$, 
\begin{equation}\label{bnd:coer}
\int_\Omega\left(|\nabla \varphi|^2-\left(\frac{\gamma\eta_\eps}{|x|^2}+h(x)\right)\varphi^2\right)\, dx\geq c\Vert\varphi\Vert_{H_{1}^2}^2.
\end{equation}
Let $G_\eps>0$ be the Green's function of $-\Delta-\left(\gamma\eta_\eps|x|^{-2}+h\right)$ on $\Omega$ with Dirichlet boundary condition. The existence follows from the coercivity and the $C^{0,\theta}$ regularity of the potential for any $\eps>0$ (see Robert \cite{rob.green}). In particular, we have that
\begin{equation}\label{eq:G:eps}
\left\{\begin{array}{ll}
-\Delta G_\eps(x,\cdot)-\left(\frac{\gamma\eta_\eps}{|\cdot|^2}+h\right)G_\eps(x,\cdot)=0&\hbox{ in }\Omega\setminus\{x\}\\
G_\eps(x,\cdot)=0&\hbox{ on }\partial\Omega
\end{array}\right.
\end{equation}
\noindent{\bf Step \ref{sec:app:c}.1: Integral bounds for $G_\eps$.} We claim that for all $\delta>0$ and $1<q<\frac{n}{n-2}$ and $\delta'\in (0,\delta)$, there exists $C(\delta,q)>0$ and $C(\delta,\delta')>0$ such that 
\begin{equation}\label{int:bnd:G}
\Vert G_\eps(x,\cdot)\Vert_{L^q(\Omega)}\leq C(\delta,q)\hbox{ and }\Vert G_\eps(x,\cdot)\Vert_{L^{\frac{2n}{n-2}}(\Omega\setminus B_{\delta'}(x))}\leq C(\delta,\delta')
\end{equation}
for all $x\in \Omega$, $|x|>\delta$. We prove the claim. We fix $f\in C^\infty_c(\Omega)$ and let $\varphi_\eps\in C^{2,\theta}(\overline{\Omega})$ be the solution to the boundary value problem
\begin{equation}
\label{eq:phi:eps}\left\{\begin{array}{ll}
-\Delta \varphi_\eps-\left(\frac{\gamma\eta_\eps}{|x|^2}+h(x)\right)\varphi_\eps= f &\hbox{ in }\Omega\\
\quad \varphi_\eps=0&\hbox{ on }\partial\Omega
\end{array}\right.
\end{equation}
Multiplying the equation by $\varphi_\eps$, integrating by parts on $\Omega$, using \eqref{bnd:coer} and H\"older's inequality, we get that 
$$\int_\Omega |\nabla\varphi_\eps|^2\, dx\leq C\Vert f\Vert_{\frac{2n}{n+2}}\Vert\varphi_\eps\Vert_{\frac{2n}{n-2}}$$
where $C>0$ is independent of $\epsilon$, $f$ and $\varphi_\epsilon$. The Sobolev inequality $\Vert\varphi\Vert_{\frac{2n}{n-2}}\leq C\Vert \nabla\varphi\Vert_2$ for $\varphi\in \huno$ then yields
\begin{equation*}
\Vert\varphi_\eps\Vert_{\frac{2n}{n-2}}\leq C\Vert f\Vert_{\frac{2n}{n+2}}
\end{equation*}
where $C>0$ is independent of $\epsilon$, $f$ and $\varphi_\epsilon$. Fix $p>n/2$ and $\delta\in (0,\delta_0)$ and $\delta_1,\delta_2>0$ such that $\delta_1+\delta_2<\delta$, and $x\in\Omega$ such that $|x|>\delta$. It follows from standard elliptic theory that 
\begin{eqnarray*}
|\varphi_\eps(x)|&\leq & \Vert \varphi_\eps\Vert_{C^0(B_{\delta_1}(x))}\\
&\leq & C\left(\Vert \varphi_\epsilon\Vert_{L^{\frac{2n}{n-2}}(B_{\delta_1+\delta_2}(x))}+\Vert f\Vert_{L^{p}(B_{\delta_1+\delta_2}(x))}\right)\\
&\leq & C\left(\Vert f\Vert_{L^{\frac{2n}{n+2}}(\Omega)}+\Vert f\Vert_{L^{p}(B_{\delta_1+\delta_2}(x))}\right)
\end{eqnarray*}
where $C>0$ depends on $p,\delta,\delta_1,\delta_2$, $\gamma$ and $\Vert h\Vert_\infty$. Therefore, Green's representation formula yields
\begin{equation}\label{rep:G:f}
\left|\int_\Omega G_\eps(x,\cdot)f\, dy\right|\leq C\left(\Vert f\Vert_{L^{\frac{2n}{n+2}}(\Omega)}+\Vert f\Vert_{L^{p}(B_{\delta_1+\delta_2}(x))}\right)
\end{equation}
for all $f\in C^\infty_c(\Omega)$.  It follows from \eqref{rep:G:f} that
$$\left|\int_\Omega G_\eps(x,\cdot)f\, dy\right|\leq C\cdot\Vert f\Vert_{L^{p}(\Omega)}$$
for all $f\in C^\infty_c(\Omega)$ where $p>n/2$. It then follows from duality arguments that for any $q\in (1, n/(n-2))$ and any $\delta>0$, there exists $C(\delta,q)>0$ such that $\Vert G_\eps(x,\cdot)\Vert_{L^q(\Omega)}\leq C(\delta,q)$ for all $\eps<\eps_0$ and $x\in \Omega\setminus B_\delta(0)$.

\medskip\noindent Let $\delta'\in (0,\delta)$ and $\delta_1,\delta_2>0$ such that $\delta_1+\delta_2<\delta'$. We get from \eqref{rep:G:f}  that 
\begin{equation}\label{rep:G:f:2}
\left|\int_\Omega G_\eps(x,\cdot)f\, dy\right|\leq C\Vert f\Vert_{L^{\frac{2n}{n+2}}(\Omega\setminus B_{\delta'}(x))}
\end{equation}
for all $f\in C^\infty_c(\Omega\setminus B_{\delta'}(x))$. Here again, a duality argument yields \eqref{int:bnd:G}, which proves the claim in  Step \ref{sec:app:c}.1. 

\smallskip\noindent Using the same method, we can get an improvement of the control, the cost being the integrability exponent $q$. When $q\in (1,n/(n-1))$, we get that $p>n$. Then, $\Vert \varphi_\eps\Vert_{C^1(B_{\delta_1}(x)\cap\Omega)}$ is controled by the $L^p$ and $L^{\frac{2n}{n+2}}$ norms. Moreover, $|\varphi_\eps(x)|\leq \Vert \varphi_\eps\Vert_{C^0(B_{\delta_1}(x)\cap\Omega)}d(x,\partial\Omega)$. The argument above then yields
\begin{equation}\label{int:bnd:G:boundary}
\Vert G_\eps(x,\cdot)\Vert_{L^q(\Omega)}\leq C(\delta,q)d(x,\partial\Omega) \hbox{ for }q\in \left(1,\frac{n}{n-1}\right).
\end{equation}

\medskip\noindent{\bf Step \ref{sec:app:c}.2: Convergence of $G_\epsilon$.} Fix $x\in \Omega\setminus\{0\}$. For $0<\eps<\eps'$, since $G_\eps(x,\cdot)$, $G_{\eps'}(x,\cdot)$ are $C^2$ outside $x$, \eqref{eq:G:eps} yields  
$$-\Delta(G_\eps(x,\cdot)-G_{\eps'}(x,\cdot))-\left(\frac{\gamma\eta_\eps}{|\cdot|^2}+h\right)(G_\eps(x,\cdot)-G_{\eps'}(x,\cdot))= \frac{\gamma(\eta_\eps-\eta_{\eps'})}{|\cdot|^2}G_{\eps'}(x,\cdot)$$
in the strong sense. The coercivity \eqref{bnd:coer} then yields $G_\eps(x,\cdot)\geq G_{\eps'}(x,\cdot)$ for $0<\eps<\eps'$ if $\gamma\geq 0$, and the reverse inequality if $\gamma<0$. It then follows from the integral bound \eqref{int:bnd:G} and elliptic regularity that there exists $G(x,\cdot)\in C^{2,\theta}(\overline{\Omega}\setminus\{0,x\})$ such that
\begin{equation}\label{lim:G:eps}
\lim_{\eps\to 0}G_\eps(x,\cdot)=G(x,\cdot)\geq 0\hbox{ in }C^{2,\theta}_{loc}(\overline{\Omega}-\{0,x\}).
\end{equation}
In particular, $G$ is symmetric and 
\begin{equation}\label{eq:G:x}
-\Delta G(x,\cdot)-\left(\frac{\gamma}{|\cdot|^2}+h\right)G(x,\cdot)=0\hbox{ in }\Omega\setminus\{x\}\hbox{ and }G(x,\cdot)=0\hbox{ on }\partial\Omega. 
\end{equation}
Moreover, passing to the limit $\eps\to 0$ in \eqref{int:bnd:G}, \eqref{int:bnd:G:boundary} and using elliptic regularity, we get that for all $\delta>0$, $1<q<\frac{n}{n-2}$ and $\delta'\in (0,\delta)$, there exist $C(\delta,q)>0$ and $C(\delta,\delta')>0$ such that for all $x\in \Omega$, $|x|>\delta$,
\begin{equation}\label{int:bnd:G:bis}
\Vert G(x,\cdot)\Vert_{L^q(\Omega)}\leq C(\delta,q)\hbox{ and }\Vert G(x,\cdot)\Vert_{L^{\frac{2n}{n-2}}(\Omega\setminus B_{\delta'}(x))}\leq C(\delta,\delta')
\end{equation}
and
\begin{equation}\label{int:bnd:G:boundary:2}
\Vert G(x,\cdot)\Vert_{L^q(\Omega)}\leq C(\delta,q)d(x,\partial\Omega) \hbox{ for }q\in \left(1,\frac{n}{n-1}\right).
\end{equation}
In particular, for any $x\in\Omega\setminus\{0\}$, $G(x,\cdot)\in L^{k}(\Omega)$ for all $1<k<n/(n-2)$ and $G(x,\cdot)\in L^{2n/(n-2)}(\Omega\setminus B_\delta(x))$ for all $\delta>0$. Moreover, for any $f\in  L^{\frac{2n}{n+2}}(\Omega)\cap L^q(\Omega\setminus B_\delta(0))$ for all $\delta>0$ with $q>n/2$, let $\varphi_\eps\in\huno$ be such that \eqref{eq:phi:eps} holds. It follows from elliptic theory that $\varphi_\eps\in C^{0,\tau}(\overline{\Omega}\setminus\{0\})$ for some $\tau\in (0,1)$ and that for all $\delta_1>0$, there exists $C(\delta_1)>0$ such that $\Vert\varphi_\eps\Vert_{C^{0,\tau}(\overline{\Omega}\setminus B_{\delta_1}(0))}\leq C(\delta_1)$. We fix $x\in \Omega\setminus \{0\}$. Passing to the limit $\eps\to 0$ in the Green identity $\varphi_\eps(x)=\int_\Omega G_\eps(x,\cdot)f\, dy$ yields
\begin{equation}\label{id:regul}
\varphi(x)=\int_\Omega G(x,\cdot)f\, dy\hbox{ for all }x\in\Omega\setminus\{0\}
\end{equation}
where $\varphi\in \huno\cap C^{0}(\overline{\Omega}\setminus\{0\})$ is the only weak solution to 
$$\left\{\begin{array}{ll}
-\Delta \varphi-\left(\frac{\gamma}{|x|^2}+h(x)\right)\varphi= f &\hbox{ in }\Omega\\
\varphi=0&\hbox{ on }\partial\Omega
\end{array}\right.$$
Since $G(x,\cdot)\geq 0$, \eqref{eq:G:x} and the strong comparison principle yield $G(x,\cdot)>0$. These points prove that $G$ is a Green's function for the operator and that (c) holds. 

\smallskip\noindent We now prove point (b). We fix $\eta\in C^\infty_c(\rn-\{x\})$ such that $\eta(y)=1$ when $y\in B_\delta(0)$ for some $\delta>0$. Then $\eta G_\eps(x,\cdot)\in C^{2,\theta}(\overline{\Omega})\cap \huno$. It follows from \eqref{eq:G:eps} and \eqref{lim:G:eps} that
$$-\Delta (\eta G_\eps(x,\cdot))-\left(\frac{\gamma\eta_\eps}{|\cdot|^2}+h\right)(\eta G_\eps(x,\cdot))={\bf 1}_{B_\delta(0)^c}f_\eps\hbox{ in }\Omega$$
where $\Vert f_\eps\Vert_{C^0(\overline{\Omega})}\leq C$ for some $C>0$ and all $\eps>0$. Therefore, with the coercivity \eqref{bnd:coer} and the convergence \eqref{lim:G:eps}, we get that
$$c\Vert \eta G_\eps(x,\cdot)\Vert_{H_1^2}^2\leq \int_{\Omega\setminus B_\delta(0)}f_\eps\eta G_\eps(x,\cdot)\, dy\leq C$$
for all $\eps>0$. Reflexivity yields convergence of $(\eta G_\eps(x,\cdot))$ in $\huno\cap L^2(\Omega)$ as $\eps\to 0$ up to extraction. The convergence in $C^2$ and uniqueness then yields $\eta G(x,\cdot)\in \huno$ and $\eta G_\eps(x,\cdot)\to \eta G(x,\cdot)$ in $\huno$ as $\eps\to 0$. The case of a general $
\eta$ is a direct consequence. This proves point (b).

\smallskip\noindent For the uniqueness, we suppose $G'$ be another Green's function. We fix $x\in \Omega$ and we define $H_x:=G_x-G'_x$. Then $H_x\in L^1(\Omega)$ and for any $f\in C^\infty_c(\Omega)$, we have that $\int_\Omega H_x f\, dy=0$.  Approximating a compactly supported function by smooth fonctions with compact support, we get that this equality holds for all  $f\in C^0_c(\Omega)$. Integration theory then yields $H_x\equiv 0$, and then $G'_x\equiv G_x$. This proves uniqueness. This finishes the proof of (a).

\smallskip\noindent This proves existence and uniqueness of the Green's function in Theorem \ref{th:green:gamma:domain}(I).

\subsection{Proof of the upper bound}\label{sec:proof}
The behavior \eqref{asymp:G} is a consequence of the classification of solutions to harmonic equations and Theorem 4.1 in Ghoussoub-Robert \cite{gr4}.

\smallskip\noindent In the proof, we will often use sub- and super-solutions to the linear problem. The following existence result is contained in Proposition 4.3 of \cite{gr4}:

\begin{proposition}\label{prop:sub:super} Let $\Omega$ be a smooth domain and $h\in C^{0}(\overline{\Omega})$ be a continuous fonction. We fix $\gamma<\frac{n^2}{4}$ and $\beta\in \{\bm,\bp\}$. Then, there exist $r>0$, and  $\overline{u}_{\beta},\underline{u}_{\beta}\in C^\infty(\overline{\Omega}\setminus\{0\})$ such that
\begin{equation}\label{ppty:ua}
\left\{\begin{array}{cc}
 \overline{u}_{\beta},\underline{u}_{\beta}=0 &\hbox{ on }\partial\Omega \cap B_r(0)\\
-\Delta \overline{u}_{\beta}-\left(\frac{\gamma}{|x|^2}+h\right)\overline{u}_{\beta}>0&\hbox{ in }\Omega\cap B_r(0)\\
-\Delta \underline{u}_{\beta}-\left(\frac{\gamma}{|x|^2}+h\right)\underline{u}_{\beta}<0&\hbox{ in }\Omega\cap B_r(0).
\end{array}\right.
\end{equation}
Moreover, for some $\tau>0$, we have that, as $x\to 0$, $x\in \Omega$,  
\begin{equation}\label{asymp:ua:plus}
\overline{u}_{\beta}(x)=\underline{u}_{\beta}(x)(1+O(|x|^\tau))=\frac{d(x,\partial\Omega)}{|x|^{\beta}}(1+O(|x|^{\tau})).
\end{equation}
\end{proposition}

\medskip\noindent{\bf Step \ref{sec:app:c}.3: Upper bound for $G(x,y)$ when one variable is far from $0$.}

\smallskip\noindent{\it Step \ref{sec:app:c}.3.1:} It follows from \eqref{eq:G:x}, elliptic theory, \eqref{int:bnd:G:boundary:2} and \eqref{int:bnd:G:bis} that for any $\delta>0$, there exists $C(\delta)>0$ such that
\begin{equation}\label{est:1:bis}
0<G(x,y)\leq C(\delta)d(y,\partial\Omega)d(x,\partial\Omega)\hbox{ for }x,y\in\Omega\hbox{ s.t. }|x|,|y|>\delta,\, |x-y|>\delta.
\end{equation}

\medskip\noindent{\it Step \ref{sec:app:c}.3.2:} 
We claim that for any $\delta>0$, there exists $C(\delta)>0$ such that
\begin{equation}\label{est:2}
|x-y|^{n-2}G(x,y)\leq C(\delta)\min\left\{1,\frac{d(x,\partial\Omega)d(y,\partial\Omega)}{|x-y|^2}\right\}\hbox{ for }x,y\in\Omega\hbox{ s.t. }|x|,|y|>\delta.
\end{equation}
Indeed, with no loss of generality, we can assume that $\delta\in (0,\delta_0)$. Let $\Omega_\delta$ be a smooth domain of $\rn$ be such that $\Omega\setminus B_{3\delta/4}(0)\subset \Omega_\delta\subset \Omega\setminus B_{\delta/2}(0)$. We fix $x\in\Omega$ such that $|x|>\delta$. Let $H_x$ be the Green's function for $-\Delta -\left(\frac{\gamma}{|x|^2}+h(x)\right)$ in $\Omega_\delta$ with Dirichlet boundary condition. Classical estimates (see \cite{rob.green}) yield the existence of $C(\delta)>0$ such that 
$$|x-y|^{n-2}H_x(y)\leq C(\delta)\min\left\{1,\frac{d(x,\partial\Omega)d(y,\partial\Omega)}{|x-y|^2}\right\}\hbox{ for all }x,y\in \Omega_\delta.$$
It is easy to check that
$$\left\{\begin{array}{ll}
-\Delta (G_x-H_x)-\left(\frac{\gamma}{|\cdot|^2}+h\right)(G_x-H_x)= 0 &\hbox{ weakly in }\Omega_\delta\\
G_x-H_x=0&\hbox{ on }\left(\partial\Omega_\delta\right)\setminus B_{3\delta/4}(0)\\
G_x-H_x=G_x&\hbox{ on }\left(\partial\Omega_\delta\right)\cap B_{3\delta/4}(0).
\end{array}\right.$$
Regularity theory then yields that $G_x-H_x\in C^{2,\theta}(\overline{\Omega_\delta})$. It follows from \eqref{est:1:bis} that $G_x(y)\leq C_1(\delta)d(y,\partial\Omega)d(x,\partial\Omega)$ on $(\partial\Omega_\delta)\cap B_{3\delta/4}(0)$ for $|x|>\delta$. The comparison principle then yields $G_x(y)-H_x(y)\leq C_1(\delta)d(y,\partial\Omega)d(x,\partial\Omega)$ for $y\in \Omega_\delta$ and $|x|>\delta$. The above bound for $H_x$ and \eqref{est:1:bis} then yields \eqref{est:2}.

\medskip\noindent{\it Step \ref{sec:app:c}.3.3:} 
We now claim that for any $0<\delta'<\delta$, there exists $C(\delta,\delta')>0$ such that
\begin{equation}\label{est:3}
|y|^{\bm}G(x,y)\leq C(\delta,\delta')d(y,\partial\Omega)d(x,\partial\Omega)\hbox{ for }x,y\in\Omega\hbox{ s.t. }|x|>\delta>\delta'>|y|.
\end{equation}
We let $\delta_1\in (0,\delta')$ that will be fixed later. We use \eqref{est:1:bis} to deduce that $G_x(y)\leq C(\delta,\delta_1)d(x,\partial\Omega)d(y,\partial\Omega)$ for all $x\in \Omega\setminus B_\delta(0)$ and $y\in \partial B_{\delta_1}(0)\cap\Omega$. Since $\delta_1<|x|$, we have that
 $$\left\{\begin{array}{ll}
-\Delta G_x-\left(\frac{\gamma}{|x|^2}+h\right)G_x= 0 &\hbox{  in }\Omega\cap B_{\delta_1}(0)\\
0\leq G_x\leq C(\delta,\delta_1)d(y,\partial\Omega)d(x,\partial\Omega)&\hbox{ on }\partial (\Omega\cap B_{\delta_1}(0))\setminus\{0\}.
\end{array}\right.$$
We choose a supersolution $\overline{u}_{\bm}$ as in \eqref{ppty:ua} of Proposition \ref{prop:sub:super}. It follows from \eqref{asymp:ua:plus} and \eqref{est:1:bis} that for $\delta_1>0$, there exists $C(\delta,\delta_1)>0$ such that $G_x(z)\leq C(\delta,\delta_1)d(x,\partial\Omega)u_{\beta_-}(z)$ for all $z\in \partial (\Omega\cap B_{\delta_1}(0))$. It then follows from the comparison principle that $G_x(y)\leq C(\delta,\delta_1)d(x,\partial\Omega)u_{\beta_-}(y)$ for all $y\in (\Omega\cap B_{\delta_1}(0))\setminus\{0\}$. Combining this with \eqref{est:1:bis} and \eqref{ppty:ua}, we obtain \eqref{est:3}.

\smallskip\noindent Note that by symmetry, we also get that for any $0<\delta'<\delta$, there exists $C(\delta,\delta')>0$ such that
\begin{equation}\label{est:4}
|x|^{\bm}G(x,y)\leq C(\delta,\delta')d(x,\partial\Omega)d(y,\partial\Omega)\hbox{ for }x,y\in\Omega\hbox{ s.t. }|y|>\delta>\delta'>|x|.
\end{equation}

\medskip\noindent{\bf Step \ref{sec:app:c}.4: Upper bound for $G(x,y)$ when both variables approach $0$.}\\

We claim first that for all $c_1,c_2,c_3>0$, there exists $C(c_1,c_2,c_3)>0$ such that for $x,y\in\Omega$ such that $c_1|x|<|y|<c_2|x|$ and $|x-y|>c_3|x|$, we have
 \begin{equation}\label{est:5:bis}
|x-y|^{n-2}G(x,y)\leq C(c_1,c_2,c_3)\frac{d(x,\partial\Omega)d(y,\partial\Omega)}{|x|^2}.\end{equation}
When one of the variables stays far from $0$, \eqref{est:5:bis} is a consequence of \eqref{est:1:bis}. We now consider a chart $\T$ at $0$ as in \eqref{def:T:bdry}. In particular, there is $\delta_0>0$, $0\in V\subset \rn$ and $\T:B_{2\delta_0}(0)\to V$ a smooth diffeomorphism such that $\T(0)=0$ and
\begin{equation}\label{chart:phi}
\T(B_{2\delta_0}(0)\cap\rnm)=\T(U)\cap\Omega\hbox{ and }\T(B_{2\delta_0}(0)\cap\partial\rnm)=\T(U)\cap\partial\Omega.
\end{equation}
Moreover, $D_{0} \mathcal{T} = \mathbb{I}_{\R^{n}}$ and
\begin{equation}\label{asymp:phi}
|\T(X)|=(1+O(|X|))|X|\hbox{ for all }X\in B_{3\delta_0/2}(0).
\end{equation} 
We fix $X\in \rnm$ such that $0<|X|<3\delta_0/2$. We define
$$H(z):=G_{\T(X)}(\T(|X|z))\hbox{ for }z\in B_{\delta_0/|X|}(0)\setminus \left\{0,\frac{X}{|X|}\right\},$$
so that 
$$-\Delta_{g_X} H-\left(\frac{\gamma}{\left(\frac{|\T(|X|z|)}{|X|}\right)^2}+|X|^2h(\T(|X|z))\right)H=0\hbox{ in }B_{\delta_0/|X|}(0)\setminus\left \{0,\frac{X}{|X|}\right\}.$$
where $g_X:=(\T^\star \eucl)_X$ is the pulled-back metric of the Euclidean metric $\eucl$ via the chart $\T$ at the point $X$. Since $H>0$, it follows from the Harnack inequality on the boundary (see Proposition 6.3 in Ghoussoub-Robert \cite{gr4}) that for all $R>0$ large enough and $r>0$ small enough, there exist $\delta_1>0$ and $C>0$ independent of $|X|<3\delta_0/2$ such that
\begin{equation*}
\frac{H(z)}{|z_1|}\leq C \frac{H(z')}{|z'_1|}\hbox{ for all }z,z'\in (B_R(0)\cap \rnm)\setminus \left(B_r(0)\cup B_r\left(\frac{X}{|X|}\right)\right),
\end{equation*}
which, via the chart $\T$, yields
\begin{equation}\label{harnack:g}
\frac{G_x(y)}{d(y,\partial\Omega)}\leq C \frac{G_x(y')}{d(y',\partial\Omega)}\hbox{ for all }y,y'\in \Omega\cap B_{R|x|/2}(0)\setminus \left(B_{2r|x|}(0)\cup B_{2r|x|}(x)\right).
\end{equation}
for all $x\in \Omega$ such that $|x|<\delta_0$. We let $W$ be a smooth domain of $\rn$ such that for some $\lambda>0$ small enough, we have
\begin{equation}\label{ppte:W.bis}
B_\lambda(0)\cap\Omega\subset W\subset B_{2\lambda}(0)\cap\Omega\hbox{ and } B_\lambda(0)\cap\partial W=B_\lambda(0)\cap\partial \Omega.
\end{equation}
We choose a subsolution $\underline{u}_{\bp}$ as in \eqref{ppty:ua} of Proposition \ref{prop:sub:super}. It follows from \eqref{asymp:ua:plus} and \eqref{est:1:bis} that for $|x|<\delta_2$ small

$$G_x(z)\geq C(R)|x|^{\bp}\left(\inf_{y\in \Omega\cap \partial B_{R|x|}(0)}\frac{G_x(y)}{d(y,\partial\Omega)}\right) \underline{u}_{\bp}(z)\hbox{ for all }z\in W \cap \partial B_{R|x|/3}(0).$$
Since $-\Delta G_x-(\gamma|\cdot|^{-2}+h)G_x=0$ outside $0$, it follows from coercivity and the comparison principle that
$$G_x(z)\geq c|x|^{\bp}\left(\inf_{y\in \Omega\cap \partial B_{R|x|}(0)}\frac{G_x(y)}{d(y,\partial\Omega)}\right) \underline{u}_{\bp}(z)\hbox{ for all }z\in W\setminus B_{R|x|/3}(0).$$
We fix $z_0\in W\setminus\{0\}$. Then for $\delta_3$ small enough, when $|x|<\delta_3$, it follows from \eqref{est:4} and the Harnack inequality \eqref{harnack:g} that there exists $C>0$ independent of $x$ such that
$$G_x(y)\leq C |x|^{-\bp-\bm}d(x,\partial\Omega)d(y,\partial\Omega)\hbox{ for all }y\in B_{R|x|}(0)\setminus \left(B_{r|x|}(0)\cup B_{r|x|}(x)\right)$$
Taking $r>0$ small enough and $R>0$ large enough, we then get \eqref{est:5:bis} for $|x|<\delta_3$. The general case for arbitrary $x\in \Omega\setminus \{0\}$ then follows from \eqref{est:2}. This completes the proof of \eqref{est:5:bis}.  

\medskip\noindent{\it Step \ref{sec:app:c}.4.2:} We claim that for all $c_1,c_2>0$, there exists $C(c_1,c_2)>0$ such that
 \begin{equation}\label{est:6}
|x-y|^{n-2}G(x,y)\leq C(c_1,c_2)\min\left\{1,\frac{d(x,\partial\Omega)d(y,\partial\Omega)}{|x-y|^2}\right\}
\end{equation}
for all $x,y\in\Omega$ s.t. $c_1|x|<|y|<c_2|x|$. To prove \eqref{est:6}, we distinguish three cases:

\smallskip\noindent{\it Case 1:} We assume that 
\begin{equation}\label{hyp:x:1}
|x|\leq C_1 d(x,\partial\Omega)\hbox{ with }C_1>1.
\end{equation}
We define 
$$H(z):=|x|^{n-2}G_x(x+|x|z)\hbox{ for }z\in B_{1/C_1}(0)\setminus \{0\}.$$
Note that this definition makes sense since for such $z$, $x+|x|z\in \Omega$. We then have that $H\in C^2(\overline{B_{1/(2C_1)}(0)}\setminus\{0\})$ and 
$$-\Delta H-\left(\frac{\gamma}{\left|\frac{x}{|x|}+z\right|^2}+|x|^2h(x+|x|z)\right)H=\delta_0\hbox{ weakly in }B_{1/(2C_1)}(0).$$
We now argue as in the proof of \eqref{est:2}. From \eqref{est:5:bis}, we have that $|H(z)|\leq C$ for all $z\in \partial B_{1/(2C_1)}(0)$ where $C$ is independent of $x\in \Omega\setminus\{0\}$ satisfying \eqref{hyp:x:1}. Let $\Gamma_0$ be the Green's function of $-\Delta -\left(\frac{\gamma}{\left|\frac{x}{|x|}+z\right|^2}+|x|^2h(x+|x|z)\right)$ at $0$ on $B_{1/(2C_1)}(0)$ with Dirichlet boundary condition. Therefore, $H-\Gamma_0\in C^2(\overline{B_{1/(2C_1)}(0)})$ and, via the comparison principle, it is bounded by its supremum on the boundary. Therefore $|z|^{n-2}H(z)\leq C$ for all $B_{1/(2C_1)}(0)\setminus\{0\}$ where $C$ is independent of $x\in \Omega\setminus\{0\}$ satisfying \eqref{hyp:x:1}. Scaling back and using \eqref{est:5:bis}, we get $|x-y|^{n-2}G_x(y)\leq C$ for all $x,y\in \Omega\setminus\{0\}$ such that $c_1|x|<|y|<c_2|x|$ and \eqref{hyp:x:1} holds. This proves \eqref{est:6} if $d(x,\partial\Omega)d(y,\partial\Omega)\geq|x-y|^2$. If $d(x,\partial\Omega)d(y,\partial\Omega)<|x-y|^2$, we get that $d(x,\partial\Omega)<2|x-y|$, and then \eqref{hyp:x:1} yields $|x|\leq 2C_1|x-y|$, and \eqref{est:6} is a consequence of \eqref{est:5:bis}.

\smallskip\noindent This ends the proof of \eqref{est:6} in Case 1. 

\smallskip\noindent{\it Case 2:} By symmetry,  \eqref{est:6} also holds when $|y|\leq C_1 d(y,\partial\Omega)$.

\smallskip\noindent{\it Case 3:} We assume that $d(x,\partial\Omega)\leq C_1^{-1} |x|$ and $d(y,\partial\Omega)\leq C_1^{-1} |y|$. We consider a chart at $0$, that is $\delta_0>0$, $0\in V\subset \rn$ and $\T:B_{2\delta_0}(0)\to V$ a smooth diffeomorphism such that $\T(0)=0$ and that \eqref{chart:phi} and \eqref{asymp:phi} hold. We fix $x'\in \rr^{n-1}$ such that $0<|x'|<3\delta_0/2$. 

\medskip\noindent We assume that $r\leq c_0|x'|$. We define
$$H_y(z):=r^{n-2}G_{\T((0,x')+r y)}(\T((0,x')+r z))\hbox{ for }y,z\in B_{\delta_0/(2r)}(0)\cap\rnm\setminus \{0\}.$$
We then have that $H_y\in C^2(\overline{B_{R_0}(0)}\cap\rnm\setminus\{0,y\})$ and 
$$-\Delta_{g_r} H_y-\left(\frac{\gamma}{\left(\frac{|\T((0,x')+rz)}{r}\right)^2}+r^2h(\T((0,x')+r z))\right)H_y=\delta_y\hbox{ weakly in }B_{R_0}(0)\cap\rnm,$$
where $g_{r}:=(\T^\star \eucl)_{(0,x')+rz}$ is the pulled-back metric of the Euclidean metric $\eucl$ via the chart $\T$ at the point $(0,x')+rz$. We now argue as in the proof of \eqref{est:2}. From \eqref{est:5}, we have that $|H_y(z)|\leq C$ for all $z\in \partial B_{R_0}(0)\cap\rnm$ where $C$ is independent of $y\in B_{R_0/2}(0)$ and $r\in (0,\delta_0/4)$. Let $\Gamma_y$ be the Green's function of $-\Delta_{g_r} -\left(\frac{\gamma}{\left(\frac{|\T((0,x')+rz)}{r}\right)^2}+r^2h(\T((0,x')+r z))\right)$ at $y$ on $B_{c_0/2}(0)\cap\rnm$ with Dirichlet boundary condition. Therefore, $H_y-\Gamma_y\in C^2(\overline{B_{c_0/2}(0)\cap\rnm})$ and, via the comparison principle, it is bounded by its supremum on the boundary. It follows from \eqref{est:5} and elliptic estimates for $\Gamma_y$ (see for instance \cite{rob.green}) that $|H_y-\Gamma_y|(z)\leq C |y_1|\cdot|z_1|$ for $z\in \partial(B_{c_0/2}(0)\cap\rnm)$ and $y\in B_{c_0/4}(0)\cap\rnm$. Applying elliptic estimates, we then get that $|H_y-\Gamma_y|(z)\leq C |y_1|\cdot|z_1|$ for $z\in B_{c_0/2}(0)\cap\rnm$ and $y\in B_{c_0/4}(0)\cap\rnm$, and since
$$\Gamma_y(z)\leq C|z-y|^{2-n}\min\left\{1,\frac{|y_1|\cdot|z_1|}{|y-z|^2}\right\}\hbox{ for all }y,z\in B_{c_0/2}(0)\cap\rnm$$
(see \cite{rob.green}), we get that
$$|z-y|^{n-2}H_y(z)\leq C\min\left\{1,\frac{|y_1|\cdot|z_1|}{|y-z|^2}\right\}\hbox{ for all }y,z\in B_{c_0/2}(0)\cap\rnm$$
where $C$ is independent of $x'\in B_{\delta_0/2}(0)\setminus\{0\}$. This yields
\begin{equation}\label{est:green:36}
|rz-ry|^{n-2}G_{\T((0,x')+r y)}(\T((0,x')+r z))\leq C\min\{1,\frac{|y_1|\cdot|z_1|}{|y-z|^2}\}
\end{equation}
for $|x'|<\delta_0/3$, $r\leq c_0|x'|$ and $|y|,|z|\leq c_0/4$.

\smallskip\noindent We now prove \eqref{est:6} in the last case. We fix $x\in\Omega\setminus \{0\}$ such that $|x|<\delta_0/3$. We assume that $d(x,\partial\Omega)\leq C_1^{-1} |x|\; , \; d(y,\partial\Omega)\leq C_1^{-1} |y|\hbox{ and }|x-y|\leq \epsilon_0|x|$.
We let $(x_1,x'), (y_1,y')\in B_{\delta_0}(0)$ be such that $x=\T(x_1,x')$ and $y=\T(y_1,y')$. Taking the norm $|(x_1,x')|=|x_1|+|x'|$, we define $r:=\max\{d(x,\partial\Omega),|x-y|\}$. Using that $|X|/2\leq |\T(X)|\leq 2|X|$ for $X\in B_{\delta_0}(0)$, up to taking $\epsilon_0>0$ small and $C_1,c_0>1$ large enough, we get that
$$\left|\frac{x_1}{r}\right|\leq \frac{c_0}{4}\; ,\; \left|\left(\frac{y_1}{r},\frac{y'-x'}{r}\right)\right|\leq \frac{c_0}{4}\hbox{ and }r\leq c_0|x'|.$$
Therefore, \eqref{est:green:36} applies and we get \eqref{est:6} in Case 3.

\medskip\noindent We are now in position to conclude. Inequality \eqref{est:6} is a consequence of Cases 1, 2, 3, \eqref{est:2} and \eqref{est:5}. This ends the proof of \eqref{est:6}.

\medskip\noindent{\it Step \ref{sec:app:c}.4.3:} We now show that there exists $C>0$ such that
 \begin{equation}\label{est:7}
|y|^{\bm}|x|^{\bp}G(x,y)\leq Cd(x,\partial\Omega) d(y,\partial\Omega)\hbox{ for }x,y\in\Omega\hbox{ such that }|y|<\frac{1}{2}|x|.
\end{equation}
The proof goes essentially as in \eqref{est:3}. For $|x|<\delta$ with $\delta>0$ small, we have that
$$-\Delta G_x-\left(\frac{\gamma}{|\cdot|^2}+h\right)G_x=0\hbox{ in }H^1(\Omega\cap B_{|x|/3}(0))\cap C^2(\Omegabar\cap  B_{|x|/3}(0)\setminus\{0\}).$$
It follows from \eqref{est:5} that $G_x(y)\leq C |x|^{-n} d(x,\partial\Omega) d(y,\partial\Omega)$ in $\Omega\cap \partial B_{|x|/3}(0)$. We choose a supersolution $\overline{u}_{\bm}$ as in \eqref{ppty:ua} of Proposition \ref{prop:sub:super}. It follows from \eqref{asymp:ua:plus} and \eqref{est:5} that there exists $C>0$ such that
$$G_x(y)\leq C |x|^{-\bp} d(x,\partial\Omega)     \overline{u}_{\bm}(y)  \hbox{ for all }y\in \Omega\cap \partial B_{|x|/3}(0).$$
The comparison principle yields that this inequality holds on $\Omega\cap  B_{|x|/3}(0)$.

\medskip\noindent{\it Step \ref{sec:app:c}.4.4:} By symmetry, we conclude that there exists $C>0$ such that
 \begin{equation}\label{est:8}
|x|^{\bm}|y|^{\bp}G(x,y)\leq Cd(x,\partial\Omega) d(y,\partial\Omega)\hbox{ for }x,y\in\Omega\hbox{ s.t. }|x|<\frac{1}{2}|y|.
\end{equation}

\medskip\noindent{\bf Step \ref{sec:app:c}.5:} Finally, it follows from \eqref{est:7}, \eqref{est:8} and \eqref{est:6} that there exists $c>0$ such that
\begin{equation}\label{G:up}
G(x,y)\leq c \left(\frac{\max\{|y|,|x|\}}{\min\{|y|,|x|\}}\right)^{\bm}|x-y|^{2-n}\min\left\{1,\frac{d(x,\partial\Omega)d(y,\partial\Omega)}{|x-y|^2}\right\}
\end{equation}
for all $x,y\in \Omega$, $x\neq y$. This proves the upper bound in \eqref{est:G:up} of Theorem \ref{th:green:gamma:asymp}. The lower-bound and the control of the gradient will be proved in Section \ref{sec:lower}.

\subsection{Behavior at infinitesimal scale}\label{sec:th:cv}
We prove three convergence results to get a comprehensive behavior of the Green's function. Throughout this subsection, we assume $\Omega$ is a smooth bounded domain of $\rn$ such that $0\in\partial\Omega$. We fix $\gamma<\frac{n^2}{4}$ and let $h\in C^{0,\theta}(\Omegabar)$ be such that $-\Delta-\gamma|x|^{-2}-h$ is coercive. We consider $G$ to be the Green's function of $-\Delta-\gamma|x|^{-2}-h$ with Dirichlet boundary condition on $\partial\Omega$. 

\begin{lemma}\label{th:cv:1}
 Let $(x_i)_i\in\Omega$ and $(r_i)_i\in (0,+\infty)$ be such that
$$\lim_{i\to +\infty}r_i=0\hbox{ and }\lim_{i\to +\infty}\frac{d(x_i,\partial\Omega)}{r_i}=+\infty.$$
Then, for all $X,Y\in \rn$ such that $X\neq Y$, we have that
\begin{equation*}
\lim_{i\to +\infty}r_i^{n-2}G(x_i+r_iX,x_i+r_iY)=\frac{1}{(n-2)\omega_{n-1}}|X-Y|^{2-n}
\end{equation*}
Moreover, the convergence holds in $C^2_{loc}((\rn)^2\setminus\hbox{Diag}(\rn))$.
\end{lemma}
To deal with the case when the points approach the boundary, we consider a chart $\T$ as in \eqref{def:T:bdry}. In particular, $D_{0} \mathcal{T} = \mathbb{I}_{\R^{n}}$.
\begin{lemma}\label{th:cv:2} 
Let $(x_i)_i\in\partial\Omega$ and $(r_i)_i\in (0,+\infty)$ and $x_0\in\partial\Omega$ be such that
$$\lim_{i\to +\infty}r_i=0,\; \lim_{i\to +\infty}x_i=x_0\in \partial\Omega\hbox{ and }\lim_{i\to +\infty}\frac{|x_i|}{r_i}=+\infty.$$
We let $\T$ be a chart at $x_0$ as in \eqref{def:T:bdry}. We define $x_i'\in\rr^{n-1}$ such that $x_i=\T(0,x_i')$. Then, for all $X,Y\in \rnm$ such that $X\neq Y$, we have that
\begin{eqnarray*}
&&\lim_{i\to +\infty}r_i^{n-2}G(\T\left((0,x'_i)+r_iX\right),\T\left((0,x'_i)+r_iY\right))\\
&&=\frac{1}{(n-2)\omega_{n-1}}\left(|X-Y|^{2-n}-|X-Y^*|^{2-n}\right)
\end{eqnarray*}
where $(Y_1,Y')^*=(-Y_1,Y')$ for $(Y_1,Y')\in \rr\times\rr^{n-1}$. Moreover, the convergence holds in $C^2_{loc}((\overline{\rnm})^2\setminus\hbox{Diag}(\overline{\rnm})\})$.

\end{lemma}

\begin{lemma}\label{th:cv:3} 
Let $(r_i)_i\in (0,+\infty)$ be such that $\lim_{i\to +\infty}r_i=0$. We let $\T$ be a chart at $0$ as in \eqref{def:T:bdry}. Then, for all $X,Y\in \overline{\rnm}\setminus\{0\}$ such that $X\neq Y$, we have that
\begin{equation*}
\lim_{i\to +\infty}r_i^{n-2}G(\T\left(r_iX\right),\T\left(r_iY\right))={\mathcal G}(X,Y)
\end{equation*}
where ${\mathcal G}(X,Y)={\mathcal G}_X(Y)$ is the Green's function for $-\Delta-\gamma|x|^{-2}$ on $\rnm$ with Dirichlet boundary condition. Moreover, the convergence holds in $C^2_{loc}((\overline{\rnm}\setminus\{0\})^2\setminus\hbox{Diag}(\overline{\rnm}\setminus\{0\}))$.

\end{lemma}

\medskip\noindent{\it Proof of Lemma \ref{th:cv:1}:} We let $(r_i)_i\in (0,+\infty)$ and $(x_i)_i\in\Omega$ as in the statement of the lemma. For any $X,Y\in\rn$, $X\neq Y$, we define
$$G_i(X,Y):=r_i^{n-2}G(x_i+r_iX,x_i+r_iY)$$
for all $i\in\nn$. Since $r_i=o(d(x_i,\partial\Omega))$ as $i\to +\infty$, for any $R>0$, there exists $i_0\in\nn$ such that this definition makes sense for any $X,Y\in B_R(0)$. Equation \eqref{eq:G:c2} yields 
\begin{equation}\label{est:eq:Gi:1}
-\Delta G_i(X,\cdot)-\left(\frac{\gamma}{\left|\frac{x_i}{r_i}+\cdot\right|^2}+r_i^2 h(x_i+r_i\cdot)\right)G_i(X,\cdot)=0\hbox{ in }B_R(0)\setminus\{X\}.
\end{equation}
The pointwise control \eqref{G:up} writes
\begin{equation}\label{est:Gi:1}
0< G_i(X,Y)\leq c \left(\frac{\max\{|x_i+r_i X|,|x_i+r_i Y|\}}{\min\{|x_i+r_i X|,|x_i+r_i Y|\}}\right)^{\bm}|X-Y|^{2-n}
\end{equation}
for all $X,Y\in B_R(0)$ such that $X\neq Y$. Since $0\in\partial\Omega$, we have that $d(x_i,\partial\Omega)\leq |x_i|$, and therefore $r_i=o(|x_i|)$ as $i\to +\infty$. Equation \eqref{est:eq:Gi:1} and inequality \eqref{est:Gi:1} yield
\begin{equation*}
-\Delta G_i(X,\cdot)+\theta_i(X,\cdot)G_i(X,\cdot)=0\hbox{ in }B_R(0)\setminus\{X\}.
\end{equation*}
where $\theta_i\to0$ uniformly in $C^0_{loc}((\rn)^2)$ and $0< G_i(X,Y)\leq c |X-Y|^{2-n}$ for all $X,Y\in B_R(0)$ such that $X\neq Y$. It then follows from standard elliptic theory that, up to a subsequence,  there exists $G_\infty(X,\cdot)\in C^2(\rn\setminus\{X\})$ such that $G_i(X,\cdot)\to G_\infty(X,\cdot)\geq 0$ in $C^2_{loc}(\rn\setminus\{X\})$ and
$$-\Delta G_\infty(X,\cdot)=0\hbox{ in }\rn\setminus\{X\}\hbox{ and }G_\infty(X,Y)\leq c |X-Y|^{2-n}\hbox{ for }X,Y\in \rn,\; X\neq Y.$$
It then follows from the classification of positive harmonic functions that there exists $\lambda>0$ such that $G_\infty(X,Y)=\lambda |X-Y|^{2-n}$ for all $X,Y\in\rn$, $X\neq Y$.

\smallskip\noindent We fix $\varphi\in C^\infty_c(\rn)$. We define $\varphi_i(x):=\varphi(r_i^{-1}(x-x_i))$ for $x\in \Omega$ (this makes sense for $i$ large enough). It follows from \eqref{eq:G:dist} that
$$\varphi_i(x_i+r_i X)=\int_\Omega G(x_i+r_i X, y)\left(-\Delta\varphi_i(y)-\left(\frac{\gamma}{|y|^2}+h(y)\right)\varphi_i(y)\right)\, dy.$$
Via a change of variable, and passing to the limit, we get that 
$$\varphi(X)=\int_{\rn} G_\infty(X,Y)\left(-\Delta\varphi(Y)\right)\, dy.$$
Since $G_\infty(X,Y)=\lambda |X-Y|^{2-n}$, we get  that $\lambda=1/((n-2)\omega_{n-1}) $. Since the limit is unique, the convergence holds without extracting a subsequence. The convergence in $C^2_{loc}((\rn)^2\setminus\hbox{Diag}(\rn))$ follows from the symmetry of $G$ and elliptic theory.\qed

\medskip\noindent{\it Proof of Lemma \ref{th:cv:2}:} The proof goes  as in the proof of lemma \ref{th:cv:1}, except that we have to take a chart due to the closeness of the boundary. We let $(r_i)_i\in (0,+\infty)$, $(x_i)_i\in\partial\Omega$ and $x_0\in \partial\Omega$ as in the statement of the lemma. We let $\T$ be a chart at $x_0$ as in \eqref{def:T:bdry} (in particular $D_{0} \mathcal{T} = \mathbb{I}_{\R^{n}}$) and we set $x'_i\in\rn$ such that $x_i=\T(0, x_i')$. In particular, $\lim_{i\to +\infty}x_i'=0$. For any $X,Y\in\overline{\rnm}$, $X\neq Y$, we define
$$G_i(X,Y):=r_i^{n-2}G(\T\left((0,x'_i)+r_iX\right),\T\left((0,x'_i)+r_iY\right))$$
for all $i\in\nn$. Here again, provided $X,Y$ remain in a given compact set, the definition of $G_i$ makes sense for large $i$. Equation \eqref{eq:G:c2} then rewrites
\begin{equation}\label{est:eq:Gi:2}
-\Delta_{g_i} G_i(X,\cdot)-\hat{\theta}_iG_i(X,\cdot)=0\hbox{ in }B_R(0)\cap\rnm\setminus\{X\}\hbox{ ; }G_i(X,\cdot)\equiv 0\hbox{ on }\partial\rnm\cap B_R(0)
\end{equation}
where 
$$\hat{\theta}_i(Y):=\frac{\gamma}{\left|\frac{\T((0,x'_i)+r_i Y)}{r_i}\right|^2}+r_i^2 h(\T((0,x'_i)+r_i Y))$$
and $g_i=\T^\star\eucl((0,x'_i)+r_i \cdot)$ is the pull-back of the Euclidean metric. In particular, since $D_{0} \mathcal{T} = \mathbb{I}_{\R^{n}}$, we get that $g_i\to\eucl$ in $C^2_{loc}(\rn)$. Since $r_i=o(|x_i|)$, we get that $r_i=o(|x'_i|)$ as $i\to +\infty$, and, using again that $D_{0} \mathcal{T} = \mathbb{I}_{\R^{n}}$, we get that $\hat{\theta}_i\to 0$ uniformly in $B_R(0)\cap\rnm$. The pointwise control \eqref{G:up} rewrite $G_i(X,Y)\leq c |X-Y|^{2-n}$ for all $X,Y\in \rnm$, $X\neq Y$. With the same arguments as above, we get that for any $X\in \overline{\rnm}$, there exists $G_\infty(X,\cdot)\in C^2(\overline{\rnm}\setminus\{X\})$ such that
$$\lim_{i\to +\infty}G_i(X,\cdot)=G_\infty(X,\cdot)\hbox{ in }C^2_{loc}(\overline{\rnm}\setminus\{X\})$$
$$\hbox{ with }\left\{\begin{array}{ll}
-\Delta G_\infty(X,\cdot)=0&\hbox{ in }\rnm\setminus\{X\}\\
G_\infty(X,\cdot)\geq 0&\\
G_\infty(X,\cdot)\equiv 0&\hbox{ on }\partial\rnm\setminus\{X\}
\end{array}\right.$$
and
$$\varphi(X)=\int_{\rnm}G_\infty(X,\cdot)(-\Delta\varphi)\, dY\hbox{ for all }\varphi\in C^\infty_c(\rnm).$$
with $0\leq G_\infty(X,Y)\leq c |X-Y|^{2-n}$ for all $X,Y\in \rnm$, $X\neq Y$. Define
$$\Gamma_{\rnm}(X,Y)=\frac{1}{(n-2)\omega_{n-1}}\left(|X-Y|^{2-n}-|X-Y^*|^{2-n}\right).$$
As one checks (see for instance \cite{rob.green}), $\Gamma_{\rnm}$ satisfies the same properties as $G_\infty$. We set $f:=G_\infty(X,\cdot)-\Gamma_{\rnm}(X,\cdot)$. As one checks, $f\in C^\infty(\overline{\rnm}\setminus\{X\})$, $-\Delta f=0$ in the distribution sense in $\rnm$, $|f|\leq C|X-\cdot|^{2-n}$  in $\rnm\setminus\{X\}$ and $f_{\partial\rnm}=0$. Hypoellipticity yields $f\in C^\infty(\overline{\rnm})$. Multiplying $-\Delta f$ by $f$ and integrating by parts, we get that $f\equiv 0$, and then $G_\infty(X,\cdot)=\Gamma_{\rnm}(X,\cdot)$.  As above, this proves the convergence without any extraction. The convergence in $C^2_{loc}((\overline{\rnm})^2\setminus\hbox{Diag}(\overline{\rnm}))$ follows from the symmetry of $G$ and elliptic theory.\qed

\medskip\noindent{\it Proof of Lemma \ref{th:cv:3}:} Here again, the proof is similar to the two preceding proofs. We let $(r_i)_i\in (0,+\infty)$ such that $\lim_{i\to +\infty}r_i=0$. We let $\T$ be a chart at $0$ as in \eqref{def:T:bdry} (in particular $D_{0} \mathcal{T} = \mathbb{I}_{\R^{n}}$). For any $X,Y\in\overline{\rnm}\setminus\{0\}$, we define
$$G_i(X,Y):=r_i^{n-2}G(\T\left(r_iX\right),\T\left(r_iY\right))$$
for all $i\in\nn$. Equation \eqref{eq:G:c2}  rewrites
\begin{equation*}
-\Delta_{g_i} G_i(X,\cdot)-\left(\frac{\gamma}{\left|\frac{\T(r_i \cdot)}{r_i}\right|^2}+r_i^2 h(\T(r_i \cdot))\right)G_i(X,\cdot)=0\hbox{ in }B_R(0)\cap\rnm\setminus\{0,X\}.
\end{equation*}
with $G_i(X,\cdot)\equiv 0$ on $B_R(0)\cap\partial\rnm$, where $g_i=\T^\star\eucl(r_i \cdot)$ is the pull-back of the Euclidean metric. In particular, since $D_{0} \mathcal{T} = \mathbb{I}_{\R^{n}}$, we get that $g_i\to\eucl$ in $C^2_{loc}(\rn)$. The pointwise control \eqref{G:up} writes
$$0\leq G_i(X,Y)\leq C \left(\frac{\max\{|X|,|Y|\}}{\min\{|X|,|Y|\}}\right)^{\bm}|X-Y|^{2-n}\hbox{ for }X,Y\in\rnm,\; X\neq Y.
$$
It then follows from elliptic theory that $G_i(X,\cdot)\to G_\infty(X,\cdot)$ in $C^2_{loc}(\overline{\rnm}\setminus\{0,X\})$. In particular, $G_\infty(X,\cdot)$ vanishes on $\partial\rnm\setminus\{0\}$ and 
\begin{equation}\label{est:G:rn}
0\leq G_\infty(X,Y)\leq C \left(\frac{\max\{|X|,|Y|\}}{\min\{|X|,|Y|\}}\right)^{\bm}|X-Y|^{2-n}\hbox{ for }X,Y\in\rnm,\; X\neq Y.
\end{equation}
Moreover, passing to the limit in Green's representation formula, we get that
$$\varphi(X)=\int_{\rnm}G_\infty(X,Y)\left(-\Delta \varphi-\frac{\gamma }{|Y|^2}\varphi\right)\, dY\hbox{ for all }\varphi\in C^\infty_c(\rnm).$$
Since $G(x,\cdot)$ is locally in $\huno$ (see (b) in Theorem \ref{th:green:gamma:domain}), we get that $(\eta G_i(X,\cdot))_i$ is uniformly bounded in $H_{1,0}^2(\rnm)$ for all $\eta\in C^\infty_c(\rn\setminus\{X\})$. Up to another extraction, we get weak convergence in $H_{1,0}^2(\rnm)$, and then $\eta G_\infty(X,\cdot)\in H_{1,0}^2(\rnm)$ for all $\eta\in C^\infty_c(\rn\setminus\{X\})$. It then follows from Theorem \ref{th:green:gamma:rn} and \eqref{est:G:rn} that $G_\infty(X,\cdot)={\mathcal G}_{X}$ is the unique Green's function of $-\Delta-\gamma|x|^{-2}$ on $\rnm$ with Dirichlet boundary condition. Here again, the convergence in $C^2$ follows from elliptic theory.\qed

\subsection{A lower bound for the Green's function}\label{sec:lower}
We let $\Omega$, $\gamma$, $h$ be as in Theorem \ref{th:green:gamma:domain}. We let $G$ be the Green's function for $-\Delta-(\gamma|x|^{-2}+h)$ on $\Omega$ with Dirichlet boundary condition. We let $(x_i), (y_i)_{i\in\nn}$ be such that $x_i,y_i\in\Omega$ and $x_i\neq y_i$ for all $i\in\nn$. We also assume that there exists $x_\infty,y_\infty\in \overline{\Omega}$ such that 
$$\lim_{i\to +\infty}x_i=x_\infty\hbox{ and }\lim_{i\to +\infty}y_i=y_\infty$$
and that there exists $c_1,c_2$ such that
\begin{equation*}
\lim_{i\to +\infty}\frac{G(x_i,y_i)}{H(x_i,y_i)}=c_1\in [0,+\infty]
\hbox{ and }\lim_{i\to +\infty}\frac{|\nabla G_{x_i}(y_i)|}{\Gamma(x_i,y_i)}
=c_2\in [0,+\infty]
\end{equation*}
where $H(x,y)$ is defined in \eqref{def:Hp:1} and 
$$\Gamma(x,y):=\left(\frac{\max\{|x|,|y|\}}{\min\{|x|,|y|\}}\right)^{\bm}|x-y|^{1-n}\min\left\{1,\frac{d(x,\partial\Omega)}{|x-y|}\right\}$$
for $x,y\in\Omega$, $x\neq y$. Note that $c_1<+\infty$ by \eqref{G:up}. We claim that
\begin{equation}\label{lim:c}
0<c_1\hbox{ and }0\leq c_2<+\infty
\end{equation}
The lower bound in \eqref{est:G:up} and the upper bound in \eqref{ineq:grad:G} both follow from \eqref{lim:c}.

\medskip\noindent This section is devoted to proving \eqref{lim:c}. We distinguish several cases:

\medskip\noindent{\bf Case 1: $x_\infty\neq y_\infty$, $x_\infty,y_\infty\in \Omega$.} As one checks, we then have that 
$$\lim_{i\to +\infty}G(x_i,y_i)=G(x_\infty,y_\infty)> 0.$$
Therefore, we get that $c_1\in (0,+\infty)$. Concerning the gradient, $\lim_{i\to +\infty}|\nabla G_{x_i}(y_i)|=|\nabla G_{x_\infty}(y_\infty)|\geq 0$ and this yields $c_2<+\infty$. This proves \eqref{lim:c} in Case 1.

\medskip\noindent{\bf Case 2: $x_\infty\in\Omega$ and $y_\infty\in \partial\Omega\setminus\{0\}$.}  Since $x_\infty,y_\infty$ are distinct and far from $0$, we have that $G(x_i,y_i)=d(y_i,\partial\Omega)\left(-\partial_{\nu}G_{x_\infty}(y_\infty)+o(1)\right)$ as $i\to +\infty$, where $\partial_\nu G_{x_\infty}(y_\infty)$ is the normal derivative of $G_{x_\infty}>0$ at the boundary point $y_\infty$. Hopf's Lemma then yields $\partial_\nu G_{x_\infty}(y_\infty)<0$. As one checks, we have that $H(x_i,y_i)=(c+o(1))d(y_i,\partial\Omega)$ as $i\to +\infty$. This then yields $0<c_1<+\infty$. Concerning the gradient, we get that $\lim_{i\to +\infty}|\nabla G_{x_i}(y_i)|=|\nabla G_{x_\infty}(y_\infty)|\geq 0$ and $\lim_{i\to +\infty}\Gamma(x_i,y_i)\in (0,+\infty)$, which yields $c_2<+\infty$. This proves \eqref{lim:c} in Case 2.

\medskip\noindent{\bf Case 3: $x_\infty\in\Omega$ and $y_\infty=0\in\partial\Omega$.} It follows from Case 2 above that there exists $c>0$ such that $G_{x_i}(y)\geq c d(y,\partial\Omega)|y|^{-\bm}$ for all $y\in \partial (\Omega\cap B_{r_0}(0))$. We take the subsolution $\underline{u}_{\bm}$ defined in Proposition \ref{prop:sub:super}. With \eqref{asymp:ua:plus}, there exists $c'>0$ such that $G_{x_i}(y)\geq c_1\underline{u}_{\bm}(y)$ for all $y\in \partial (\Omega\cap B_{r_0}(0))$. Since $G_{x_i}$ is locally in $H_{1,0}^2$ around $0$, the comparison principle and \eqref{asymp:ua:plus} yields $G_{x_i}(y)\geq c" d(y,\partial\Omega)|y|^{-\bm}$ for all $y\in\Omega\cap B_{r_0}(0)$. This yields $c_1>0$. 

\medskip\noindent We deal with the gradient. We let $\T$ be a chart at $0$ as in \eqref{def:T:bdry} and we define
$$G_i(y):=r_i^{\bm-1}G_{x_i}(\T(r_i y))\hbox{ for }y\in \rnm\cap B_2(0)$$
with $r_i\to 0$. It follows from \eqref{G:up}  that $G_i(y)\leq C |y_1|\cdot |y|^{-\bm}$ for all $y\in \rnm\cap B_2(0)$. It follows from \eqref{eq:G:c2} that $-\Delta_{g_i}G_i-\left(\gamma|\cdot|^2+o(1)\right)G_i=0$ in $\rnm\cap B_2(0)$ where $g_i:=\T^\star\eucl(r_i \cdot)$ and $o(1)\to 0$ in $L^\infty_{loc}(\rn)$. Elliptic regularity then yields $|\nabla G_i(y)|\leq C$ for $y\in \rnm\cap B_{3/2}(0)$. We now let $r_i:=|\tilde{y}_i|$ where $y_i:=\T(\tilde{y}_i)$, so that $r_i\to 0$. We then have that $|\nabla G_i(\tilde{y_i}/r_i)|\leq C$, which rewrites $|\nabla G_{x_i}(y_i)|\leq C|y_i|^{-\bm}$. By estimating $\Gamma(x_i,y_i)$, we then get that $c_2<+\infty$. This proves \eqref{lim:c} in Case 3.

\medskip\noindent{\bf Case 4:  $x_\infty\neq y_\infty$, $x_\infty,y_\infty\in \partial\Omega\setminus\{0\}$.} Since $x_\infty,y_\infty$ are distinct and far from $0$, we have that $G(x_i,y_i)=d(y_i,\partial\Omega)d(x_i,\partial\Omega)\left(\partial_{\nu_x}\partial_{\nu_y}G_{x_\infty}(y_\infty)+o(1)\right)$ as $i\to +\infty$, where $\partial_{\nu_x}$ is the normal derivative along the first coordinate,  and $\partial_{\nu_y}$ is the normal derivative along the second coordinate. Since $y\mapsto G_x(y)$ is positive for $x,y\in\Omega$, $x\neq y$, and solves \eqref{eq:G:c2}, Hopf's maximum principle yields $-\partial_{\nu_y}G(x,y_\infty)>0$ for $x\in \Omega$. Moreover, it follows from the symmetry of $G$ that $-\partial_{\nu_y}G(x,y_\infty)>0$ solves also \eqref{eq:G:c2}. Another application of Hopf's principle yields $\partial_{\nu_x}\partial_{\nu_y}G_{x_\infty}(y_\infty)>0$. Estimating independently $H(x_i,y_i)$, we get that $0<c_1<+\infty$.

\smallskip\noindent We deal with the gradient. We have that $|\nabla_y G_{x_i}(y_i)|=|\nabla_y (G_{x_i}-G_{\tilde{x_i}})(y_i)|$ where $\tilde{x}_i\in\partial\Omega$ is the projection of $x_i$ on $\partial\Omega$. The $C^2-$control then yields $|\nabla_y G_{x_i}(y_i)|\leq Cd(x_i,\partial\Omega)$. Estimating independently $\Gamma(x_i,y_i)$, we get that $c_2<+\infty$. This proves \eqref{lim:c} in Case 4.

\medskip\noindent{\bf Case 5:  $x_\infty\neq y_\infty$, $x_\infty\in \partial\Omega\setminus\{0\}$ and $y_\infty=0$.} It follows from Cases 2 and 4 that $G_{x_i}(y)\geq C d(x_i,\partial\Omega)d(y_i,\partial\Omega)$ for all $y\in  \partial(B_{|x_\infty|/2}(0)\cap\Omega)$. Using a sub-solution as in Case 3, we get that 
$G_{x_i}(y)\geq c d(x_i,\partial\Omega)d(y,\partial\Omega)|y|^{-\bm}$ for all $y\in  \partial(B_{|x_\infty|/2}(0)\cap\Omega)$. This yields $0<c_1$.

\medskip\noindent For the gradient estimate, we choose a chart $\T$ around $y_\infty=0$ as in \eqref{def:T:bdry}, and we let $r_i:=|\tilde{y}_i|\to 0$ where $y_i=\T(\tilde{y}_i)$we define $G_i(y):=r_i^{\bm-1}G_{x_i}(\T(r_i y))/d(x_i,\partial\Omega)$ for $y\in \rnm\cap B_2(0)$ where $r_i\to 0$ . The pointwise control \eqref{G:up} and equation \eqref{eq:G:c2} yields the convergence of $(G_i)$ in $C^1_{loc}(\overline{\rnm}\cap B_2(0)\setminus\{0\})$ as $i\to +\infty$. The boundedness of $|\nabla G_i|$ yields $c_2<+\infty$. This proves \eqref{lim:c} in Case 5.

\medskip\noindent Since $G$ is symmetric, it follows from Cases 1 to 5 that \eqref{lim:c} holds when $x_\infty\neq y_\infty$.

\medskip\noindent We now deal with the case $x_\infty=y_\infty$, which rewrites $\lim_{i\to +\infty}|x_i-y_i|=0$. Via a rescaling, we are essentially back to the case $x_\infty\neq y_\infty$ via the convergence Theorems \ref{th:cv:1}, \ref{th:cv:2} and \ref{th:cv:3}. 

\medskip\noindent{\bf Case 6:  $|x_i-y_i|=o(d(x_i,\partial\Omega))$ as $i\to +\infty$.} We set $r_i:=|x_i-y_i|\to 0$ as $i\to +\infty$ and we define
$$G_i(Y):=r_i^{n-2}G(x_i, x_i+r_iY)\hbox{ for }Y\in \frac{\Omega-x_i}{r_i}\setminus\{0\}.$$
It follows from Theorem \ref{th:cv:1} that $G_i\to c_n|\cdot|^{2-n}$ in $C^2_{loc}(\rn\setminus\{0\})$ as $i\to +\infty$, with $c_n:=((n-2)\omega_{n-1})^{-1}$. We define $Y_i:=\frac{y_i-x_i}{|y_i-x_i|}$, and we then get that $|y_i-x_i|^{n-2}G(x_i,y_i)=G_i(Y_i)\to c_n$ as $i\to +\infty$. Estimating $H(x_i,y_i)$ (and noting that $d(x_i,\partial\Omega)\leq |x_i-0|=|x_i|$), we get that $0<c_1<+\infty$.

\smallskip\noindent The convergence of the gradient yields $|\nabla G_i(Y_i)|\leq C$ for all $i$. With the original function $G$ and points $x_i$, $y_i$, this yields $c_2<+\infty$. This proves \eqref{lim:c} in Case 6.

\medskip\noindent{\bf Case 7:  $d(x_i,\partial\Omega)=O(|x_i-y_i|)$ and $|x_i-y_i|=o(|x_i|)$ as $i\to +\infty$.} Then $\lim_{i\to +\infty}x_i=x_\infty\in \partial\Omega$. We let $\T$ be a chart at $x_\infty$ as in \eqref{def:T:bdry}, in particular $D_{0} \mathcal{T} = \mathbb{I}_{\R^{n}}$. We let $x_i=\T(x_{i,1}, x_i')$ and $y_i=\T(y_{i,1}, y_i')$ where $(x_{i,1}, x_i'),(y_{i,1}, y_i')\in (-\infty, 0)\times \rr^{n-1}$ are going to $0$ as $i\to +\infty$. In particular $d(x_i,\partial\Omega)=(1+o(1))|x_{i,1}|$ and $d(y_i,\partial\Omega)=(1+o(1))|y_{i,1}|$  as $i\to +\infty$. We define $r_i:=|(y_{i,1}, y_i')-(x_{i,1}, x_i')|$. In particular $r_i=(1+o(1))|x_i-y_i|$ as $i\to +\infty$. The hypothesis of Case 7 rewrite $x_{i,1}=O(r_i)$ and $r_i=o(|(x_{i,1}, x_i')|)$. Consequently, we have that $y_{i,1}=O(r_i)$ and $r_i=o(|x_{i}'|)$ as $i\to +\infty$. We define
$$G_i(X,Y):=r_i^{n-2}G(\T\left((0,x'_i)+r_iX\right),\T\left((0,x'_i)+r_iY\right))$$
for $X,Y\in \rnm$ such that $X\neq Y$. It follows from Theorem \ref{th:cv:2} that
$$\lim_{i\to +\infty}G_i(X,Y)=c_n\left(|X-Y|^{2-n}-|X-Y^*|^{2-n}\right):=\Psi(X,Y)$$
for all $X,Y\in \overline{\rnm}$, $X\neq Y$, and this convergence holds in $C^2_{loc}$. We define $X_i:=(r_i^{-1}x_{i,1},0)$ and $Y_i:=(r_i^{-1}y_{i,1}, r_i^{-1}(y_i'-x_i'))$: the definition of $r_i$ yields $X_i\to X_\infty\in \overline{\rnm}$ and $Y_i\to Y_\infty\in \overline{\rnm}$ as $i\to +\infty$. Therefore, we get that
$$|x_i-y_i|^{n-2}G(x_i,y_i)=(1+o(1))G_i(X_i,Y_i)\to \Psi(X_\infty,Y_\infty)$$
as $i\to +\infty$, and 
\begin{equation}\label{lim:bdy}
|X_{\infty,1}|=\lim_{i\to +\infty}\frac{|x_{i,1}|}{r_i}=\lim_{i\to +\infty}\frac{d(x_i,\partial\Omega)}{r_i}.
\end{equation}
\noindent{\it Case 7.1: $X_{\infty,1}\neq 0$ and $Y_{\infty,1}\neq 0$.} We then get that $\lim_{i\to +\infty}|x_i-y_i|^{n-2}G(x_i,y_i)=\Psi(X_\infty,Y_\infty)>0$. Moreover, it follows from \eqref{lim:bdy} that $d(x_i,\partial\Omega)d(y_i,\partial\Omega)=(c+o(1))|x_i-y_i|^2$ as $i\to +\infty$ for some $c>0$. Since $|x_i|=(1+o(1))|y_i|$ as $i\to +\infty$ (this follows from the assumption of Case 7), we get that $\lim_{i\to +\infty}|x_i-y_i|^{n-2}H(x_i,y_i)\in (0,+\infty)$. Then $0<c_1<+\infty$.

\noindent{\it Case 7.2: $X_{\infty,1}\neq 0$ and $Y_{\infty,1}=0$.} Then $Y_{i,1}\to 0$ as $i\to +\infty$, and then, there exists $(\tau_i)_i\in (0,1)$ such that $G_i(X_i,Y_i)=Y_{i,1}\partial_{Y_1}G_i(X_i, (\tau_i Y_{i,1}, Y'_i))$. Letting $i\to +\infty$ and using the convergence of $G_i$ in $C^1$, we get that
\begin{eqnarray*}
|x_i-y_i|^{n-2}G(x_i,y_i)&=&(1+o(1))G_i(X_i,Y_i)=Y_{i,1}\partial_{Y_1}G_i(X_i, \tau_i Y_i)\\
&=&\frac{d(y_i,\partial\Omega)}{|x_i-y_i|}\left(-\partial_{Y_1}\Psi(X_\infty,Y_\infty)+o(1)\right)
\end{eqnarray*}
as $i\to +\infty$. As one checks, $\partial_{Y_1}\Psi(X_\infty,Y_\infty)<0$. Arguing as in Case 7.1, we get that $0<c_1<+\infty$.

\noindent{\it Case 7.3: $X_{\infty,1}=Y_{\infty,1}=0$.} As in Case 7.2, there exists $(\tau_i)_i,(\sigma_i)_i\in (0,1)$ such that
$G_i(X_i,Y_i)=Y_{i,1}X_{i,1}\partial_{Y_1}\partial_{X_1}G_i((\sigma_i X_{i,1}, X'_i)X_i, (\tau_i Y_{i,1}, Y'_i))$. We conclude as above, noting that $\partial_{Y_1}\partial_{X_1}\Psi(X_\infty,Y_\infty)>0$. Then $0<c_1<+\infty$.

\medskip\noindent The gradient estimate is proved as in Cases 1 to 6. This proves \eqref{lim:c} in Case 7.

\medskip\noindent{\bf Case 8:  $d(x_i,\partial\Omega)=O(|x_i-y_i|)$, $|x_i|=O(|x_i-y_i|)$ and $|y_i|=O(|x_i-y_i|)$ as $i\to +\infty$.} In particular, $x_\infty=y_\infty=0$. We take a chart at $0$ as in Case 7, and we define $(x_{i,1},x_i'),(y_{i,1},y_i')$ similarly. We define $r_i:=|(y_{i,1}, y_i')-(x_{i,1}, x_i')|=(1+o(1))|x_i-y_i|$ as $i\to +\infty$. We define
$$G_i(X,Y):=r_i^{n-2}G(\T\left(r_iX\right),\T\left(r_iY\right))$$
for $X,Y\in \rnm$. It follows from Theorem \ref{th:cv:3} that $G_i\to {\mathcal G}$ in $C^2_{loc}((\overline{\rnm}\setminus\{0\})^2\setminus\hbox{Diag}(\overline{\rnm}\setminus\{0\}))$, where $\mathcal G$ is the Green's function for $-\Delta-\gamma|\cdot|^{-2}$ in $\rnm$. Then$$|x_i-y_i|^{n-2}G(x_i,y_i)=(1+o(1))G_i(X_i,Y_i)={\mathcal G}(X_\infty,Y_\infty)+o(1)$$
as $i\to +\infty$. 

\medskip\noindent{\it Case 8.1: We assume that $X_{\infty,1}\neq 0$ and $Y_{\infty,1}\neq 0$.} Then we get $0<c_1<+\infty$ as in Case 7.1.

\medskip\noindent{\it Case 8.2: We assume that $X_\infty\in \rnm$ and $Y_\infty\in \partial\rnm\setminus\{0\}$ or $X_\infty,Y_\infty\in \partial\rnm\setminus\{0\}$}. Then we argue as in Cases 7.2 and 7.3 to get $0<c_1<+\infty$ provided $\{\partial_{Y_1}{\mathcal G}(X_\infty,Y_\infty)<0\hbox{ if }X_\infty\in \rnm\hbox{ and }Y_\infty\in \partial\rnm\}$ and $\{\partial_{Y_1}\partial_{X_1}{\mathcal G}(X_\infty,Y_\infty)>0\hbox{ if }X_\infty,Y_\infty\in \partial\rnm\}$.  So we are just left with proving these two inequalities. 

\medskip\noindent We assume that $X_\infty\in \rnm$. It follows from Theorem \ref{th:green:gamma:rn} below that ${\mathcal G}(X_\infty,\cdot)>0$ is a solution to $(-\Delta-\gamma|\cdot|^{-2}){\mathcal G}(X_\infty,\cdot)=0$ in $\rnm-\{X_\infty\}$, vanishing on $\partial\rnm\setminus\{0\}$. Hopf's maximum principle then yields $-\partial_{Y_1}{\mathcal G}(X_\infty,Y_\infty)>0$ for $Y_\infty\in \partial\rnm\setminus\{0\}$.

\medskip\noindent We fix $Y_\infty\in \partial\rnm\setminus\{0\}$. For $X\in\rnm$, we then define $H(X):=-\partial_{Y_1}{\mathcal G}(X,Y_\infty)>0$ by the above argument. Moreover, $(-\Delta-\gamma|\cdot|^{-2})H=0$ in $\rnm$, vanishing on $\partial\rnm\setminus\{0, Y_\infty\}$. Hopf's maximum principle yields $-\partial_{X_1}H(X_\infty)=\partial_{Y_1}\partial_{X_1}{\mathcal G}(X_\infty,Y_\infty)>0$ for $X_\infty,Y_\infty\in \partial\rnm\setminus\{0\}$

 \medskip\noindent{\it Case 8.3: we assume that $X_\infty=0$ or $Y_\infty=0$.} Since $|X_\infty-Y_\infty|=1$, without loss of generality, we can assume that $X_\infty\neq 0$. It follows from Cases 8.1 and 8.2 that there exists $C>0$ such that 
 \begin{equation}\label{low:G:34}
 C^{-1}\frac{d(x_i,\partial\Omega)}{|x_i|^{n-\bm}}\frac{d(y,\partial\Omega)}{|y|^{\bm}}\leq G_{x_i}(y)\leq  C\frac{d(x_i,\partial\Omega)}{|x_i|^{n-\bm}}\frac{d(y,\partial\Omega)}{|y|^{\bm}}
 \end{equation}
 for all $y\in \partial (B_{|x_i|/2}(0)\cap\Omega)$. We let $\underline{u}_{\bm}$ be the sub-solution given by Proposition \ref{prop:sub:super}. Arguing as in Case 3, it then follows from the comparison principle that \eqref{low:G:34} holds for $y\in B_{|x_i|/2}(0)\cap\Omega$. Since $|y_i|=o(|x_i|)$, we then get that  \eqref{low:G:34} holds with $y:=y_i$. Estimating $H(x_i,y_i)$, we then get that $0<c_1<+\infty$.

 \medskip\noindent The gradient estimate is proved as in Cases 1 to 6. This proves \eqref{lim:c} in Case 8.

\medskip\noindent Since $G$ is symmetric, it follows from Cases 7 and 8 that \eqref{lim:c} holds when $x_\infty=y_\infty$.

\medskip\noindent In conclusion, we get that \eqref{lim:c} holds, which proves the initial claim. As noted previously, both the lower bound in \eqref{est:G:up} and the upper bound in \eqref{ineq:grad:G} follow from these results. 

\smallskip\noindent We are now left with proving \eqref{est:G:4}. We let $(\txi)_i, (\tyi)_i\in \Omega$ be such that
$$\tyi=o(|\txi|)\hbox{ and }\txi=o(1)\hbox{ as }i\to +\infty,$$
and $(h_i)_i\in C^{0,\theta}(\Omega)$ such that $\lim \limits_{i\to+\infty}h_i=h$ in $C^{0,\theta}$. It follows from \eqref{est:G:up} that, up to extraction, there exists $l>0$ such that
\begin{equation}\label{est:G:pf}
G_{h_i}(\txi,\tyi)=(l+o(1))\frac{d(\txi,\partial\Omega)}{|\txi|^{\bp}}\frac{d(\tyi,\partial\Omega)}{|\tyi|^{\bm}}
\end{equation}
From now on, to avoid unnecessary notations, the extraction is \underline{fixed}. We define 
$$r_i:=|\txi|\, ;\, s_i:=|\tyi|\, ;\, \tau_i:=s_i^{-1}\T^{-1}(\tyi)\in\rnm\hbox{ and }\theta_i:=r_i^{-1}\T^{-1}(\txi)\in\rnm,$$
and $\theta_\infty,\tau_\infty\in \overline{\rnm}$ such that
\begin{equation}\label{def:theta}
\txi=\T(r_i\theta_i)\, ;\, \tyi=\T(s_i\tau_i)\, ; \, \theta_i\to \theta_\infty\neq 0\hbox{ and }\tau_i\to\tau_\infty\neq 0\hbox{ as }i\to +\infty.
\end{equation}

\begin{step}\label{step:G:pf:1} We fix $R>0$. We claim that
\begin{equation}\label{cv:G:pf:1}
G_{h_i}(\txi,y)=(l+o(1))\frac{d(\txi,\partial\Omega)}{|\txi|^{\bp}}\frac{d(y,\partial\Omega)}{|y|^{\bm}}\hbox{ as }i\to +\infty
\end{equation}
uniformly for $y\in \Omega\cap \T(B_{Rs_i}\setminus B_{R^{-1}s_i})$.
\end{step}

\smallskip\noindent{\it Proof of Step \ref{step:G:pf:1}:} For $z\in B_{2R}\setminus B_{(2R)^{-1}}$, we define
$$G_i(z):=\frac{s_i^{\bm-1}|\txi|^{\bp}}{d(\txi,\partial\Omega)}G_{h_i}(\txi, \T(s_i z)).$$
As one checks, \eqref{cv:G:pf:1} is equivalent to prove that
\begin{equation}\label{cv:G:pf:1:equiv}
G_i(y)=(l+o(1))\frac{|y_1|}{|y|^{\bm}}\hbox{ uniformly for }y\in B_{R}(0)\setminus B_{R^{-1}}(0)
\end{equation}
Since $s_i=o(|\txi|)$ and \eqref{rem:T:bdry} holds, it follows from the control \eqref{est:G:up} that there exists $C>0$ such that
\begin{equation}\label{pf:G:1}
\frac{1}{C}\cdot \frac{|z_1|}{|z|^{\bm}}\leq G_i(z)\leq C \cdot \frac{|z_1|}{|z|^{\bm}}\hbox{ for all }z\in\rnm\cap  B_{2R}(0)\setminus B_{(2R)^{-1}}(0).
\end{equation}
As for \eqref{est:eq:Gi:2}, it follows from \eqref{eq:G:c2} that
\begin{equation}\label{pf:G:2}
-\Delta_{g_i} G_i-\left(\frac{\gamma s_i^2}{|\T(s_i \cdot)|^2}+O(s_i^2)\right)G_i=0\hbox{ in }B_R(0)\cap\rnm\hbox{ ; }G_i\equiv 0\hbox{ on }\partial\rnm\cap B_R(0)\setminus\{0\}.
\end{equation}
It follows from \eqref{pf:G:1}, \eqref{pf:G:2} and standard elliptic theory that there exists $G\in C^2(\overline{\rnm}\setminus\{0\})$ such that, up to a subsequence, 
\begin{equation}\label{lim:Gi:pf}
\lim_{i\to +\infty}G_i=G\hbox{ in }C^2_{loc}(\overline{\rnm}\setminus\{0\})
\end{equation}
with
$$-\Delta G-\frac{\gamma}{|x|^2}G=0\hbox{ in }\overline{\rnm}\setminus\{0\}\, ;\, G=0\hbox{ on }\partial\rnm\setminus\{0\}\, ;$$
$$\frac{1}{C}\cdot \frac{|z_1|}{|z|^{\bm}}\leq G(z)\leq C \cdot \frac{|z_1|}{|z|^{\bm}}\hbox{ for all }z\in \rnm\setminus\{0\}.$$
It the follows from Proposition 6.4 in \cite{gr4} that there exists $\lambda>0$ such that
\begin{equation}\label{G:pf:explicit}
G(z)=\lambda \cdot \frac{|z_1|}{|z|^{\bm}}\hbox{ for all }z\in \rnm.
\end{equation}
We claim that $\lambda=l$. We prove the claim. It follows from \eqref{est:G:pf} and the definition \eqref{def:theta} of $\tau_i$ that
\begin{equation}\label{eq:G:l}
G_i(\tau_i)=(l+o(1))\frac{|\tau_{i,1}|}{|\tau_i|^{\bm}}\hbox{ and }\tau_i\to\tau_\infty\neq 0\hbox{ as }i\to +\infty.
\end{equation}

\smallskip\noindent{\it Case 1:} we assume that $\tau_\infty\in\rnm\setminus\{0\}$, that is $\tau_{\infty,1}\neq 0$. Passing to the limit in \eqref{eq:G:l}, using the convergence \eqref{lim:Gi:pf} and the explicit form \eqref{G:pf:explicit}, we get that
$$l\frac{|\tau_{\infty,1}|}{|\tau_{\infty}|^{\bm}}=\lambda\frac{|\tau_{\infty,1}|}{|\tau_\infty|^{\bm}},$$
and therefore, since $\tau_{\infty,1}\neq 0$, we get that $\lambda=l$.

\smallskip\noindent{\it Case 2:} we assume that $\tau_\infty\in \partial\rnm\setminus\{0\}$, that is $\tau_{i,1}\to 0$ as $i\to +\infty$. With a Taylor expansion, we get that there exists a sequence $(t_i)_{i\in\nn}\in (0,1)$ such that $G_i(\tau_i)=\partial_1G_i(t_i \tau_{i,1},\theta_i')\tau_{i,1}$ for all $i\in \nn$. With the convergence \eqref{lim:Gi:pf} of $G_i$ to $G$  in $C^1$, we get that 
$$G_i(\tau_i)=\left(\partial_1G(\tau_\infty)+o(1)\right)\cdot\tau_{i,1}=\left(\frac{\lambda}{|\tau_\infty|^{\bm}}+o(1)\right)\cdot |\tau_{i,1}|.$$
Since $\tau_{i,1}\neq 0$ for all $i\in\nn$, it follows form \eqref{eq:G:l} that $\lambda=l$.

\smallskip\noindent Therefore, in both cases, we have proved that $\lambda=l$. It follows from this uniqueness that the convergence of $G_i$ holds with no extraction. 

\smallskip\noindent We now prove \eqref{cv:G:pf:1:equiv}. We let $(z_i)_i\in \rnm\setminus\{0\}$ be such that $z_i\to z_\infty\in \overline{\rnm}\setminus\{0\}$. Then $G_i(z_i)\to G(z_\infty)$ as $i\to +\infty$. Therefore, if $z_{\infty,1}\neq 0$, we get that $G_i(z_i)=(1+o(1)) G(z_i)$ as $i\to +\infty$. We now assume that $z_{\infty,1}=0$, that is $z_{i,1}\to 0$ as $i\to +\infty$. We use the $C^1-$convergence of $(G_i)$ and argue as in Case 2 above to get that $\lim_{i\to +\infty}|z_{i,1}|^{-1}G_i(z_i)=-\partial_1G(z_\infty)\neq 0$. As one checks, this yields also $G_i(z_i)=(1+o(1)) G(z_i)$ as $i\to +\infty$. As noticed above, this proves \eqref{cv:G:pf:1} and ends Step \ref{step:G:pf:1}.\qed

\begin{step}\label{step:G:pf:2} We fix $R>0$. We claim that
\begin{equation}\label{cv:G:pf:2}
G_{h_i}(\txi,y)=(l+o(1))\frac{d(\txi,\partial\Omega)}{|\txi|^{\bp}}\frac{d(y,\partial\Omega)}{|y|^{\bm}}\hbox{ as }i\to +\infty
\end{equation}
uniformly for $y\in \Omega\cap \T(B_{Rs_i}(0))$.
\end{step}

\smallskip\noindent{\it Proof of Step \ref{step:G:pf:2}:} For $r>0$ small, we choose $\bar{u}_{\bm}\in C^2(\Omega\cap B_r(0))$ a supersolution to $-\Delta \bar{u}_{\bm}-(\gamma|x|^{-2}+h_i)\bar{u}_{\bm}>0$ as in \eqref{ppty:ua} and \eqref{asymp:ua:plus}. Note that, due to the convergence of $(h_i)$ to $h$ in $C^0$, the choice of $\bar{u}_{\bm}$ can be made independently of $i$. We fix $\eps>0$. It follows from the convergence \eqref{cv:G:pf:1} of Step \ref{step:G:pf:1} and \eqref{asymp:ua:plus} that there exists $i_0\in\nn$
\begin{equation}\label{est:sup:G:pf}
G_{h_i}(\txi,y)\leq (l+\eps)\frac{d(\txi,\partial\Omega)}{|\txi|^{\bp}}\bar{u}_{\bm}(y)\hbox{ for all }y\in \partial \left(\Omega\cap \T(B_{Rs_i}(0))\right)\hbox{ for all }i\geq i_0.
\end{equation}
Note that $G_{h_i}(\txi,\cdot), \bar{u}_{\bm}\in H_1^2\left(\Omega\cap \T(B_{Rs_i}(0))\right)$ (these are variational super- or sub-solutions) and that the operator $-\Delta -(\gamma|x|^{-2}+h_i)$ is coercive. Since $G_{h_i}(\txi,\cdot)$ is a solution and $\bar{u}_{\bm}$ is a supersolution to $-\Delta u-(\gamma|x|^{-2}+h_i)u=0$, it follows from the comparison principle that \eqref{est:sup:G:pf} holds for $y\in \Omega\cap \T(B_{Rs_i}(0))$. With \eqref{asymp:ua:plus}, we get that there exists $i_1\in\nn$ such that
\begin{equation}\label{est:sup:G:pf:2}
G_{h_i}(\txi,y)\leq (l+2\eps)\frac{d(\txi,\partial\Omega)}{|\txi|^{\bp}}\frac{d(y,\partial\Omega)}{|y|^{\bm}}\hbox{ for all }y\in \Omega\cap \T(B_{Rs_i}(0))\hbox{ for all }i\geq i_1.
\end{equation}
Using a subsolution $\underline{u}_{\bm}$ as in \eqref{ppty:ua} and \eqref{asymp:ua:plus} and arguing as above, we get that
\begin{equation}\label{est:sup:G:pf:3}
G_{h_i}(\txi,y)\geq (l-2\eps)\frac{d(\txi,\partial\Omega)}{|\txi|^{\bp}}\frac{d(y,\partial\Omega)}{|y|^{\bm}}\hbox{ for all }y\in \Omega\cap \T(B_{Rs_i}(0))\hbox{ for all }i\geq i_2.
\end{equation}
The inequalities \eqref{est:sup:G:pf:2} and \eqref{est:sup:G:pf:3} put together yield \eqref{cv:G:pf:2}. This ends Step \ref{step:G:pf:2}.\qed

\medskip\noindent We now vary the $x-$variable. 

\begin{step}\label{step:G:pf:3} We fix $R,R'>0$. We claim that
\begin{eqnarray}
&&G_{h_i}(\txi,y)=(l+o(1))\frac{d(x,\partial\Omega)}{|x|^{\bp}}\frac{d(y,\partial\Omega)}{|y|^{\bm}}\hbox{ as }i\to +\infty\label{cv:G:pf:3}\\
&&\hbox{ uniformly for }y\in \Omega\cap \T(B_{Rs_i}(0))\hbox{ and }x\in \Omega\cap \T(B_{R'r_i}(0)\setminus B_{(R')^{-1}r_i}(0)) .\nonumber
\end{eqnarray}
\end{step}

\smallskip\noindent{\it Proof of Step \ref{step:G:pf:3}:} We fix a sequence $(y_i)_i\in \Omega$ such that $y_i\in \T(B_{Rs_i}(0))$ for all $i\in \nn$. For $z\in B_{2R'}\setminus B_{(2R')^{-1}}$, we define
$$\tilde{G}_i(z):=\frac{|y_i|^{\bm}r_i^{\bp-1}}{d(y_i,\partial\Omega)}G_{h_i}(\T(s_i z), y_i).$$
As one checks, \eqref{cv:G:pf:3} is equivalent to prove that
\begin{equation}\label{cv:G:pf:3:equiv}
\tilde{G}_i(x)=(l+o(1))\frac{|x_1|}{|x|^{\bp}}\hbox{ uniformly for }x\in B_{R'}(0)\setminus B_{(R')^{-1}}(0)
\end{equation}
Since $|y_i|=o(r_i)$ as $i\to +\infty$ and \eqref{rem:T:bdry} holds, it follows from the control \eqref{est:G:up} that there exists $C>0$ such that
\begin{equation}\label{pf:G:3}
\frac{1}{C}\cdot \frac{|z_1|}{|z|^{\bp}}\leq \tilde{G}_i(z)\leq C \cdot \frac{|z_1|}{|z|^{\bp}}\hbox{ for all }z\in\rnm\cap  B_{2R'}\setminus B_{(2R')^{-1}}.
\end{equation}
As for \eqref{est:eq:Gi:2}, it follows from \eqref{eq:G:c2} that
\begin{equation}\label{pf:G:4}
-\Delta_{g_i} \tilde{G}_i-\left(\frac{\gamma r_i^2}{|\T(r_i \cdot)|^2}+O(r_i^2)\right)\tilde{G}_i=0\hbox{ in }B_{2R'}(0)\cap\rnm\hbox{ ; }\tilde{G}_i\equiv 0\hbox{ on }\partial\rnm\cap B_{2R'}(0)\setminus\{0\}.
\end{equation}
It follows from \eqref{pf:G:3}, \eqref{pf:G:4} and standard elliptic theory that there exists $\tilde{G}\in C^2(\overline{\rnm}\setminus\{0\})$ such that, up to a subsequence, 
\begin{equation}\label{lim:Gi:pf:bis}
\lim_{i\to +\infty}\tilde{G}_i=\tilde{G}\hbox{ in }C^2_{loc}(\overline{\rnm}\setminus\{0\})
\end{equation}
with
$$-\Delta \tilde{G}-\frac{\gamma}{|x|^2}\tilde{G}=0\hbox{ in }\overline{\rnm}\setminus\{0\}\, ;\, \tilde{G}=0\hbox{ on }\partial\rnm\setminus\{0\}\, ;$$
$$\frac{1}{C}\cdot \frac{|z_1|}{|z|^{\bp}}\leq \tilde{G}(z)\leq C \cdot \frac{|z_1|}{|z|^{\bp}}\hbox{ for all }z\in \rnm\setminus\{0\}.$$
It the follows from Proposition 6.4 in \cite{gr4} that there exists $\mu>0$ such that
\begin{equation}\label{G:pf:explicit:bis}
\tilde{G}(z)=\mu \cdot \frac{|z_1|}{|z|^{\bp}}\hbox{ for all }z\in \rnm.
\end{equation}
We claim that $\mu=l$. We prove the claim. It follows from \eqref{cv:G:pf:2} and the definition \eqref{def:theta} of $\theta_i$ that
\begin{equation}\label{eq:G:l:2}
\tilde{G}_i(\theta_i)=(l+o(1))\frac{|\theta_{i,1}|}{|\theta_i|^{\bp}}\hbox{ and }\theta_i\to\theta_\infty\neq 0\hbox{ as }i\to +\infty.
\end{equation}

\smallskip\noindent{\it Case 1:} we assume that $\theta_\infty\in\rnm\setminus\{0\}$, that is $\theta_{\infty,1}\neq 0$. Passing to the limit in \eqref{eq:G:l:2}, using the convergence \eqref{lim:Gi:pf:bis} and the explicit form \eqref{G:pf:explicit:bis}, as in Case 1 of Step \ref{step:G:pf:1}, we get that $l|\theta_{\infty,1}|\cdot |\theta_{\infty}|^{-\bm}=\mu|\theta_{\infty,1}|\cdot|\theta_\infty|^{-\bm}$, and therefore, since $\theta_{\infty,1}\neq 0$, we get that $\mu=l$.

\smallskip\noindent{\it Case 2:} we assume that $\theta_\infty\in \partial\rnm\setminus\{0\}$, that is $\theta_{i,1}\to 0$ as $i\to +\infty$. With a Taylor expansion, we get that there exists a sequence $(\tilde{t}_i)_{i\in\nn}\in (0,1)$ such that $\tilde{G}_i(\theta_i)=\partial_1\tilde{G}_i(\tilde{t}_i \theta_{i,1},\theta_i')\theta_{i,1}$ for all $i\in \nn$. With the convergence \eqref{lim:Gi:pf:bis} of $\tilde{G}_i$ to $\tilde{G}$  in $C^1$, we get that 
$$\tilde{G}_i(\theta_i)=\left(\partial_1\tilde{G}(\theta_\infty)+o(1)\right)\cdot\theta_{i,1}=\left(\frac{\mu}{|\theta_\infty|^{\bp}}+o(1)\right)\cdot |\theta_{i,1}|.$$
Since $\theta_{i,1}\neq 0$ for all $i\in\nn$, it follows form \eqref{eq:G:l:2} that $\mu=l$.

\smallskip\noindent Therefore, in both cases, we have proved that $\mu=l$. It follows from this uniqueness that the convergence of $\tilde{G}_i$ holds with no extraction. As for Step \ref{step:G:pf:1}, we get \eqref{cv:G:pf:1}. This ends Step \ref{step:G:pf:3}.\qed

\begin{step}\label{step:G:pf:4} We fix $R,R'>0$. We claim that
\begin{equation}
G_{h_i}(x,y)=(l+o(1)+O(|x|^{\bp-\bm}))\frac{d(x,\partial\Omega)}{|x|^{\bp}}\frac{d(y,\partial\Omega)}{|y|^{\bm}}\hbox{ as }i\to +\infty\label{cv:G:pf:4}
\end{equation}
uniformly for $y\in \Omega\cap \T(B_{Rs_i}(0))$ and $x\in \Omega\setminus \T(B_{(R')^{-1}r_i}(0))$.
\end{step}

\smallskip\noindent{\it Proof of Step \ref{step:G:pf:4}:} The differs from Step \ref{step:G:pf:2} since one works on domains exteriors to the ball of radius $r_i$. Here again, we choose $(y_i)_i$ such that $y_i\in \T(B_{Rs_i}(0))$. For $r>0$ small, we choose $\bar{u}_{\bp}\in C^2(\Omega\cap B_r(0))$ a supersolution to $-\Delta \bar{u}_{\bp}-(\gamma|x|^{-2}+h_i)\bar{u}_{\bp}>0$ as in \eqref{ppty:ua} and \eqref{asymp:ua:plus}. Note that, due to the convergence of $(h_i)$ to $h$ in $C^0$, the choice of $\bar{u}_{\bm}$ can be made independently of $i$. We fix $\eps>0$. It follows from the convergence \eqref{cv:G:pf:3} of Step \ref{step:G:pf:3} and \eqref{asymp:ua:plus} that there exists $i_0\in\nn$
\begin{equation}\label{est:sup:G:pf:4}
G_{h_i}(x,y_i)\leq (l+\eps)\frac{d(y_i,\partial\Omega)}{|y_i|^{\bm}}\bar{u}_{\bp}(x)\hbox{ for all }x\in \Omega\cap\partial  \T(B_{R'r_i}(0))\hbox{ for all }i\geq i_0.
\end{equation}
We fix $\delta>0$ such that $\delta<r$. We choose a supersolution $\bar{u}_{\bm}$ as in \eqref{ppty:ua} and \eqref{asymp:ua:plus}. It follows from the upper bound \eqref{est:G:up} that for some $i_1\in\nn$, there exists $C>0$ such that
\begin{equation}\label{est:sup:G:pf:5}
G_{h_i}(x,y_i)\leq C\frac{d(y_i,\partial\Omega)}{|y_i|^{\bm}}\bar{u}_{\bm}(x)\hbox{ for all }x\in \Omega\cap\partial B_\delta(0)\hbox{ for all }i\geq i_1.
\end{equation}
Therefore,
\begin{equation}\label{est:G:up:34}
G_{h_i}(x,y_i)\leq w_i(x)\hbox{ for all }x\in \partial\left(\Omega\cap \T(B_\delta(0)\setminus B_{(R')^{-1}r_i}(0))\right)
\end{equation}
where 
$$w_i:=\frac{d(y_i,\partial\Omega)}{|y_i|^{\bm}}\left((l+\eps)\bar{u}_{\bp}+C\bar{u}_{\bm}\right)$$
and, since $\bar{u}_{\bp},\bar{u}_{\bm}$ are supersolution, 
$$-\Delta w_i-\left(\frac{\gamma}{|x|^2}+h_i\right)w_i\geq 0\hbox{ in }\Omega\cap \T(B_\delta(0)\setminus B_{(R')^{-1}r_i}(0)).$$
Since $-\Delta -\left(\gamma|x|^{-2}+h_i\right)$ is coercive, the maximum principle holds and \eqref{est:G:up:34} holds on $\Omega\cap \T(B_\delta(0)\setminus B_{(R')^{-1}r_i}(0))$. With \eqref{asymp:ua:plus}, we get that there exists $i_2\in\nn$ such that
\begin{equation}\label{est:sup:G:pf:45}
G_{h_i}(x,y_i)\leq \left(l+2\eps+C|x|^{\bp-\bm}\right)\frac{d(x,\partial\Omega)}{|x|^{\bp}}\frac{d(y,\partial\Omega)}{|y|^{\bm}}
\end{equation}
for all $x\in \Omega\cap \T(B_\delta(0)\setminus B_{(R')^{-1}r_i}(0))$ for all $i\geq i_2$. Using subsolutions and arguing as above, we get that for some $i_3\in\nn$
\begin{equation}\label{est:sup:G:pf:55}
G_{h_i}(x,y_i)\geq \left(l-2\eps-C|x|^{\bp-\bm}\right)\frac{d(x,\partial\Omega)}{|x|^{\bp}}\frac{d(y,\partial\Omega)}{|y|^{\bm}}
\end{equation}
for all $x\in \Omega\cap \T(B_\delta(0)\setminus B_{(R')^{-1}r_i}(0))$ for all $i\geq i_3$. The inequalities \eqref{est:sup:G:pf:45} and \eqref{est:sup:G:pf:55} put together yield \eqref{cv:G:pf:4}. This ends Step \ref{step:G:pf:4}.\qed

\begin{step}\label{step:G:pf:5} We let $(X_i)_i,(Y_i)_i\in \Omega$ such that $|Y_i|=o(|X_i|)$ and $X_i=o(1)$ as $i\to +\infty$. We assume that there exists $l'>0$ such that
$$G_{h_i}(X_i,Y_i)=(l'+o(1))\frac{d(X_i,\partial\Omega)}{|X_i|^{\bp}}\frac{d(Y_i,\partial\Omega)}{|Y_i|^{\bm}}\hbox{ as }i\to +\infty.$$
Then $l'=l$.
\end{step}

\smallskip\noindent{\it Proof of Step \ref{step:G:pf:5}:} We define
$$\sigma_i:=\min\{|\tyi|, |Y_i|\}\hbox{ and }\rho_i:=\max\{|\txi|, |X_i|\}.$$
We let $(z_i)_i, (t_i)_i\in\Omega$ such that $c_1\sigma_i\leq |z_i|\leq c_2\sigma_i$ and $c_1\rho_i\leq |t_i|\leq c_2\rho_i$ for all $i\in\nn$. Since $|z_i|=O(s_i)$, $r_i=O(|t_i|)$ and $t_i\to 0$ as $i\to +\infty$, it follows from \eqref{cv:G:pf:4} that
$$G_{h_i}(z_i,t_i)=(l+o(1))\frac{d(z_i,\partial\Omega)}{|z_i|^{\bm}}\frac{d(t_i,\partial\Omega)}{|t_i|^{\bp}}\hbox{ as }i\to +\infty.$$
In addition, since $|z_i|=O(|Y|_i)$, $|X_i|=O(|t_i|)$ and $t_i\to 0$ as $i\to +\infty$, it follows from \eqref{cv:G:pf:4} that
$$G_{h_i}(z_i,t_i)=(l'+o(1))\frac{d(z_i,\partial\Omega)}{|z_i|^{\bm}}\frac{d(t_i,\partial\Omega)}{|t_i|^{\bp}}\hbox{ as }i\to +\infty.$$
Therefore, we get that $l'=l$. This ends Step \ref{step:G:pf:5}.\qed

\begin{step}\label{step:G:pf:6} We let $(X_i)_i,(Y_i)_i\in \Omega$ such that $|Y_i|=o(|X_i|)$ and $X_i=o(1)$ as $i\to +\infty$. Then
$$G_{h_i}(X_i,Y_i)=(l+o(1))\frac{d(X_i,\partial\Omega)}{|X_i|^{\bp}}\frac{d(Y_i,\partial\Omega)}{|Y_i|^{\bm}}\hbox{ as }i\to +\infty.$$
\end{step}

\smallskip\noindent{\it Proof of Step \ref{step:G:pf:6}:} We argue by contradiction and we assume that there exists $\eps_0>0$ and a subsequences $(\varphi(i))_i$ such that $|U_{\varphi(i)}-l|\geq\eps_0$ for all $i\in\nn$ where
$$U_i:=\frac{G_{h_i}(X_i,Y_i)|Y_i|^{\bm}|X_i|^{\bp}}{d(X_i,\partial\Omega)d(Y_i,\partial\Omega)}.$$
Since $(U_{\varphi(i)})$ is bounded, up to another extraction, there exists $l^{\prime\prime}>0$ such that $U_{\varphi(i)}\to l^{\prime\prime}$ as $i\to +\infty$. Therefore, $|l-l^{\prime\prime}|\geq\eps_0$ and $l^{\prime\prime}\neq l$. Since \eqref{est:G:pf} holds for the subfamily $(\varphi(i))$, it then follows from Step \ref{step:G:pf:5} that $l^{\prime\prime}=l$, contradicting $l^{\prime\prime}\neq l$. This ends Step \ref{step:G:pf:6}.

\medskip\noindent We are now in position to prove \eqref{est:G:4}, that is the convergence with no extraction of subsequence. It follows from \eqref{est:G:pf} and Step \ref{step:G:pf:5} applied to $(h_i)_i$ and to the null function that there exists a subsequence $(h_{\varphi(i)})$ and $l,L_{\gamma,\Omega}>0$ such that for any $(x_i)_i,(y_i)_i\in\Omega$ such that $|y_i|=o(|x_i|)$ and $x_i=o(1)$ as $i\to +\infty$, then
\begin{equation}\label{G:l1}
G_{h_{\varphi(i)}}(x_i,y_i)=(l+o(1))\frac{d(x_i,\partial\Omega)}{|x_i|^{\bp}}\frac{d(y_i,\partial\Omega)}{|y_i|^{\bm}},
\end{equation}
and
\begin{equation}\label{G:l2}
G_{0}(x_i,y_i)=(L_{\gamma,\Omega}+o(1))\frac{d(x_i,\partial\Omega)}{|x_i|^{\bp}}\frac{d(y_i,\partial\Omega)}{|y_i|^{\bm}}
\end{equation}
as $i\to +\infty$. We fix a sequence $(x_i)_i\in\Omega$ such that $x_i\to 0$ and $d(x_i,\partial\Omega)\geq |x_i|/2$ as $i\to +\infty$. In the distribution sense, we have that
$$-\Delta(G_{h_{\varphi(i)}}(x_i,\cdot)-G_{0}(x_i,\cdot))+h_{\varphi(i)} (G_{h_{\varphi(i)}}(x_i,\cdot)-G_{0}(x_i,\cdot))=(0-h_{\varphi(i)})G_{0}(x_i,\cdot)\hbox{ in }\Omega$$
in the distribution sense and $G_{h_{\varphi(i)}}(x_i,\cdot)-G_{0}(x_i,\cdot)=0$ on $\partial\Omega$. It follows from \eqref{est:G:up} that for any $1<p<\frac{n}{n-2}$, we have that $\Vert G_{0}(x_i,\cdot)\Vert_p\leq C(p)$ for all $i\in\nn$. It then follows from elliptic theory that $G_{h_{\varphi(i)}}(x_i,\cdot)-G_{h}(x_i,\cdot)\in W^{2,p}(\Omega)$ and that 
$$\Vert G_{h_{\varphi(i)}}(x_i,\cdot)-G_{0}(x_i,\cdot)\Vert_{W^{2,p}}\leq C \Vert h_{\varphi(i)}\Vert_\infty$$
For $1<p<\min\{n/2;n/(n-2)\}$, we define $q:=\frac{np}{n-2p}$. Sobolev embeddings then yield
$$\Vert G_{h_{\varphi(i)}}(x_i,\cdot)-G_{0}(x_i,\cdot)\Vert_{L^q(\Omega)}\leq C \Vert h_{\varphi(i)}\Vert_\infty.$$
We let $(\eps_i)_>0$ such that $\eps_i\to 0$ as $i\to +\infty$. We define $\alpha_i:=\eps_i|x_i|$ so that $\alpha_i=o(|x_i|)$ as $i\to +\infty$. We have that
$$\int_{B_{\alpha_i}(0)}\left|G_{h_{\varphi(i)}}(x_i,y)-G_{0}(x_i,y)\right|^q\, dy\leq C \Vert h_{\varphi(i)}\Vert_\infty^q.$$
It then follows from \eqref{G:l1}, \eqref{G:l2} and the boundedness of $(h_i)$ in $C^0$ that
$$\int_{B_{\alpha_i}(0)}\left|(l-L_{\gamma,\Omega}+o(1))\frac{d(x_i,\partial\Omega)}{|x_i|^{\bp}}\frac{d(y,\partial\Omega)}{|y|^{\bm}}\right|^q\, dy\leq C.$$
We assume by contradiction that $l\neq L_{\gamma,\Omega}$, so that
$$\frac{d(x_i,\partial\Omega)}{|x_i|^{\bp}}\left(\int_{B_{\alpha_i}(0)}\left|\frac{d(y,\partial\Omega)}{|y|^{\bm}}\right|^q\, dy\right)^{1/q}\leq C .$$
If $n\leq q(1-\bm)$, then the integral is infinite. This is a contradiction. Therefore $n>q(1-\bm)$. Estimating the integral and using that $|x_i|\leq 2d(x_i,\partial\Omega)$, we get that
$$|x_i|^{1-\bp}\alpha_i^{1-\bm+\frac{n}{q}}\leq C .$$
With $\alpha_i=\eps_i|x|_i$, $\bm+\bp=n$ and the definition of $q$, we get that 
$$|x_i|^{-n\left(1-\frac{1}{p}\right)}\eps_i^{1-\bm+\frac{n}{q}}\leq C.$$
Since $|x_i|\to 0$, with a suitable choice of $\eps_i\to 0$, we get a contradiction.

\smallskip\noindent Therefore $l=L_{\gamma,\Omega}$ that is independent of the choice of the sequence $(h_i)$. This proves \eqref{est:G:4} and ends the proof of Theorem \ref{th:green:gamma:domain}.

\section[Appendix E: Green's function on $\rnm$]{Appendix E: Green's function for the Hardy-Schr\"odinger operator on $\rnm$}\label{sec:G:rnm}
In this section, we prove the following:
\begin{theorem}\label{th:green:gamma:rn} Fix $\gamma<\frac{n^2}{4}$. For all $p\in\rnm\setminus\{0\}$, there exists $G_p\in L^1(\rnm)$ such that 

\smallskip\noindent{\bf (i)} $\eta G_p\in H_{1,0}^2(\rnm)$ for all $\eta\in C^\infty_c(\rn-\{p\})$,

\smallskip\noindent{\bf (ii)} For all $\varphi\in C^\infty_c(\rnm)$, we have that
\begin{equation}\label{eq:kernel}
\varphi(p)=\int_{\rnm}G_p(x)\left(-\Delta \varphi-\frac{\gamma }{|x|^2}\varphi\right)\, dx,
\end{equation}
\medskip\noindent Moreover, if $G_p,G'_p$ satisfy $(i)$ and $(ii)$ and are positive, then there exists $C\in\rr$ such that $G_p(x)-G'_p(x)=C|x_1|\cdot|x|^{-\bm}$ for all $x\in\rnm\setminus\{0,p\}$.

\medskip\noindent In particular, there exists one and only one function ${\mathcal G}_p={\mathcal G}(p,\cdot)>0$ such that (i) and (ii) hold with $G_p={\mathcal G}_p$ and

\smallskip\noindent{\bf (iii) } ${\mathcal G}_p(x)=O\left(\frac{|x_1|}{|x|^{\bp}}\right)\hbox{ as }|x|\to+\infty. $

We then say that ${\mathcal G}$ is the Green's function for $-\Delta-\gamma|x|^{-2}$ on $\rnm$ with Dirichlet boundary condition.

\medskip\noindent In addition, ${\mathcal G}$ satisfies the following properties:

\smallskip\noindent{\bf (iv) }  For all $p\in\rnp$, there exists $c_0(p),c_\infty(p)>0$ such that
\begin{equation}\label{eq:23}
{\mathcal G}_p(x)\sim_{x\to 0} \frac{c_0(p)|x_1|}{|x|^{\bm}}\hbox{ and }{\mathcal G}_p(x)\sim_{x\to \infty} \frac{c_\infty(p)|x_1|}{|x|^{\bp}}
\end{equation}
and \begin{equation}\label{eq:24}
{\mathcal G}_p(x)\sim_{x\to p}\frac{1}{(n-2)\omega_{n-1}|x-p|^{n-2}}.
\end{equation}
\smallskip\noindent{\bf (v) }  There exists $c>0$ independent of $p$ such that
\begin{equation}\label{est:G:glob}
c^{-1} {\mathcal H}_p(x)\leq {\mathcal G}_p(x)\leq c {\mathcal H}_p(x)
\end{equation}
where
\begin{equation}\label{def:Hp}
{\mathcal H}_p(x):=\left(\frac{\max\{|p|,|x|\}}{\min\{|p|,|x|\}}\right)^{\bm}|x-p|^{2-n}\min\left\{1,\frac{|x_1|\cdot|p_1|}{|x-p|^2}\right\}
\end{equation}
\end{theorem}

\medskip\noindent{\it Proof of Theorem \ref{th:green:gamma:rn}:} 
We shall again proceed with several steps. 

\smallskip\noindent{\bf Step \ref{sec:G:rnm}.1: Construction of a positive kernel at a given point:} For a fixed $p_0\in\rn\setminus\{0\}$, we show that there exists $G_{p_0}\in C^2(\overline{\rnm}\setminus\{0,p_0\})$ such that
\begin{equation}\label{lim:G:bis}
\left\{\begin{array}{ll}
-\Delta G_{p_0}-\frac{\gamma}{|x|^2}G_{p_0}=0&\hbox{ in }\rnm\setminus\{0,p_0\}\\
G_{p_0}>0&\\G_{p_0}\in L^{\frac{2n}{n-2}}(B_\delta(0)\cap\rnm)&\hbox{ with }\delta:=|p_0|/4\\
G_{p_0}\hbox{ satisfies }(ii)\hbox{ with }p=p_0.
\end{array}\right.
\end{equation}
Indeed, let $\tilde{\eta}\in C^\infty(\rr)$ be a nondecreasing function such that $0\leq \tilde{\eta}\leq 1$, $\tilde{\eta}(t)=0$ for all $t\leq 1$ and $\tilde{\eta}(t)=1$ for all $t\geq 2$. For $\eps>0$, set $\eta_\eps(x):=\tilde{\eta}\left(\frac{|x|}{\eps}\right)$ for all $x\in \rn$.

\medskip\noindent We let $\Omega_1$ be a smooth bounded domain of $\rn$ such that $\rnm\cap B_1(0)\subset \Omega_1\subset \rnm\cap B_3(0)$. We define $\Omega_R:=R\cdot \Omega_1$ so that
$\rnm\cap B_R(0)\subset \Omega_R\subset \rnm\cap B_{3R}(0).$
We argue as in the proof of \eqref{bnd:coer} to deduce that the operator $-\Delta-\frac{\gamma\eta_\eps}{|x|^2}$ is coercive on $\Omega_R$ and that there exists $c>0$ independent of $R,\eps>0$ such that
\begin{equation*}
\int_{\Omega_R}\left(|\nabla \varphi|^2-\frac{\gamma\eta_\eps}{|x|^2}\varphi^2\right)\, dx\geq c\int_{\Omega_R}|\nabla \varphi|^2\, dx \quad \hbox{for all $\varphi\in C^\infty_c(\Omega_R)$.}
\end{equation*}

\medskip\noindent Consider $R,\eps>0$ such that $R>2|p_0|$ and $\eps<\frac{|p_0|}{6}$, and let $G_{R,\eps}$ be the Green's function of $-\Delta-\frac{\gamma\eta_\eps}{|x|^2}$ in $\Omega_R$ with Dirichlet boundary condition. We have that $G_{R,\eps}>0$ since the operator is coercive.

\medskip\noindent Fix $R_0>0$ and $q'\in (1,\frac{n}{n-2})$, then by arguing as in the proof of \eqref{int:bnd:G}, we get that there exists $C=C(\gamma,p_0, q', R_0)$ such that 
 \begin{equation}\label{bnd:G:1}
 \Vert G_{R,\eps}(p_0,\cdot)\Vert_{L^{q'}(B_{R_0}(0)\cap\rnm)}\leq C\hbox{ for all }R>R_0\hbox{ and }0<\eps<\frac{|p_0|}{6}, 
 \end{equation} 
 and
 \begin{equation}\label{bnd:G:0}
 \Vert G_{R,\eps}(p_0,\cdot)\Vert_{L^{\frac{2n}{n-2}}(B_{\delta_0}(0)\cap\rnm)}\leq C\hbox{ for all }R>R_0\hbox{ and }0<\eps<\frac{|p_0|}{6},
 \end{equation} 
where $\delta:=|p_0|/4$. Arguing again as in Step \ref{sec:app:c}.2 of the proof of Theorem \ref{th:green:gamma:domain}, there exists $G_{p_0}\in C^2(\overline{\rnm}\setminus\{0,p_0\})$ such that
\begin{equation}\label{lim:G}
\left\{\begin{array}{ll}
G_{R,\eps}(p_0,\cdot)\to G_{p_0}\geq 0&\hbox{ in }C^2_{loc}(\overline{\rnm}\setminus\{0,p_0\})\hbox{ as }R\to +\infty,\; \eps\to 0\\
-\Delta G_{p_0}-\frac{\gamma}{|x|^2}G_{p_0}=0&\hbox{ in }\rnm\setminus\{0,p_0\}\\
G_{p_0}\equiv 0\hbox{ on }\partial\rnm\setminus\{0\}\\
G_{p_0}\in L^{\frac{2n}{n-2}}(B_\delta(0)\cap \rnm)
\end{array}\right.
\end{equation}
and $\eta G_{p_0}\in H_{1,0}^2(\rnm)$ for all $\eta\in C^\infty_c(\rn\setminus\{p_0\})$. Fix $\varphi\in C^\infty_c(\rnm)$. For $R>0$ large enough, we have that $\varphi(p_0)=\int_{\rnm}G_{R,\eps}(p_0,\cdot)(-\Delta\varphi-\gamma\eta_\eps|x|^{-2}\varphi)\, dx$. The integral bounds above yield $x\mapsto G_{p_0}(x) |x|^{-2}\in L^1_{loc}(\rnm)$. Therefore, we get 
\begin{equation}\label{formula:G}
\varphi(p_0)=\int_{\rnm}G_{p_0}(x)\left(-\Delta \varphi-\frac{\gamma }{|x|^2}\varphi\right)\, dx\hbox{ for all }\varphi\in C^\infty_c(\rnm).
\end{equation}
As a consequence, $G_{p_0}>0$.

\medskip\noindent{\bf Step \ref{sec:G:rnm}.2: Asymptotic behavior at $0$ and $p_0$ for solutions to \eqref{lim:G:bis}.}  It  follows from Theorem 6.1 in Ghoussoub-Robert \cite{gr4} that either $G_{p_0}$ behaves like $|x_1|\cdot|x|^{-\bm}$ or $|x_1|\cdot|x|^{-\bp}$ at $0$. Since $G_{p_0}\in L^{\frac{2n}{n-2}}(B_\delta(0)\cap\rnm)$ for some small $\delta>0$ and $\bm<\frac{n}{2}<\bp$, we  get that there exists $c_0>0$ such that 
\begin{equation}\label{asymp:G:0}
\lim_{x\to 0}\frac{G_{p_0}(x)}{|x_1|\cdot |x|^{-\bm}}=c_0.
\end{equation}

\noindent Since $G_{p_0}$ is positive and smooth in a neighborhood of $p_0$, it follows from \eqref{formula:G} and the classification of 
solutions to harmonic equations that
\begin{equation}\label{asymp:p}
G_{p_0}(x)\sim_{x\to p_0}\frac{1}{(n-2)\omega_{n-1}|x-p_0|^{n-2}}.
\end{equation}

\medskip\noindent{\bf Step \ref{sec:G:rnm}.3: Asymptotic behavior at $\infty$ for solutions to \eqref{lim:G:bis}:} We let 
$$\tilde{G}_{p_0}(x):=\frac{1}{|x|^{n-2}}G_{p_0}\left(\frac{x}{|x|^2}\right)\hbox{ for all }x\in\rnm\setminus\left\{0, \frac{p_0}{|p_0|^2}\right\},$$ 
be the Kelvin's transform of $G$. We have that
 $$-\Delta \tilde{G}_{p_0}-\frac{\gamma}{|x|^2}\tilde{G}_{p_0}=0\hbox{ in }\rnm\setminus\left\{0, \frac{p_0}{|p_0|^2}\right\}\; ;\; \tilde{G}\equiv 0\hbox{ on }\partial\rnm\setminus\{p_0\}.$$
Since $\tilde{G}_{p_0}>0$, it follows from Theorem 6.1 in \cite{gr4} that there exists $c_1>0$ such that
\begin{equation*}
\hbox{either }\tilde{G}_{p_0}(x)\sim_{x\to 0}c_1\frac{|x_1|}{|x|^{\bm}}\hbox{ or }\tilde{G}_{p_0}(x)\sim_{x\to 0}c_1\frac{|x_1|}{|x|^{\bp}}.
\end{equation*}
Coming back to $G_{p_0}$, we get that
\begin{equation}\label{choice}
\hbox{either }G_{p_0}(x)\sim_{|x|\to \infty}c_1\frac{|x_1|}{|x|^{\bp}}\hbox{ or }G_{p_0}(x)\sim_{|x|\to \infty}c_1\frac{|x_1|}{|x|^{\bm}}.
\end{equation}
Assuming we are in the second case, for any $c\leq c_1$, we define
$$\bar{G}_c(x):=G_{p_0}(x)-c\frac{|x_1|}{|x|^{\bm}}\hbox{ in }\rnm\setminus\{0, p_0\},$$
which satisfy $-\Delta \bar{G}_c-\frac{\gamma}{|x|^2}\bar{G}_c=0$ in $\rnm\setminus\{0, p_0\}$. It follows from \eqref{choice} and \eqref{asymp:p} that for $c<c_1$, $\bar{G}_c>0$ around $p_0$ and $\infty$. Using that $\eta \bar{G}_c\in H_{1,0}^2(\rnm)$ for all $\eta\in C^\infty_c(\rn\setminus\{p_0\})$, it follows from the coercivity of $-\Delta-\gamma|x|^{-2}$ that $\bar{G}_c>0$ in $\rnm\setminus\{0, p_0\}$ for $c<c_1$. Letting $c\to c_1$ yields $\bar{G}_{c_1}\geq 0$, and then $\bar{G}_{c_1}> 0$. Since $\bar{G}_{c_1}(x)=o(|x_1|\cdot|x|^{-\bm})$ as $|x|\to \infty$, another Kelvin transform and Theorem 6.1 in \cite{gr4} yield  $|x_1|^{-1}|x|^{\bp}\bar{G}_{c_1}(x)\to c_2>0$ as $|x|\to \infty$ for some $c_2>0$. Then there exists $c_3>0$ such that

\begin{equation}\label{again}
\lim_{x\to 0}\frac{\bar{G}_{c_1}(x)}{|x_1|\cdot |x|^{-\bm}}=c_3>0\hbox{ and }\lim_{x\to \infty}\frac{\bar{G}_{c_1}(x)}{|x_1|\cdot |x|^{-\bp}}=c_2.
\end{equation}
Since $x\mapsto |x_1|\cdot|x|^{-\bm}\in H_{1,loc}^2(\rn)$, we get that 
$$\varphi(p)=\int_{\rnm}\bar{G}_{c_1}(x)\left(-\Delta \varphi-\frac{\gamma }{|x|^2}\varphi\right)\, dx\hbox{ for all }\varphi\in C^\infty_c(\rnm).$$

\medskip\noindent{\bf Step \ref{sec:G:rnm}.4: Uniqueness:} Let $G_1,G_2>0$ be 2 functions such that $(i),(ii)$ hold for $p:=p_0$, and set $H:=G_1-G_2$. It follows from Steps 2 and 3 that there exists $c\in\rr$ such that $H'(x):=H(x)-c|x_1|\cdot|x|^{-\bm}$ satisfies
\begin{equation}\label{bnd:H:prime}
H'(x)=_{x\to 0}O\left(|x_1|\cdot|x|^{-\bm}\right)\hbox{ and }H'(x)=_{|x|\to \infty}O\left(|x_1|\cdot|x|^{-\bp}\right).
\end{equation}
We then have that $\eta H'\in H_{1,0}^2(\rnm)$ for all $\eta\in C^\infty_c(\rn\setminus\{p_0\})$ and
$$\int_{\rnm}H'(x)\left(-\Delta\varphi-\frac{\gamma}{|x|^2}\varphi\right)\, dx=0 \quad \hbox{for all $\varphi\in C^\infty_c(\rnm)$.}$$
The ellipticity of the Laplacian then yields $H'\in C^\infty(\overline{\rnm}\setminus\{0\})$. The pointwise bounds \eqref{bnd:H:prime} yield that $H'\in H_{1,0}^2(\rnm)$. Multiplying $-\Delta H'-\frac{\gamma}{|x|^2}H'=0$ by $H'$, integrating by parts and the coercivity yield $H'\equiv 0$, and therefore, $(G_1-G_2)(x)=c|x_1|\cdot|x|^{-\bm}$ for all $x\in\rnm$. This proves uniqueness.

\medskip\noindent{\bf Step \ref{sec:G:rnm}.5: Existence.}  It follows from Steps 2 and 3 that, up to substracting a multiple of $x\mapsto |x_1|\cdot|x|^{-\bm}$, there exists a unique function ${\mathcal G}_{p_0}>0$ satisfying (i), (ii) and the pointwise control (iii). Moreover,  \eqref{asymp:G:0}, \eqref{asymp:p} and \eqref{again} yield \eqref{eq:23} and \eqref{eq:24}. As a consequence, \eqref{est:G:glob} holds with $p=p_0$.

\smallskip\noindent For $p\in\rnp$, consider $\rho_p: \rnm\to\rnm$ a linear isometry fixing $\rnm$ such that $\rho_p(\frac{p_0}{|p_0|})=\frac{p}{|p|}$, and define
$${\mathcal G}_p(x):=\left(\frac{|p_0|}{|p|}\right)^{n-2}{\mathcal G}_{p_0}\left(\left(\rho_p^{-1}\left(\frac{|p_0|}{|p|}x\right)\right)\right)\hbox{ for all }x\in \rn\setminus\{0,p\}.$$
As one checks, ${\mathcal G}_p>0$ satisfies (i), (ii), (iii), \eqref{eq:23}, \eqref{eq:24} and \eqref{est:G:glob}.

\smallskip\noindent The definition of ${\mathcal G}_p$ is independent of the choice of $\rho_p$. Indeed, for any linear isometry $\rho_{p_0}: \rnm\to\rnm$ fixing $p_0$ and $\rnm$, ${\mathcal G}_{p_0}\circ \rho_{p_0}^{-1}$ satisfies  (i), (ii), (iii), and therefore ${\mathcal G}_{p_0}\circ \rho_{p_0}^{-1}={\mathcal G}_{p_0}$. The argument goes similarly of any isometry fixing $p$.

\end{document}